\documentclass[11pt]{article}
\usepackage{amsmath,amssymb} 
\usepackage{natbib} 
\usepackage[english]{babel}
\usepackage{graphicx,color}
\usepackage[dvipsnames]{xcolor}
\usepackage[font=small,labelfont=bf]{caption}
\usepackage{hyperref}
\usepackage{geometry}	
\usepackage{lscape} 
\usepackage[shortcuts]{extdash}
\usepackage{empheq}
\usepackage{enumitem}
\usepackage{slashbox}
\usepackage{multirow}
\usepackage{hhline} 
\usepackage{placeins} 
\definecolor{darkblue}{rgb}{0,0,1}
\definecolor{matlab_blue}{rgb}{0,0.447,0.741}
\definecolor{matlab_red}{rgb}{0.85,0.325,0.098}
\definecolor{matlab_yellow}{rgb}{0.929,0.694,0.125}

\hypersetup{colorlinks=true, breaklinks=true, linkcolor=darkblue, menucolor=darkblue, citecolor=darkblue, urlcolor=darkblue}

\newcommand{\bitm}{\begin{itemize}}
\newcommand{\eitm}{\end{itemize}}
\newcommand{\bnumr}{\begin{enumerate}}
\newcommand{\enumr}{\end{enumerate}}

\newcommand{\mcalF}{\mathcal{F}}
\newcommand{\mcalG}{\mathcal{G}}
\newcommand{\mcalI}{\mathcal{I}}
\newcommand{\mcalL}{\mathcal{L}}

\newcommand{\mcalP}{\mathcal{P}}
\newcommand{\mcalS}{\mathcal{S}}
\newcommand{\mcalT}{\mathcal{T}}

\newcommand{\bcalG}{\boldsymbol{\mathcal{G}}}

\newcommand{\bcalN}{\boldsymbol{\mathcal{N}}}

\newcommand{\bcalT}{\boldsymbol{\mathcal{T}}}

\newcommand{\mrA}{\mathrm{A}}

\newcommand{\mrN}{\mathrm{N}}
\newcommand{\mrT}{\mathrm{T}}

\newcommand {\eqb}[1]{\begin{equation}\begin{array}{#1}}
\newcommand {\eqe}{\end{array}\end{equation}}

\newcommand {\esb}[1]{\begin{equation*}\begin{array}{#1}}
\newcommand {\ese}{\end{array}\end{equation*}}

\newcommand {\tr}{\mathrm{tr}\,}

\newcommand {\grad}{\mathrm{grad}\,}

\newcommand {\divz}{\mathrm{div}\,}

\newcommand {\II}{{I\kern-.3em I}}
\newcommand {\III}{{I\kern-.3em I\kern-.3em I}}

\newcommand {\mrd}{\mathrm{d}}
\newcommand {\mre}{\mathrm{e}}
\newcommand {\mrf}{\mathrm{f}}

\newcommand {\mrn}{\mathrm{n}}

\newcommand {\mrp}{\mathrm{p}}

\newcommand {\mrr}{\mathrm{r}}

\newcommand {\mrt}{\mathrm{t}}

\newcommand {\mrv}{\mathrm{v}}
\newcommand {\mrw}{\mathrm{w}}

\newcommand {\mrF}{\mathrm{F}}

\newcommand {\md}{\mathbf{d}}

\newcommand {\mf}{\mathbf{f}}

\newcommand {\mt}{\mathbf{t}}

\newcommand {\mv}{\mathbf{v}}
\newcommand {\mw}{\mathbf{w}}
\newcommand {\mx}{\mathbf{x}}

\newcommand {\bb}{\boldsymbol{b}}

\newcommand {\be}{\boldsymbol{e}}

\newcommand {\bi}{\boldsymbol{i}}

\newcommand {\bn}{\boldsymbol{n}}

\newcommand {\br}{\boldsymbol{r}}

\newcommand {\bt}{\boldsymbol{t}}

\newcommand {\bv}{\boldsymbol{v}}

\newcommand {\bx}{\boldsymbol{x}}
\newcommand {\by}{\boldsymbol{y}}

\newcommand {\bzero}{\boldsymbol{0}}

\newcommand {\bxi}{\mbox{\boldmath$\xi$}}
\newcommand {\bet}{\mbox{\boldmath$\eta$}}

\newcommand {\bome}{\mbox{\boldmath$\omega$}}

\newcommand {\bsig}{\mbox{\boldmath$\sigma$}}

\newcommand {\bXi}{\mbox{\boldmath$\Xi$}}

\newcommand {\mA}{\mathbf{A}}

\newcommand {\mG}{\mathbf{G}}

\newcommand {\mN}{\mathbf{N}}

\newcommand {\mT}{\mathbf{T}}

\newcommand {\bD}{\boldsymbol{D}}

\newcommand {\bG}{\boldsymbol{G}}

\newcommand {\bT}{\boldsymbol{T}}

\newcommand {\bone}{\mathbf{1}}

\newcommand {\IR}{{\rm\kern.24em
   \vrule width.02em height1.53ex depth-.05ex
   \kern-.3em R}}
\newcommand {\ic}{{\rm\kern.20em
   \vrule width.02em height1.0ex depth-.05ex
   \kern-.22em c}}
\newcommand {\ia}{{\rm\kern.20em
   \vrule width.02em height1.05ex depth-.0ex
   \kern-.25em a}}
\newcommand {\IC}{{\rm\kern.24em
   \vrule width.02em height1.4ex depth-.05ex
   \kern-.26em C}}
\newcommand {\ID}{{\rm\kern.34em
   \vrule width.02em height1.5ex depth-.05ex
   \kern-.36em D}}
\newcommand {\IS}{{\rm\kern.24em
   \vrule width.02em height1.6ex depth.05ex
   \kern-.26em S}}
\newcommand {\IT}{{\rm\kern.50em
   \vrule width.02em height1.55ex depth-.05ex
   \kern-.52em T}}

\newcommand {\IE}{{\rm\kern.24em
   \vrule width.02em height1.55ex depth-.05ex
   \kern-.33em E}}
\newcommand {\IEa}{{\rm\kern.24em
   \vrule width.02em height1.55ex depth-.05ex
   \kern-.33em E}^{1}_{ijkl}}
\newcommand {\IEb}{{\rm\kern.24em
   \vrule width.02em height1.55ex depth-.05ex
   \kern-.33em E}^{2}_{ijkl}}

\newcommand {\sS}{\mathcal{S}}

\newcommand {\Ass}[2]{\kern 0.9ex \vrule width0.45em height0.2ex depth0ex \kern -2.1ex \bigwedge_{#1}^{#2}}
\newcommand {\ASS}[2]{\kern 1.45ex \vrule width0.5em height0.2ex depth0ex \kern -2.65ex \bigwedge_{#1}^{#2}}

\newenvironment{packeditemize}{
\begin{itemize}
  \setlength{\listparindent}{-10pt}
  \setlength{\itemsep}{2pt}
  \setlength{\parskip}{0pt}
  \setlength{\parsep}{0pt}
}{\end{itemize}} 

\newcommand{\hh}{\hspace{4mm}}
\newcommand{\sforall} {\hh\forall\,}  
\newcommand{\wbox}[4]{\put({#1},{#2}){\colorbox{white}{\makebox({#3},{#4}){}}}}   
\newcommand{\twhite}[1]{\textcolor{white}{ #1}}						
\newcommand{\tmred}[1]{\textcolor{matlab_red}{ #1}}						
\newcommand{\tmblue}[1]{\textcolor{matlab_blue}{ #1}}						
\newcommand{ \wput}[3]{\put({#1},{#2}){\colorbox{white}{\makebox{#3}}}} 
\newcommand{\vwput}[3]{\put({#1},{#2}){\rotatebox{90}{\colorbox{white}{\makebox{#3}}}}} 

\newcommand{\dof}{n_\mathrm{dof}}
\newcommand{\nel}{n_\mathrm{el}}
\newcommand{\nno}{n_\mathrm{no}}

\begin{document}


\begin{center}
\Large{\bf{New hybrid quadrature schemes for weakly singular kernels applied to isogeometric boundary elements for 3D~Stokes flow}}\\

\end{center}

\begin{center}
\large{Maximilian $\mathrm{Harmel}^a$ and Roger A.~$\mathrm{Sauer}^{a,b,c,}
\renewcommand{\thefootnote}{\fnsymbol{footnote}}\footnote[1]{corresponding author, email: roger.sauer@pg.edu.pl, sauer@aices.rwth-aachen.de}\renewcommand{\thefootnote}{\arabic{footnote}}$	}\\
\vspace{4mm}

\small{\textit{
$\mathrm{ }^a$Aachen Institute for Advanced Study in Computational Engineering Science (AICES), \\ 
RWTH Aachen University, Templergraben 55, 52062 Aachen, Germany\\
$\mathrm{ }^b$Faculty of Civil and Environmental Engineering, Gda\'nsk University of Technology,\\
ul. Narutowicza 11/12, 80-233 Gda\'nsk, Poland\\
$\mathrm{ }^c$Dept.~of Mechanical Engineering, Indian Institute of Technology Guwahati, Assam 781039, India
}}

\end{center}

\vspace{-4mm}

\begin{center}

\small{Published\footnote{This pdf is the personal version of an article whose journal version is available at \href{https://doi.org/10.1016/j.enganabound.2023.04.037}{\!www.sciencedirect.com}}
in \textit{Engineering Analysis with Boundary Elements}, \href{https://doi.org/10.1016/j.enganabound.2023.04.037}{DOI: 10.1016/j.enganabound.2023.04.037} \\
Submitted on 17 October 2022; Received in revised form on 18 April 2023; Accepted on 26 April 2023} 

\end{center}

\vspace{-9mm}

\rule{\linewidth}{.15mm}
{\bf Abstract}

This work proposes four novel hybrid quadrature schemes for the efficient and accurate evaluation of weakly singular boundary integrals ($1/r$ kernel) on arbitrary smooth surfaces. Such integrals appear in boundary element analysis for several partial differential equations including the Stokes equation for viscous flow and the Helmholtz equation for acoustics. The proposed quadrature schemes apply a Duffy transform-based quadrature rule \citep{Duffy82} to surface elements containing the singularity and classical Gaussian quadrature to the remaining elements. Two of the four schemes additionally consider a special treatment for elements near to the singularity, where refined Gaussian quadrature and a new moment-fitting quadrature rule are used.
\\The hybrid quadrature schemes are systematically studied on flat B-spline patches and on NURBS spheres considering two different sphere discretizations: An exact single-patch sphere with degenerate control points at the poles and an approximate discretization that consist of six patches with regular elements. The efficiency of the quadrature schemes is further demonstrated in boundary element analysis for Stokes flow, where steady problems with rotating and translating curved objects are investigated in convergence studies for both, mesh and quadrature refinement. Much higher convergence rates are observed for the proposed new schemes in comparison to classical schemes.

{\bf Keywords:} Boundary element analysis, Duffy quadrature, isogeometric analysis, moment-fitting quadrature, singular integrals, Stokes flow.

\vspace{-5mm}
\rule{\linewidth}{.15mm}


%

\section{Introduction}\label{sec:intro}
%
One of the main difficulties of the boundary element method (BEM) is the efficient approximation of singular integrals that appear in the boundary integral equation (BIE). One approach to overcome this challenge is to superimpose known solutions to the unknown fields such that singularities are removed \citep{Cruse74,Liu99,Liu00,Klaseboer09}. This approach is referred to as nonsingular, regularized or desingularized BEM. It is applied to linear elasticity \citep{Scott13,Taus19}, to Stokes flow \citep{Taus16,Harmel18} and to the Helmholtz equation \citep{Simpson14, Peake15}, among others. \cite{Klaseboer12} further apply nonsingular BEM to fluid mechanics by considering Stokes equations for viscous flow, Laplace equation for potential flow and Helmholtz equation for free-streamline flow.
\\\\The nonsingular BEM avoids singular integrals and is thus commonly used in recent papers, but it requires additional integrals and knowledge about analytical solutions and is further disadvantageous in efficiency and implementation \citep{Khayat05}. An alternative approach considered is to use the classical BIE without regularization and to approximate the singular boundary integrals with special quadrature rules. Such an approach is considered in this paper for weakly singular boundary integrals (kernel proportional to $1/r$)\footnote{$r$ denotes the distance between field point and the singularity} that are required and sufficient for steady and three-dimensional BIEs of several linear problems such as Stokes flow, potential flow, elastostatics, heat conduction and acoustics. Existing quadrature rules for surface elements containing the singularity (denoted \textit{singular elements}) are discussed in Sec.~\ref{sec:intro_sing} and for the adjacent elements near to the singularity (denoted \textit{near singular elements}) in Sec.~\ref{sec:intro_near_sing}. 
The application of these quadrature rules to various boundary element~(BE) problems is then discussed in Sec.~\ref{sec:intro_app}.

\subsection{Approximation of weakly singular integrals}\label{sec:intro_sing}
Considerable progress has been made with quadrature rules based on variable transformations that map the physical domain to a parent domain so that the singularity is removed through the introduction of the Jacobian of the transformation, e.g.~see \cite{Schwartz69,Takahasi73,Lean85,Khayat05} and \cite{Cano15}. 
The Duffy transformation from a square to a triangle \citep{Fairweather79,Duffy82} is of particular interest here. In finite element~(FE) analysis for fracture mechanic, numerical quadrature based on the Duffy transformation~(\textit{Duffy quadrature} in short) is used to integrate singular shape function derivatives \citep{Tracey71, Stern78} and in the context of the extended FE method \citep{Laborde05, Bechet05, Lv18}. 
Duffy quadrature is further used to integrate rational bubble functions with multiple singularities in FE analysis for Stokes flow \citep{Schneier15}. More recently, \cite{Tan19} propose the Duffy-Distance transformation that includes a further mapping that has to be adjusted depending on the shape of the element. Numerical investigations show that weakly singular integrals on spheres and cylinders can be approximated mostly with high accuracy, whereas the location of the singularity and the aspect ratio of the elements have a strong negative impact. Other important quadrature rules based on variable transformations are the tanh rule \citep{Haber77} and the sinh-tanh rule \citep{Borwein06}, the quadrature rule from \cite{Telles87} that is based on a non-linear coordinate transformation and a polar coordinate transformation that is also applicable to elements with collapsed edges \citep{Taus16}.
The application of the quadrature rules from Sec.~\ref{sec:intro_sing}, including Duffy quadrature, to various BE problems is discussed in more detail in Sec.~\ref{sec:intro_app}. 
\\\\In addition to variable transformation methods, there are many other approaches for the approximation of weakly singular boundary integrals: \cite{Guiggiani92} propose the approximation of singular integrals by singularity subtraction, \cite{Dautray85} and \cite{Atkinson95} use piecewise polynomial functions to approximate surface and integrals, while \cite{Niu05} propose semi-analytical integration of the singular kernel.
Recently, \cite{Ochiai22} propose the triple-reciprocity BEM that allows a direct integration on Lagrange elements and \cite{Velazquez22} present a method that constructs efficient quadrature rules for the approximation of singular integrals from the finite part of known integrals including the shape functions.

\subsection{Approximation of weakly near singular integrals}\label{sec:intro_near_sing}
Weakly singular integrals are approximated efficiently with adaptive quadrature rules \citep{Lachat76,Gao00} by dividing near singular elements into subelements and applying Gauss-Legendre quadrature to the subelements. The number of subelements and the number of quadrature points for each subelement are determined individually based upon analytical \citep{Stroud66_book} or numerical \citep{Bu95} criteria. Similarly as on singular elements (see Sec.~\ref{sec:intro_sing}), quadrature rules based on variable transformation can also be used on nearly singular elements to approximate weakly singular integrals. The quadrature rule from \cite{Telles87} is also applicable to near singular elements, while Duffy quadrature \citep{Duffy82} is only accurate on the singular element. Other variable transformation considered for the approximation of weakly near singular integrals include the radial variable transformation \citep{Hayami94}, the distance transformation \citep{Ma02,Xie13} and the sinh transformation \citep{Johnston07,Xie21}. 
The quadrature scheme presented in \cite{Gong20} combines the benefits of the sinh transformation method and adaptive methods such it is capable of integral kernels of the type $1/r$, $1/r^2$ and $1/r^3$ on near singular elements. Many other quadrature rules for weakly singular integrals on nearly singular elements can be found in the literature as for example the parameterized Gaussian quadrature \citep{Lutz92}, line integral approaches \citep{Krishnasamy94,Liu98} and analytical \citep{Padhi01,Zhou08} and semi-analytical methods \citep{Sladek01, Niu05, Han22}. 
\\\\An alternative approach is to construct quadrature rules by the moment fitting method, i.e.~solve a variant formulation of the moment fitting equations for the quadrature weights and abscissas. First, moment fitting quadrature rules for the exact integration of polynomial functions on triangles \citep{Lyness75} and on quadrilaterals \citep{Dunavant85} were presented. Moment fitting quadrature rules for polynomials are formalized and generalized by \cite{Wandzurat03}, while a numerical algorithm for polynomials on triangles, squares and cubes is given by \cite{Xiao10}. In all these approaches the integral on the right hand side of the moment fitting equations is computed analytically. Moment fitting quadrature for discontinuous functions is proposed for crack propagation analysis with the extended FE method \citep{Mousavi10b,Zhang18}, where the moment fitting integrals are evaluated partially analytically and partially numerically applying the method of \cite{Lasserre98,Lasserre99}. \cite{Joulaian16} and \cite{Hubrich17} manipulate the moment fitting integrals so they can be computed by Gaussian quadrature numerically. This quadrature rule approximates discontinuous functions on domains of arbitrary geometry and topology  highly efficiently. Moment fitting quadrature rules are further used for the construction of highly efficient quadrature rules in FE analysis for shells \citep{Zou21}.
The present paper proposes a new Gauss-Legendre quadrature rule with adjusted weights that are determined by the moment fitting method. The proposed quadrature rule determines weakly near singular integrals exactly on plane surfaces with regular quadrilateral elements and accurately on curved surfaces. It requires only a moderate number of quadrature points and uses the same abscissas as standard Gauss-Legendre quadrature. The Gauss-Legendre quadrature rule with adjusted weights is thus simple to implement and computationally efficient.

\subsection{Application to BE analysis}\label{sec:intro_app}
Solving boundary value problems using BE analysis requires to approximate the BIE on the entire surface including singular, near singular and the remaining regular elements. Numerous approaches using the quadrature rules from Sec.~\ref{sec:intro_sing} and Sec.~\ref{sec:intro_near_sing} to approximate the boundary integrals efficiently can be found in the literature, so the following enumeration of application is not exhaustive.
\\\\Duffy quadrature \citep{Duffy82} is applied to the singular element for Stokes flow problems \citep{Varnhorn89,Johnson89,Barakat18} and for acoustic problems \citep{Amini90}. In these approaches, the integrals on the remaining elements are approximated by standard Gauss-Legendre quadrature or by simple trapezoidal rules. \cite{Venaas20} additionally consider local refinement of Gauss-Legendre quadrature on near singular elements for acoustic problems without investigating its influence in detail.
The quadrature rule from \cite{Telles87} is combined with Gauss-Legendre quadrature for linear elastic problems considering homogeneous material \citep{Karam88,Simpson12} and liquid inclusions \citep{Dai21} and for FSI problems considering Stokes Flow \citep{Patino21}. The adaptive integral method for near singular integrals \citep{Gao00} is used to solve Laplace's equation by isogeometric BE analysis for potential problems \citep{Gong17}, thermal problems \citep{Gong18} and thermoelasticity problems \citep{Gong20}, whereas the singularity subtraction method \citep{Guiggiani92} is applied in \cite{Jarvenpaa06} to electromagnetics and in \cite{Peng17} to fracture mechanics. \cite{Ata18} consider the tanh-sinh quadrature rule \citep{Borwein06} for BE analysis of Stokes flow problems and \cite{Keuchel17} apply the sinh-quadrature \citep{Johnston07} to nearly singular integrals in the Burton-Miller formulation of the Helmholtz equation for acoustic problems. The quadrature rule based on polar coordinate transformation \citep{Taus16} is extended in
\cite{Taus19} to handle elements with large curvatures and aspect ratios efficiently and to improve the accuracy on near singular elements. The extended quadrature rule is applied to linear elasticity considering Gaussian quadrature on the regular elements.
\cite{Giuliani18} provide a generalized parallel implementation for adaptive, geometry aware, and high order boundary element methods using quadrature rules from \cite{Lachat76}, \cite{Duffy82} and \cite{Telles87} without investigating their influence on the accuracy of the BE results.
\\\\The above mentioned approaches are limited in application (geometry and discretization) or lack systematic investigation.   
This paper, therefore, proposes four hybrid quadrature schemes which
\begin{packeditemize}
	\item approximate weakly singular integrals on singular, near singular and regular elements accurately,
	\item adaptively combine a new moment fitting quadrature rule with Duffy quadrature \citep{Duffy82} and rings of Gauss-Legendre quadrature with varying quadrature density,
	\item are systematically studied for discretization and quadrature refinement,
	\item show fast convergence for three-dimensional BE examples of steady Stokes flow with flat and curved isogeometric surfaces,
	\item {also apply to other BE formulations including those for Helmholtz' equation and linear elasticity},
	\item are applicable to various discretizations (including isogeometric, Lagrange and Hermite shape functions) using quadrilateral elements of arbitrary genus zero surfaces.
\end{packeditemize}
This paper focus on Stokes flow, but the presented quadrature schemes are also applicable to other BE applications like heat transfer \citep{Mera02,Zang21}, acoustics \citep{Amini90,Venaas20}, elastostatics \citep{Simpson12, Taus19}, mechanical contact \citep{Yac70,Zirakashvili20} and electromagnetics \citep{Rajski19,Takahashi22}.
The quadrature schemes are further applicable to coupled BE and surface FE formulations that can be used to study droplets \citep{Brown80, Sauer14_IJNMF}, bubbles \citep{Wang03,Boedec17}, shells \citep{Heltai17, Maestre17} and wetting \citep{Osman15, Luginsland17}.
They can be used to solve unsteady BE problems \citep{Falletta15,Aimi19,Aimi22} since the singular integrals investigated here also occur there.
\\\\The remainder of this paper is organized as follows: Sec.~\ref{sec:flow} presents an overview of the underlying BE theory for incompressible Stokes flow. Boundary quadrature of weakly singular integrals on singular and near singular elements is discussed in Sec.~\ref{sec:quad}. In addition to known rules, the new moment fitting quadrature rule is also presented in this chapter. Sec.~\ref{sec:hybrid} presents the four new hybrid quadrature scheme and investigates their accuracy and efficiency in detail. The hybrid quadrature schemes are applied to BE analysis for three-dimensional Stokes flow problems in Sec.~\ref{sec:examples}. The paper then concludes with Sec.~\ref{sec:conclusion}.
%
\section{Incompressible Stokes flow}\label{sec:flow}
Incompressible steady-state Stokes flow is briefly outlined in this section: The constitutive relation for an incompressible Newtonian fluid is introduced in Sec.~\ref{sec:flow_consti}, while the governing equations, namely the linear momentum balance for Stokes flow and the mass balance for incompresible flow, are presented in Sec.~\ref{sec:flow_equil}. The singular Green's functions are discussed and investigated in Sec.~\ref{sec:flow_Greens}, and the boundary integral equation (BIE) for Stokes flow is then presented in Sec.~\ref{sec:flow_BIE}.

\subsection{Fluid constitution}\label{sec:flow_consti}
The three-dimensional Cauchy stress tensor for an incompressible Newtonian fluid with dynamic viscosity $\eta$ is given by
\begin{equation}\label{eq:flow_consti}
\bsig = -p\, \bone +2\eta\, \bD~,
\end{equation}
where $\bD:= \left(\grad\bv + (\grad\bv)^\mrT\right)/2$ is the symmetric part of the velocity gradient and $p=-\tr\bsig/3$ is the fluid pressure in domain $\mcalF$. The incompressibility of the fluid is enforced by the continuity equation that is given in absence of mass sources or sinks by
\begin{equation}\label{eq:flow_continuity}
\divz \bv = 0 \hh\mathrm{in}~\mcalF~,
\end{equation}
where $\bv$ denotes the velocity.

\subsection{Fluid equilibrium}\label{sec:flow_equil}

The steady-state motion of a viscous\footnote{Reynolds number $\mathrm{Re}<1$ \citep{Pozrikidis02_book}} fluid flow is governed by the steady Stokes equation
\begin{equation}\label{eq:flow_stokes}
\divz \bsig  = -\rho \,\bar{\bb} \hh\mathrm{in}~\mcalF~,
\end{equation}
together with the continuity equation \eqref{eq:flow_continuity} and the Dirichlet and Neumann boundary conditions
%
\begin{equation}\label{eq:flow_bc}
\begin{aligned}
  \bv(\bx) &= \bar \bv \hh\sforall \bx \in \partial_\mrd \mcalF\,,\\
 \bt(\bx)  &= \bar \bt \hh\sforall \bx \in \partial_\mrn \mcalF\,,
\end{aligned}
\end{equation}
with velocity $\bv$, surface traction $\bt := \bsig \, \bn$ and outward unit normal vector $\bn$. Dirichlet and Neumann boundary regions are denoted by $\partial_\mrd \mcalF$  and $\partial_\mrn \mcalF$, respectively and the entire boundary of the fluid domain is denoted by $\partial \mcalF:=\partial_\mrn \mcalF \cup \partial_\mrd \mcalF$, which coincides here with the surface $\sS$ of an immersed body, as there are no other boundaries considered. Apart from \eqref{eq:flow_bc}, a condition on pressure $p$ is needed on $\mcalS$ to define the flow field. \footnote{This paper considers incompressible flow with no slip between surface and fluid. In this case, the fluid pressure on $\mcalS$ is given by $p(\bx)= - \bt \cdot \bn$, since $\bn \cdot \bD\, \bn = 0$, c.f.~Eq.\eqref{eq:flow_consti}.}
%
\\\\This work focuses on the efficient quadrature of boundary integrals for boundary element~(BE) analysis of steady Stokes flow. However, all following results apply to arbitrary weakly singular boundary integrals including BE for the Helmholtz equation and linear elasticity. The equilibrium equation for linear elasticity is even mathematically equivalent to the steady Stokes equation \eqref{eq:flow_stokes}.\footnote{For linear elasticity, $\bv$ in \eqref{eq:flow_consti} corresponds to the displacement field.} The quadrature schemes presented in this work are further applicable to unsteady BE problems since the singular integrals investigated here also appear there.



\subsection{Boundary representation of the fluid flow}\label{sec:flow_boundary}
This section discusses a representation of steady Stokes equation \eqref{eq:flow_stokes} that only lives on the fluid boundary $\mcalS$. The transformation of Stokes equation into the boundary representation based on the Green's functions is presented in Sec.~\ref{sec:flow_Greens} and discussed in Sec.~\ref{sec:flow_BIE}. This boundary representation, called boundary integral equation (BIE), is the most important equation in boundary element analysis, which is summarized in Sec.~\ref{sec:BE}. The singular nature of the Green's functions leading to the singular integral kernels in the BIE is investigated in Sec.~\ref{sec:flow_Greens_singu}.

\subsubsection{Green's functions for steady Stokes flow}\label{sec:flow_Greens}
Considering a given point load $\bar\bb^\infty(\bx-\by)$, that is applied at the fixed source point $\by$, leads to the singularly forced steady Stokes equation
\begin{equation}\label{eq:flow_stokes_singu}
	\divz \bsig =  - \rho \,\bar\bb^\infty \hh\hh\mathrm{in} ~\mcalF~.
\end{equation}
The distance between an arbitrary field point $\bx \in \mcalF$ and source point $\by$ is defined as $r:=\| \br \|$ with $\br:=\bx-\by$.
The point load can then be written as 
\begin{equation}
	\bar\bb^\infty(r) =\delta(r)\; \bb^\infty ~,
\end{equation}
where $\delta(r)$ denotes the Dirac delta function and $\bb^\infty$ denotes a constant force vector that describes magnitude and orientation of the point load. The solutions for velocity and stress of a singularly forced Stokes flow \eqref{eq:flow_stokes_singu} are given by the components
\begin{equation}\label{eq:flow_infty}
\begin{aligned}
v_i^\infty(\bx)&= \frac{1}{8\pi\eta} G_{ij}(\br)\, b^\infty_j~, \\
\sigma_{ij}^\infty(\bx)&=  \frac{1}{8\pi} T_{ijk}(\br)\, b^\infty_k~,
\end{aligned}
\end{equation}
where Latin indices $i,j$ and $k$ (from one to three) indicate components in Cartesian coordinates following Einstein's summation convention.
The Green's functions for velocity and stress are given by
\begin{equation}\label{eq:flow_green}
\begin{aligned}
G_{ij}(\br)&=  \frac{\delta_{ij}+\bar r_i \bar r_j}{r}~, \\[0.4em]
T_{ijk}(\br)&=  -6 \,\frac{\bar r_i \bar r_j \bar r_k}{r^2}~,
\end{aligned}
\end{equation}
where $\bar r_i$ is the component of $\bar{\br}:={\br}/{r}$ and $\delta_{ij}$ denotes Kronecker's delta. From \eqref{eq:flow_green} follows that the indices of both Green's functions are arbitrarily interchangable.\footnote{The Green's function tensors are symmetric such that $G_{ij}(\br) \!=\! G_{ji}(\br)$ and $T_{ijk}(\br)\!=\!  T_{ikj}(\br)\!=\!  T_{jik}(\br)\!=\!T_{jki}(\br)\!=\!T_{kij}(\br)\!=\!T_{kji}(\br)$.} Further, the Green's functions are  symmetric with respect to $\br$ such that $G_{ij}(\br)\!=\!G_{ij}(-\br)$ and $T_{ijk}(\br)\!=\!-T_{ijk}(-\br)$.
\\\\
Integrating the Green's functions over the fluid domain $\mcalF$ and applying the divergence theorem yield the important integral identities for the velocity Green's function
\begin{equation}\label{eq:identity_SL}
 \int_{\partial \mcalF}  G_{ij}(\br)\, n_j(\bx) \, \mathrm da_x =	0 \hh \by \in \mcalS~,
\end{equation}
and the stress Green's function
\begin{equation}\label{eq:identity_DL}
 -\frac{1}{4 \pi}\int_{\partial \mcalF}^\mathrm{pv}  T_{ijk}(\br)\, n_k(\bx) \, \mathrm da_x = \,\delta_{ij} \hh \by \in \mcalS~,
\end{equation}
where $\int^\mathrm{pv} \mrd \square$ denotes the Cauchy principal value integral.

\subsubsection{Boundary integral equation}\label{sec:flow_BIE}
The governing equations for incompressible Stokes flow, \eqref{eq:flow_consti}, \eqref{eq:flow_continuity} and \eqref{eq:flow_stokes}, are transformed into a boundary representation that is referred to as the boundary integral equation (BIE). The BIE allows to determine the fluid velocity at any point $\by\in\mcalF$ by solving only surface integrals defined on $\mcalS$.
To derive the BIE, a general incompressible Stokes flow ($\bv, \bsig$) is related to a singularly forced incompressible Stokes flow ($\bv^\infty,\bsig^\infty$) by use of the introduced Green's functions \eqref{eq:flow_infty} and Lorentz reciprocal theorem \citep{Lorentz1896}. Integrating the resulting identity over the fluid domain and applying the divergence theorem yields the BIE
\begin{equation}\label{eq:flow_BIE_domain}
v_i(\by)= - \frac{1}{8\pi\,\eta} \int_{\mcalS} G_{ij}(\bx-\by)\, {t_j}(\bx) \,\mrd a_x  + \frac{1}{8\pi}\int_{\mcalS} v_j (\bx) \, T_{ijk}(\bx-\by)\, n_k (\bx)  \, \mrd a_x ~,
\end{equation}
where $\by\in\mcalF$ denotes a source point within the fluid domain. Considering the limiting behavior of the singular integrals yields the BIE for a point on the surface ($\by \in \mcalS$)
\begin{equation}\label{eq:flow_BIE_surf}
v_i(\by)= - \frac{1}{2\varphi\,\eta} \int_{\mcalS} G_{ij}(\bx-\by)\, {t_j}(\bx) \,\mrd a_x  + \frac{1}{2\varphi}\int_{\mcalS}^\mathrm{pv} v_j(\bx) \, T_{ijk}(\bx-\by)\, n_k(\bx)  \, \mrd a_x ~,
\end{equation}
where $\varphi$ denotes the solid angle of the enclosed domain.\footnote{For a smooth surface (locally at least $C^1$-continuous) the solid angle becomes $\varphi= 2\pi$.}
Note that \eqref{eq:flow_BIE_domain} and \eqref{eq:flow_BIE_surf} consider exterior flow problems where an object or particle enclosed by surface $\mcalS$ is surrounded by the infinite fluid domain $\mcalF$. Similar BIEs for interior and two-sided problems 
and their derivations can be found, see for instance \cite{Pozrikidis92_book} and \cite{Harmel21_thesis}.
\\\\BIEs for unsteady Stokes flow and corresponding Green's functions can be found for instance in \cite{Power93} and \cite{Jiang12}. Those BIEs involve additional complexities like trivariate domain integrals and Laplace or Fourier transforms. However, the quadrature schemes proposed in Sec.~\ref{sec:hybrid} are also applicable to BIEs for unsteady Stokes, since those also contain the boundary integrals from \eqref{eq:flow_BIE_surf}. The quadrature schemes can further be used to approximate singular boundary integrals arising in other unsteady BE formulations including those for potential problems \citep{Falletta15} and elastodynamics \citep{Aimi19,Aimi22}.

\subsubsection{Boundary discretization and collocation}\label{sec:flow_BE}
The surface geometry $\mcalS$ and the continuous BIE \eqref{eq:flow_BIE_surf} are discretized into $\nel$ elements and $\nno$ nodes (control points in isogeometric analysis). The discretized BIE is then collocated at $\nno$ source points $\by_A \in \mcalS$, with $A = 1,\ldots,\nno$, to obtain a square boundary element~(BE) system. The source points $\by_A$ are thus referred to as collocation points in the following. Note that integration over the whole surface has to be performed once for each collocation point. Explicit expressions for the discretized BIE and the collocated BE system can be found in Appendix~\ref{sec:BE}.
%
%
\\\\This paper considers isogeometric discretizations, where the control points are not necessarily located on the surface, which makes them unsuitable for collocation. Instead, the location of the collocation points are determined by the Greville abscissae (see e.g.~\cite{Greville64}, \cite{Johnson05} and \cite{Aurrichio10}).

\subsubsection{Singular nature of the Green's functions}\label{sec:flow_Greens_singu}
\begin{figure}[h]
\unitlength\linewidth
\begin{picture}(1,0.46)
\put(0.11,0.28){\includegraphics[trim = 100 150 100 50 ,clip, width=.24\linewidth]{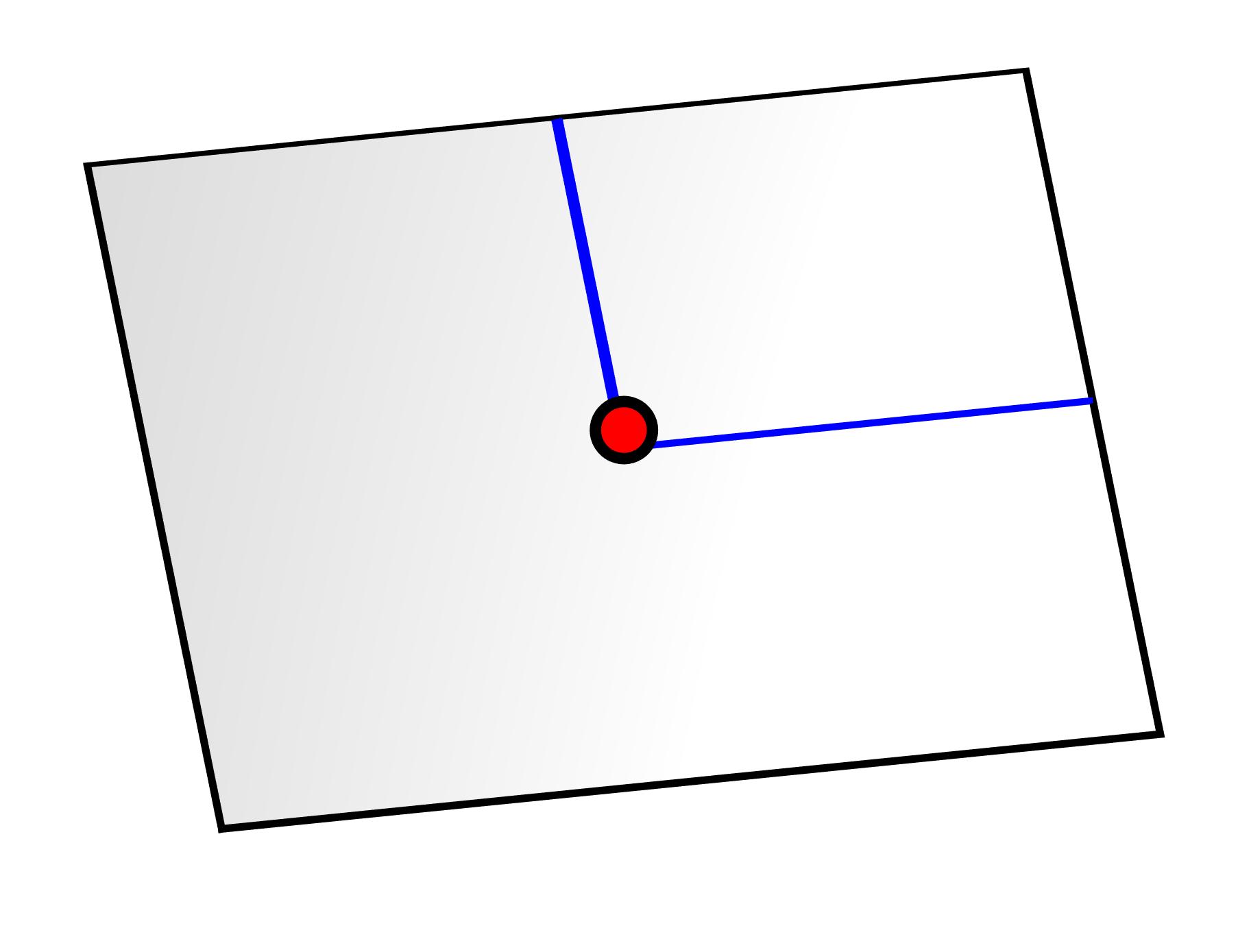}}
\put(0.44,0){\includegraphics[trim = 30 0 60 30 ,clip, width=.55\linewidth]{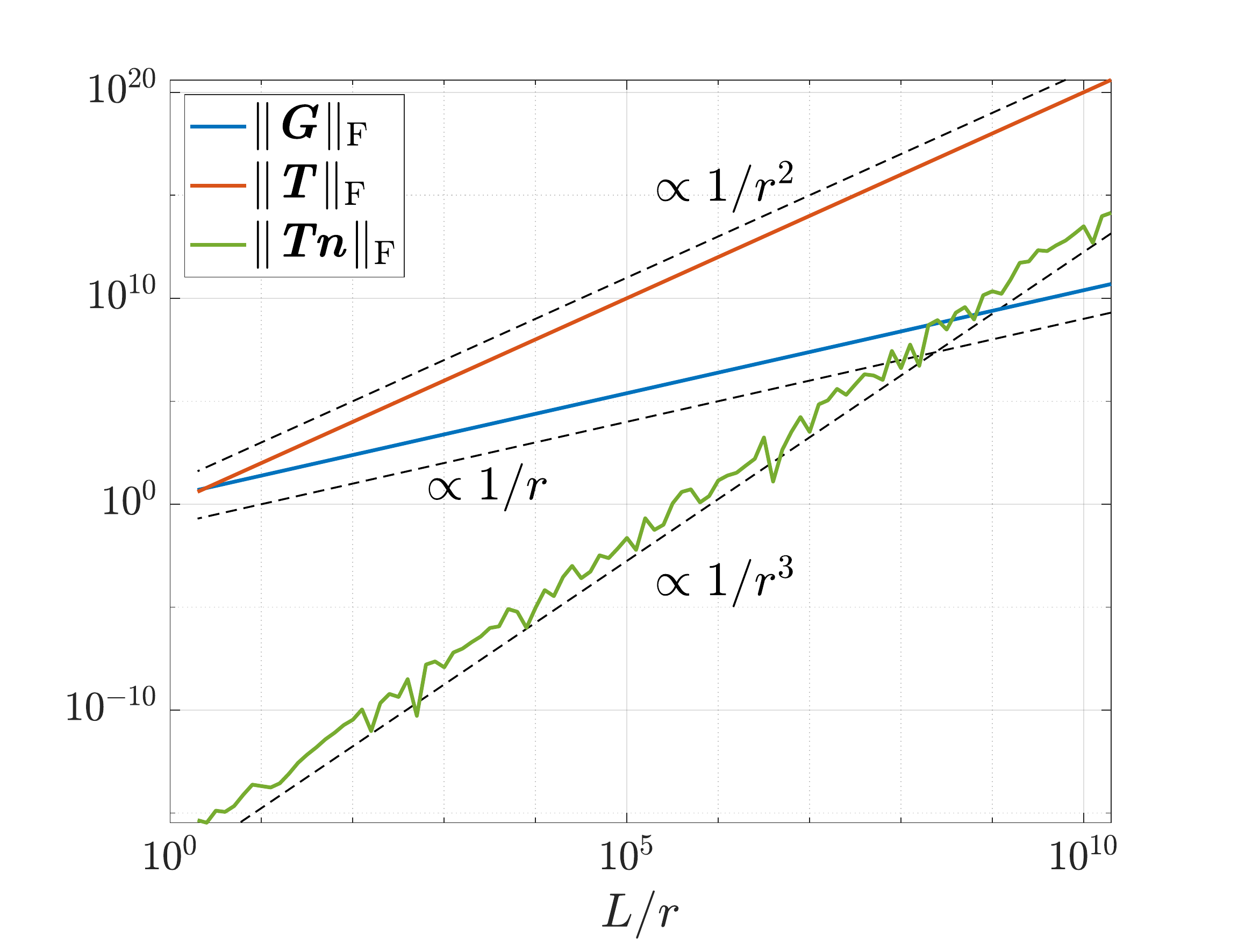}}
\put(0.02,0){\includegraphics[trim = 0 0 0 0,clip, width=.4\linewidth]{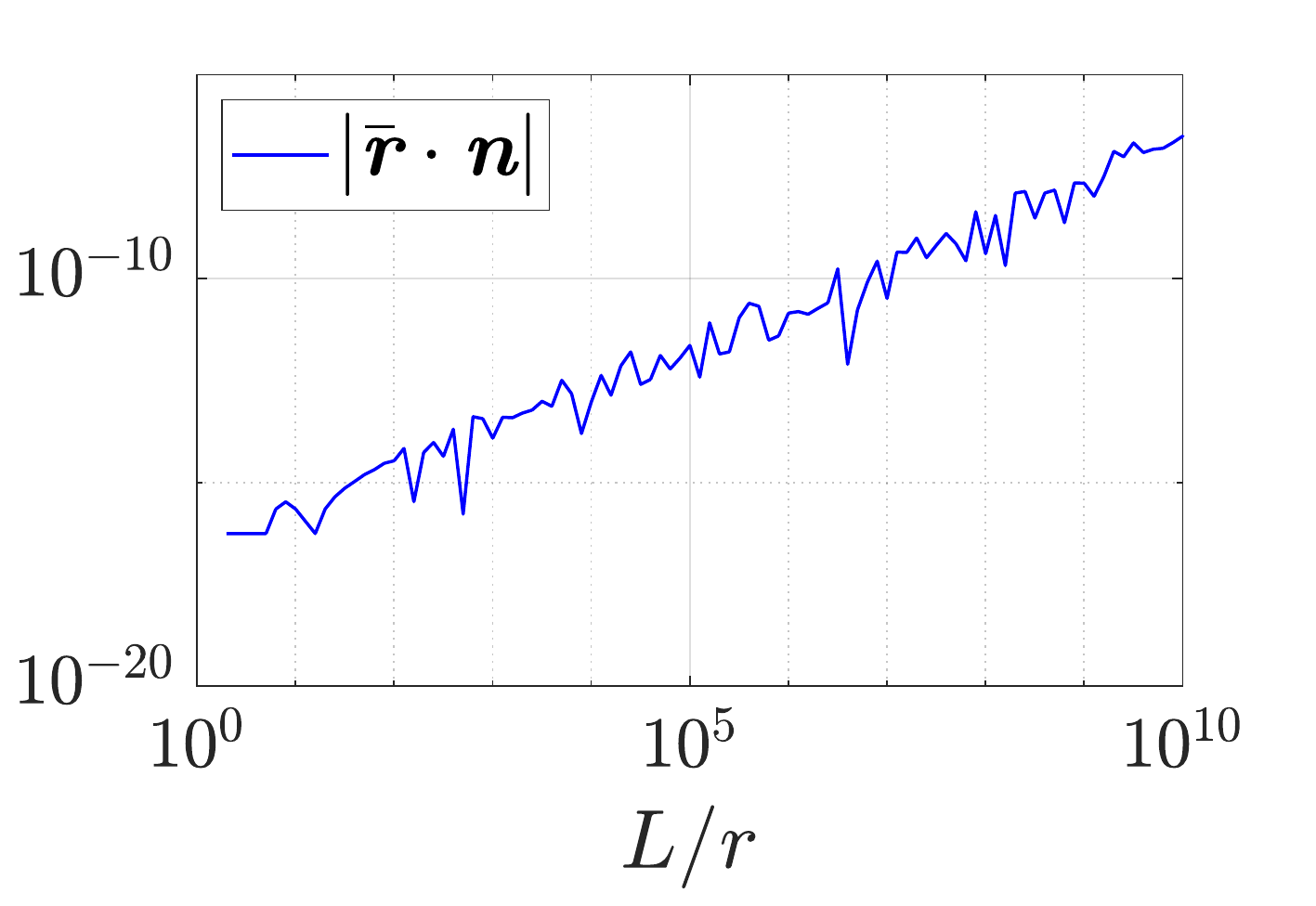}}
\put(0.11,0.29){a.}\put(0.03,0.01){b.} \put(0.45,0.01){c.} \put(0.2,0.36){\large $\by$} 
\end{picture}
\caption{\textit{Singular nature of Green's functions}: a. Tilted planar surface (size $L\times L$) with singularity at $\by$; b.~Limiting behavior of $|\bar \br \cdot \bn|$ for $r:=\|\br\|$ approaching 0; c.~Limiting behavior of the Green's functions.}\label{fig:GreensFcn_2}
\end{figure}
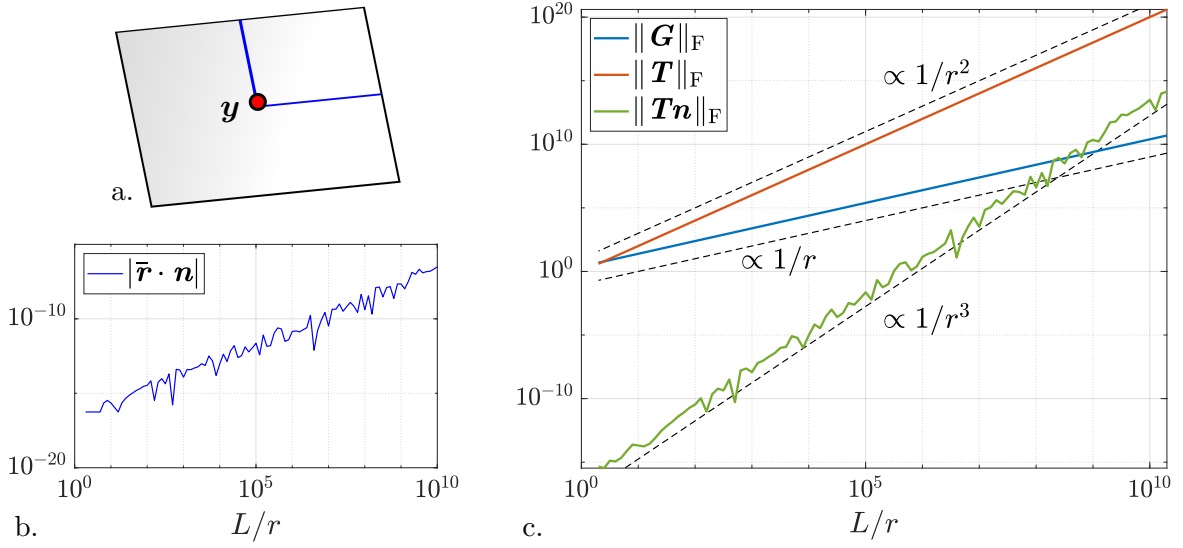

The Green's functions \eqref{eq:flow_green} become singular for $r$ approaching to zero. Source point $\by$ is therefore a singularity of both Green's functions. According to their definitions, the velocity Green's function $G_{ij}(\br)$ is proportional to $1/r$, while the stress Green's function $T_{ijk}(\br)$  is proportional to $1/r^2$. This kind of function behavior is referred to as \textit{weakly singular} and \textit{strongly singular}, respectively.
Since these singularities are an essential component of boundary element (BE) analysis, their asymptotic behavior for $r\rightarrow 0$ is numerically investigated on spherical and on flat surfaces in the following. The Green's functions behavior is evaluated by a scalar invariant, the Frobenius norm, that is defined for the second order tensor $\bG(\br)=G_{ij}(\br)\, \be_i \otimes \be_j$~by
\begin{equation}\label{eq:green1_inv}
	\|\bG(\br)\|_\mrF= \sqrt{G_{ij}(\br)\,G_{ij}(\br)}~,
\end{equation}
and for the third order tensor $\bT(\br)=T_{ijk}(\br)\, \be_i \otimes \be_j \otimes \be_k$ by
\begin{equation}\label{eq:green2_inv}
	\|\bT(\br)\|_\mrF= \sqrt{T_{ijk}(\br)\,T_{ijk}(\br)}~,
\end{equation}
where \eqref{eq:green1_inv} and \eqref{eq:green2_inv} again follow Einstein's summation convention. The second-order tensor $\bT\bn$ is part of an integral kernel of the BIE in Sec.~\ref{sec:flow_BIE} and therefore of particular interest for BE analysis. Its asymptotic behavior is additionally investigated here by determining $\|\bT\bn\|_\mrF$ according to \eqref{eq:green1_inv}.
%
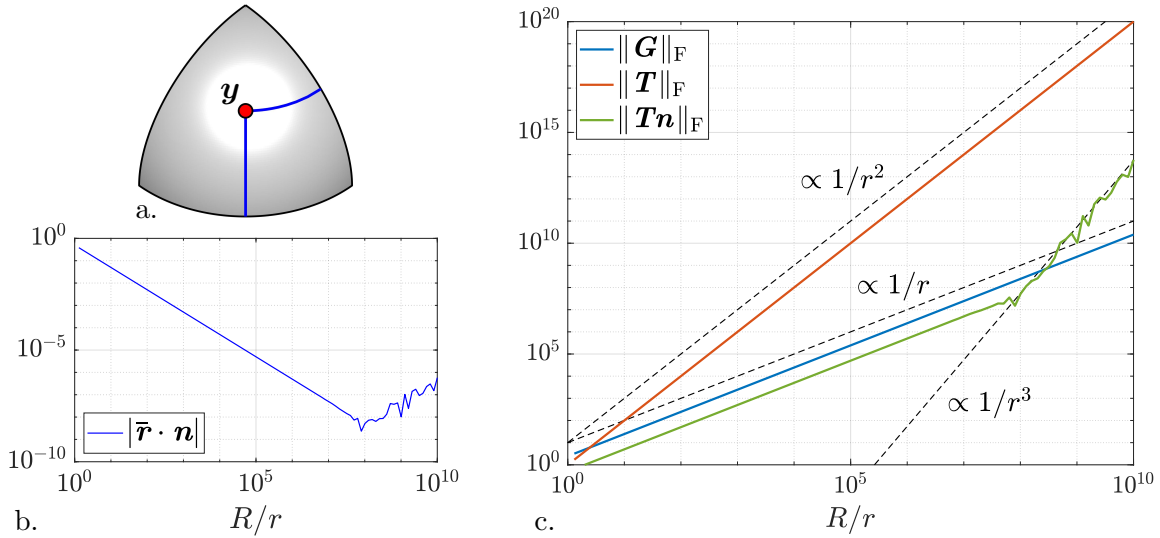
\begin{figure}[h]
\unitlength\linewidth
\begin{picture}(1,0.46)
	\put(0.02,0){\includegraphics[trim = 0 0 0 0,clip, width=.4\linewidth]{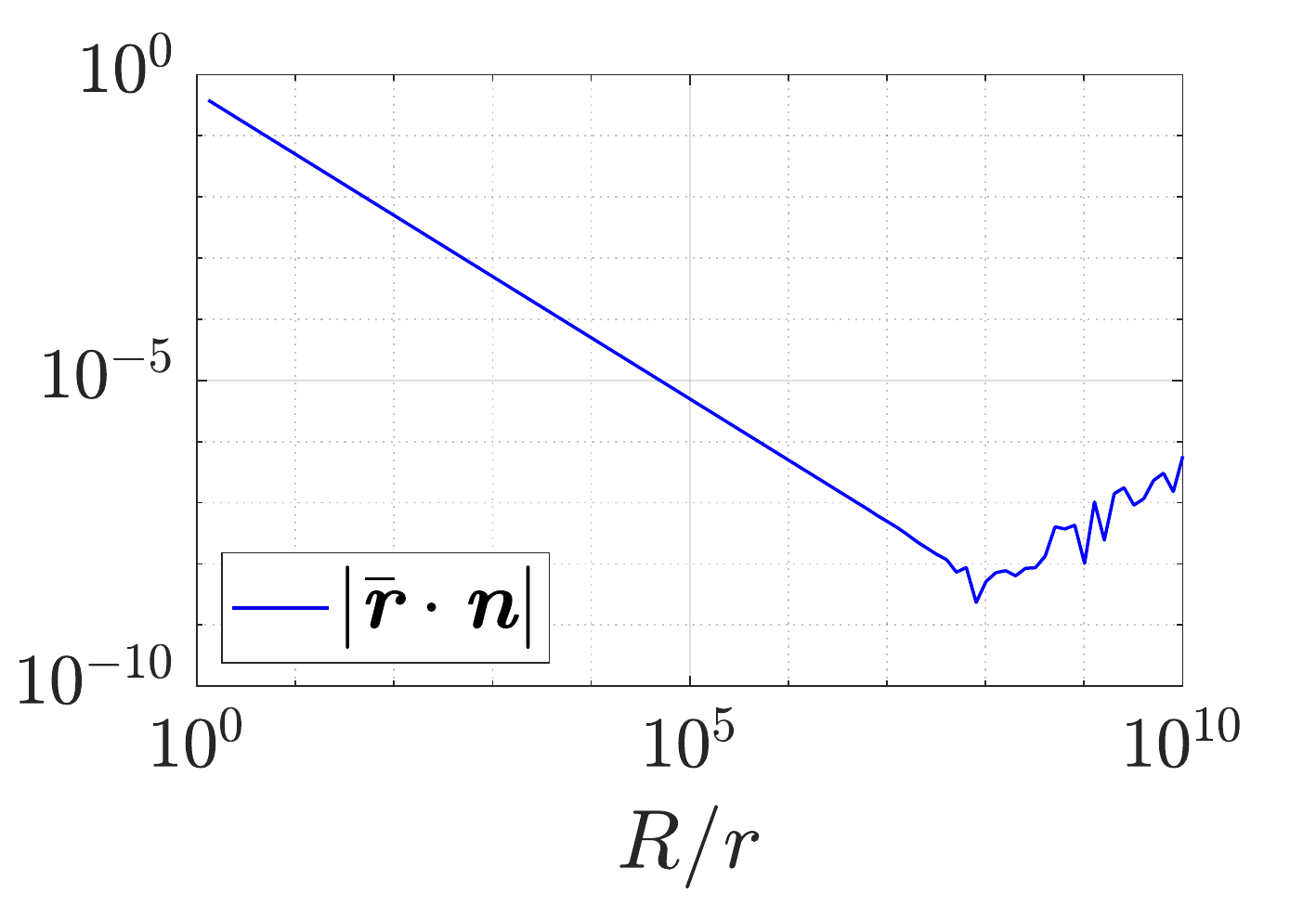}}
	\put(0.13,0.27){\includegraphics[trim = 250 40 240 40 ,clip, width=.18\linewidth]{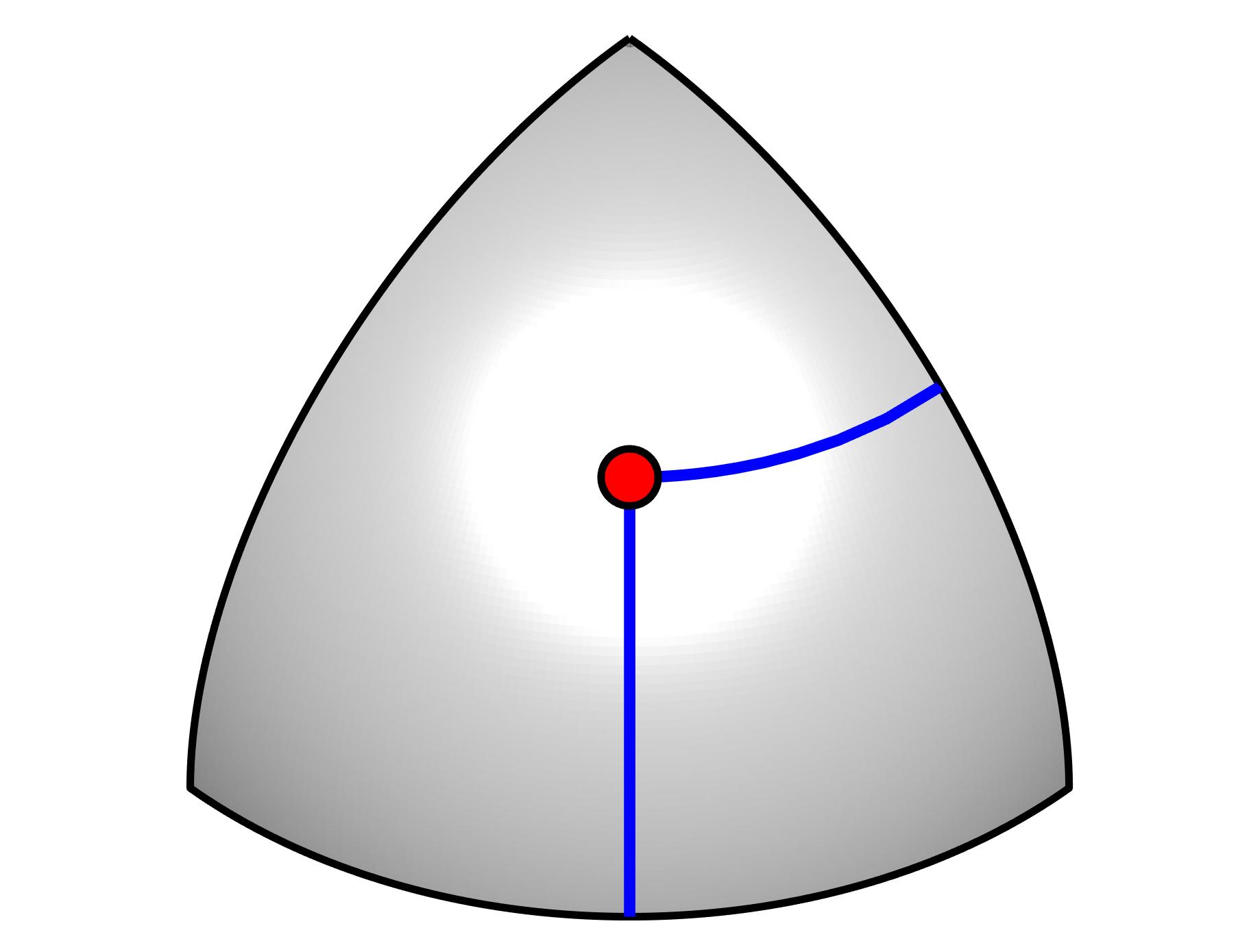}}
	\put(0.44,0){\includegraphics[trim = 40 0 30 20 ,clip, width=.55\linewidth]{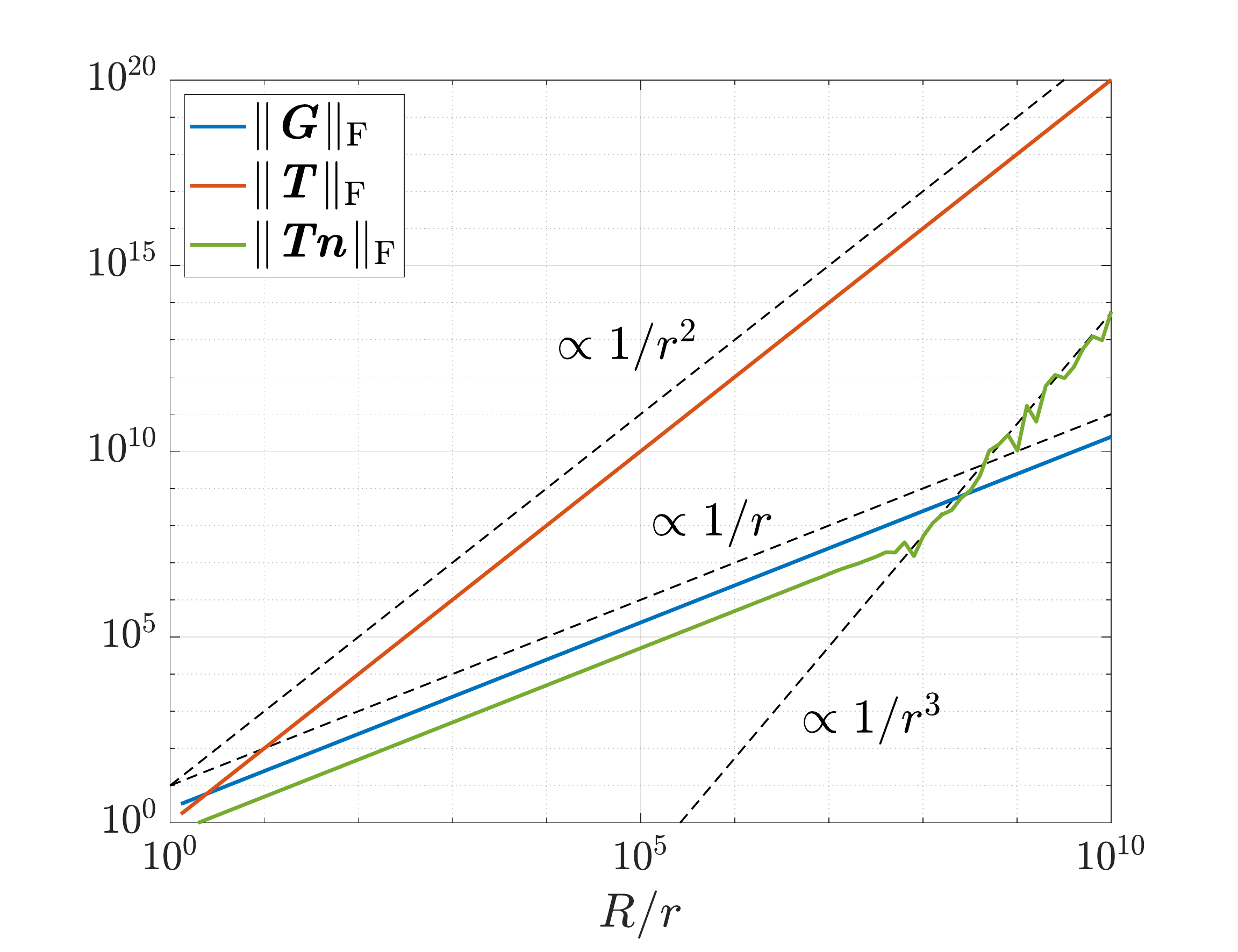}}
	\put(0.13,0.27){a.}\put(0.03,0.01){b.} \put(0.46,0.01){c.} \put(0.2,0.37){$\large \by$} 	
\end{picture}
\caption{\textit{Singular nature of Green's functions}: a.~Spherical surface (radius $R$) with singularity at $\by$; b.~Limiting behavior of $|\bar\br \cdot \bn|$ for $r:=\|\br\|$ approaching 0; c.~Limiting behavior of the Green's functions.}\label{fig:GreensFcn_1}
\end{figure}
\\\\A spherical surface with radius $R$ is considered first, see Fig.~\ref{fig:GreensFcn_1}a. Given source point $\by$, the Green's functions behavior is investigated as field variable $\bx$ approaches $\by$ along the blue curve. Fig.~\ref{fig:GreensFcn_1}c shows that $\|\bG\|_\mrF$ is proportional to $1/r$ (\textit{weakly singular}), while $\|\bT\|_\mrF$ is proportional to $1/r^2$ (\textit{strongly singular}). The $\|\bG\|_\mrF$ and $\|\bT\|_\mrF$ show the asymptotic behavior that is expected from the Green's functions definition \eqref{eq:flow_green}. Fig.~\ref{fig:GreensFcn_1}b further shows that $\|\bT\bn\|_\mrF$ is at first only \textit{weakly singular}, although $\|\bT\|_\mrF$ is \textit{strongly singular}. Let us take a closer look on the expression to substantiate this numerical result with analytical findings: From \eqref{eq:flow_green} follows that $\|\bT\bn\,\|_\mrF \propto  \bar \br\cdot \bn /r^2$. For a smooth surface with continuously varying normal vector, the unit vectors $\bar \br$ and $\bn$ become more and more orthogonal as $\bx$ approaches $\by$, i.e.~$r\rightarrow 0$. The dot product $\bar \br\cdot \bn$ thus changes linearly with respect to $r$ (see Fig.~\ref{fig:GreensFcn_1}b) such that the behavior of $\|\bT\bn\|_\mrF$ for $r\rightarrow 0$ is mildened to be only \textit{weakly singular}.
However, for $r<10^{-8}\,R$ the behavior of $\|\bT\bn\|_\mrF$ vs.~$1/r$ changes suddenly: From a perfectly linear relation ($\|\bT\bn\,\|_\mrF \propto 1/r$) to a distorted cubic relation ($\|\bT\bn\,\|_\mrF \propto 1/r^3$). This behavior results from the inaccurate numerical computation of the dot product $\bar \br\cdot \bn$ for very small $r$ that is shown in Fig.~\ref{fig:GreensFcn_1}b. Within the numerical quadrature, $\bT\bn$ is typically not required to be evaluated for such small values of $r$. The discussed change in the convergence behavior of $\bT\bn$ has therefore no effect on the approximations of the singular integrals in Sec.~\ref{sec:flow_BIE}.
\\\\Second, a source point on a tilted planar surface of side length $L$ is considered (see Fig.~\ref{fig:GreensFcn_2}a). It can be seen in Fig.~\ref{fig:GreensFcn_2}c that the Green's functions show the same asymptotic behavior here as for spherical surfaces: $\|\bG\|_\mrF$ is \textit{weakly singular} and $\|\bT\|_\mrF$ is \textit{strongly singular} for $r\rightarrow 0$.
However, the numerical evaluation of $\|\bT\,\bn\|_\mrF$ yields an overall cubic asymptotic behavior deviating from theoretical expectations: For a source point located on a plane surface, the vectors $\bn$ and $\bar\br$ are orthogonal. Analytically their dot product thus yields $\bn \cdot \bar\br=0$ and the expression $\bT\bn$ consequently vanishes. The deviating behavior of $\|\bT\,\bn\|_\mrF$ is caused by the diverging numerical result for $r\rightarrow 0$ on plane surfaces that is shown in Fig.~\ref{fig:GreensFcn_2}b. On flat surfaces it is therefore advisable to use the known value $\bT\bn=\bzero$ in the BE computation instead of evaluating it numerically. Nevertheless, it should be noted that moderate values of $r$ result in very small $\|\bT\,\bn\|_\mrF$ (smaller than $10^{-10}$ for $r>10^{-2}L$).
%

\section{Boundary quadrature of singularities}\label{sec:quad}
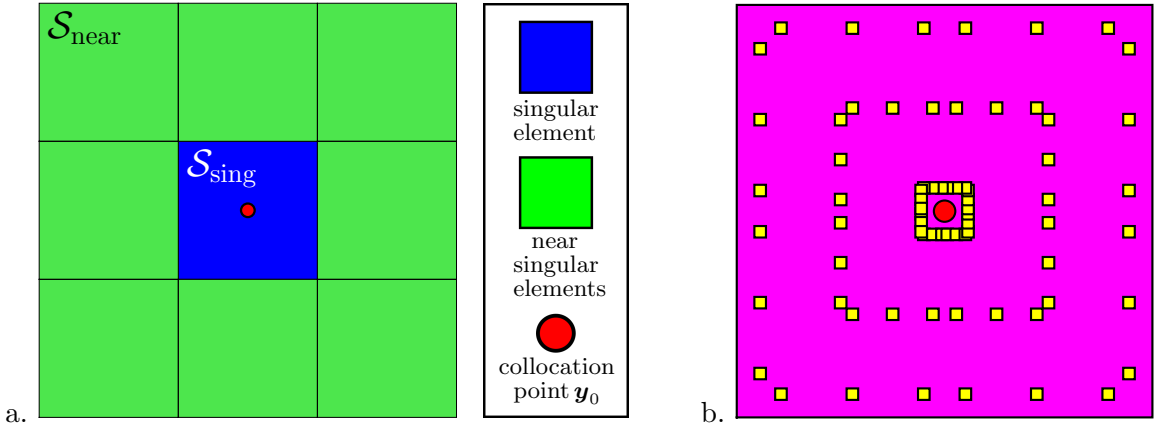
\begin{figure}[h]
\unitlength\linewidth
\begin{picture}(1,.35)
	%
    \put(0.05,0){\includegraphics[trim = 130 50 110 35, clip, width=.35\linewidth ]{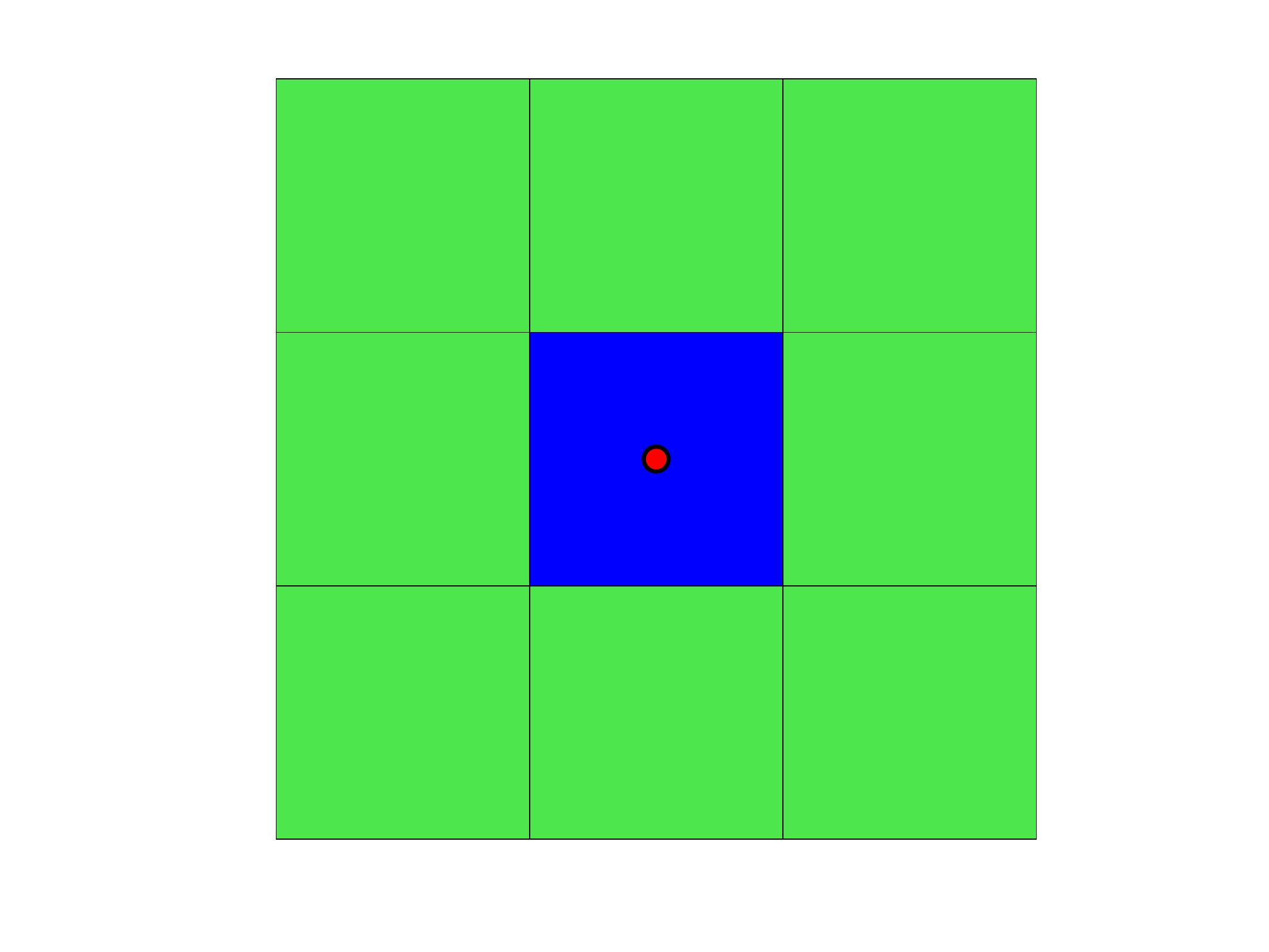}}
    \put(0.42,0.004){\includegraphics[trim = 0 0 0 0, clip, height=.345\linewidth ]{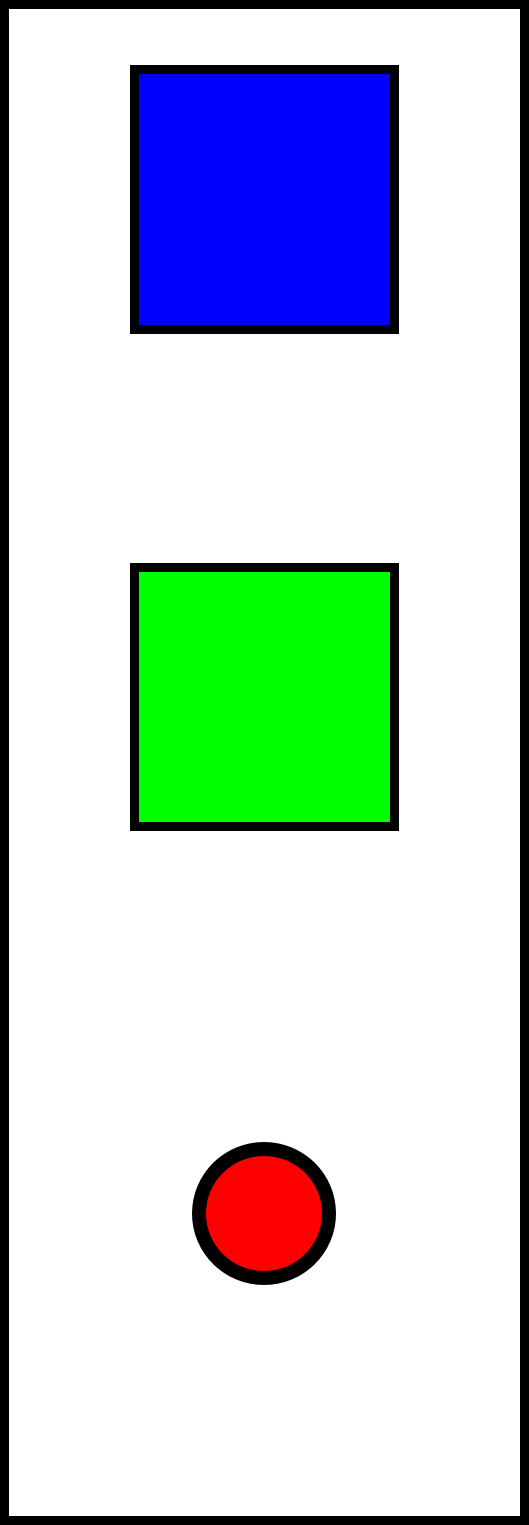}}
	\put(0.445,0.255){\footnotesize singular} \put(0.445,0.235){\footnotesize element}
    \put(0.46,0.145){\footnotesize near} \put(0.445,0.125){\footnotesize singular}  \put(0.445,0.105){\footnotesize elements}
    \put(0.433,0.04){\footnotesize collocation} \put(0.445,0.02){\footnotesize point\,$\by_0$}
    \put(0.625,0){\includegraphics[trim = 65 65 50 50, clip, width=.35\linewidth]{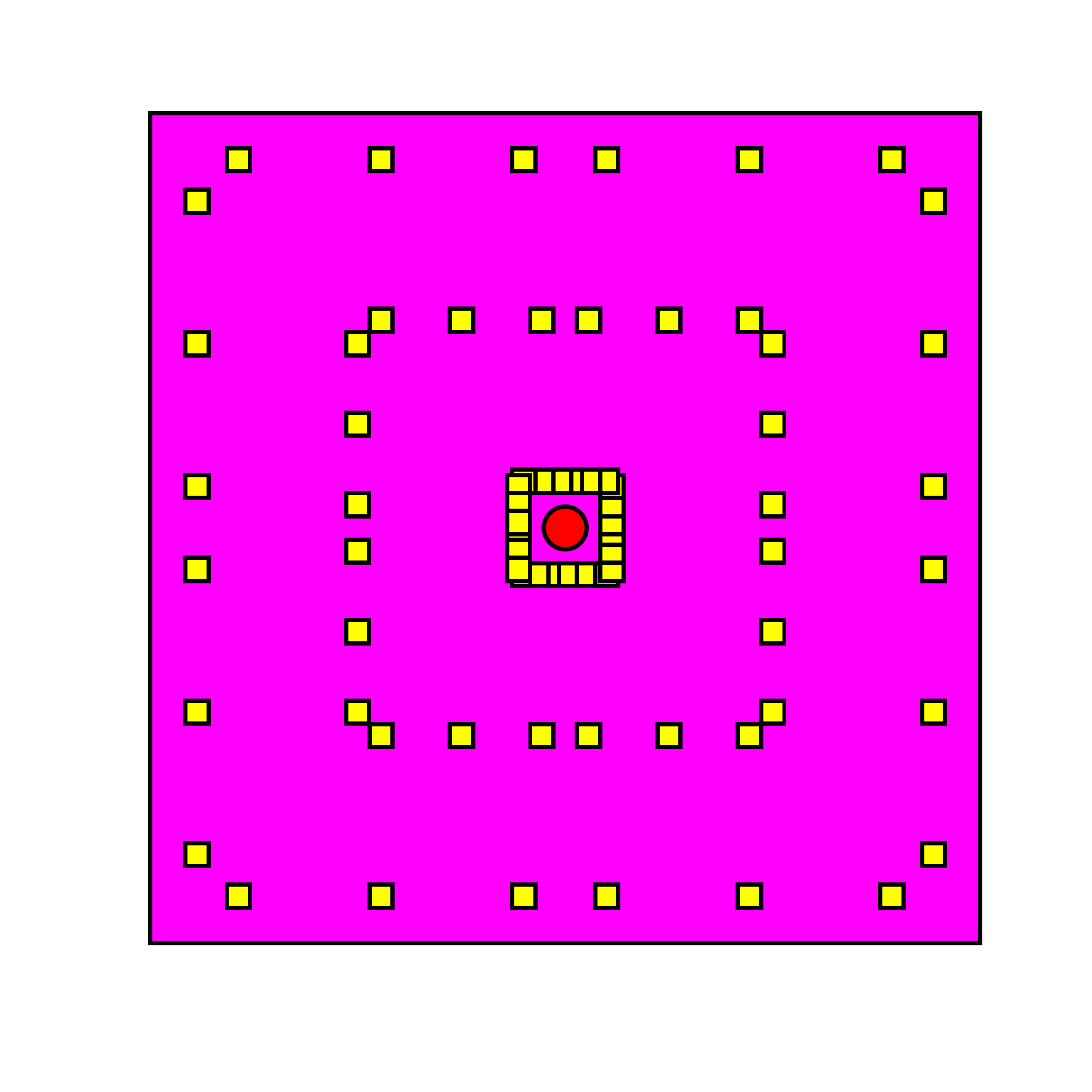}}
	\put(0.025,0){a.}\put(0.6,0){b.}
	\put(0.06,0.32){\Large $\mcalS_\mathrm{near}$}\put(0.175,0.205){\Large $\twhite{\mcalS_\mathrm{sing}}$}
	\end{picture}
\caption{\textit{Boundary quadrature of singularities}: a.~Singular element (blue, $\mcalS_\mathrm{sing}$) with collocation point $\by_0$ and the surrounding near singular elements~(green, $\mcalS_\mathrm{near}$). b.~Duffy quadrature points on the singular element for $\tilde n_\mathrm{gp}=3$ (see Table~\ref{tab:nqp} and \eqref{app_eq:quad_Duffy_nqp}) and collocation point $\by_0$.\protect\footnotemark~Duffy quadrature points for different collocation point positions are shown in Fig.~\ref{fig:quad_Duff2}.}\label{fig:quad_intro}
\end{figure}

This section investigates singular boundary quadrature on the plane sheet from Fig.~\ref{fig:quad_intro}a that consists of $3\times 3$ plane biquadratic B-spline elements with side length $L^e$. The depicted collocation point $\by_0$ is located on the center element which is therefore referred to as singular element (blue, surface $\mcalS_\mathrm{sing}$), whereas the adjacent elements are referred to as near singular elements (green, surface $\mcalS_\mathrm{near}$).
The weakly singular integral
\begin{equation}\label{eq:quad_int}
\mcalI_\square :=\int_{\mcalS_\square} \frac{1}{\|\bx-\by_0\|} \,\mrd a
\end{equation}
that is representative of the behavior of the two integrals in the BIE \eqref{eq:flow_BIE_surf} is approximated on $\mcalS_\mathrm{sing}$ and on $\mcalS_\mathrm{near}$ using various quadrature rules: Sec.~\ref{sec:quad_singular} studies Duffy quadrature and classical and modified Gauss-Legendre quadrature on the singular element and Sec.~\ref{sec:quad_nearly} investigates classical Gauss-Legendre quadrature and a new Gauss-Legendre quadrature rule with adjusted weights on the near singular elements. The quadrature accuracy is evaluated by means of the relative quadrature error
\begin{equation}\label{eq:quad_err}
e^\square_\mathrm{rel}:=  \frac{\|\mcalI_\square^h - \mcalI_\square\|}{\mcalI_\square}~,
\end{equation}
where $\mcalI^h_\square$ denotes the numerical approximation of integral \eqref{eq:quad_int} on surface $\mcalS_\square$.\footnote{Analytical solutions: 
$\mcalI_\mathrm{sing} =\log\left(12\sqrt{2}+17\right) L^e$ and
$\mcalI_\mathrm{near} = \left( 6\,\log(\sqrt{2} + 1) + \log\left(2\sqrt{2} + 3\right) \right) L^e$}

\footnotetext{The color scheme from Fig.~\ref{fig:quad_intro} is used throughout the paper to visualize the quadrature schemes: Blue, green and white are used to distinguish between singular, near singular and regular elements (see Fig.~\ref{fig:quad_sheet_discr}a for a 4x4 mesh). Various other colors, including magenta for Duffy quadrature, are used to illustrate the quadrature rules applied to these elements (see Fig.~\ref{fig:rules_cmp} for the same 4x4 mesh).}
\subsection{Boundary quadrature on singular elements}\label{sec:quad_singular}
\begin{figure}[h]
\unitlength\linewidth
\begin{picture}(1,0.45)
\put(0,0){\includegraphics[trim = 40 0 40 20, clip, width=.5\linewidth]{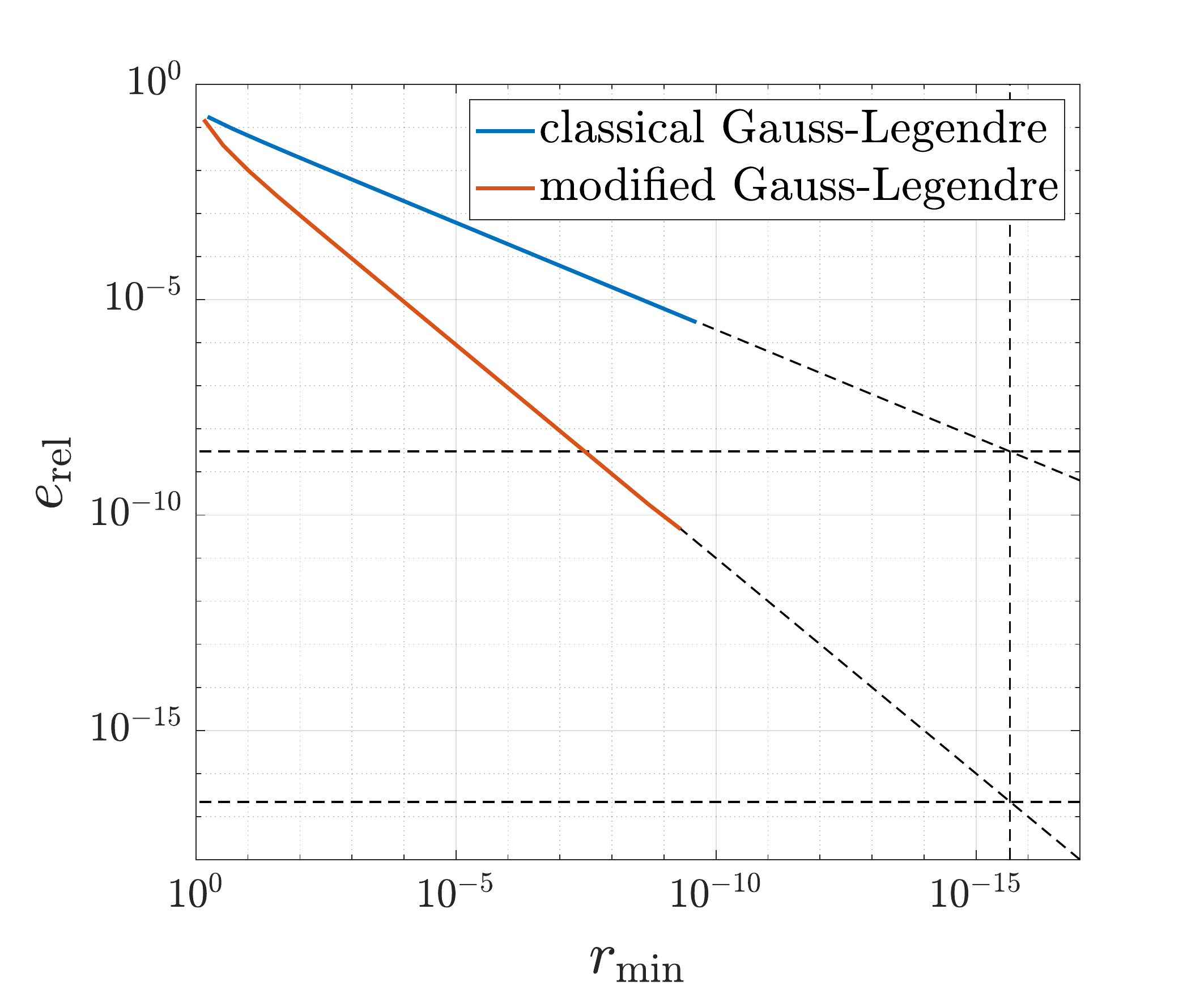}}
\put(0.0,0){a.} 
\vwput{0.005}{0.25}{\large $e_\mathrm{rel}^\mathrm{sing}$}
\wput{0.23}{0.015}{\large $r_\mathrm{min}/L^\mre$}
\put(0.21,0.21){$\propto 1/n_\mathrm{qp}$} \put(0.31,0.325){$\propto 1/\sqrt{n_\mathrm{qp}}$}
\put(.425,0.14){\rotatebox{90}{\large $\varepsilon_\mathrm{machine}$}}
\put(.08,0.27){$e_\mathrm{class}$}
\put(.08,0.102){$e_\mathrm{mod}$}
\put(0.5,0){\includegraphics[trim = 40 0 40 20, clip, width=.5\linewidth]{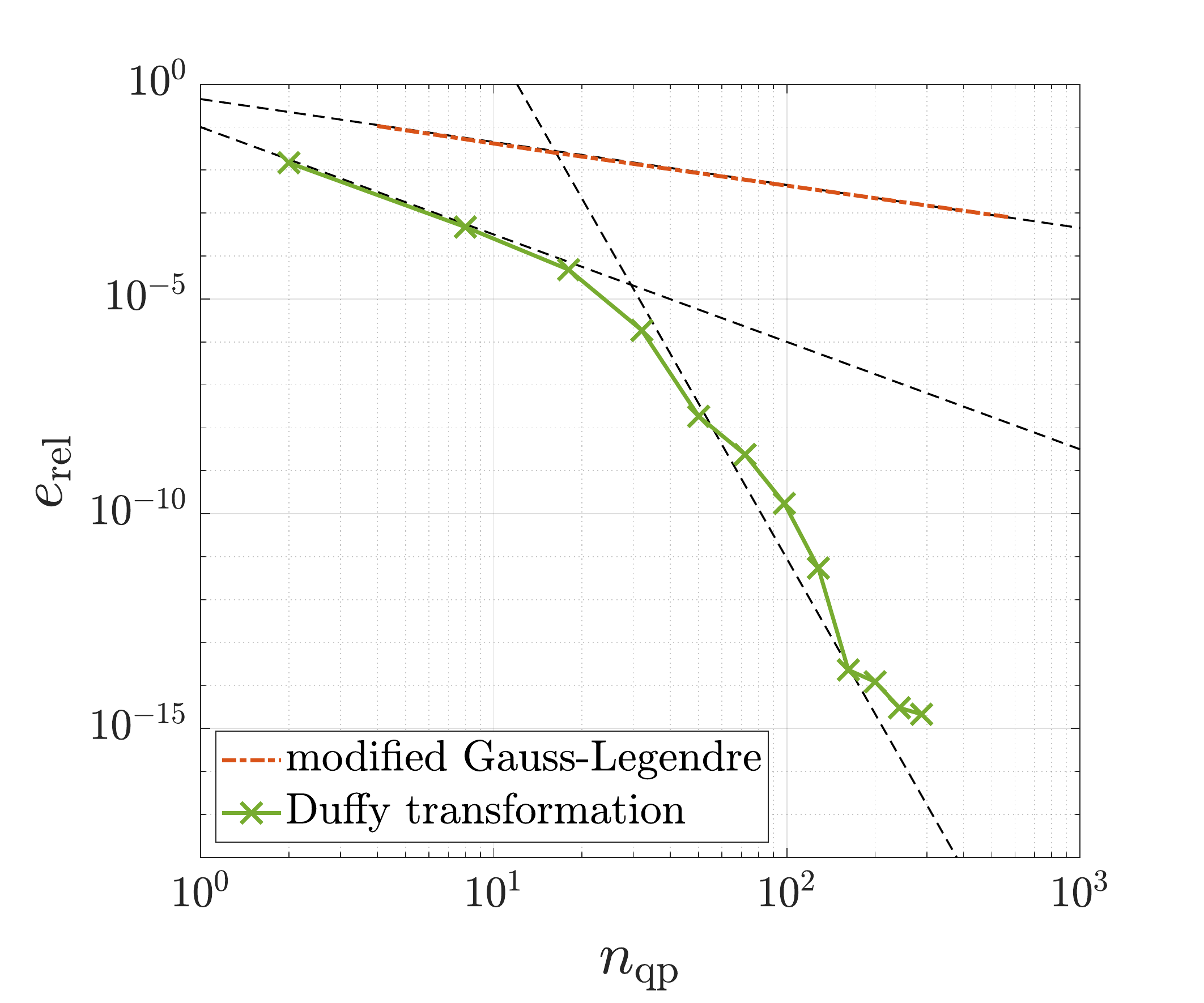}}
 \put(0.85,0.4){$\propto 1/n_\mathrm{qp}$}
\put(0.84,0.32){$\propto 1/n_\mathrm{qp}^{2.5}$} \put(0.74,0.2){$\propto 1/n_\mathrm{qp}^{12}$}
\put(0.51,0){b.} 
\end{picture}
\caption{\textit{Boundary quadrature on singular elements}: a.~{Relative quadrature error} $e_\mathrm{rel}^\mathrm{sing}$ for Gauss-Legendre quadrature vs.~the minimum distance between quadrature point and collocation point $r_\mathrm{min}$ \eqref{app_eq:quad_dist}, the vertical dashed line shows the smallest admissible distance $\varepsilon_\mathrm{machine} \approx {2.2\times 10^{-16}}$. The horizontal lines show that classical Gauss-Legendre quadrature has an inherent error bound $e_\mathrm{class}\approx3\times10^{-9}$, while modified Gauss-Legendre quadrature theoretically meets $e_\mathrm{mod}\approx \varepsilon_\mathrm{machine}$. b.~Relative quadrature error $e_\mathrm{rel}^\mathrm{sing}$ for modified Gauss-Legendre quadrature and Duffy quadrature vs.~the number of quadrature points $n_\mathrm{qp}$.}\label{fig:quad_res}
\end{figure}
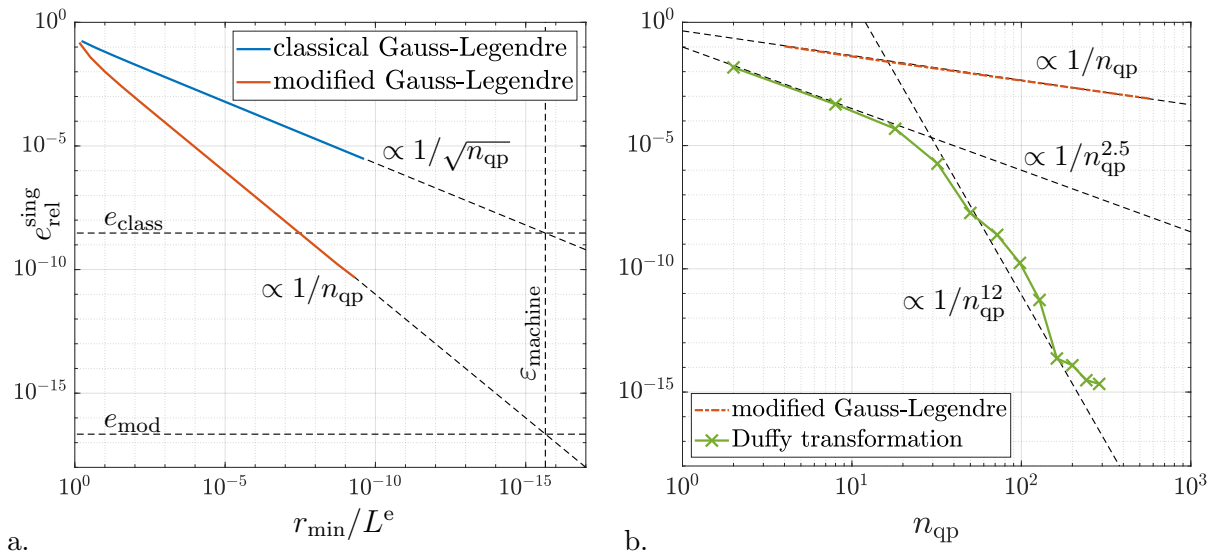

%
Classical Gauss-Legendre quadrature \citep{Gauss1815}, modified Gauss-Legendre quadrature (\cite{Heltai14} for biquadratic NURBS) and Duffy quadrature \citep{Fairweather79,Duffy82} are investigated in Appendix~\ref{app:quad} with respect to their suitability for approximation of weakly singular integrals. The resulting advantages and deficiencies of these three quadrature rules are summarized here.
\\\\Fig.~\ref{fig:quad_res} shows that the relative quadrature error \eqref{eq:quad_err} decreases with increasing number of quadrature points $n_\mathrm{qp}$ \eqref{eq:nqp} for all three quadrature rules. However, their convergence behavior is very different: The error for classical and modified Gauss-Legendre quadrature with respect to $n_\mathrm{qp}$ yield the rather low convergence rates of $1/2$ and $1$, respectively.\footnote{The minimum distance between quadrature points and collocation points \eqref{app_eq:quad_dist} is inversely proportional to the number of quadrature points, i.e.~$r_\mathrm{min}\propto n_\mathrm{qp}^{-1}$.}
Fig.~\ref{fig:quad_res}a shows that the relative quadrature error for modified Gauss-Legendre quadrature can be reduced to the range of machine precision by increasing the number of quadrature points to $n_\mathrm{qp}>10^{15}$, while classical Gauss-Legendre quadrature cannot reach machine precision as it has an inherent error bound of $e_\mathrm{rel}^\mathrm{sing} >e_\mathrm{class}\approx3\times10^{-9}$.
\\Duffy quadrature leads to a much lower error for small $n_\mathrm{qp}$ and, most importantly,  to a significantly improved convergence behavior: Fig.~\ref{fig:quad_res}b shows that the relative error converges to the analytical solution with $e_\mathrm{rel} \propto 1/{n_\mathrm{qp}}^{2.5}$ for small $n_\mathrm{qp}\leq 18$ and with $e_\mathrm{rel} \propto 1/{n_\mathrm{qp}}^{12}$ for higher $n_\mathrm{qp}$ such that is in the range of machine precision for $n_\mathrm{qp}>200$. Hence, Duffy quadrature is several orders of magnitude more accurate than modified Gauss-Legendre quadrature, not to mention classical Gauss-Legendre quadrature.
\\\\Further, collocation points and classical Gauss-Legendre quadrature points may coincide, in which case the integral kernel becomes infinite and useless. Classical Gauss-Legendre quadrature is therefore less robust than modified Gauss-Legendre quadrature and Duffy quadrature, where the coincidence of collocation points and quadrature points is generally prevented.\footnote{The quadrature points for modified Gauss-Legendre quadrature and for Duffy quadrature are obtained by combining quadrature rules for subelements such that the collocation points do not coincide with the quadrature points (see Fig.~\ref{fig:quad_GL}b and c for modified Gauss-Legendre quadrature and Fig.~\ref{fig:quad_Duff2} for Duffy quadrature).}
\\\\For classical and modified Gauss-Legendre quadrature, it is sufficient to evaluate the shape functions only once on the master element for all collocation points and elements.\footnote{Isogeometric shape function values are determined from universal Bernstein polynomials and element specific B\'ezier extraction operators \citep{Borden11}, see Appendix~\ref{sec:BE}.} Gauss-Legendre quadrature is therefore simple to implement. Duffy quadrature, on the other hand, needs to be modified depending on the location of the collocation points (see Fig.~\ref{fig:quad_intro}b) and it is thus more difficult to implement. The above findings are summarized in Table~\ref{tab:existing}.
\begin{table}[h]
\centering
\begin{tabular}{|c||c|c|c|} 
 \hline
  &  efficiency & robustness & implementation \\ 
 \hline
   classical Gauss-Legendre &  poor & low & simple \\ 
   modified Gauss-Legendre  &  low & high & simple  \\
   Duffy quadrature &  high & high & involved \\\hline
\end{tabular}
\caption{\textit{Boundary quadrature on singular elements}: Suitability of classical and modified Gauss-Legendre quadrature and Duffy quadrature for singular boundary quadrature (see Appendix \ref{app:quad}). } \label{tab:existing}
\end{table}
\\In conclusion, Duffy quadrature is highly accurate and thus most recommend for singular integral approximation.
Modified Gauss-Legendre quadrature is a robust alternative that is simple to implement but less accurate, while classical Gauss-Legendre quadrature is not recommend for the quadrature of singular kernels since it lacks robustness and efficiency.

\subsection{Boundary quadrature on near singular elements}\label{sec:quad_nearly}
This section deals with the approximation of singular integrals on so-called near singular elements. The collocation point is located near these elements, but not on within them (green elements in Fig.~\ref{fig:quad_intro}a). The accuracy of Gauss-Legendre quadrature on near singular elements is investigated in Sec.~\ref{sec:quad_nearly_gauss}, while a new and efficient quadrature rule for these elements is presented in Sec.~\ref{sec:quad_nearly_adjusted}. The quadrature rules discussed in Sec.~\ref{sec:quad_singular}, modified Gauss-Legendre quadrature and Duffy quadrature, are only beneficial on singular elements and thus not considered here.

\subsubsection{Classical Gauss-Legendre quadrature}\label{sec:quad_nearly_gauss}
\begin{figure}[h]
\unitlength\linewidth
\begin{picture}(1,.5)
	\put(0.02,0.015){\includegraphics[trim = 0 20 50 20, clip, width=.575\linewidth]{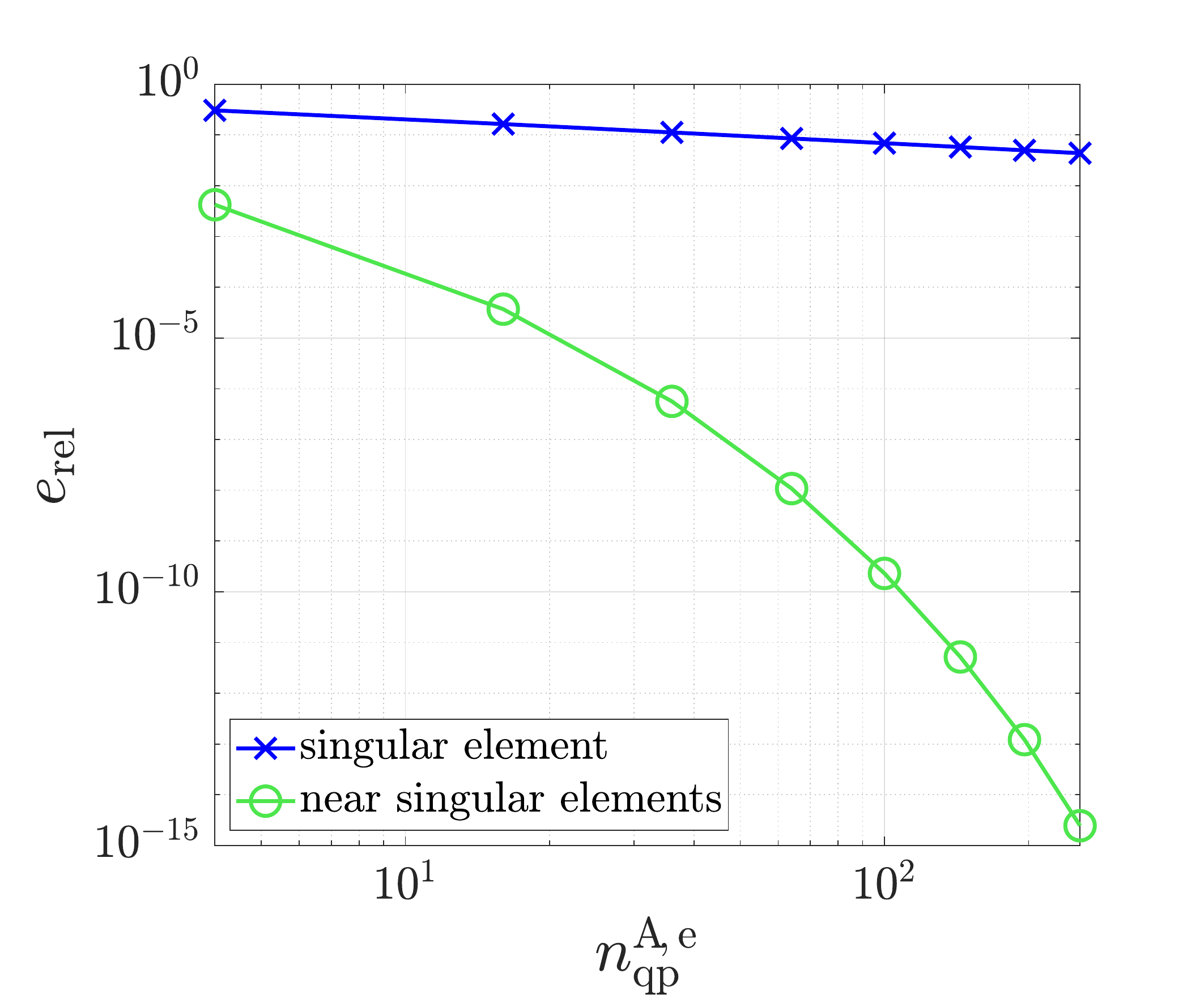}}
	\wput{0.26}{0.02}{$n_\mathrm{qp}^e:= \tilde n_\mathrm{qp} \times \tilde n_\mathrm{qp}$} 
	\put(0.654,0.23){\includegraphics[trim = 100 20 30 30, clip, width=.307\linewidth]{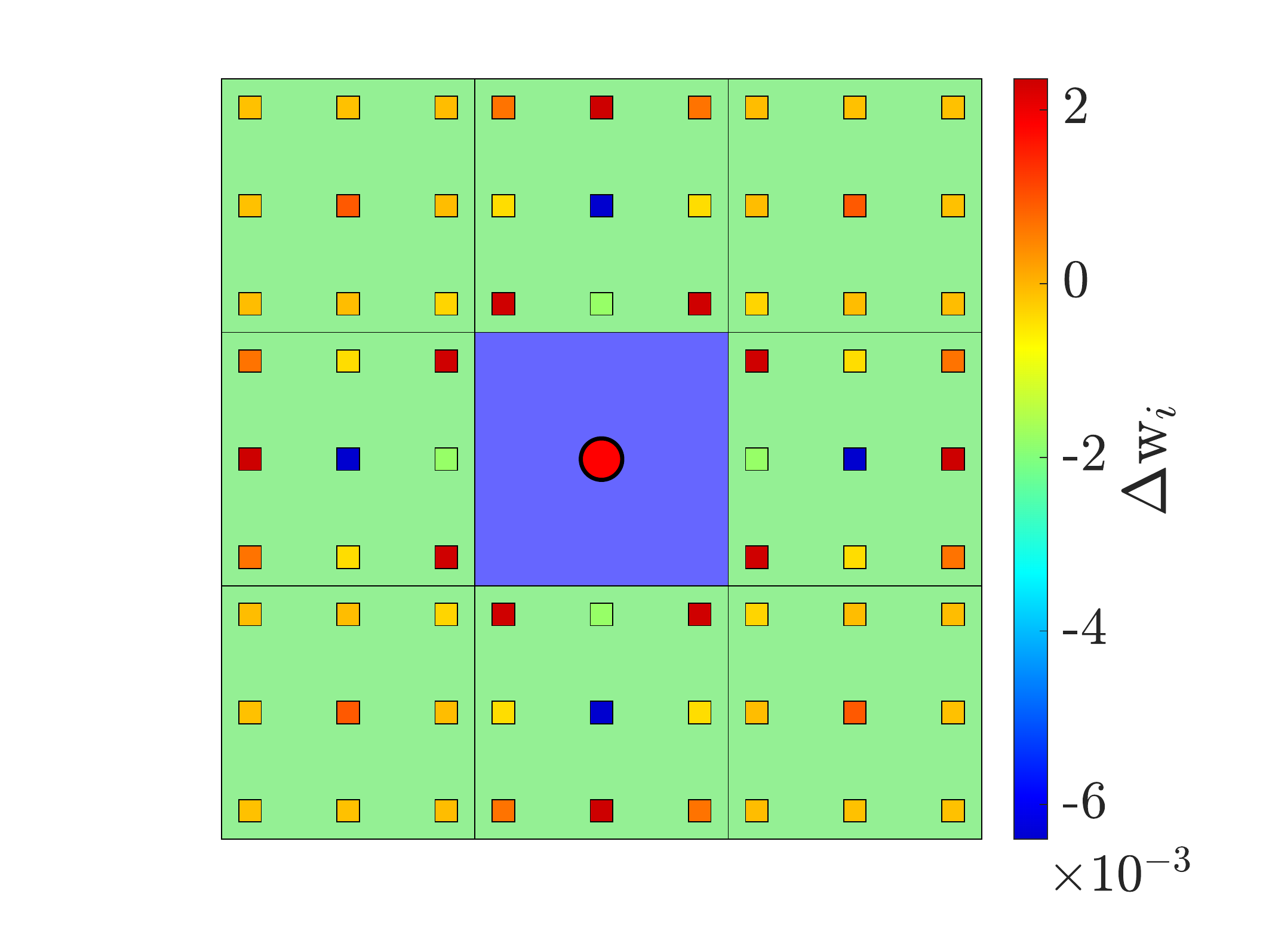}}
	\put(0.65,0){\includegraphics[trim = 120 50 100 30, clip, width=.25\linewidth]{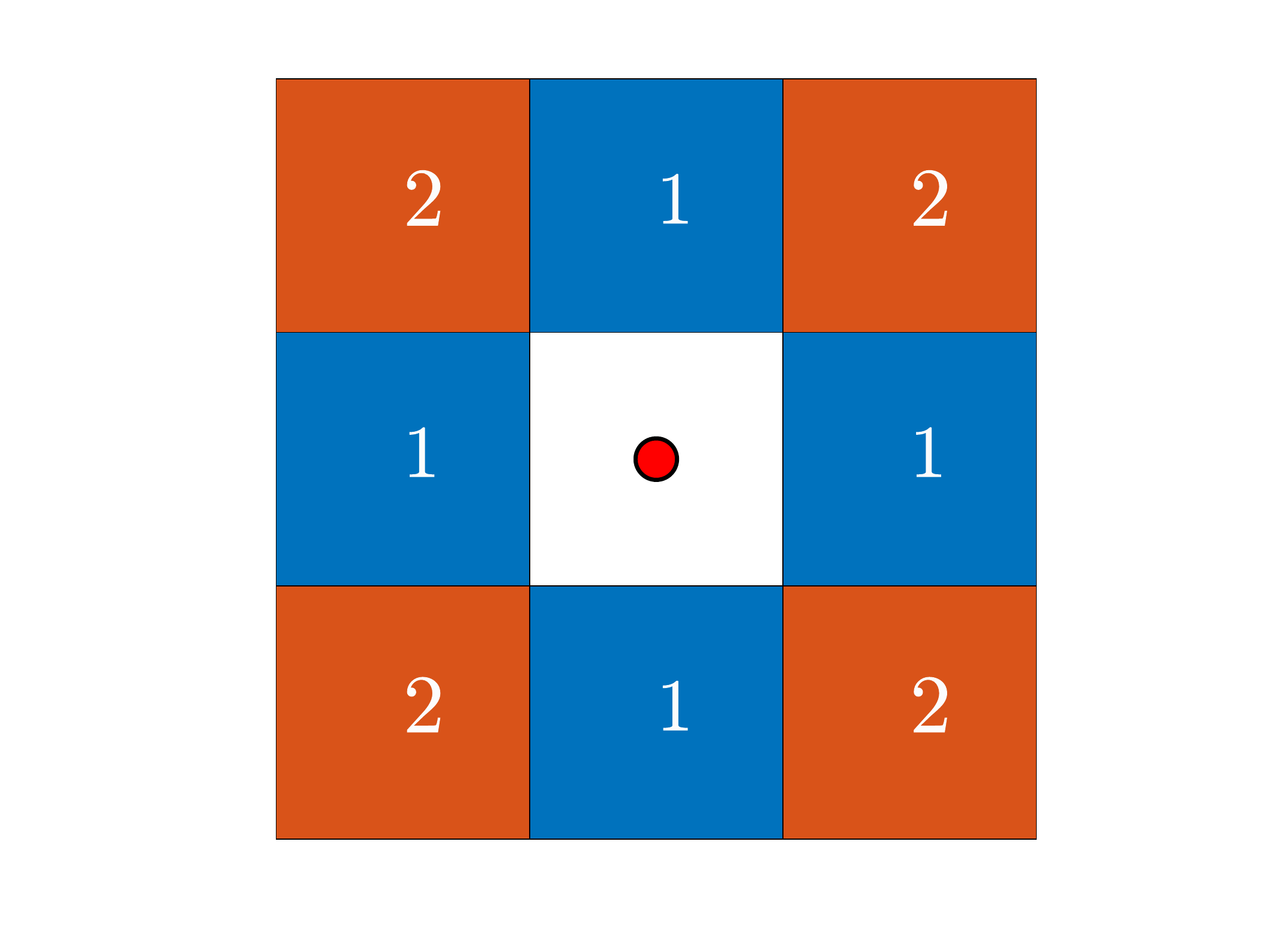}}
	\put(0,0){a.}\put(0.63,0.25){b.}\put(0.63,0){c.}
	\end{picture}
\caption{\textit{Boundary quadrature on singular and near singular elements}: a.~Relative error for classical Gauss-Legendre quadrature with $\tilde n_\mathrm{qp}=2,4,6,\ldots,16$ on the singular element and on the near singular elements from Fig.~\ref{fig:quad_intro}a vs.~the number of quadrature points per element $n_\mathrm{qp}^e$. b.~Adjusted weights for the exact calculation of singular integrals on $3\times3$ plane quadrilateral elements with $\tilde n_\mathrm{qp}=3$.
~c.~Two element types that have the same quadrature values but different orientations.}\label{fig:quad_nearly}
\end{figure}
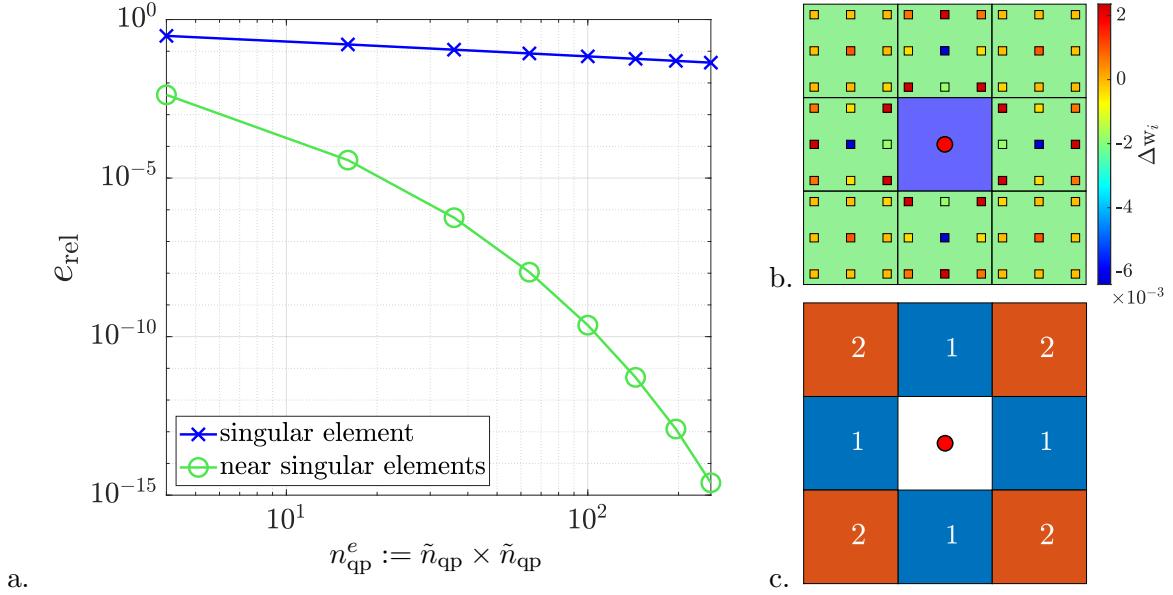
Fig.~\ref{fig:quad_nearly}a shows that the relative quadrature error \eqref{eq:quad_err} for classical Gauss-Legendre quadrature on near singular elements decreases with increasing number of classical Gauss-Legendre quadrature points until it within machine precision for $n_\mathrm{qp}^e= 16 \times 16$ quadrature points per element ($\tilde n_\mathrm{qp}=16$).\footnote{$\tilde n_\mathrm{qp}$ is introduced for the systematic application of quadrature rules in Sec.~\ref{sec:hybrid}. It denotes the number of univariate quadrature points for type Gauss-Legendre quadrature, i.e.~$ n_\mathrm{qp}^e= \tilde n_\mathrm{qp} \times \tilde n_\mathrm{qp}$ (see Fig.~\ref{fig:quad_GL} for $\tilde n_\mathrm{qp}=4$). For Duffy quadrature, $ n_\mathrm{qp}^e$ additionally depends on the position of the collocation point (see Fig.~\ref{fig:quad_Duff2} for $\tilde n_\mathrm{qp}=3$).} The refinement of classical Gauss-Legendre quadrature is on near singular elements much more effective than on singular elements. However, determining an appropriate number of quadrature points per singular element requires a careful trade-off between accuracy and computational effort. Four variables that refer to the number of quadrature points are used in this paper. An overview of these variable is given in Table~\ref{tab:nqp} and more details can be found in Appendix~\ref{app:quad}.
\begin{table}[h]
\centering
\renewcommand*{\arraystretch}{1.2}
\begin{tabular}{|c|c|c|} 
 \hline
 quantity & definition & description \\\hline
 $\tilde n_\mathrm{qp}$ & -- & univariate number of quadrature points\\ 
 $n_\mathrm{qp}^e$ & Eqs.~\eqref{eq:quad_adjusted_n}, \eqref{app_eq:quad_Gauss_nqp}, \eqref{app_eq:quad_Duffy_nqp} and  \eqref{app_eq:quad_GL2_n} & number of quadrature points per element  \\
 $n_\mathrm{qp}$ &  Eq.~\eqref{eq:nqp} & total number of quadrature points \\
 $\bar n_\mathrm{qp}$ & Table \ref{tab:n_qp} & mean number of quadrature points\\\hline
\end{tabular}
\caption{\textit{Boundary quadrature of singularities}: Variables that refer to the number of quadrature points.}\label{tab:nqp}
\end{table}
\subsubsection{Gauss-Legendre quadrature with adjusted weights}\label{sec:quad_nearly_adjusted}
As discussed in Sec.~\ref{sec:quad_nearly_gauss}, a large number of Gauss-Legendre quadrature points is required to obtain accurate approximations of singular integrals on near singular elements. This section presents a new quadrature rule that yields high accuracy by adjusting the quadrature weights for singular kernels instead of increasing the number of quadrature points. The new quadrature rule is therefore referred to as Gauss-Legendre quadrature with adjusted weights.
\\\\ The main idea is to use the locations of $n_\mathrm{qp}^e$ classical Gauss-Legendre quadrature points in the parameter space, denoted by $\bxi_i$ (see Appendix \ref{app:quad_Gauss1}), and determine the corresponding quadrature weights $\mrw_i^e=\mrw^e(\bxi_i)$ by solving a linear moment fitting problem such that $1/r$~kernels are integrated exactly. Quadrature rules based on moment fitting have been introduced by \cite{Joulaian16}, \cite{Thiagarajan16} and \cite{Hubrich17} and applied to surface integration by \cite{Hubrich19} and \cite{Zou21}. Here, the quadrature weights on near singular element $e$ have to fulfill
\begin{equation}\label{eq:quad_adjusted_condition}
	\sum_{i=1}^{n_\mathrm{qp}^e} \frac{J(\bxi_i)\,\mrN_A(\bxi_i)}{r(\bxi_i)}\,  \mrw_i^e \stackrel{!}{=} \int_{\Omega^e} \frac{N_A}{r} \,\mrd a := \mrf_A
\end{equation}
for each shape function $\mrN_A$ with $A=1,\,\ldots,\,(p+1)(q+1)$. The local surface stretch at $\bxi_i$ is denoted by $J(\bxi_i)$ such that $\mrd a(\bxi_i) \approx J(\bxi_i)\, \mrw_i^e$.\footnote{ The reader is referred to Appendices \ref{sec:BE_discr} and \ref{sec:BE_quad} for more details on BE discretization and on the mapping to the parameter space, respectively.} The integral kernel $N_A/r$ is proportional to the discrete BE kernels \eqref{eq:BE_GreensInt}, since $\|\bG\| \propto 1/r$ and  $\|\bT\,\bn \|\propto 1/r$ as shown in Fig.~\ref{fig:GreensFcn_1}c. Quadrature weights that fulfill \eqref{eq:quad_adjusted_condition} on element $e$ thus determine the singular BE integrals on this element exactly. The right hand side of \eqref{eq:quad_adjusted_condition} is determined to 
machine precision by classical Gauss-Legendre quadrature with $16\times 16$ quadrature points (see Fig.~\ref{fig:quad_nearly}a).
\\\\From~\eqref{eq:quad_adjusted_condition} follows that 
\begin{equation}\label{eq:quad_adjusted_n}
n_\mathrm{qp}^e=(p+1)\times (q+1)
\end{equation}
quadrature points are required to obtain the following square moment fitting system of equations for near singular element $e$
	\begin{equation}\label{eq:quad_adjusted_system}
		\mA^e\,\mw^e = \mf^e~,
	\end{equation}
with weight vector $\mw^e= [\mrw_1^e,\ldots,\mrw_{n_\mathrm{qp}^e}^e]^\mrT$, solution vector $\mf^e=[\mrf_1^e,\ldots, \mrf_{n_\mathrm{qp}^e}^e]^\mrT$ and matrix components 
\begin{equation}\label{eq:quad_adjusted_A}
	\mrA_{AB}^e :=  \frac{1}{r(\bxi_B)}\, \mrN_A(\bxi_B) \, J(\bxi_B)~,
\end{equation}
with subscripts $A=1,\ldots,n_\mathrm{qp}^e$ and $B=1,\ldots,n_\mathrm{qp}^e$. The system of equations \eqref{eq:quad_adjusted_system} is solved on each near singular element for the respective weight vector $\mw^e$.
\\\\The proposed Gauss-Legendre quadrature rule with adjusted weights is illustrated by applying it to the biquadratic B-spline sheet from Fig.~\ref{fig:quad_intro}a. The adjusted quadrature weights are determined by solving \eqref{eq:quad_adjusted_system} with $n_\mathrm{qp}^e=9$ for each of the eight near singular elements. Fig.~\ref{fig:quad_nearly}b shows the difference between adjusted and classical Gauss-Legendre weights, defined on singular element $e$ by
\begin{equation}\label{eq:quad_adjusted_dw}
\Delta \mrw_i^e  = \mrw_i^e - \mrw^\mathrm{GL}_i~.
\end{equation}
It can be seen that small adjustments of the weights are sufficient to obtain integral approximations on near singular elements that are exact to machine precision. Due to symmetry, only two different sets of quadrature weight values occur on near singular elements: Larger weight differences on elements sharing one edge with the singular element (element type 1) and smaller weight differences on elements sharing one point with the singular element (element type 2). The near singular elements in Fig.~\ref{fig:quad_nearly}c are colored to indicate elements of type~1 (blue) and elements of type~2 (red). The element type, and thus the quadrature weights, are represented in this way repeatedly throughout the remainder of this paper. The classical Gauss-Legendre weights ($\mrw^\mathrm{GL}_1,\ldots,\mrw_9^\mathrm{GL}$) and the quadrature weight differences for element type 1 ($\Delta \mrw^1_1,\ldots,\Delta\mrw^1_9$) and for element type 2 ($\Delta\mrw^2_1,\ldots,\Delta\mrw^2_9$) are given in Table~\ref{tab:Gauss}.
\\\\The same weight values are applied to all elements of the same element type, but the orientation of the weights may differ. The application of Gauss quadrature with adjusted weights therefore requires some additional bookkeeping. However, this is the only additional effort compared to classical Gauss-Legendre quadrature, since both quadrature rules use the same quadrature point locations. Shape functions and their derivatives can simply be determined a priori at the known locations $\bxi^\mathrm{GL}$ in the parameter domain as discussed in Appendix \ref{sec:BE}. This property makes Gauss quadrature with adjusted weights highly suitable for hybrid quadrature based on classical Gauss quadrature, as will be seen in Sec.~\ref{sec:hybrid}. Both set of weights are derived for a very specific case here, but they yield exact (to the range of machine precision) integrals on regular, quadrilateral and plane NURBS elements of biquadratic order and arbitrary size, as will be seen in Sec.~\ref{sec:hybrid_sheet}. Moreover, Gauss quadrature with adjusted weights is also beneficial on curved surface, where it improves the quadrature accuracy compared to classical Gauss quadrature (see Sec.~\ref{sec:hybrid_sphere}).
\\\\The adjusted weights are determined once on the near singular elements of a plane and regular mesh (green elements in Fig.~\ref{fig:quad_intro}a). These weights (see Tables \ref{tab:Gauss} and \ref{tab:weights}) are then applied to other surfaces. This approach is very efficient in cpu time as bookkeeping is the only additional computational effort compared to classical Gauss-Legendre quadrature. The BE integrals from \eqref{eq:flow_BIE_surf} are determined exactly on plane and regular near singular elements, but only approximately on other, e.g.~curved, near singular elements.
An alternative approach, left for future work, would be to determine adjusted weights for different geometries and discretizations. This approach would be especially efficient for (rotational) symmetric geometries with regular mesh such as cylindrical surfaces. The adjusted weights depend on the curvature and thus on the element size. However, they only need to be evaluated once and can then be stored and re-used for the near singular elements of the same shape. On completely arbitrary surface geometries, it is advisable to use Gaussian quadrature with $16\times 16$ quadrature points directly on near singular elements, rather than computing adjusted weights first.


\section{Hybrid quadrature schemes}\label{sec:hybrid}
The quadrature rules for singular and near singular elements from Sec.~\ref{sec:quad} are combined here to obtain hybrid quadrature schemes capable of efficient singular integral approximation on the entire surface. The hybrid quadrature schemes are then applied to plane B-spline surfaces in Sec.~\ref{sec:hybrid_sheet} and to curved NURBS surfaces in Sec.~\ref{sec:hybrid_sphere}.\footnote{The numerical computations in Secs.~\ref{sec:hybrid} and \ref{sec:examples} are carried out in MATLAB using IEEE~754 double precision ({53~bit} with machine precision ${\varepsilon_\mathrm{machine}\approx 2.2\times 10^{-16}}$).} The findings on hybrid quadrature schemes are summarized in Sec.~\ref{sec:hybrid_conclusion}.
The following four hybrid quadrature schemes are investigated:
\begin{packeditemize}
 \item[a.] \textbf{hybrid Gauss-Legendre quadrature~(G)}
	\begin{packeditemize}
\setlength{\itemindent}{-5pt}
		\item singular element: modified Gauss-Legendre quadrature with $\tilde n_\mathrm{qp}= n_0$
		\item near singular elements: classical Gauss-Legendre quadrature with $\tilde n_\mathrm{qp}=n_0$
	\end{packeditemize}
 \item[b.] \textbf{hybrid Duffy-Gauss quadrature (DG)}:
	\begin{packeditemize}
\setlength{\itemindent}{-5pt}
		\item singular element: Duffy quadrature with $\tilde n_\mathrm{qp}=2 n_0$\footnotemark
		\item near singular elements: classical Gauss-Legendre quadrature with $\tilde n_\mathrm{qp}=n_0$
	\end{packeditemize}
 \item[c.] \textbf{hybrid Duffy-Gauss quadrature with progressive refinement (DGr)}
	\begin{packeditemize} 
\setlength{\itemindent}{-5pt}
		\item singular element: Duffy quadrature with $\tilde n_\mathrm{qp}=2 n_0$
		\item near singular elements: classical Gauss-Legendre quadrature with $\tilde n_\mathrm{qp}=2 n_0$
	\end{packeditemize}
 \item[d.] \textbf{hybrid Duffy-Gauss quadrature with adjusted weights (DGw)}
	\begin{packeditemize} 
\setlength{\itemindent}{-5pt}
		\item singular element: Duffy quadrature with $\tilde n_\mathrm{qp}=2 n_0$
		\item near singular elements: Gauss-Legendre quadrature with adjusted weights ($\tilde n_\mathrm{qp}=3$)
	\end{packeditemize}
\end{packeditemize}
\footnotetext{Fig.~\ref{fig:quad_res} shows that increasing the number of Duffy quadrature points drastically reduces the quadrature error on the singular element. However, the overall error for DG, DGr and DGw is dominated by the remaining elements where classical Gauss-Legendre quadrature is applied (see Fig.~\ref{fig:err_cmp} for a B-spline sheet). A higher number of Duffy quadrature points would increase the computational effort without increasing the quadrature error and is therefore not considered here.}
In all cases, classical Gauss-Legendre quadrature with $\tilde n_\mathrm{qp}= n_0$ is used on the regular elements (see Fig.~\ref{fig:quad_sheet_discr} for singular, near singular and regular elements on a plane sheet). The quadrature density $n_0$ controls the total number of quadrature points $n_\mathrm{qp}$ for the hybrid quadrature schemes.\footnote{Variables defining the number of quadrature points for an element or for the entire discretization are summarized in Table~\ref{tab:nqp}. More details can be found in Appendix~\ref{app:quad}.} It should be noted that even though Duffy quadrature can be applied to the entire surface to evaluate the singular integrals efficiently, this only gives accurate results for the entire surface but not for individual elements, as is required in BE analysis.
\\\\As described in Sec.~\ref{sec:intro}, hybrid Duffy-Gauss quadrature~(DG) has been used in BE analysis for Stokes flow problems \citep{Barakat18}. Progressive refinement of Gauss-Legendre quadrature is additionally considered for acoustic problems  \citep{Venaas20}. This approach resembles hybrid Duffy-Gauss quadrature with progressive refinement (DGr), but \cite{Venaas20} do neither describe the refinement strategy nor investigate its influence on the BE result. Hybrid Duffy-Gauss quadrature with adjusted weights (DGw) is an entirely new scheme that is presented in this paper for the first time.
\\\\Each of the four schemes presented is applicable to arbitrary surface discretizations with quadrilateral elements.\footnote{The schemes can also be adapted to triangular elements, since quadrature points and weights on triangles are well-known for Duffy quadrature \cite{Duffy82} and for classical Gauss-Legendre quadrature \citep{Abramowitz64}. The latter can be used to determine the adjusted weights on triangular elements by proceeding exactly as described in Section 3.2.2.} However, the investigation here focus on biquadratic isogeometric discretizations: The quadrature schemes are applied to a biquadratic B-spline sheet in Sec.~\ref{sec:hybrid_sheet} and to a biquadratic NURBS sphere in Sec.~\ref{sec:hybrid_sphere}.

\subsection{Hybrid quadrature on a B-spline sheet}\label{sec:hybrid_sheet}
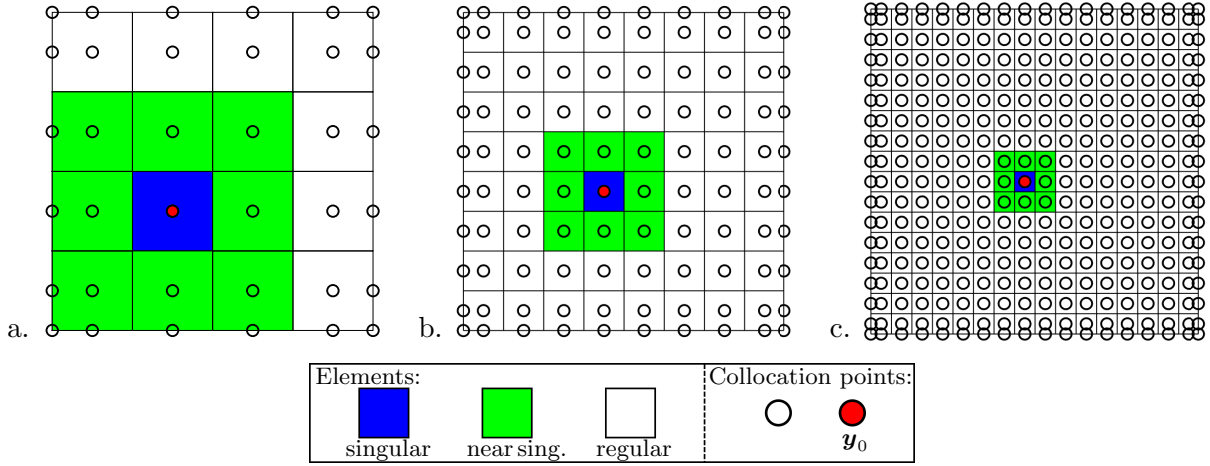
\begin{figure}[h]
\unitlength\linewidth
\begin{picture}(1,0.37)
	\put(0.02,0.08){\includegraphics[trim = 100 50 80 40, clip, width=.3\linewidth]{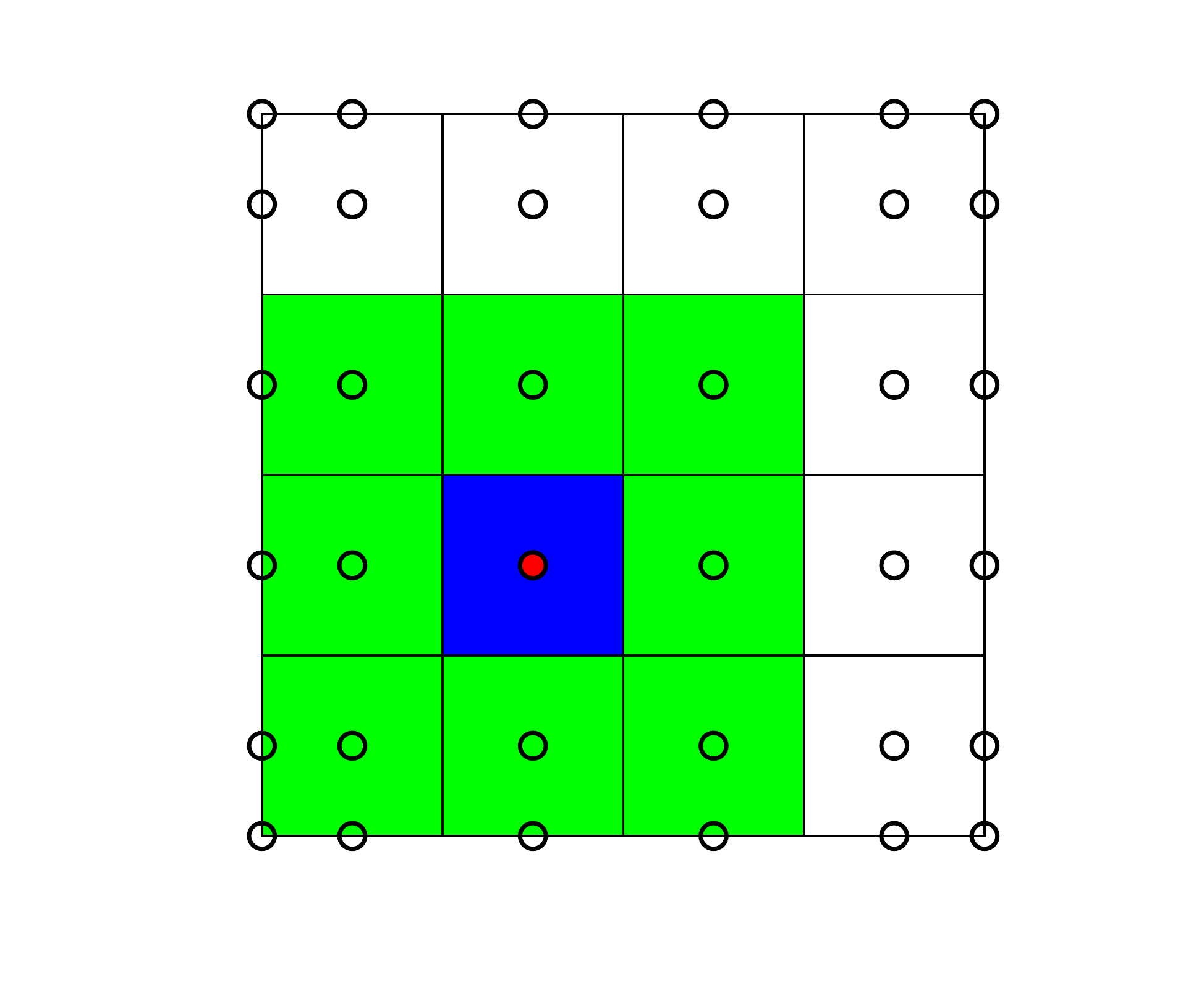}}
	\put(0.36,0.08){\includegraphics[trim = 100 50 80 40, clip, width=.3\linewidth]{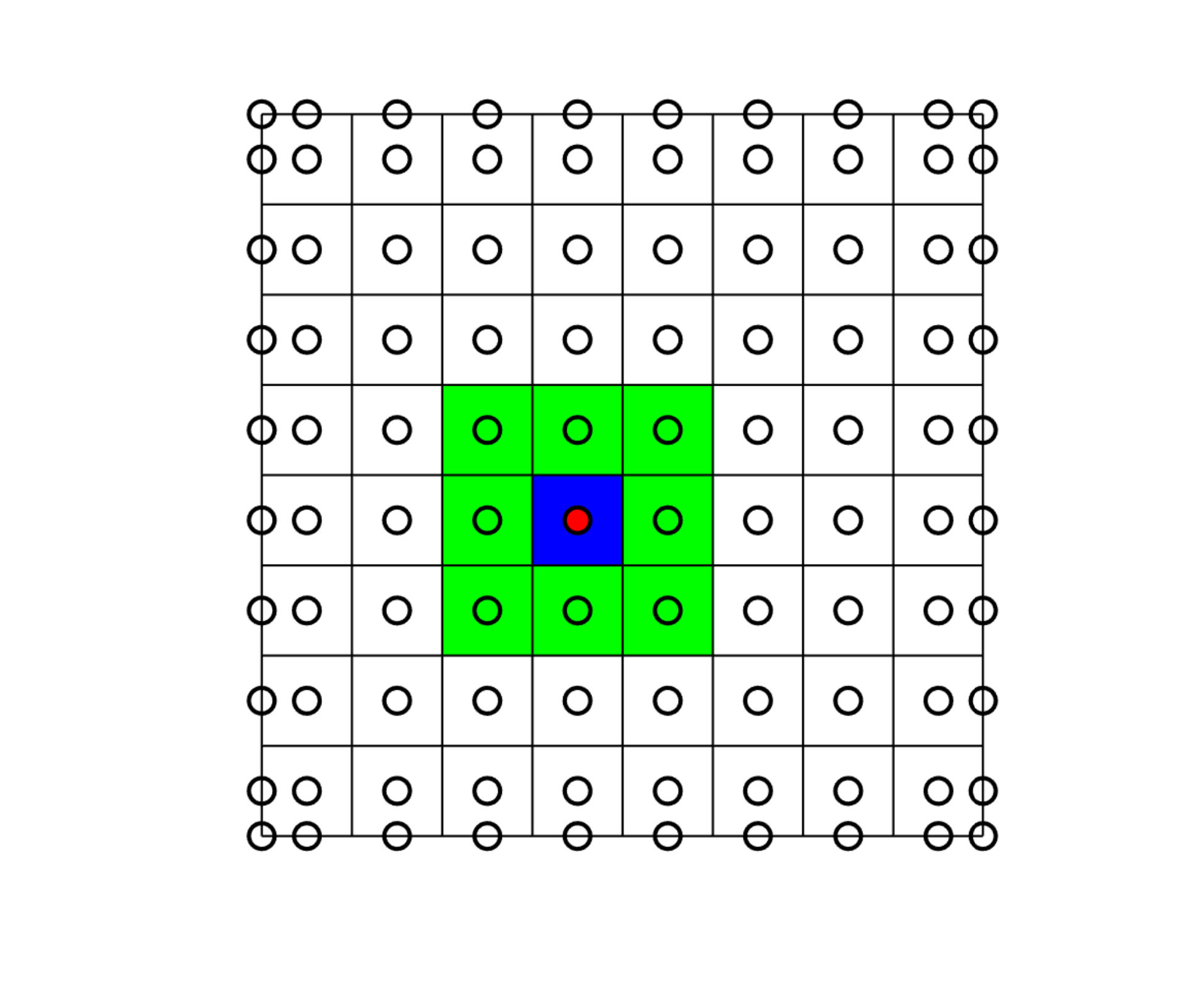}}
	\put(0.7,0.08){\includegraphics[trim = 100 50 80 40, clip, width=.3\linewidth]{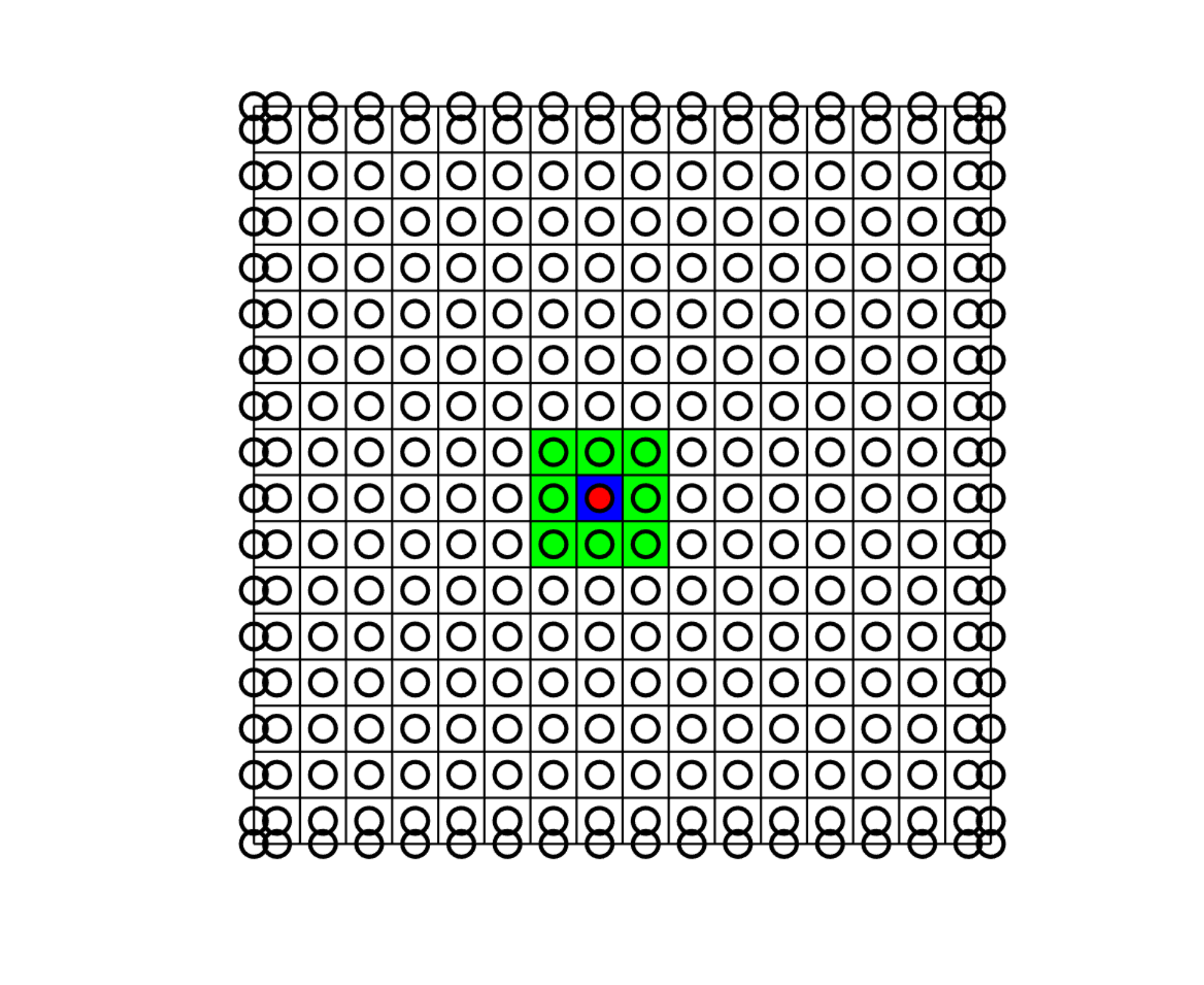}}
	\put(0.25,-0.015){\includegraphics[trim = 0 0 0 0, clip, width=.5\linewidth]{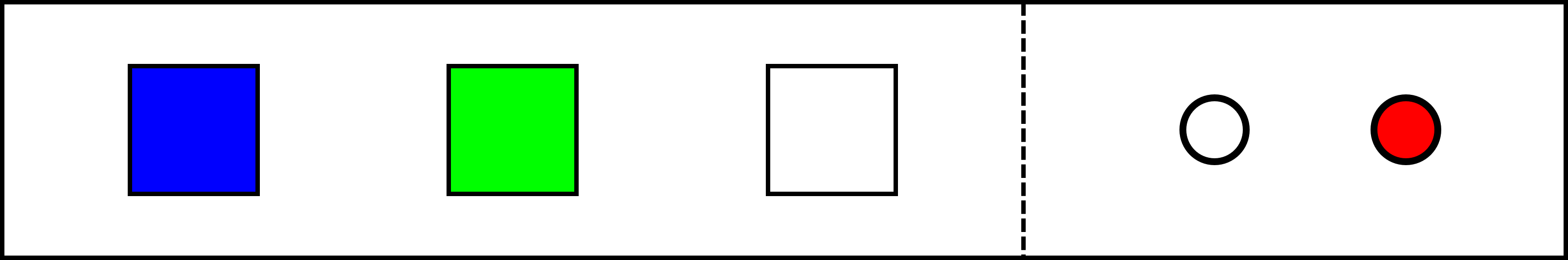}}
	\put(0.255,0.051){\footnotesize Elements:}
	\put(0.28,-0.009){\footnotesize singular}  \put(0.38,-0.009){\footnotesize near\,sing.} \put(0.487,-0.009){\footnotesize regular}
	\put(0.58,0.051){\footnotesize Collocation points:} 
	\put(0.69,0){\footnotesize $\by_0$}
	\put(0,0.09){a.}\put(0.34,0.09){b.}\put(0.68,0.09){c.}
\end{picture}
\caption{\textit{Hybrid quadrature on a B-spline sheet}: Singular, near singular and regular elements for the quadrature of a singularity at collocation point $\by_0$ on a biquadratic sheet with $4\times4$ elements ($\ell=1$, a.), with $8\times8$ elements ($\ell=2$, b.) and with $16\times16$ elements ($\ell=3$, c.).} \label{fig:quad_sheet_discr}
\end{figure}

The introduced hybrid quadrature schemes are applied to a biquadratic B-spline sheet in this section. The coarsest discretization level $\ell=1$ with $4\times 4$ elements is shown in Fig.~\ref{fig:quad_sheet_discr}a, while Figs.~\ref{fig:quad_sheet_discr}b and c show the subsequent levels $\ell=2,3$ created by successive knot insertion. Black circles represent the $n_\mathrm{no}$ collocation points on the sheet. Among those, $\by_0$ (red face) is chosen to illustrate the introduced quadrature schemes. The quadrature rules applied to the B-spline elements are discussed in Sec.~\ref{sec:hybrid_sheet1}. The accuracy of the hybrid quadrature schemes is then compared in Sec.~\ref{sec:hybrid_sheet2}.

\subsubsection{Elemental quadrature}\label{sec:hybrid_sheet1}
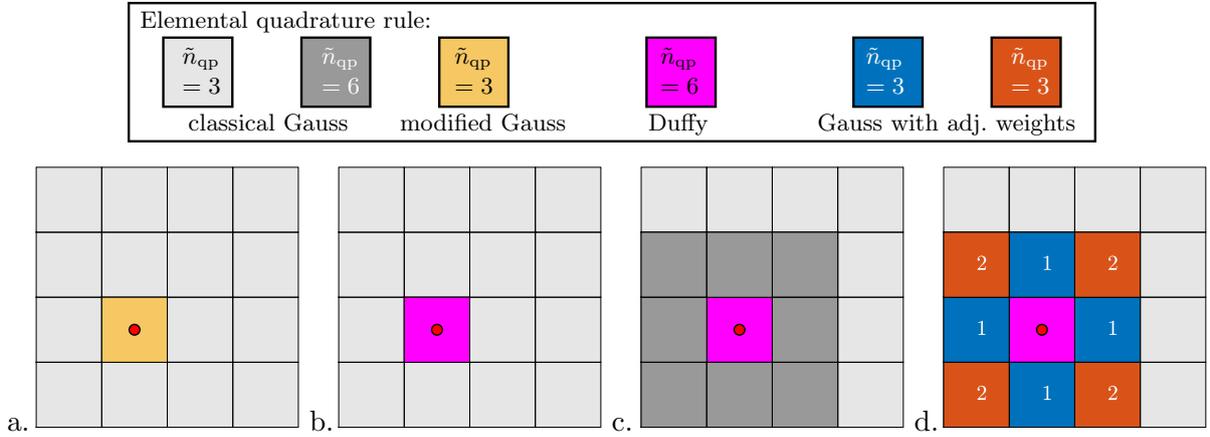
\begin{figure}[t]
\unitlength\linewidth
\begin{picture}(1,0.36)
	\put(0.1,.24){\includegraphics[trim = 0 0 0 0, clip, width=.8\linewidth]{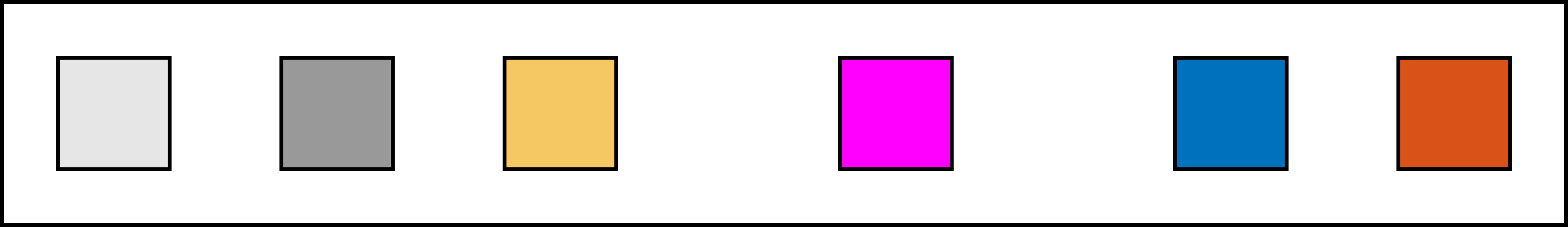}} 
	\put(0.11,.335){\footnotesize Elemental quadrature rule:}	
	\put(0.53,.25){\footnotesize Duffy}	\put(0.67,.25){\footnotesize Gauss with adj.~weights}
	\put(0.15,.25){\footnotesize classical~Gauss}	\put(0.325,.25){\footnotesize modified~Gauss}
	\put(0.145,.305){\footnotesize $\tilde n_\mathrm{qp}$}	\put(0.145,.28){\footnotesize $=3$}
	\put(0.26,.305){\footnotesize \twhite{$\tilde n_\mathrm{qp}$}}	\put(0.26,.28){\footnotesize \twhite{$=6$}}
	\put(0.37,.305){\footnotesize $\tilde n_\mathrm{qp}$}	\put(0.37,.28){\footnotesize $=3$}
	\put(0.54,.305){\footnotesize $\tilde n_\mathrm{qp}$}	\put(0.54,.28){\footnotesize $=6$}
	\put(0.71,.305){\footnotesize \twhite{$\tilde n_\mathrm{qp}$}}	\put(0.71,.28){\footnotesize \twhite{$=3$}}
	\put(0.83,.305){\footnotesize \twhite{$\tilde n_\mathrm{qp}$}}	\put(0.83,.28){\footnotesize \twhite{$=3$}}
	\put(0.02,0){\includegraphics[trim = 100 50 80 30, clip, width=.225\linewidth]{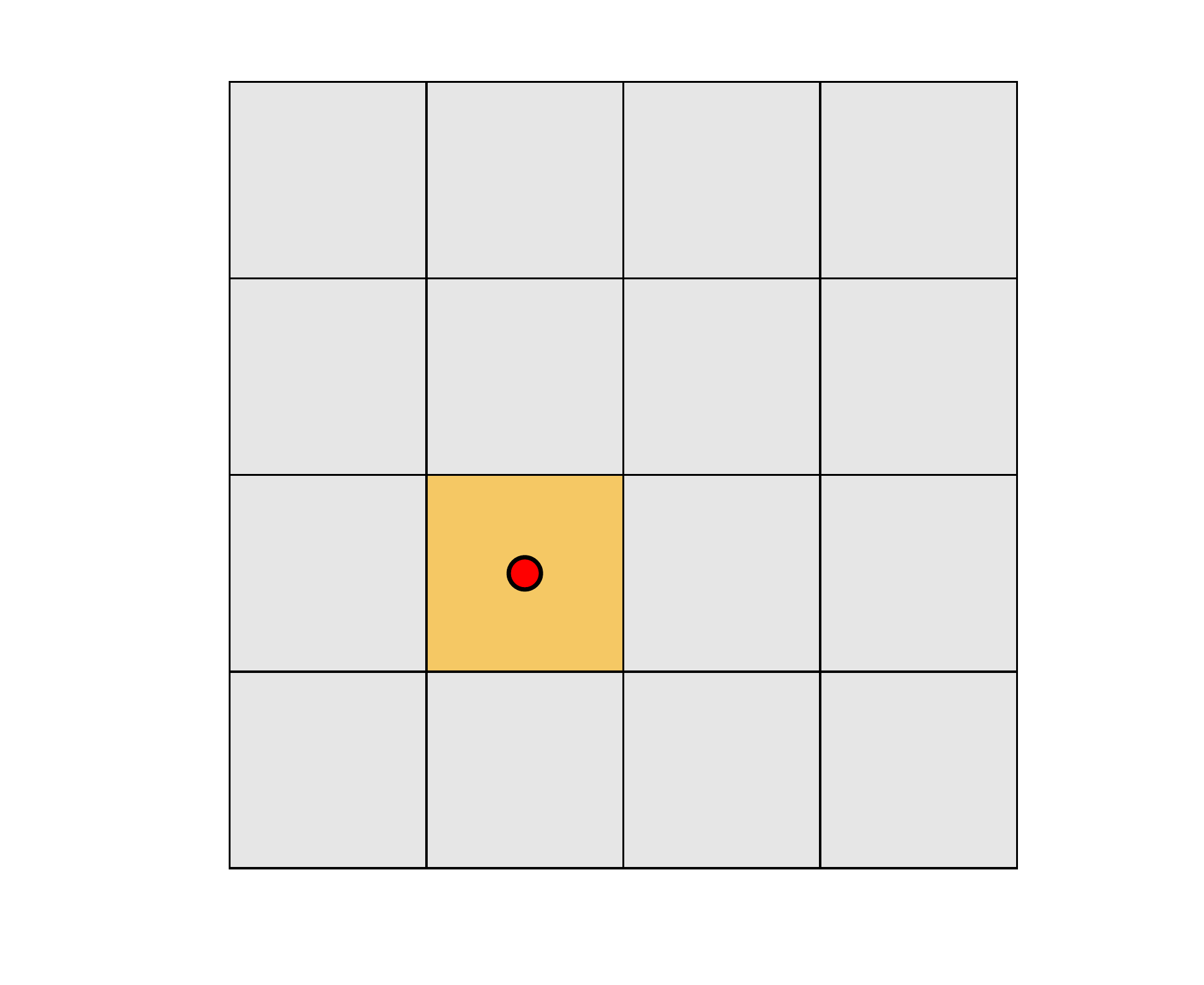}}
	\put(0.27,0){\includegraphics[trim = 100 50 80 30, clip, width=.225\linewidth]{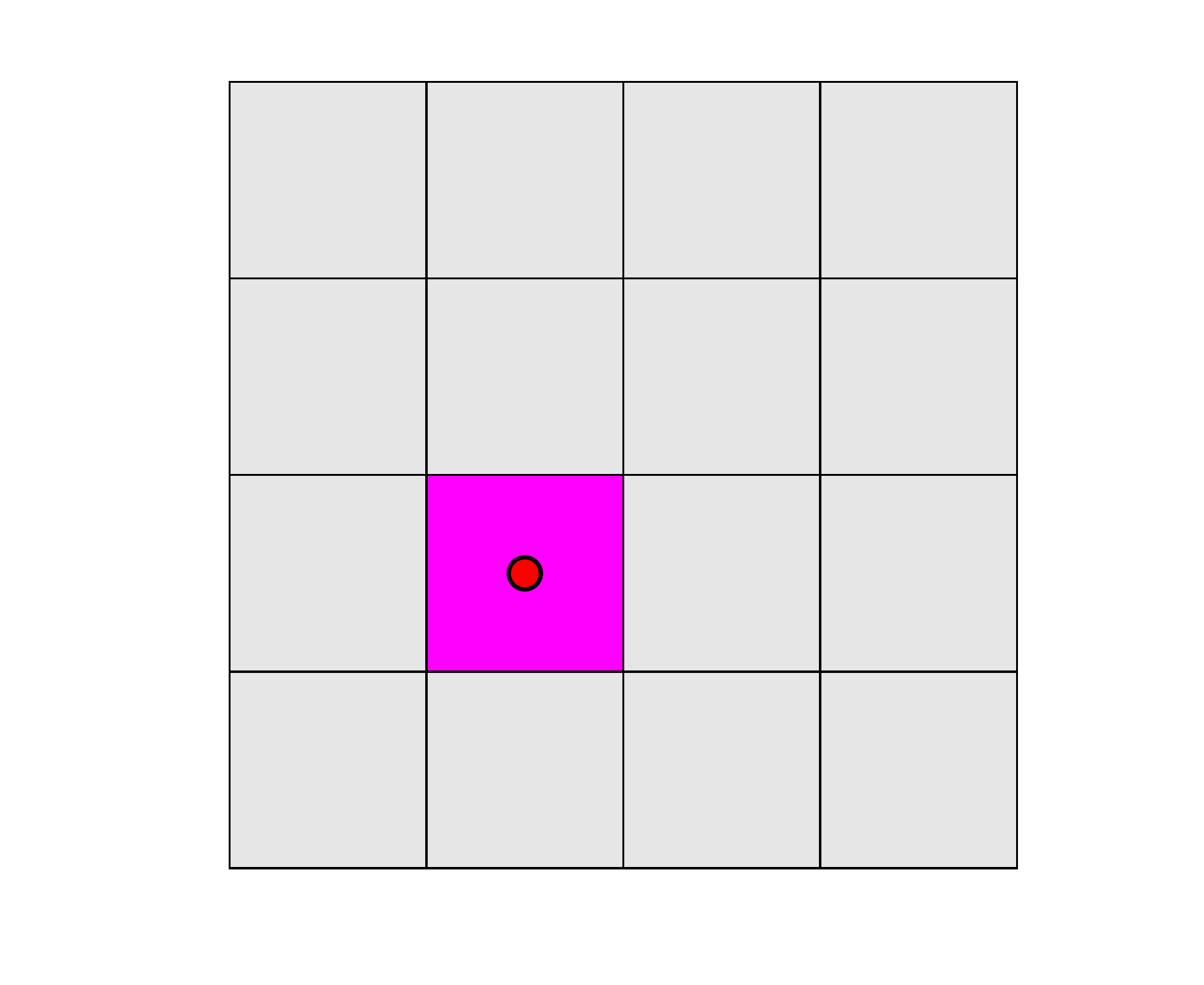}}
	\put(0.52,0){\includegraphics[trim = 100 50 80 30, clip, width=.225\linewidth]{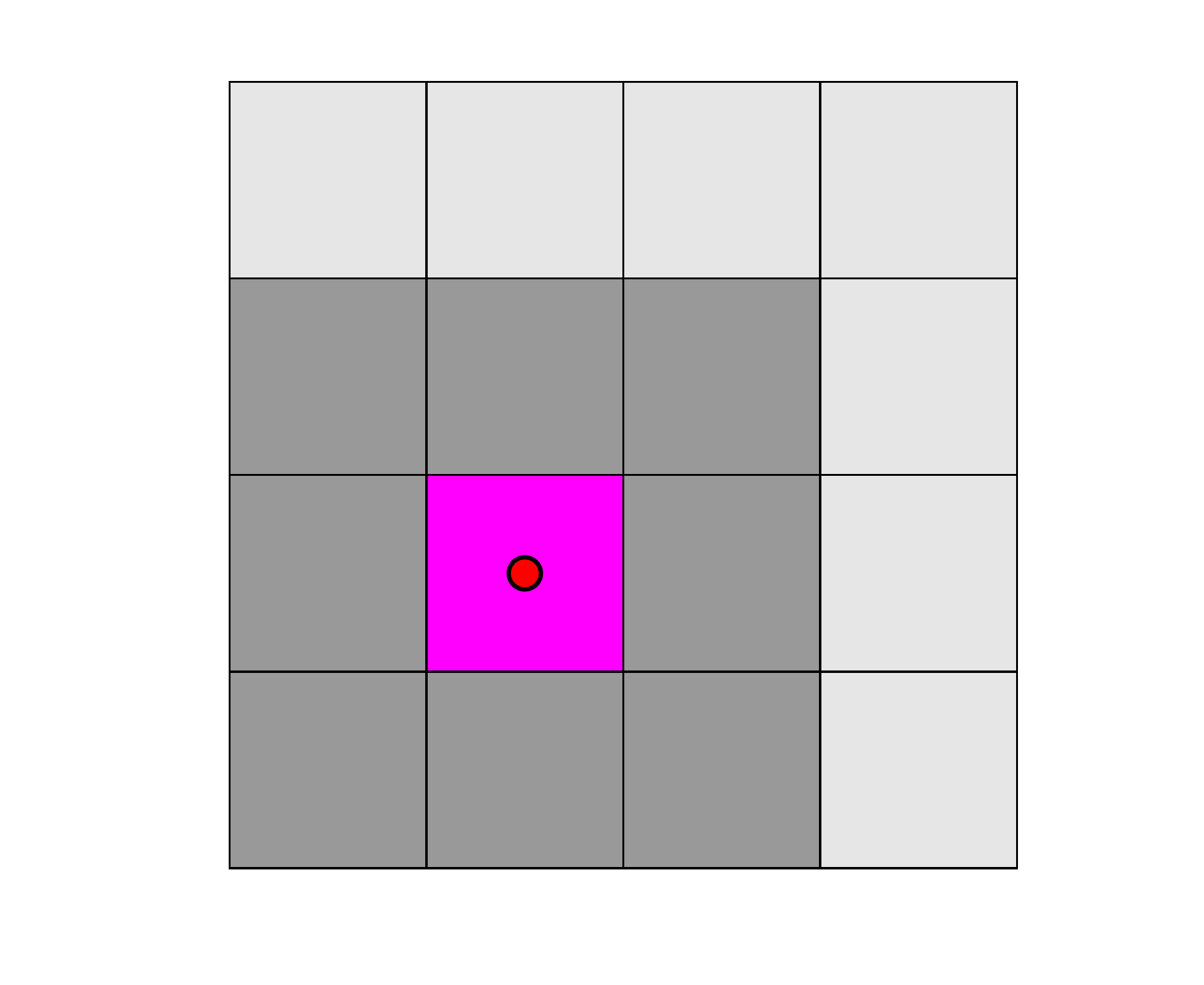}}
	\put(0.77,0){\includegraphics[trim = 100 50 80 30, clip, width=.225\linewidth]{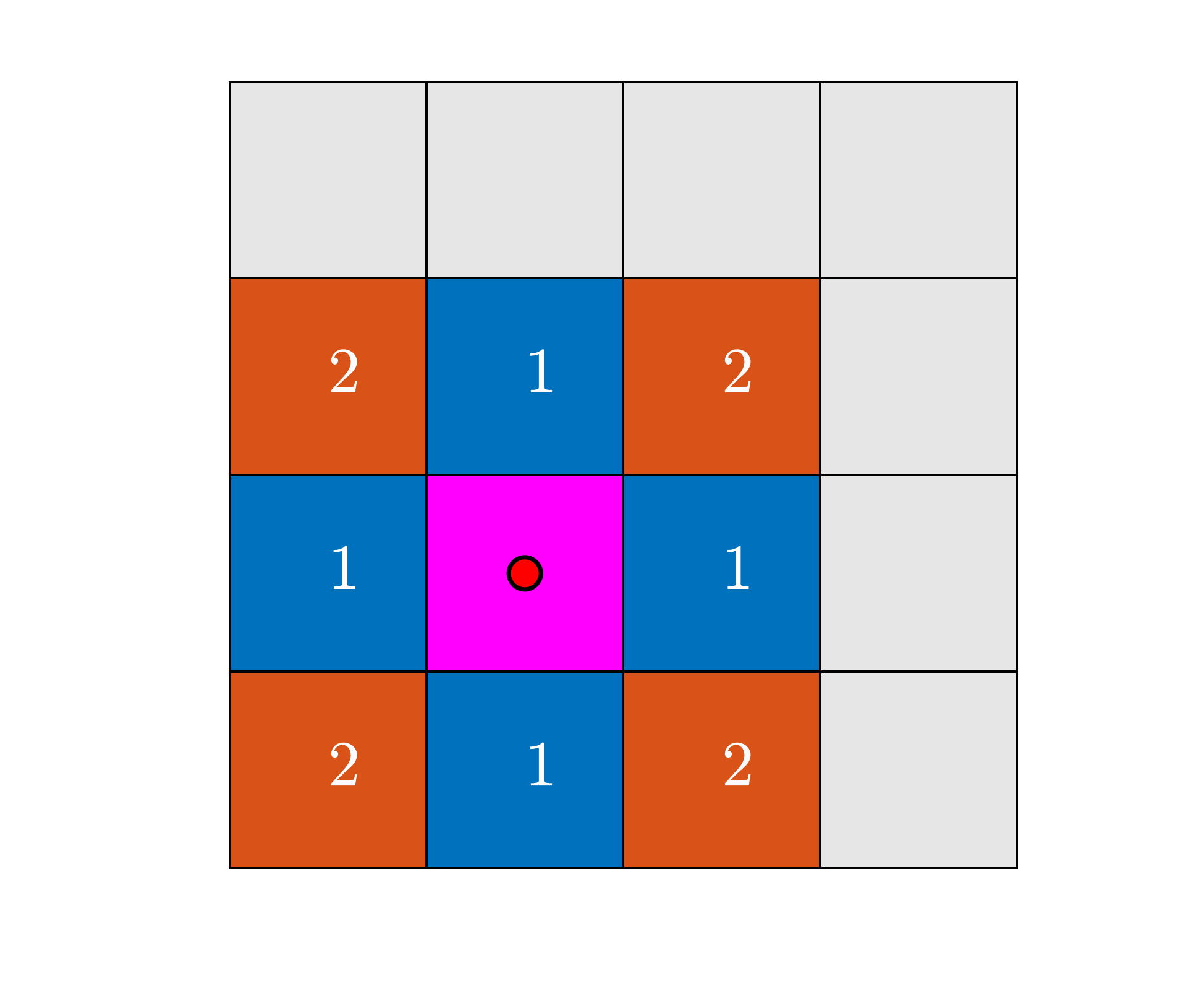}}
	\put(0,0){a.}\put(.25,0){b.}\put(.5,0){c.}\put(0.75,0){d.}
	\end{picture}
\caption{\textit{Hybrid quadrature on a B-spline sheet}: Quadrature rules used for Gauss-Legendre quadrature~(G,~a.), for hybrid Duffy-Gauss quadrature (DG,~b.), for Duffy-Gauss quadrature with progressive refinement~(DGr,~c.) and for hybrid Duffy-Gauss quadrature with adjusted weights~(DGw,~d.) considering collocation point $\by_0$ and quadrature density $n_0=3$.}\label{fig:rules_cmp}
\end{figure}
\begin{figure}[h]
\unitlength\linewidth
\begin{picture}(1,.58)
	\put(0.075,0.48){\includegraphics[trim = 0 0 0 0, clip, width=.8\linewidth ]{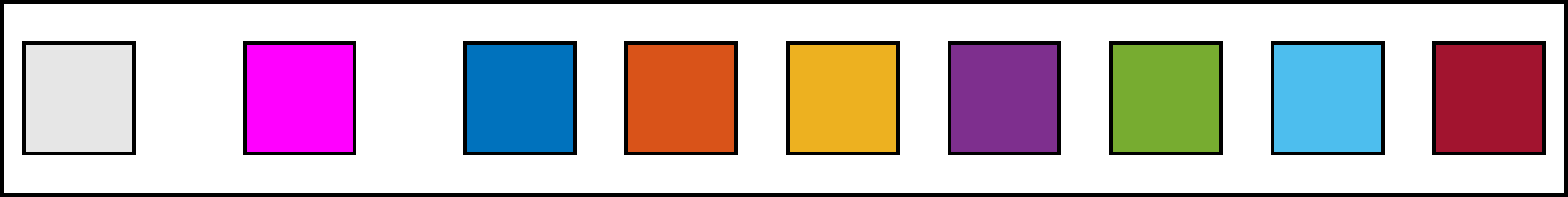}}
	\put(0.08,0.565){\scriptsize Elemental quadrature rule:}
	\put(0.09,0.4865){\scriptsize Gauss} \put(0.205,0.4865){\scriptsize Duffy} \put(0.47,0.4865){\scriptsize Gauss with adjusted weights}
	\put(0.3330,0.523){\twhite 1} \put(0.417,0.523){\twhite 2} \put(0.5000,0.523){\twhite 3} \put(0.5810,0.523){\twhite 4} \put(0.6650,0.523){\twhite 5} \put(0.7480,0.523){\twhite 6} \put(0.8320,0.523){\twhite 7}	
	\put(0,0.24){\includegraphics[trim = 90 40 70 20, clip, width=.237\linewidth ]{figures/quad/scheme3/m1_nReg_3_nDuff_6/y15_rules.pdf}}	
	\put(0.24,0.24){\includegraphics[trim = 90 40 70 20, clip, width=.237\linewidth]{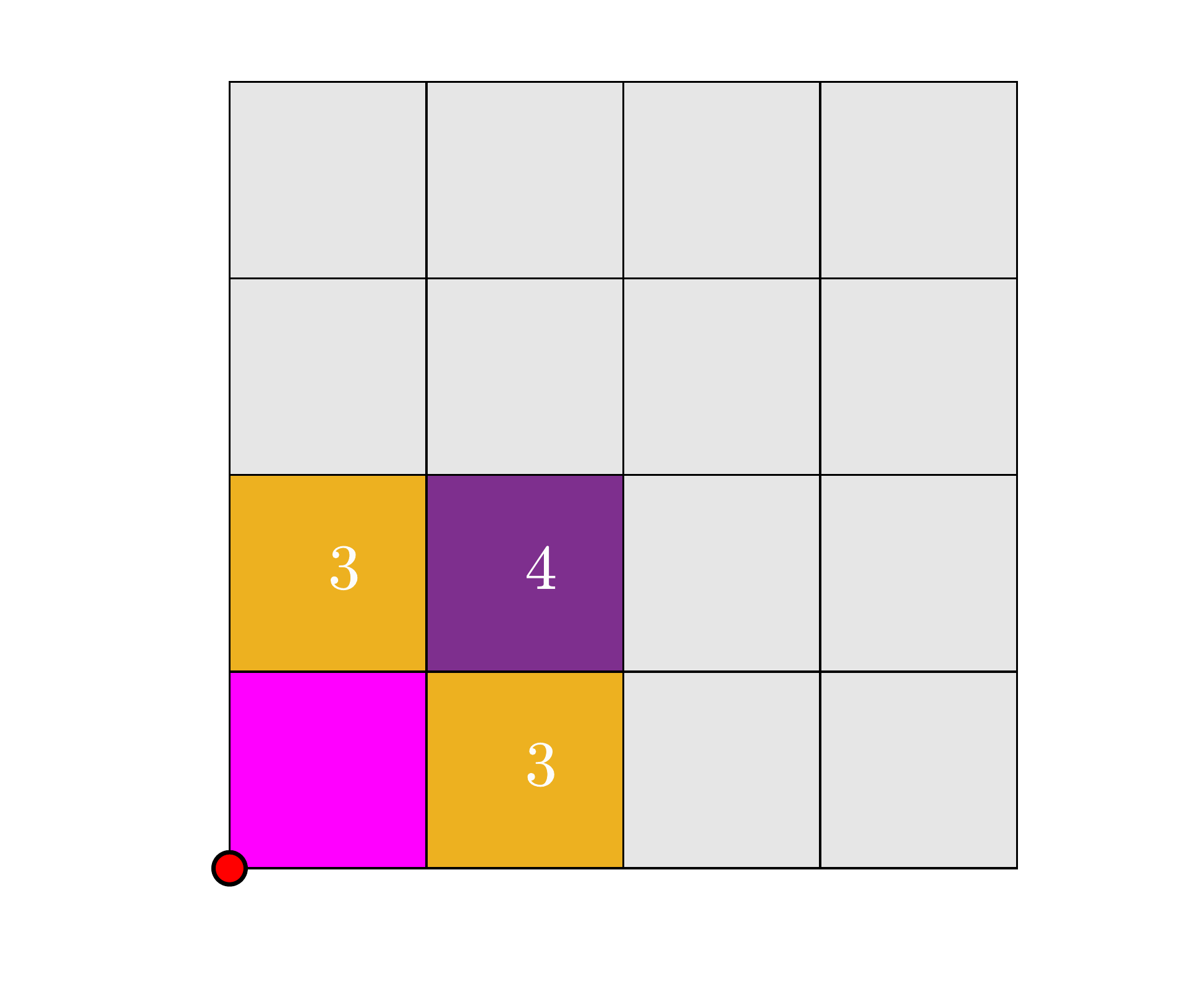}}
	\put(0.48,0.24){\includegraphics[trim = 90 40 70 20, clip, width=.237\linewidth ]{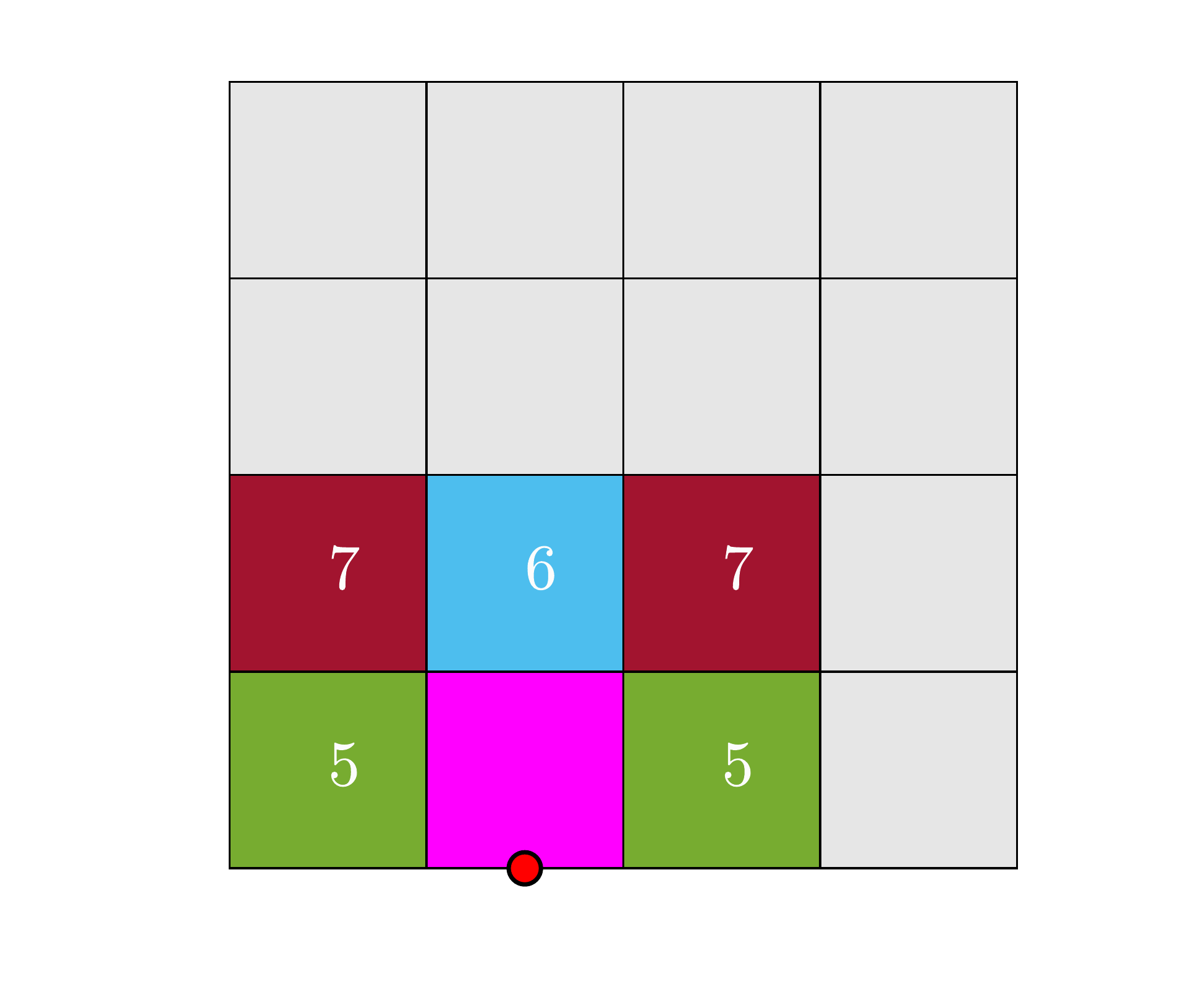}}
	\put(0.72,0.24){\includegraphics[trim = 90 40 70 20, clip, width=.237\linewidth ]{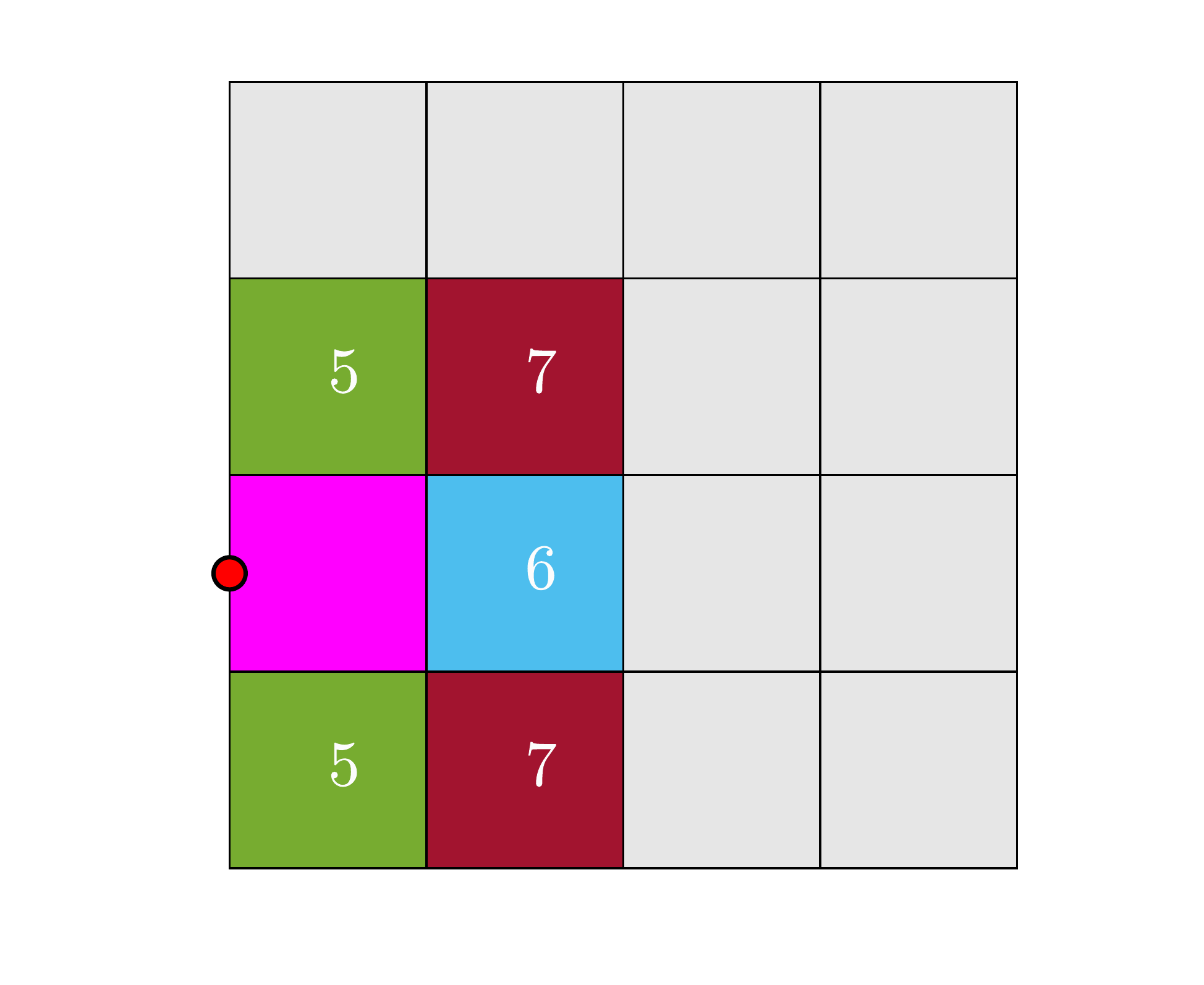}}
	\put(0,0){\includegraphics[trim = 90 40 130 30, clip, width=.2322\linewidth ]{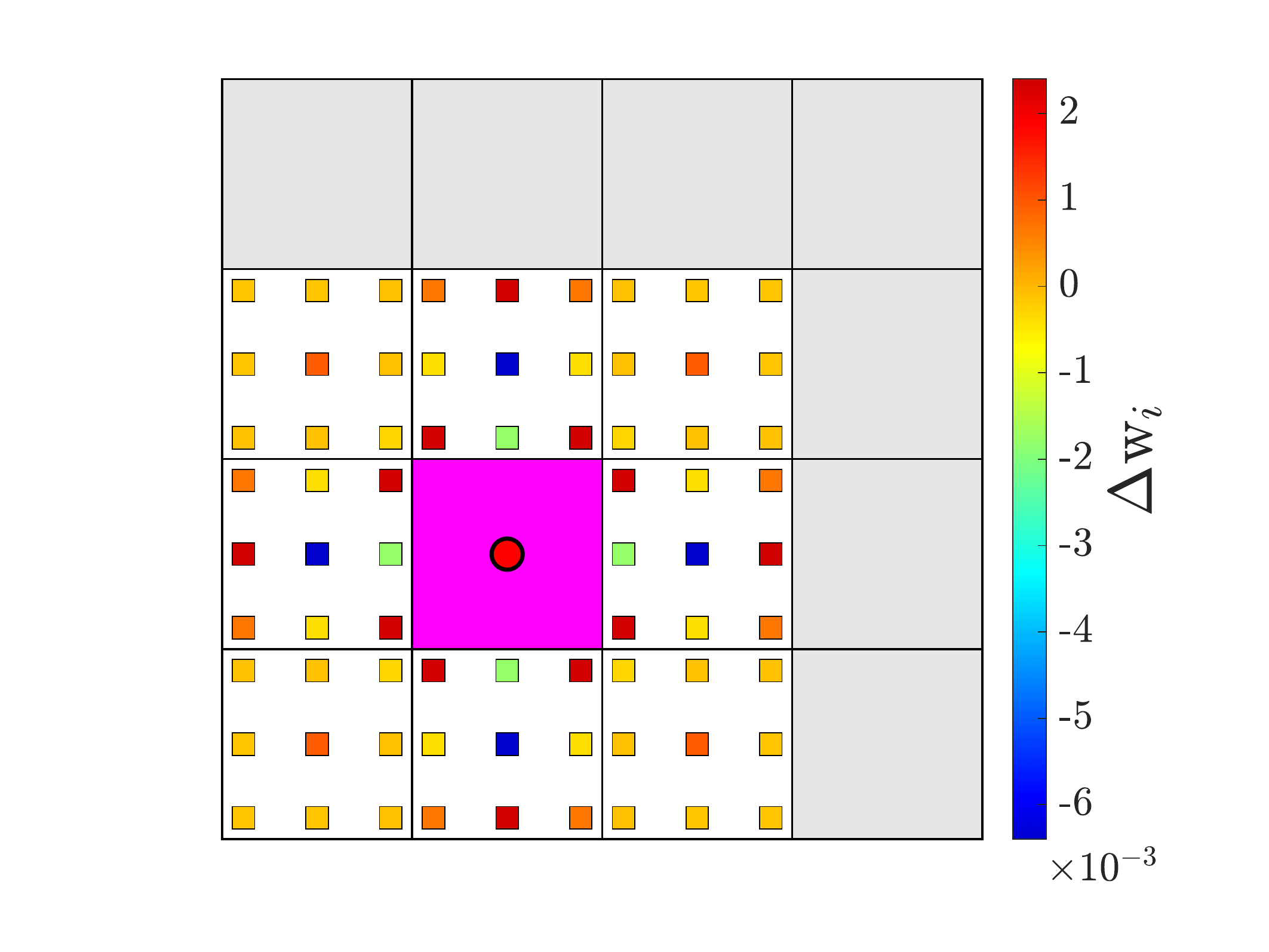}}
	\put(0.24,0){\includegraphics[trim = 90 40 130 30, clip, width=.2322\linewidth]{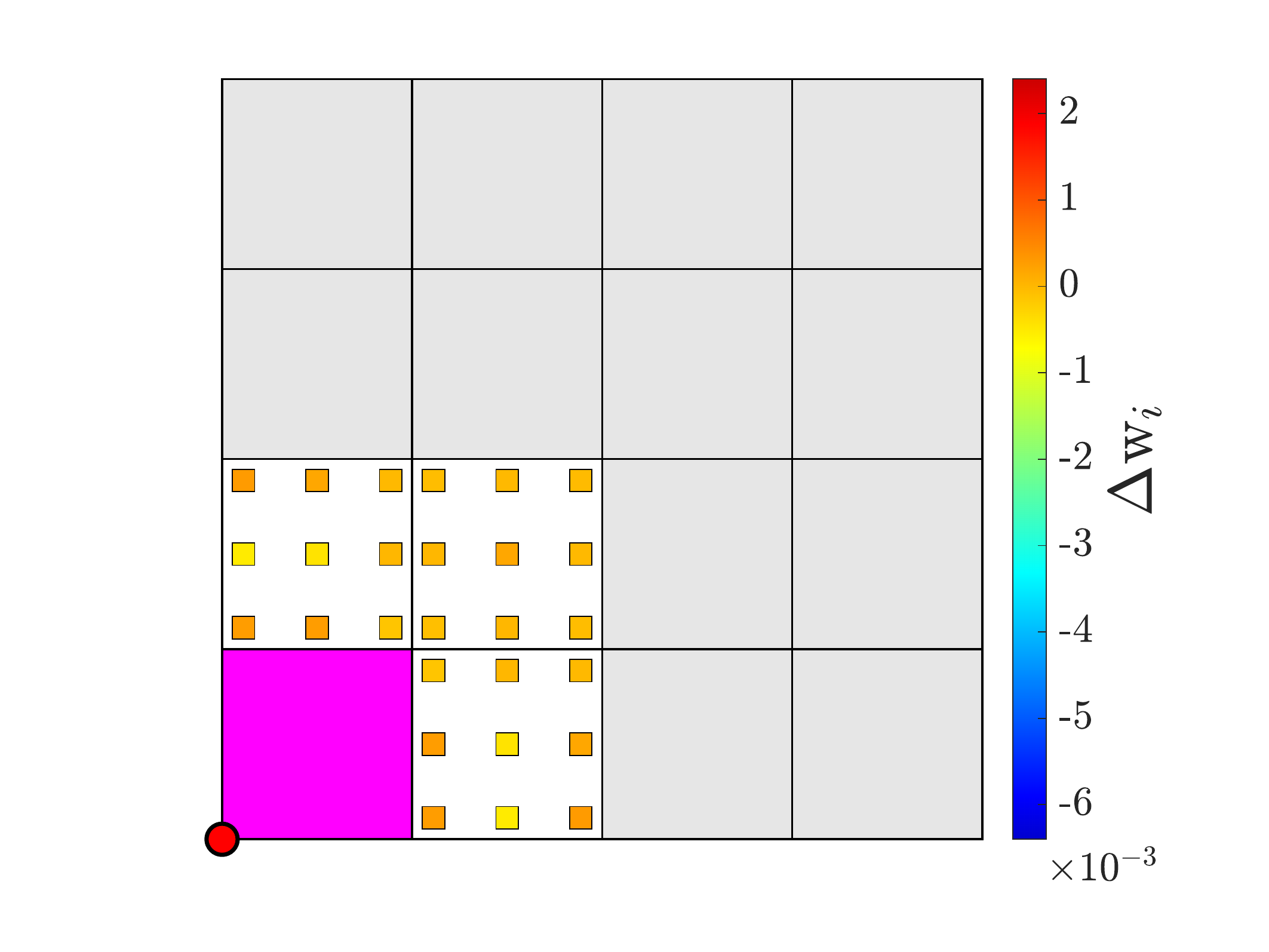}}
	\put(0.48,0){\includegraphics[trim = 90 40 130 30, clip, width=.2322\linewidth ]{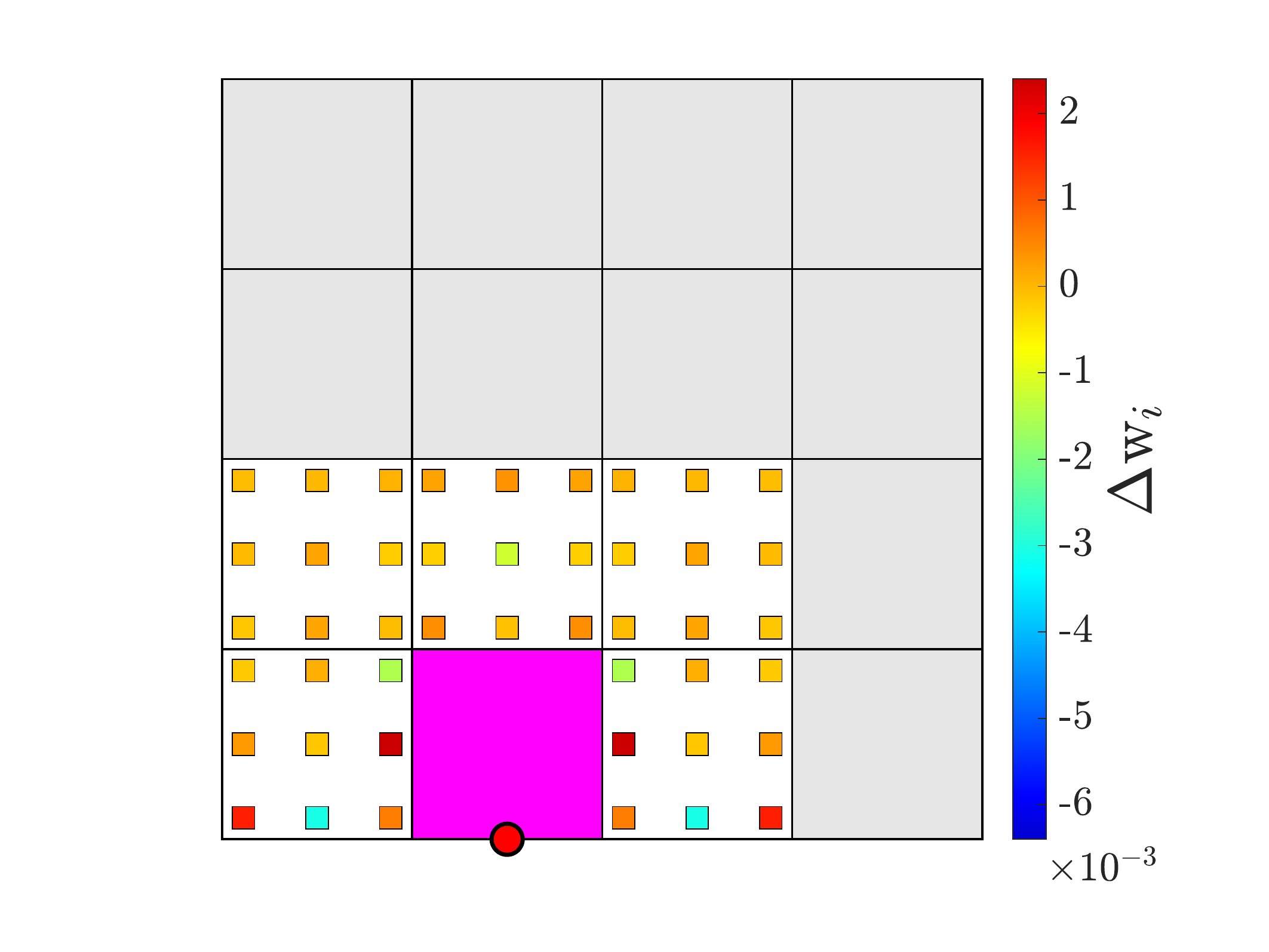}}
 	\put(0.72,0){\includegraphics[trim = 90 40 130 30, clip, width=.2322\linewidth ]{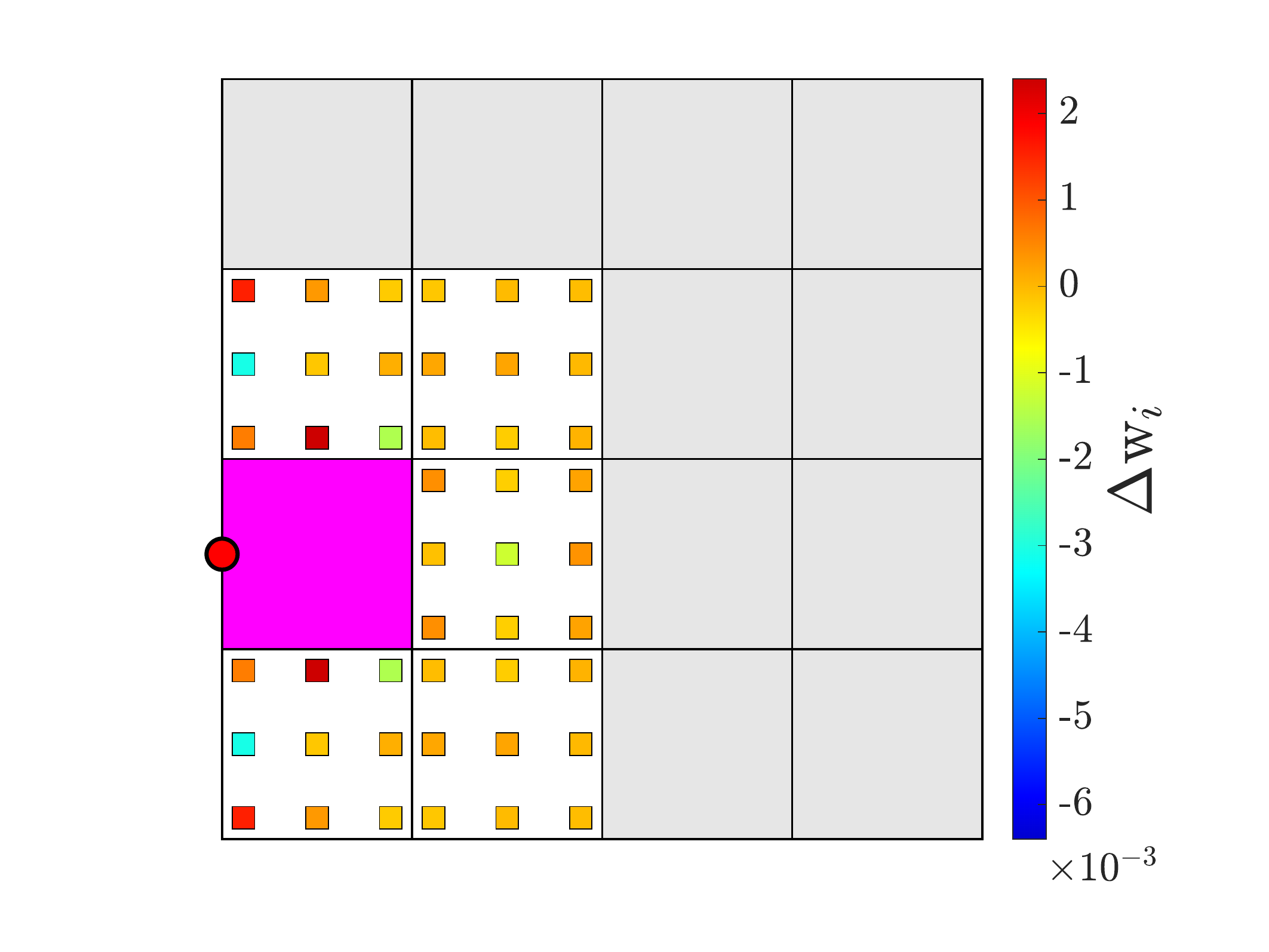}}
	\put(0.95,-0.01){\includegraphics[trim = 475 20 50 30, clip, height=.24\linewidth ]{figures/quad/scheme3/m1_nReg_3_nDuff_6/y13_dw2_abs.pdf}}
	\end{picture}
\caption{\textit{Hybrid quadrature on a B-spline sheet}: Element type (above) and the corresponding weight differences (below) for quadrature scheme DGw and various exemplary collocation points (red dots).}\label{fig:DGw}
\end{figure}

The four hybrid quadrature schemes are applied to the B-spline sheets from Fig.~\ref{fig:quad_sheet_discr}. Fig.~\ref{fig:rules_cmp} shows the elemental quadrature rules used for quadrature schemes G~(a), DG~(b), DGr~(c) and DGw~(d) on the B-spline sheet from Fig.~\ref{fig:quad_sheet_discr}a ($\ell=1$) considering collocation point $\by_0$ and quadrature density $n_0=3$. The elemental quadrature rules on the sheets from Fig.~\ref{fig:quad_sheet_discr}b ($\ell=2$) and c ($\ell=3$) are shown in Fig.~\ref{fig:app_rules}. The corresponding quadrature point locations and weights are given in Fig.~\ref{fig:app_w} for $\ell=1,2,3$.
\\\\In BE analysis, the integration over the whole surface has to be performed for each collocation point once (see Appendix~\ref{sec:BE_collo}). The quadrature weights for Gauss-Legendre quadrature with adjusted weights as well as the weights and locations of the quadrature points for modified Gauss-Legendre quadrature and for Duffy quadrature depend on the location of the collocation point. The introduced quadrature schemes thus require small modifications based on the location of the collocation point. The collocation points on a biquadratic NURBS sheet (black circles in Fig.~\ref{fig:quad_sheet_discr}) can be divided into three types based on their location:
\begin{packeditemize}
	\item \textbf{type C}: at a corner point of an element
	\item \textbf{type E}: at the midpoint of an elemental edge
	\item \textbf{type M}: at the midpoint of an element (e.g.~$\by_0$)
\end{packeditemize}
The location of the Duffy quadrature points are shown in Fig.~\ref{fig:quad_intro}b for collocation points of type M and in Fig.~\ref{fig:quad_Duff2} for collocation points of type C and E. Quadrature points for modified Gauss-Legendre quadrature on a biquadratic element are shown in Fig.~\ref{fig:quad_GL}b. This modification is beneficial in terms of robustness and efficiency for collocation points of type E and M, while classical Gauss-Legendre quadrature (Fig.~\ref{fig:quad_GL}a) is more efficient for type C collocation points.
\\\\Gauss-Legendre quadrature with adjusted weights has been presented and discussed in Sec.~\ref{sec:quad_nearly_adjusted} for collocation points of type M. The quadrature weights for collocation points of other types are obtained in the same manner. Fig.~\ref{fig:DGw} shows the element types for DGw and exemplary collocation points of type C, E and M. It can be seen that only seven sets of quadrature weights are required for the integration over all elements and collocation points. These seven sets, which are required and sufficient for regular biquadratic B-spline sheets of arbitrary refinement, are given in Tables~\ref{tab:Gauss} and~\ref{tab:weights}. The largest weight differences occur by far for element type 1 and 5. Therefore, a simpler but still promising approach is to apply adjusted weight quadrature only to elements of type 1 and 5 and classical Gauss quadrature to the remaining near singular elements.
%

\subsubsection{Quadrature accuracy}\label{sec:hybrid_sheet2}
\begin{figure}[h!]
\unitlength\linewidth
\begin{picture}(1,1.27)
	\put(0.01,0.93){\includegraphics[trim = 90 50 130 0, clip, width=.3\linewidth]{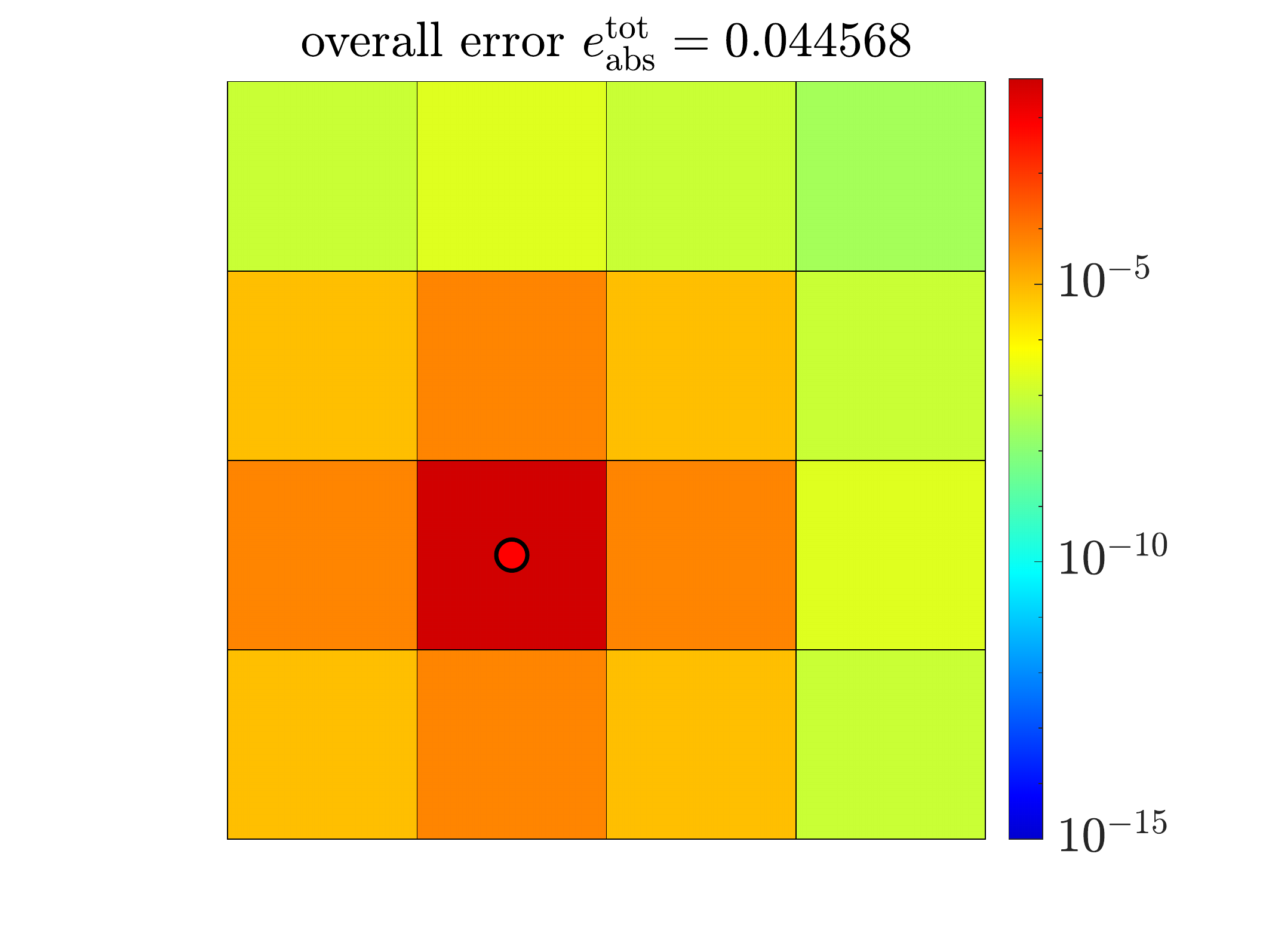}}
	\put(0.31,0.93){\includegraphics[trim = 90 50 130 0, clip, width=.3\linewidth ]{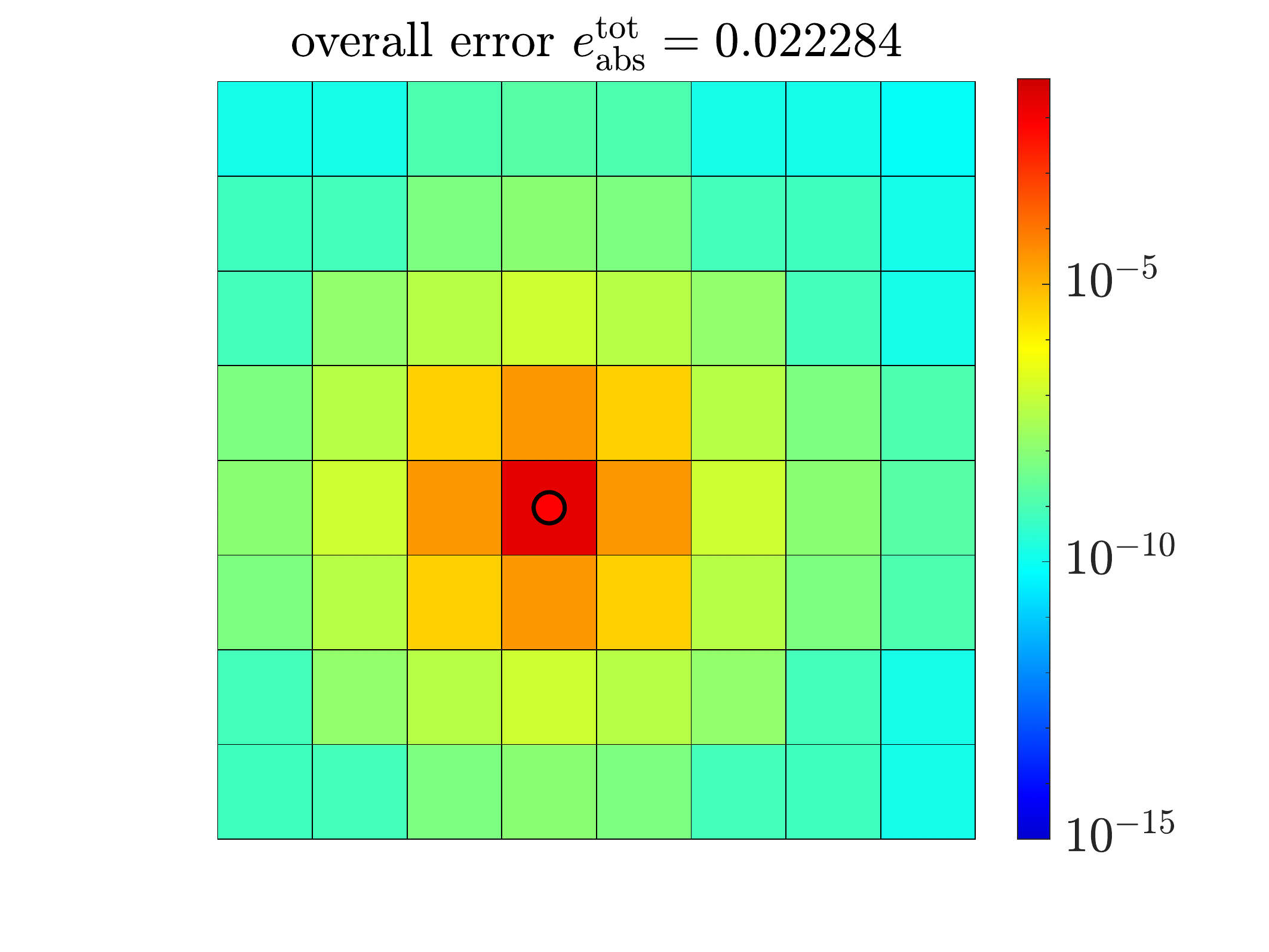}}
	\put(0.61,0.93){\includegraphics[trim = 90 50 130 0, clip, width=.3\linewidth ]{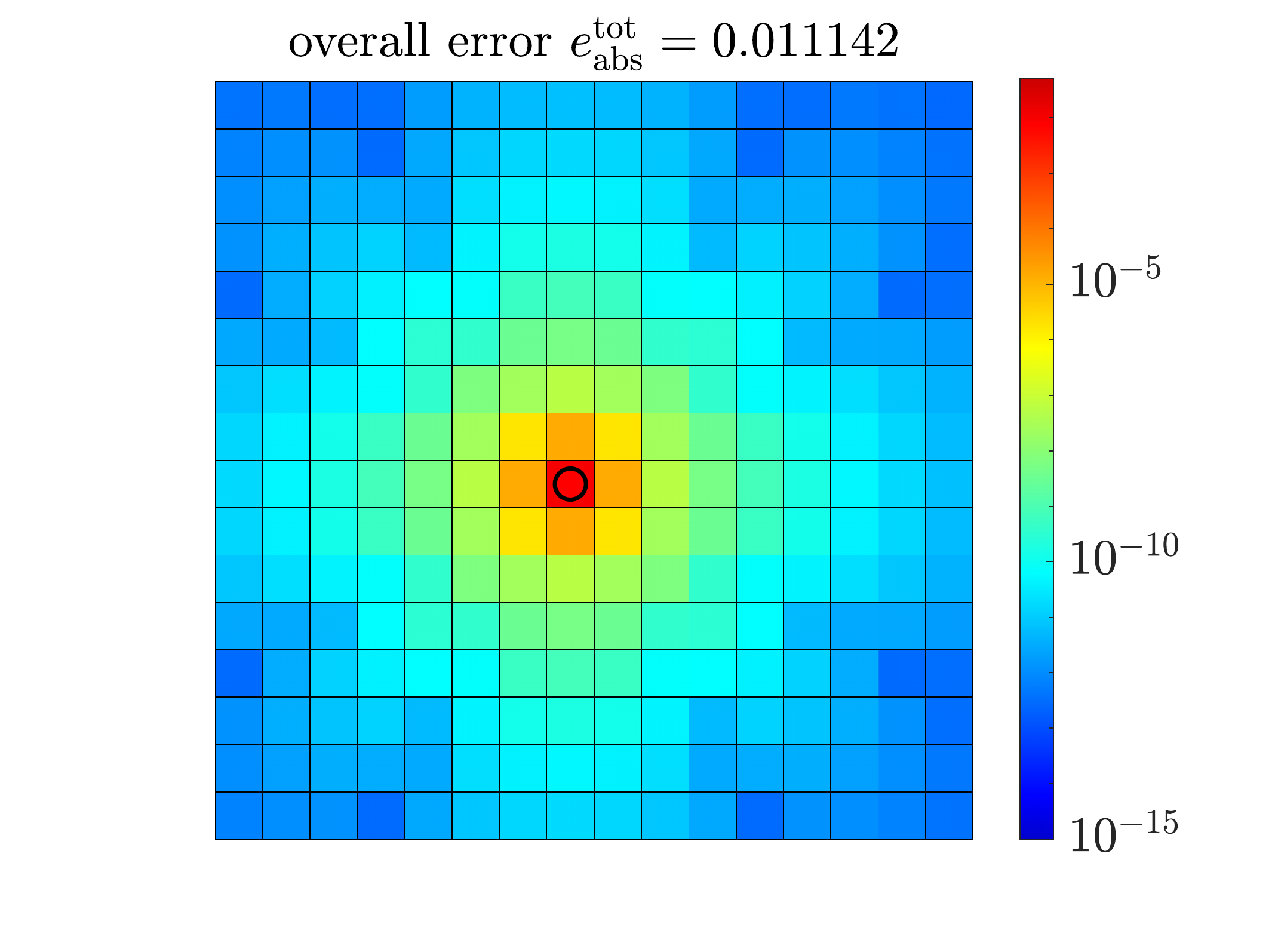}}
	\put(0.01,0.62){\includegraphics[trim = 90 50 130 0, clip, width=.3\linewidth]{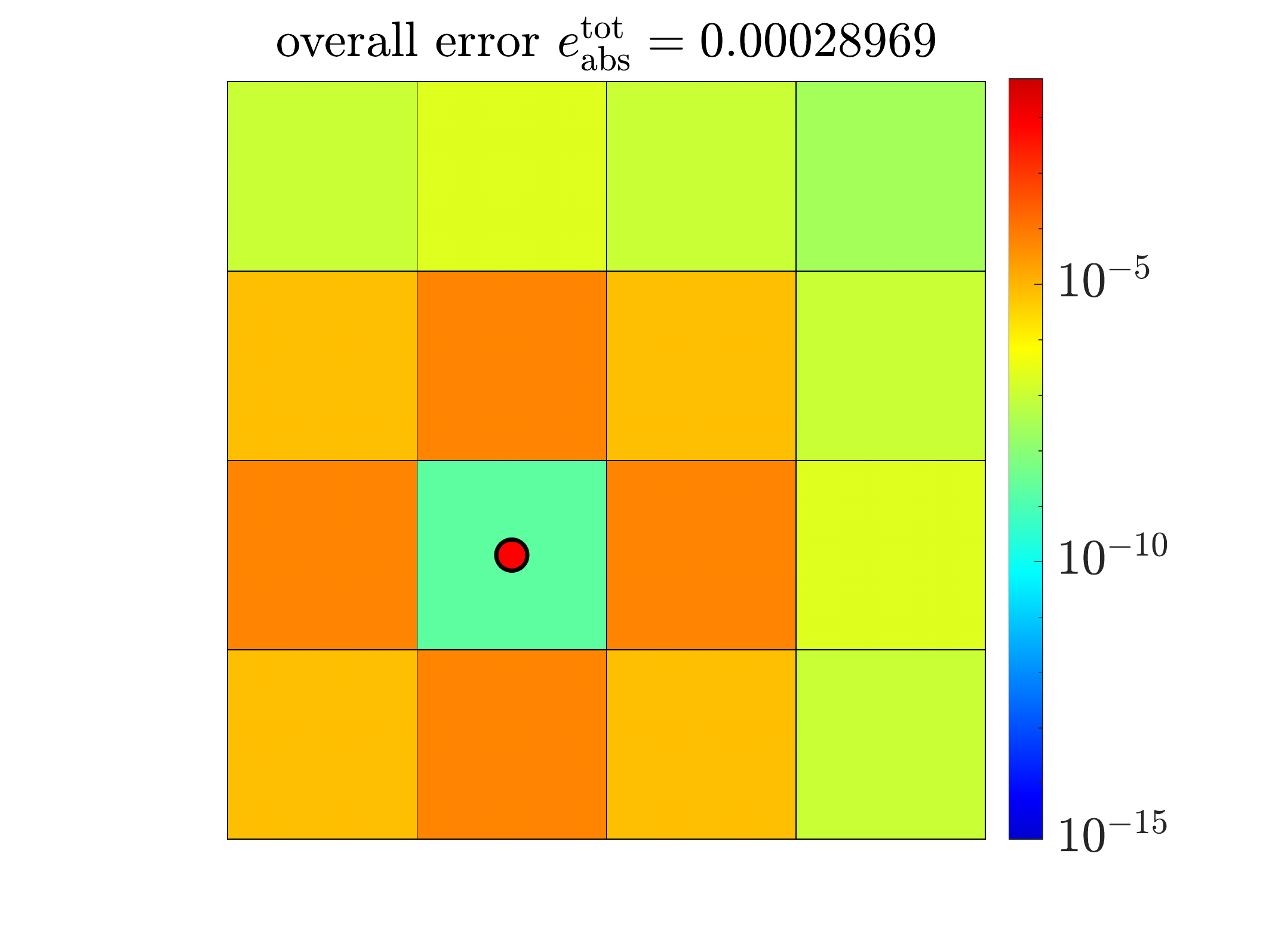}}
	\put(0.31,0.62){\includegraphics[trim = 90 50 130 0, clip, width=.3\linewidth ]{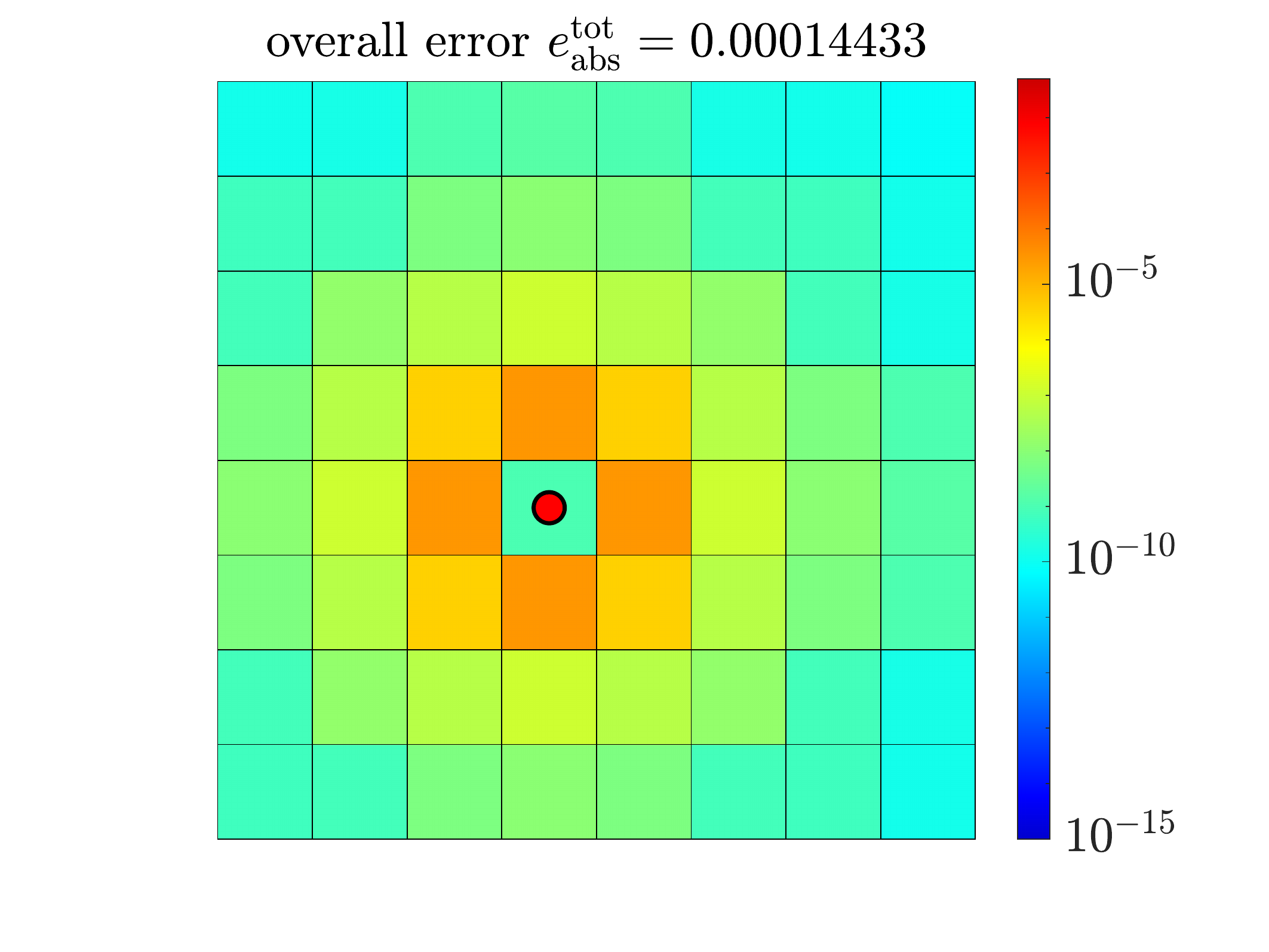}}
	\put(0.61,0.62){\includegraphics[trim = 90 50 130 0, clip, width=.3\linewidth ]{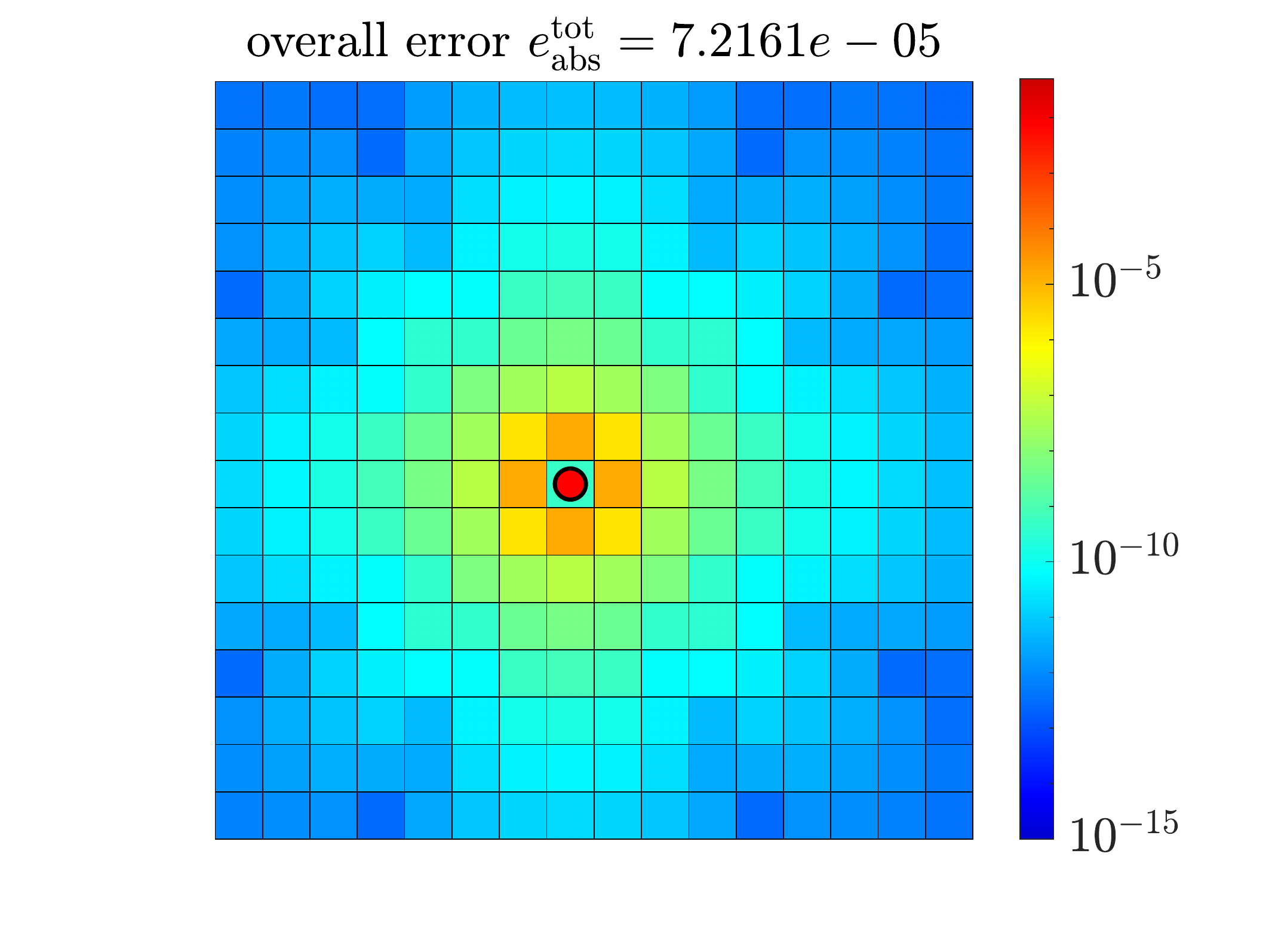}}
	\put(0.01,0.31){\includegraphics[trim = 90 50 130 0, clip, width=.3\linewidth]{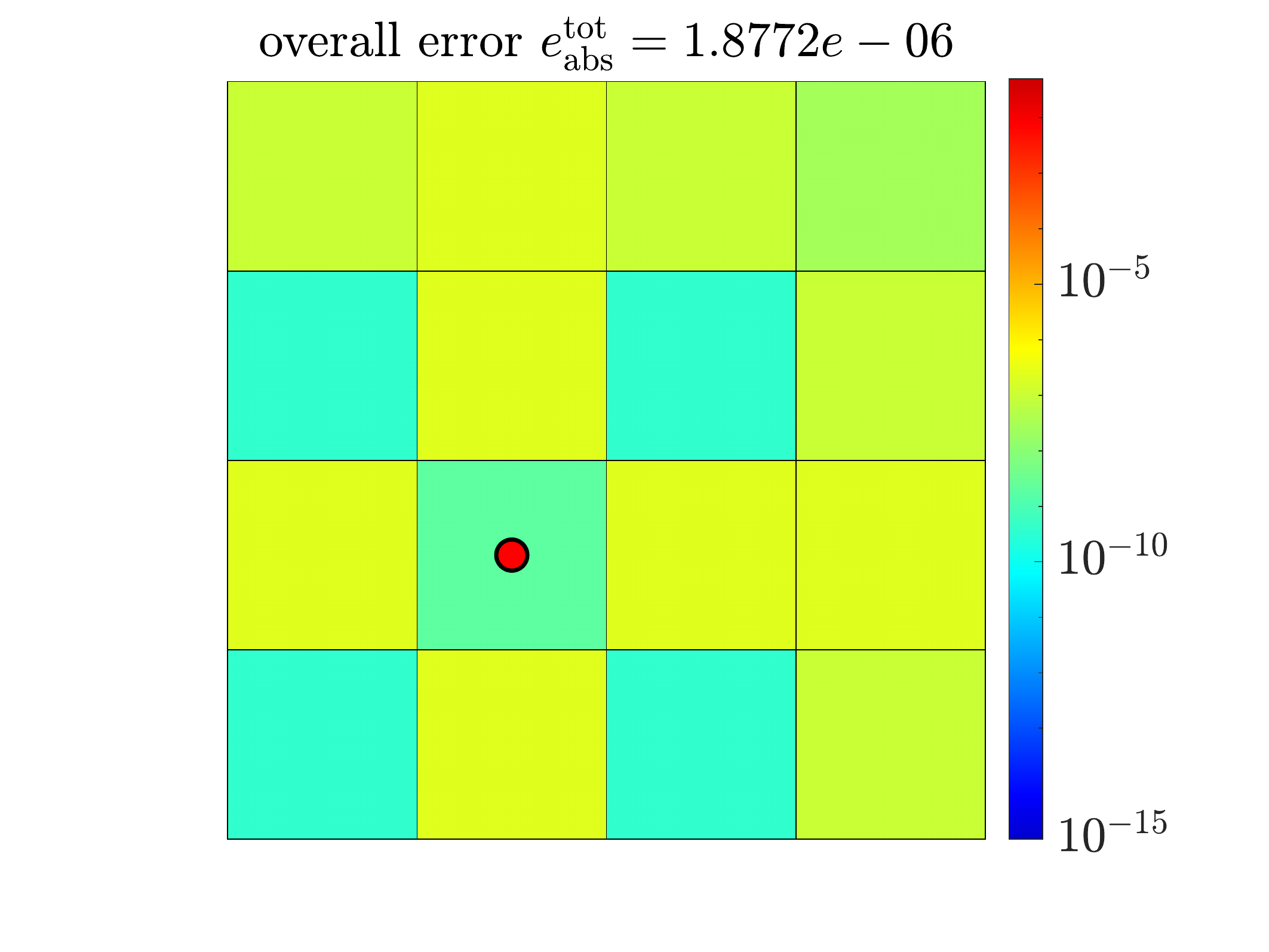}}
	\put(0.31,0.31){\includegraphics[trim = 90 50 130 0, clip, width=.3\linewidth ]{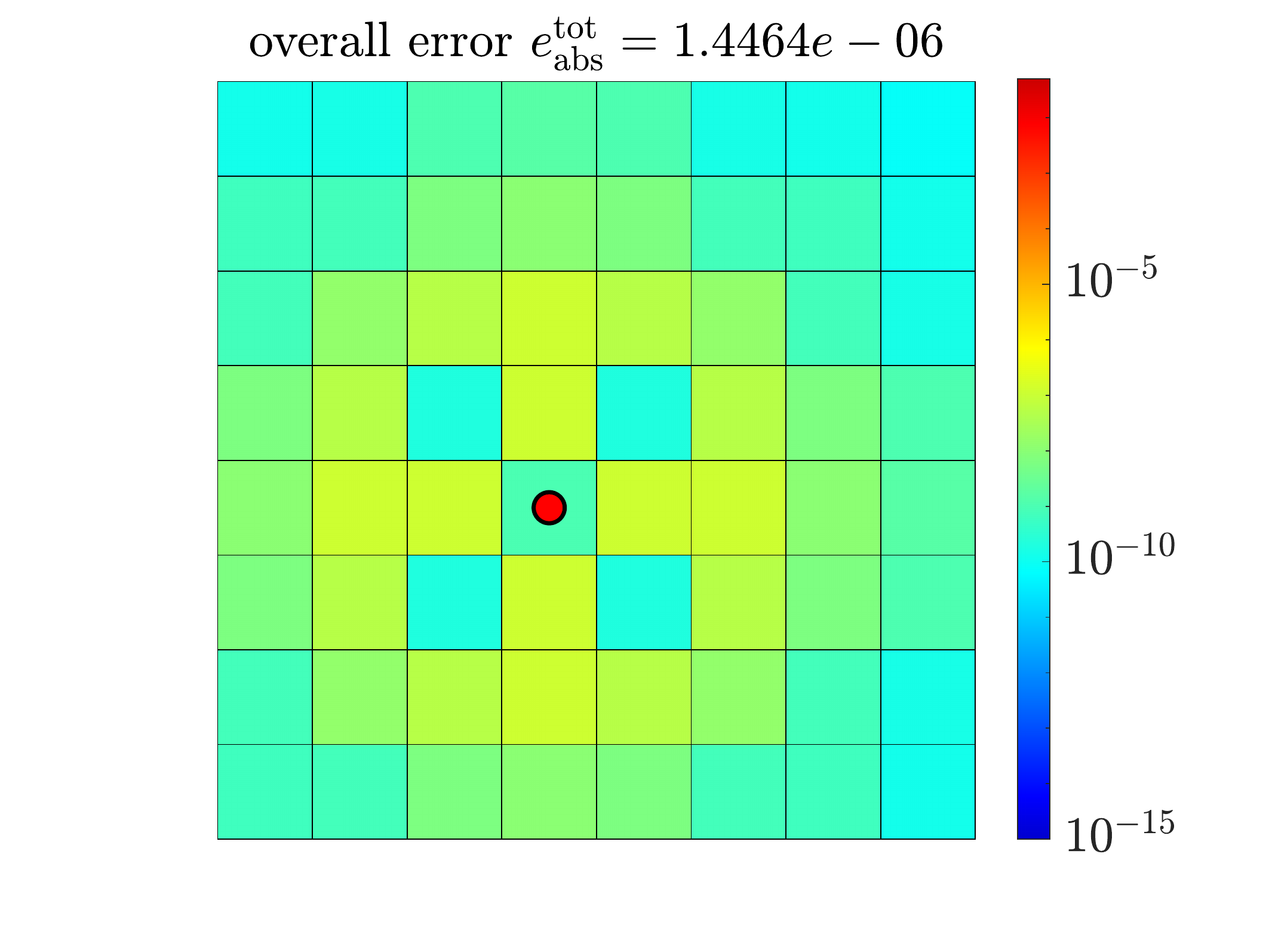}}
	\put(0.61,0.31){\includegraphics[trim = 90 50 130 0, clip, width=.3\linewidth ]{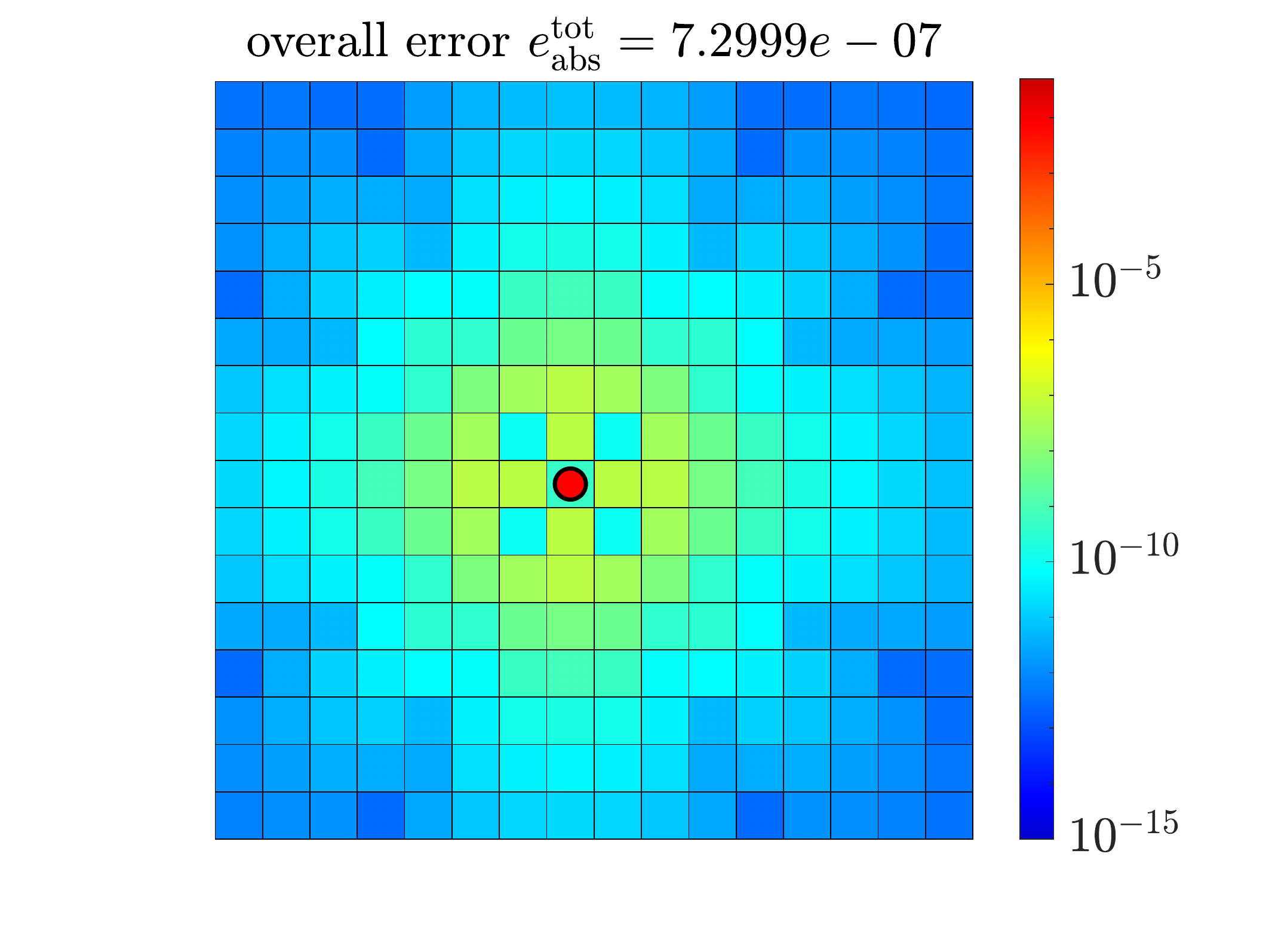}}
 \put(0.01,0){\includegraphics[trim = 90 50 130 0, clip, width=.3\linewidth]{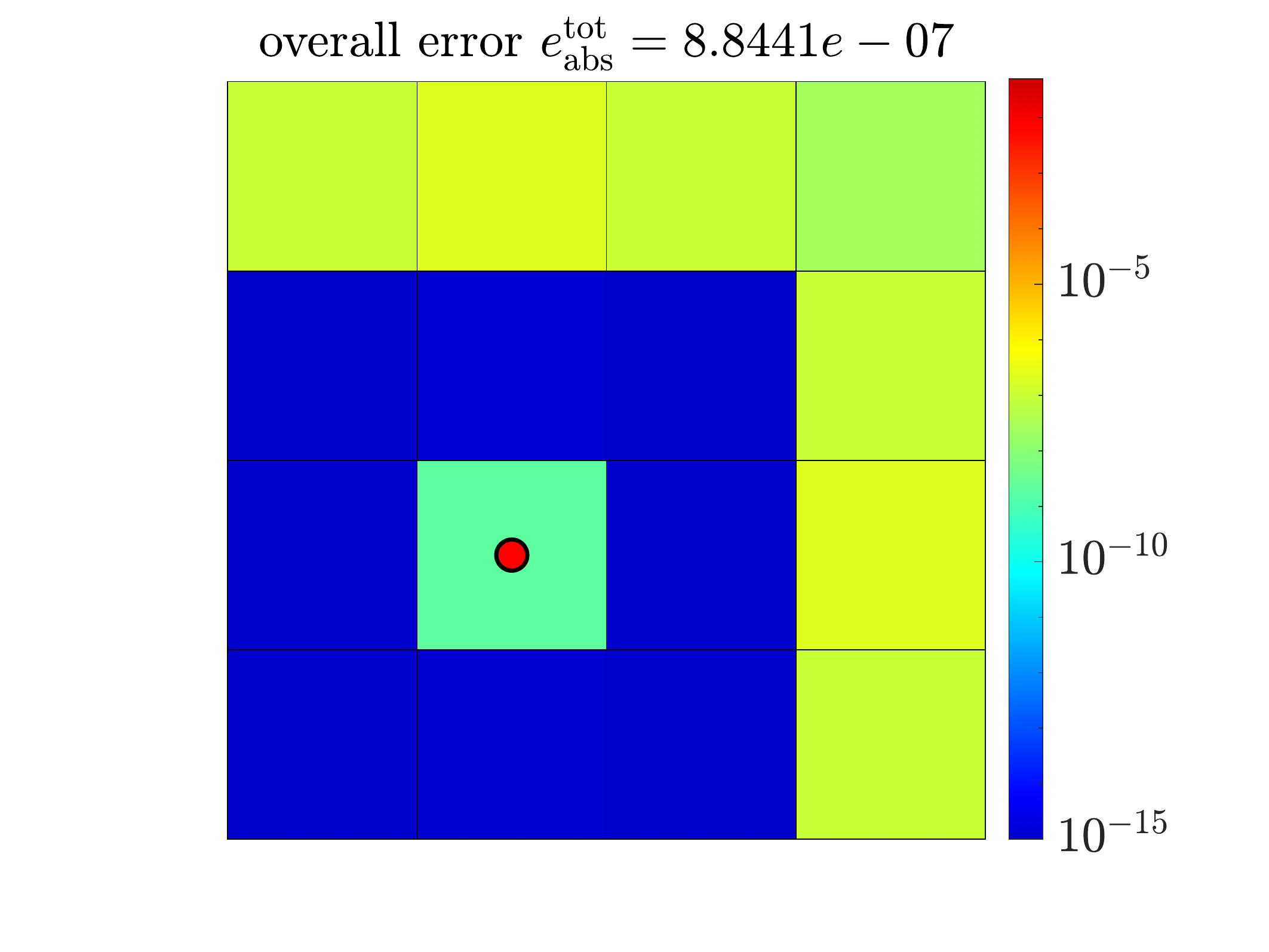}}
 	\put(0.31,0){\includegraphics[trim = 90 50 130 0, clip, width=.3\linewidth ]{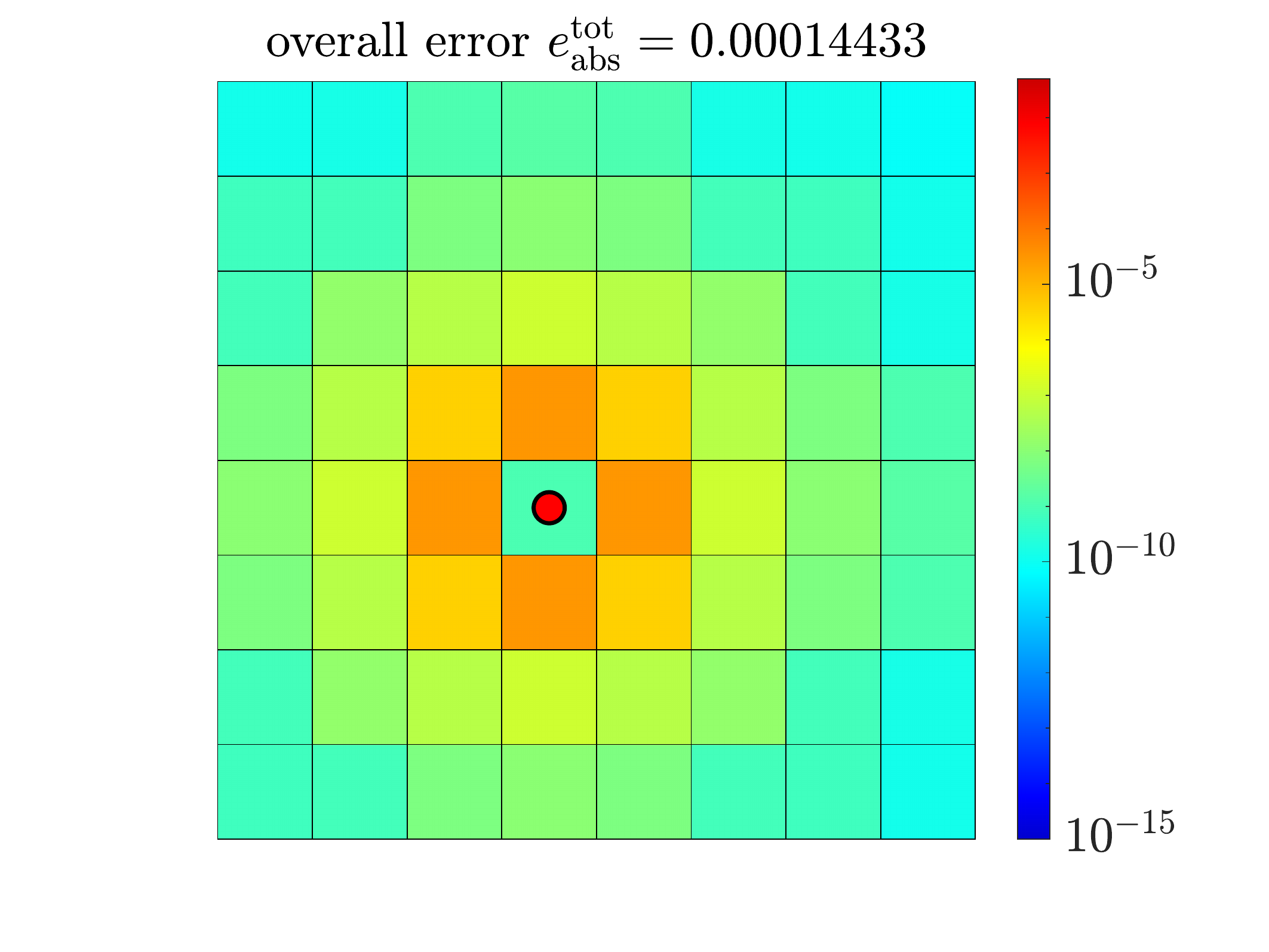}}
	\put(0.61,0){\includegraphics[trim = 90 50 130 0, clip, width=.3\linewidth ]{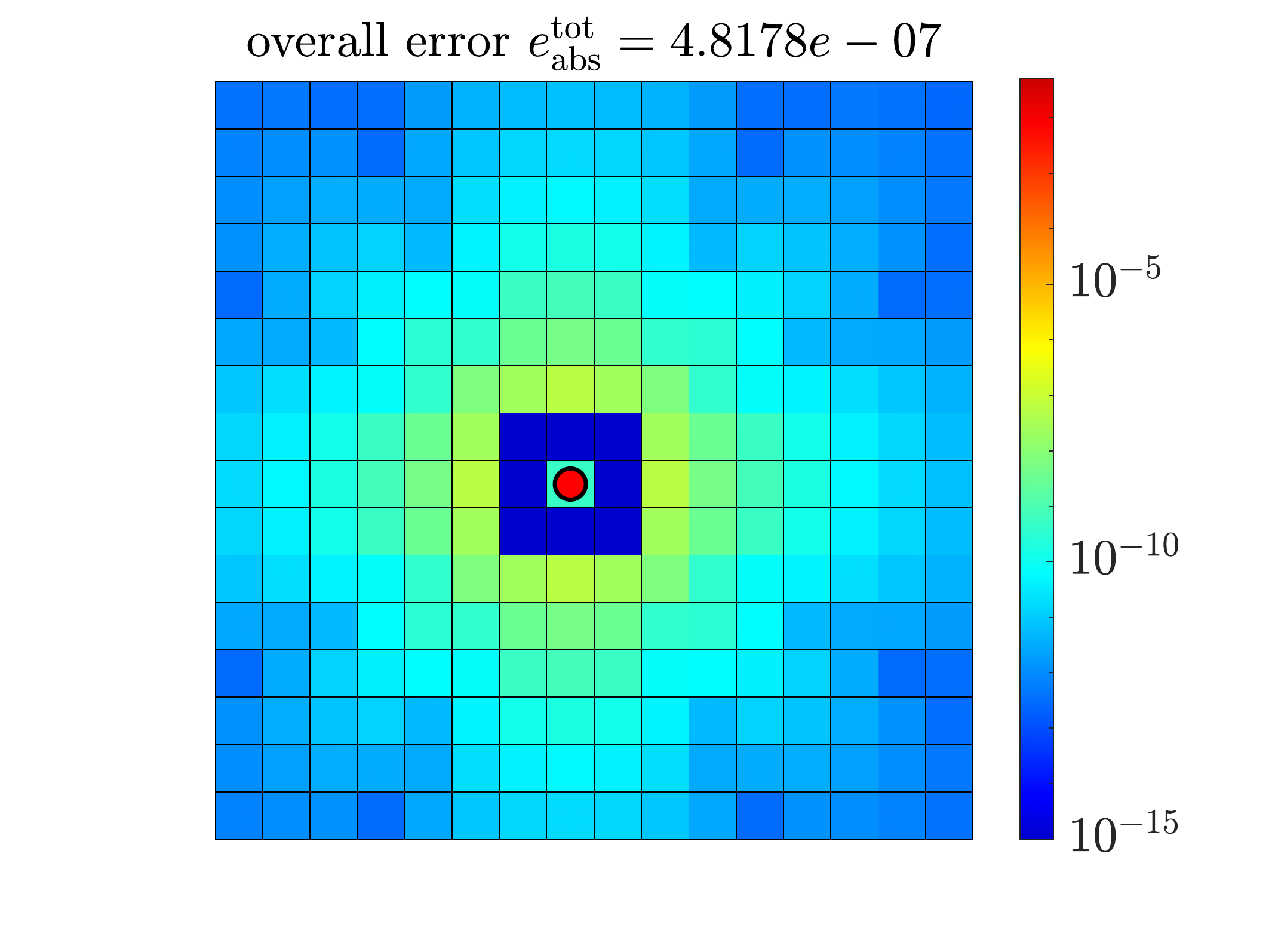}}
	\put(0.91,-0.005){\includegraphics[trim = 475 40 30 0, clip, width=.08\linewidth ]{figures/quad/scheme3/m3_nReg_3_nDuff_6/y153_error2.pdf}}
	\put(0.91,0.305){\includegraphics[trim = 475 40 30 0, clip, width=.08\linewidth ]{figures/quad/scheme3/m3_nReg_3_nDuff_6/y153_error2.pdf}}
	\put(0.91,0.615){\includegraphics[trim = 475 40 30 0, clip, width=.08\linewidth ]{figures/quad/scheme3/m3_nReg_3_nDuff_6/y153_error2.pdf}}
 \put(0.91,0.925){\includegraphics[trim = 475 40 30 0, clip, width=.08\linewidth ]{figures/quad/scheme3/m3_nReg_3_nDuff_6/y153_error2.pdf}}
	\put(0,.93){a.}\put(0,.62){b.}\put(0,.31){c.}\put(0,0){d.}
	\put(0.14,1.25){$\ell=1$}\put(0.43,1.25){$\ell=2$}\put(0.735,1.25){$\ell=3$}
	\vwput{0.96}{0.13}{$e_\mathrm{abs}^e$}\vwput{0.96}{0.435}{$e_\mathrm{abs}^e$}\vwput{0.96}{0.75}{$e_\mathrm{abs}^e$}\vwput{0.96}{1.06}{$e_\mathrm{abs}^e$}
	\end{picture}
\caption{\textit{Hybrid quadrature on a B-spline sheet}: Absolute quadrature error for G~(a.), DG~(b.), DGr~(c.) and for DGw~(d.) considering collocation point $\by_0$ and $n_0=3$.}\label{fig:err_cmp}
\end{figure}
\begin{figure}[h]
\unitlength\linewidth
\begin{picture}(1,.44)
\put(0,0){\includegraphics[trim = 10 0 40 20, clip, width=.5\linewidth]{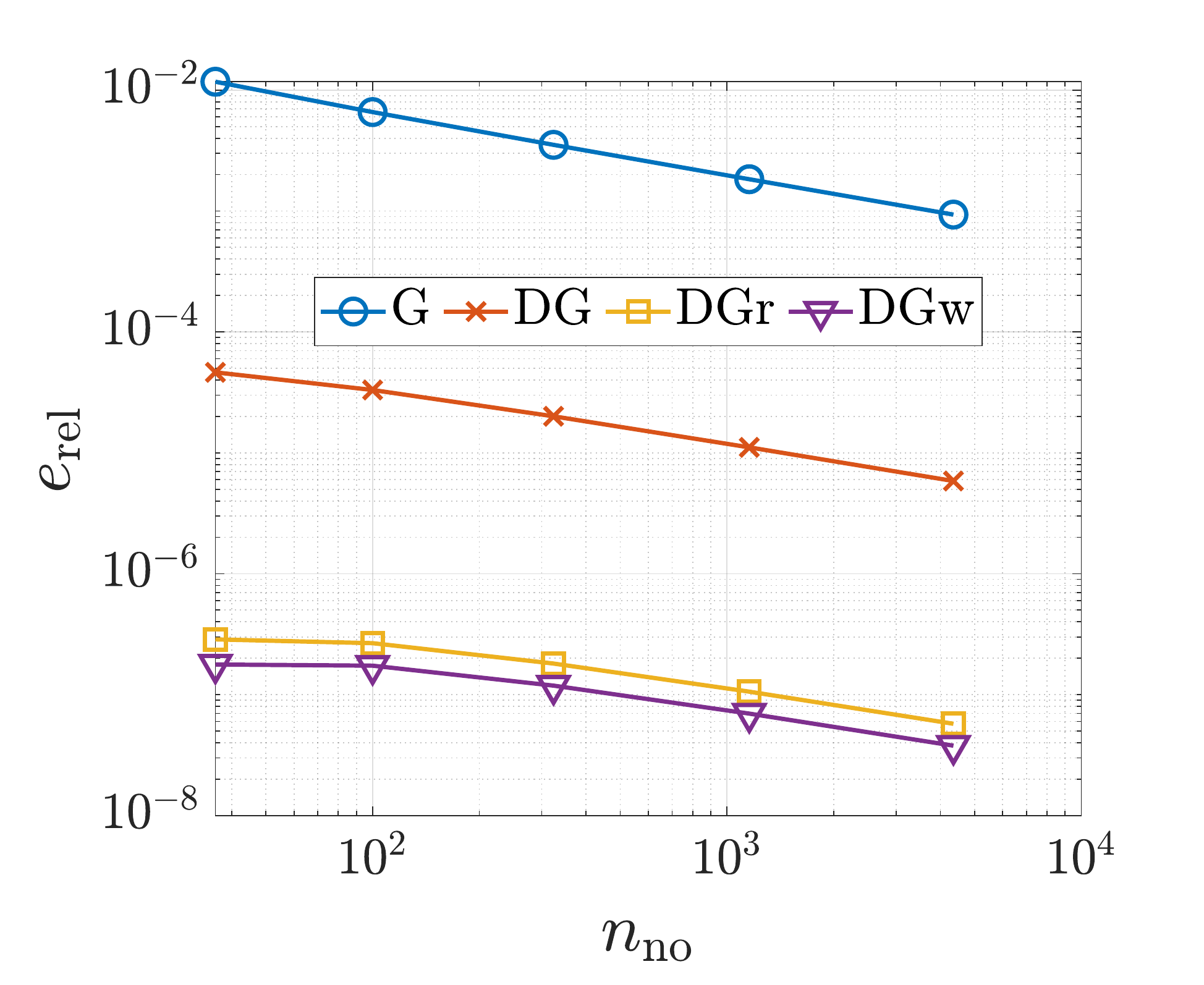}}
\put(0.5,0){\includegraphics[trim = 10 0 40 20, clip, width=.5\linewidth]{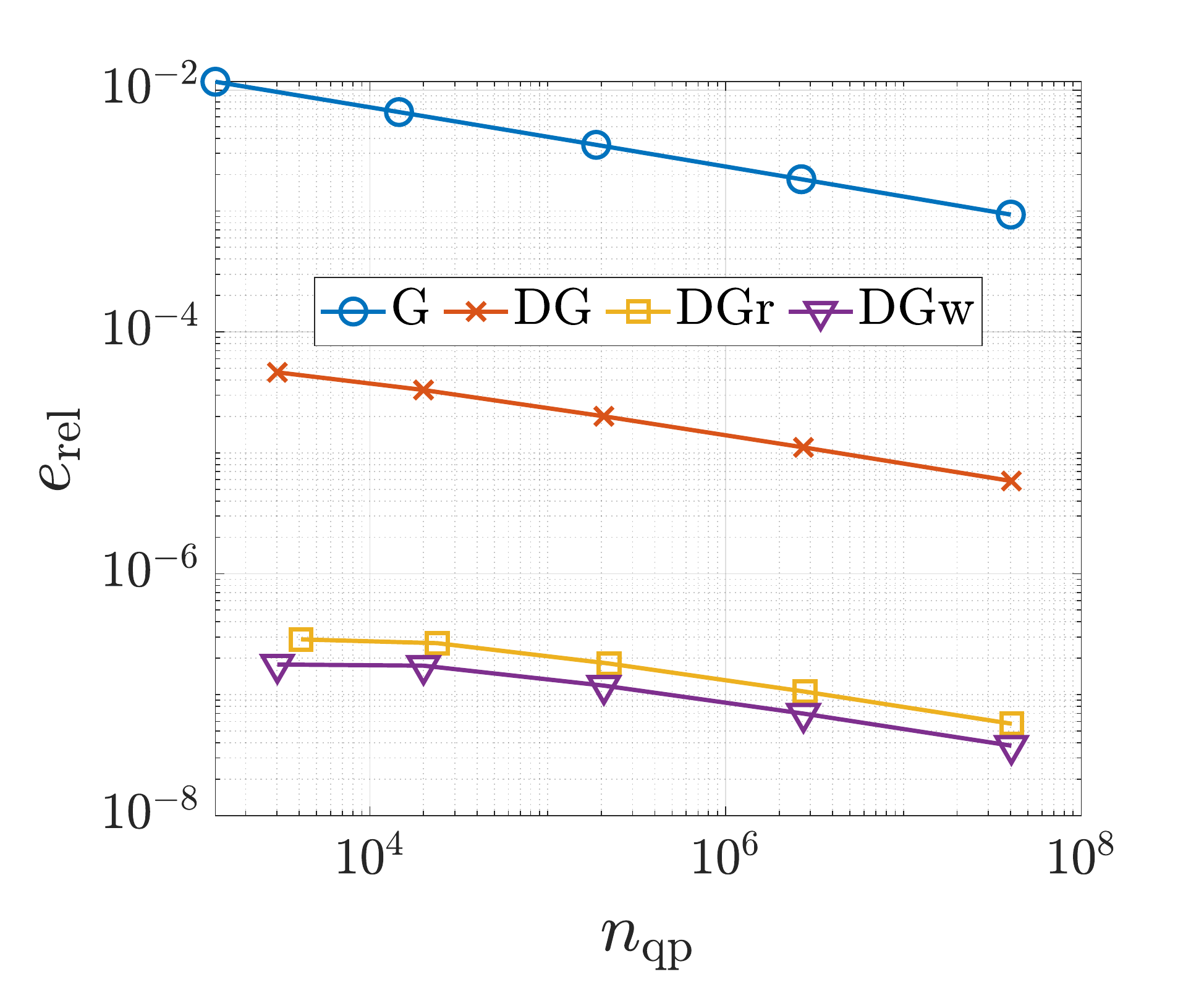}}
	\end{picture}
\caption{\textit{Hybrid quadrature on a B-spline sheet}: Quadrature error for $n_0=3$ vs.~the number of nodes~(a.) and vs.~the total number of quadrature points~(b.) for $\ell=1,\ldots,5$.}\label{fig:convergence}
\end{figure}

\noindent The absolute quadrature error on element $e$ due to evaluating the singular integral\footnote{The singular integral \eqref{eq:quad_int_elem} is representative of the behavior of the integrals in the BIE \eqref{eq:flow_BIE_surf}.}
\begin{equation}\label{eq:quad_int_elem}
\mcalI_e =\int_{\Omega^e} \frac{1}{\|\bx-\by_0\|} \,\mrd a_x
\end{equation}
at collocation point $\by_0$ is defined as
\begin{equation}\label{eq:quad_err0_elem}
e_\mathrm{abs}^{e} = |\mcalI_{e}^h - \mcalI_{e}|~,
\end{equation}
where $\mcalI_{e}^h$ denotes the numerical approximation of integral \eqref{fig:err_cmp}. Fig.~\ref{fig:err_cmp} shows the absolute quadrature error on B-spline discretizations of refinement level $\ell=1,2,3$ with $n_0=3$. It can be seen that the maximum elemental error decreases with increasing $\ell$ and decreases from scheme a.~(G) to d.~(DGw). The total quadrature error for $\by_0$, defined as
\begin{equation}\label{eq:quad_err_elem}
e_\mathrm{abs}^\mathrm{tot} =  |\mcalI^h - \mcalI|
\end{equation}
with the elemental summations $\mcalI= \sum_e \mcalI_e$ and $\mcalI^h= \sum_e \mcalI_e^h$, also decreases with increasing $\ell$ and from scheme a.~to d.~as Fig.~\ref{fig:err_cmp} shows. The only exception occurs for scheme DGw from  $\ell=1$ to $\ell=2$ (see Fig.~\ref{fig:err_cmp}d). The reason for the increasing error is to be found in the most inner layer of classical Gauss-Legendre elements (third layer in total) that yields the highest elemental errors: The sheets for $\ell=2,3$ contain the entire third layer, while the $\ell=1$ sheet only contains half of that layer. The same effect can be observed for DGw, where the overall error decreases only by factor 1.3 from $\ell=1$ to $\ell=2$. For all other cases, $e_\mathrm{abs}^0 $ is halved with each refinement step, i.e.~linear convergence with convergence rate $\mu=1/2$.
\\\\Fig.~\ref{fig:convergence} depicts the mean relative quadrature error 
\begin{equation}\label{eq:quad_mean_err}
e_\mathrm{rel}:= \frac{1}{n_\mathrm{no}}\,\sum_{A=1}^{n_\mathrm{no}}{ \frac{\|\mcalI_A^h - \mcalI_A\|}{\mcalI_A}}~,
\end{equation}
on the investigated B-spline. Fig.~\ref{fig:convergence}a shows linear convergence with rate $\mu=1/2$ for each quadrature scheme. However, the error for quadrature scheme DG is about two orders of magnitude lower than for scheme G, while Duffy-Gauss quadrature scheme with progressive refinement~(DGr) reduces the error by two further orders of magnitude. Duffy-Gauss quadrature with adjusted weights (DGw) is even more accurate and yields the lowest error among the investigated schemes.
\\Fig.~\ref{fig:convergence}b shows the mean relative error vs.~the total number of quadrature points that is given by
\begin{equation}\label{eq:nqp}
n_\mathrm{qp} := \sum_{A=1}^{n_\mathrm{no}} \sum_{e=1}^{n_\mathrm{el}} n_\mathrm{qp}^{A,\,e}~,
\end{equation}
where $n_\mathrm{qp}^{A,\,e}$ denotes the number of quadrature points for collocation point $\by_A$ on element $e$. The figure shows the same convergence behavior as noted in Fig.~\ref{fig:convergence}a. It also shows that DGr uses slightly more quadrature points per refinement level than schemes DGw, DG and G. However, the gain in accuracy clearly outweighs the additional quadrature effort as Fig.~\ref{fig:convergence}b shows. Therefore DGr and especially DGw are the most efficient quadrature schemes for singular integral approximation on biquadratic B-spline sheets.
\\Fig.~\ref{fig:convergence} further shows that increasing the number of quadrature points or adjusting the quadrature weights in one ring of elements (DGr or DGw) leads to a beneficial jump in accuracy (c.f.~DG). Using further rings, additional, albeit smaller jumps in accuracy can be expected.

\subsection{Hybrid quadrature on a NURBS sphere}\label{sec:hybrid_sphere}
\begin{figure}[h]
\unitlength\linewidth
\begin{picture}(1,0.5)	
	\put(0,0.28){\includegraphics[trim = 20 20 0 0, clip, width=.23\linewidth]{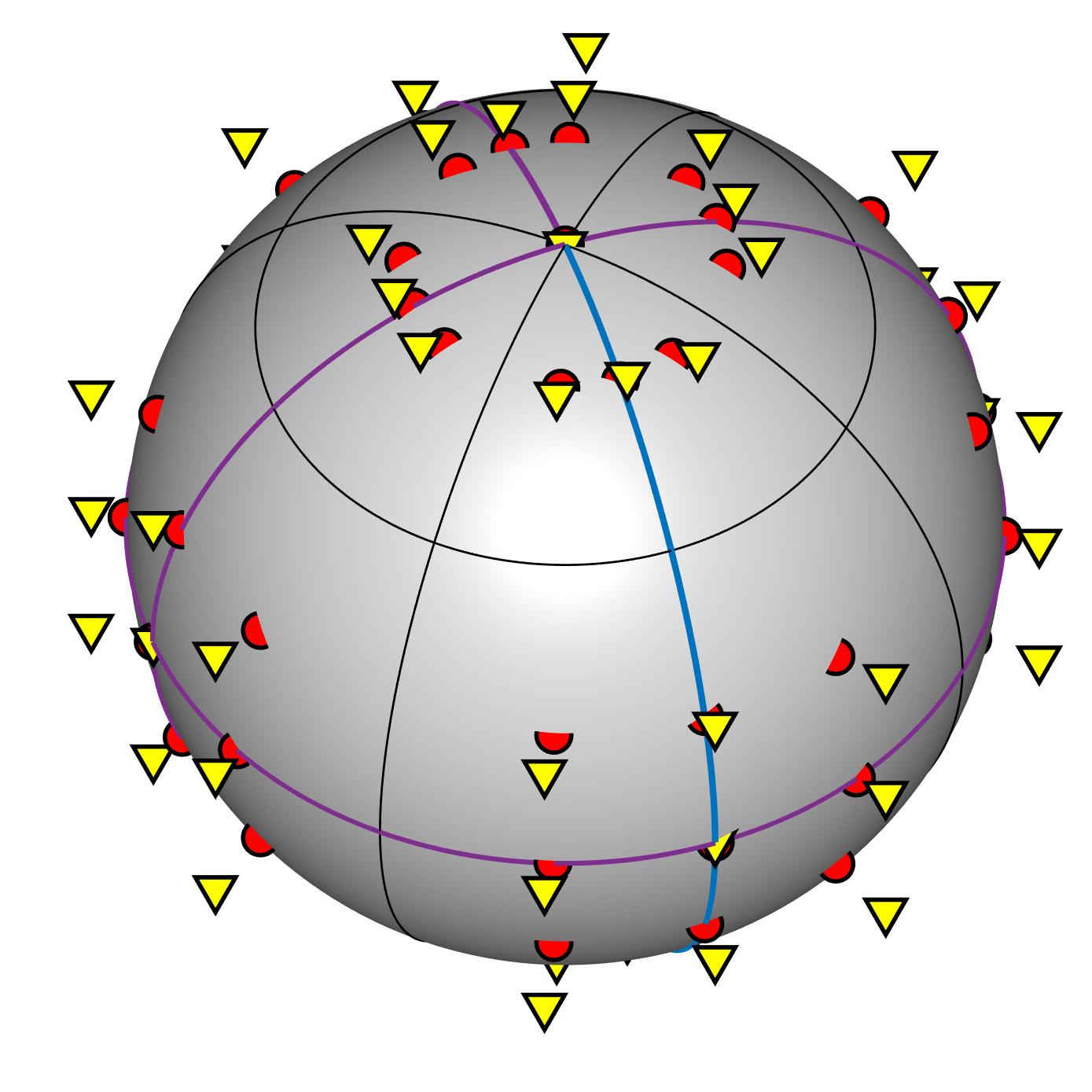}}
	\put(0.23,0.28){\includegraphics[trim = 20 20 0 0, clip, width=.23\linewidth]{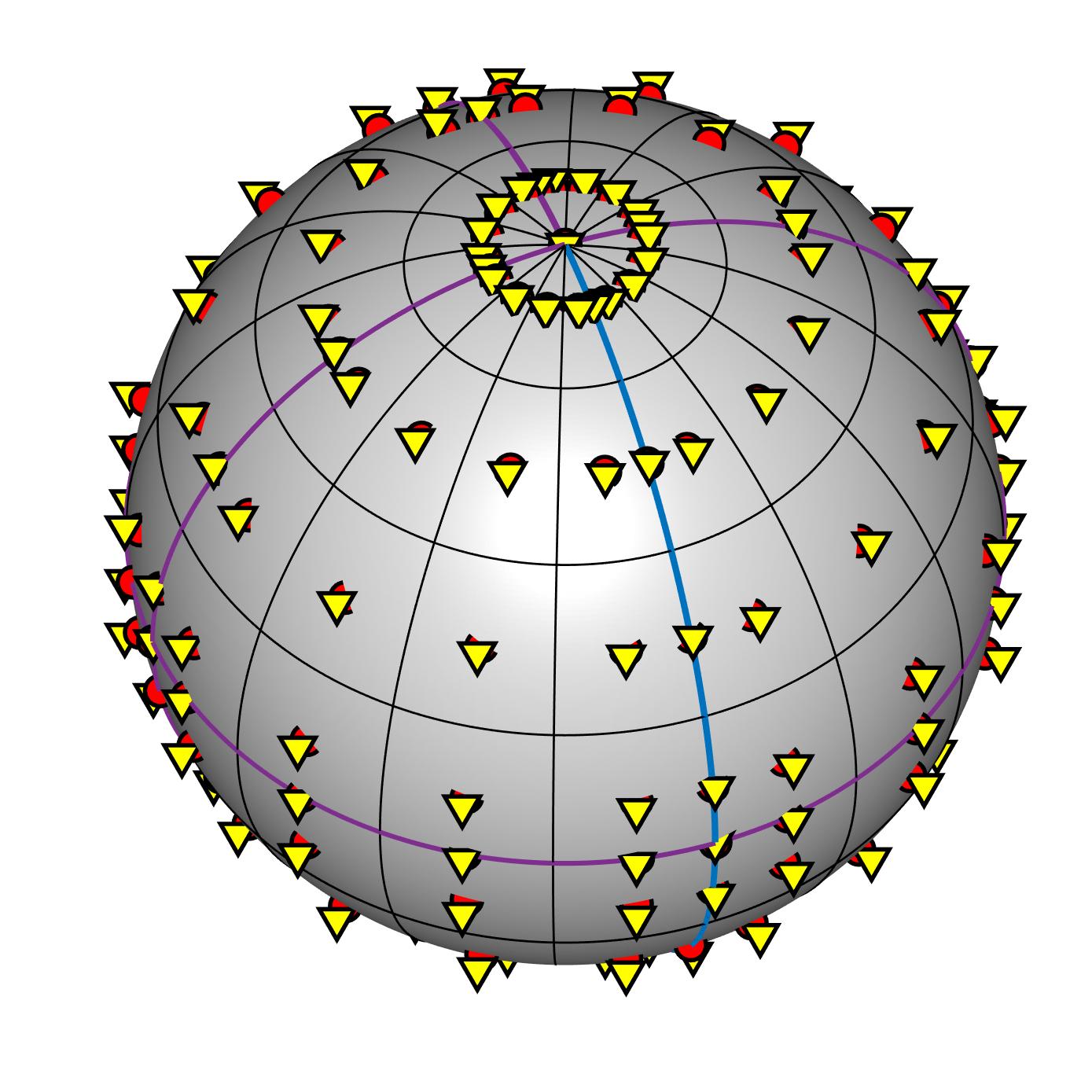}}
	\put(0,0.06){\includegraphics[trim = 20 20 0 0, clip, width=.23\linewidth]{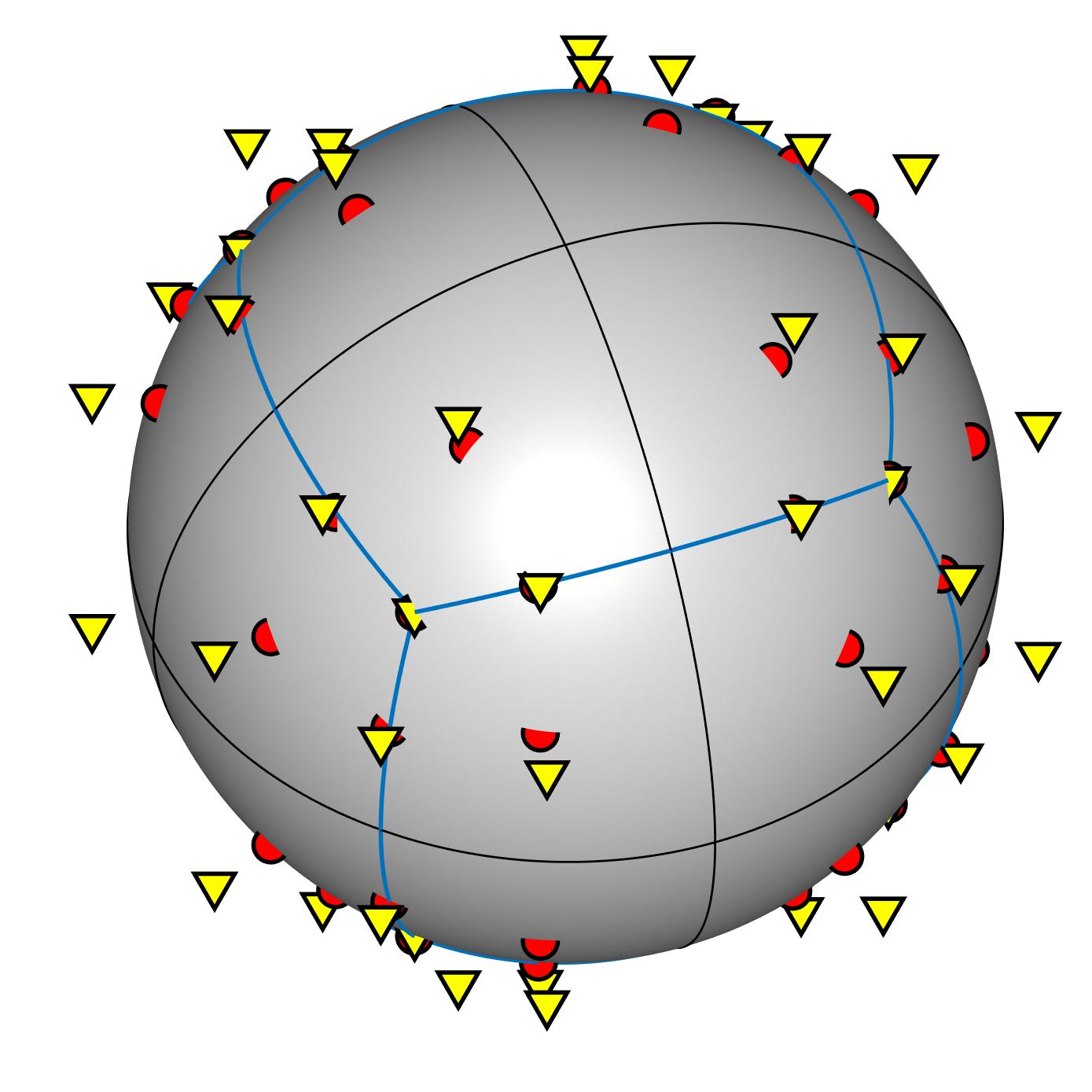}}
	\put(0.23,0.06){\includegraphics[trim = 20 20 0 0, clip, width=.23\linewidth]{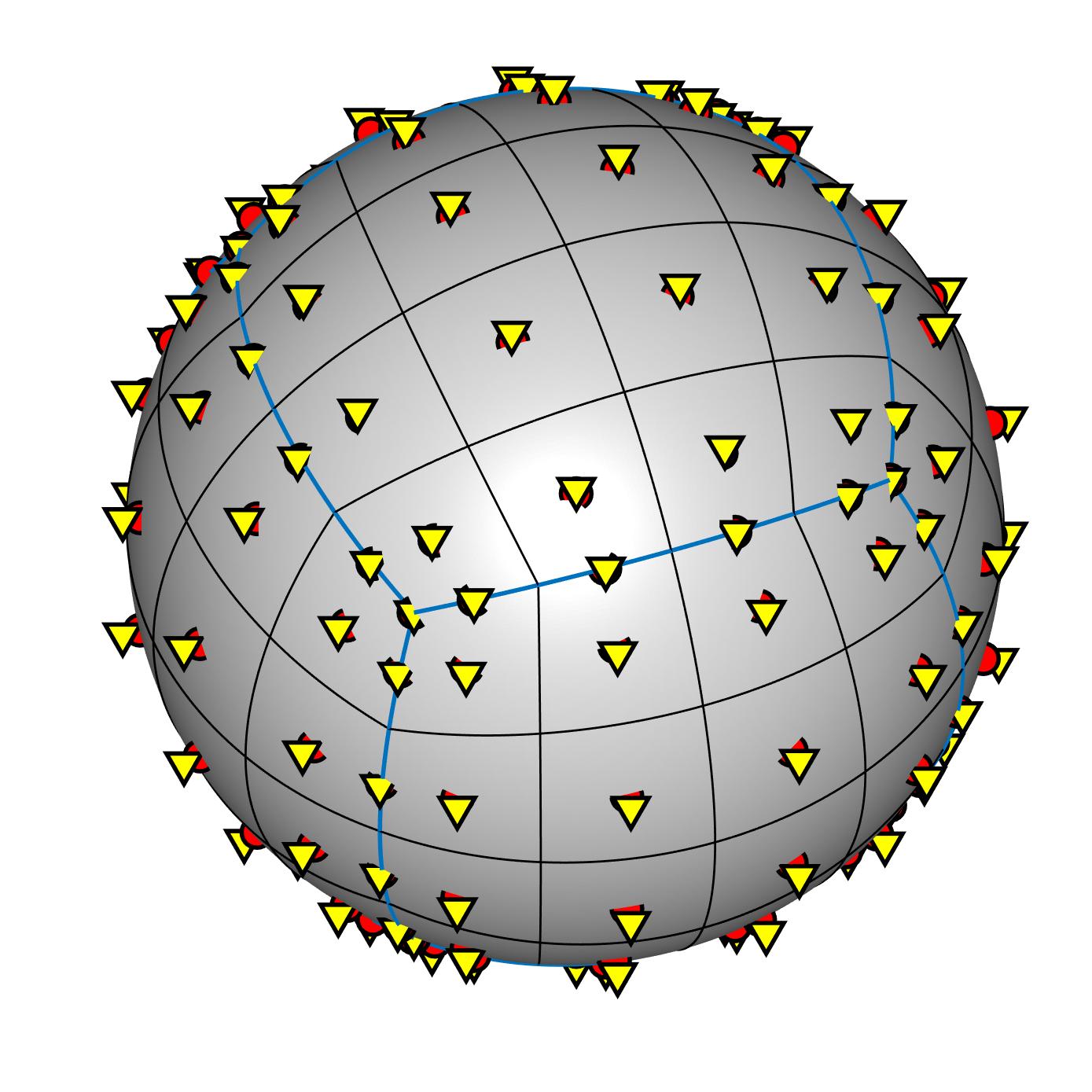}}
	\put(0,0){\includegraphics[trim = 120 45 40 10, clip, width=.24\linewidth]{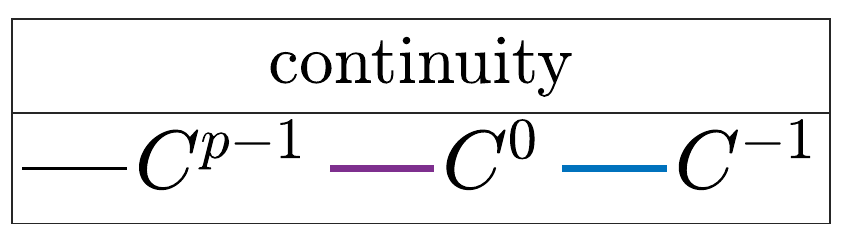}}	
	\put(0.25,0){\includegraphics[trim = 70 34 30 10, clip, width=.2\linewidth]{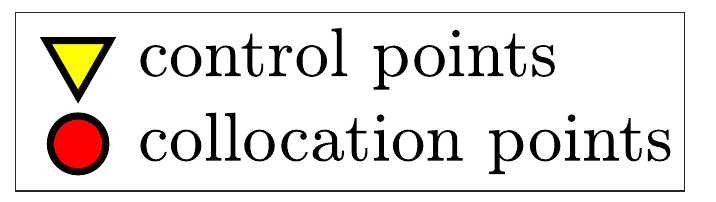}}	
	\put(0.46,0){\includegraphics[trim = 0 0 35 35, clip, width=.54\linewidth]{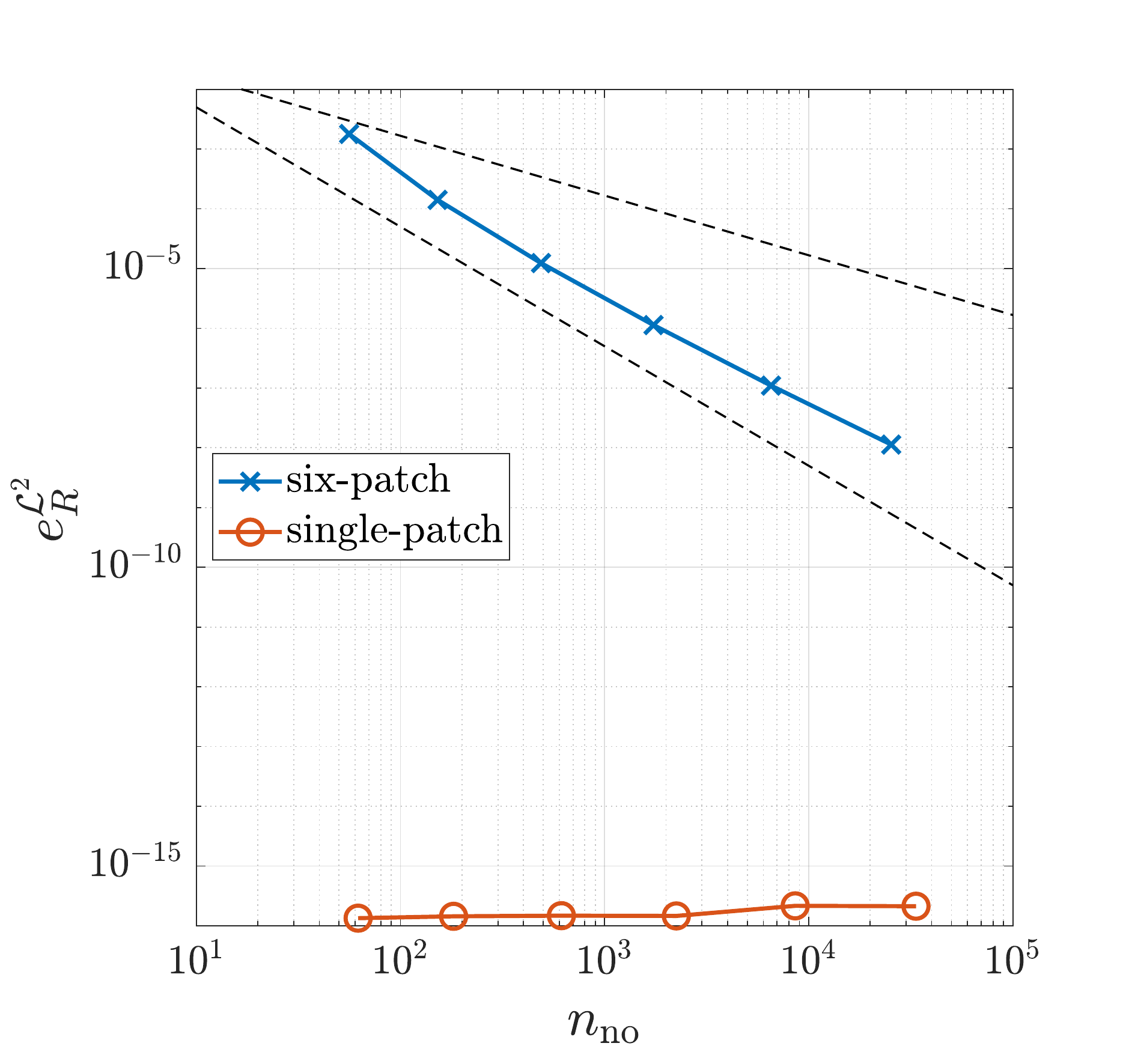}}
	\put(0.01,0.3){\small a.} \put(0.24,0.3){\small b.} \put(0.01,0.08){\small c.} \put(0.24,0.08){\small d.}
	\put(0.48,0.01){\small e.} 
	\put(0.87,0.275){\rotatebox{-30}{\small$\propto 1/{n_\mathrm{no}}^2$}} 
	\put(0.87,0.425){\rotatebox{-15}{\small$\propto 1/n_\mathrm{no}$}} 
\end{picture}
\caption{\textit{Hybrid quadrature on a NURBS sphere}: Single-patch NURBS sphere of discretization level $\ell=1$~(a) and $\ell=2$~(b) and six-patch NURBS sphere of level $\ell=1$~(c) and $\ell=2$~(d). Collocation points are located on surface $\mcalS^h=\mcalS$, while control points are not. e.~$\mcalL^2$ norm of the radius error for $\ell=1,\ldots,6$.} \label{fig:quad_sphere_discretization}
\end{figure}
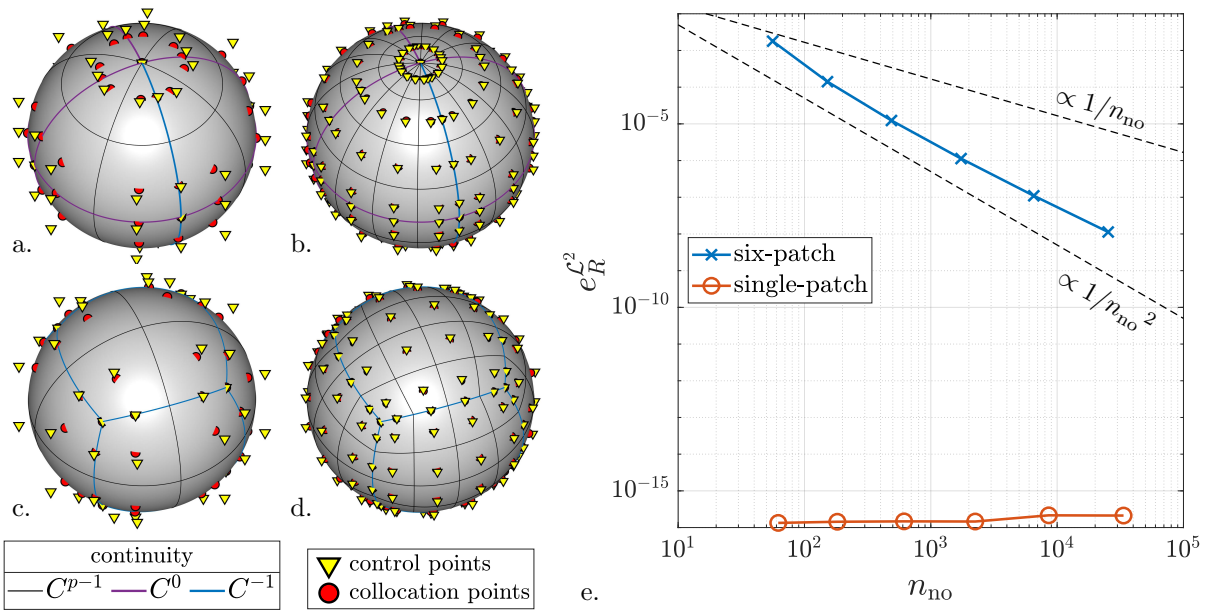

Next, the four hybrid quadrature schemes are applied to spherical surfaces to investigate their suitability for curved NURBS surfaces. The procedure advocated here is visualized and investigated by considering biquadratic NURBS spheres, but is applicable to any closed and open surface defined by a mapping from a planar parameter domain. It should be further noted that any genus zero surface can be mapped conformally onto a sphere \citep{Gu04}, so this example has far reaching applications.
Two different spherical discretizations with biquadratic NURBS elements are considered here: The single-patch NURBS sphere from Fig.~\ref{fig:quad_sphere_discretization}a and b is investigated in Sec.~\ref{sec:hybrid_sphere_single}, while the six-patch NURBS sphere from Fig.~\ref{fig:quad_sphere_discretization}c and d is investigated in Sec.~\ref{sec:hybrid_sphere_six}. Fig.~\ref{fig:quad_sphere_discretization}e shows that the single-patch discretization is exactly spherical independent of $\ell$, whereas the six-patch sphere is only approximately spherical. For the six-patch sphere, the $\mcalL^2$~norm of the radius error\footnote{The $\mcalL^2$~norm of the radius error $e_\mathrm{R} = (R^h-R)/R$ is defined as $e_R^{\mathcal L^2} := \sqrt{ \int_\mcalS  e_R^2 \,\mrd a / A_\mcalS}$, where $A_\mcalS=4 \pi R^2$ denotes the surface area of a sphere.} shows almost quadratic convergence w.r.t.~the number of control points ($e_R^{\mathcal L^2}  \propto 1/{n_\mathrm{no}}^2$).
\\The number of elements and control points are given in Table~\ref{tab:NURBS} for the single-patch and the six patch sphere, both of refinement level $\ell=1,\ldots,6$. The corresponding number of quadrature points are given in Table~\ref{tab:n_qp} for the four hybrid quadrature schemes with quadrature density $n_0=3$. Those numbers also apply to the numerical examples in Sec.~\ref{sec:examples}, which consider spherical and spheroidal geometries.
\begin{table}[h]
\centering
\begin{tabular}{|c||c|c|c|c|c|c|} 
 \hline
  $\ell$ & 1 & 2 & 3& 4 & 5 & 6 \\ 
 \hline
  \multirow{2}{*}{$\nel$} &  \tmred{32} & \tmred{128} &   \tmred{512} & \tmred{2,048} &  \tmred{8,192} &  \tmred{32,768} \\  
  & \tmblue{24} &  \tmblue{96} & \tmblue{384} &   \tmblue{1,536} & \tmblue{6,144} &  \tmblue{24,576}  \\\hline
  \multirow{2}{*}{$\nno$} &  \tmred{62} & \tmred{182} & \tmred{614} & \tmred{2,246} &  \tmred{8,582} &  \tmred{33,542} \\
   & \tmblue{56} &  \tmblue{152} & \tmblue{488} &   \tmblue{1,736} & \tmblue{6,536} &  \tmblue{25,352}	 \\
  \hline
\end{tabular}
\caption{\textit{Hybrid quadrature on a NURBS sphere}: Number of elements $\nel$ and control points $\nno$ required for NURBS spheres of discretization level $\ell=1,\ldots,6$ (red font: \tmred{single-patch sphere}, blue font: \tmblue{six-patch sphere}). The number of collocation points is also $\nno$.} \label{tab:NURBS}
\end{table}
\begin{table}[h]
%
%
\centering
\begin{tabular}{|c||c|c|c|c|c|c|} 
 \hline 
   \multirow{2}{*}{$\bar n_\mathrm{qp}$} & \multicolumn{6}{c|}{$\ell$} \\ \cline{2-7}
   & 1 & 2 & 3& 4 & 5 & 6 \\\hhline{|=#=|=|=|=|=|=|}
  \multirow{2}{*}{G} & \tmred{297} &  \tmred{1,161} & \tmred{4,616} &   \tmred{18,440} & \tmred{73,735} &  \tmred{294,919}\\
	                 & \tmblue{225} &  \tmblue{873} & \tmblue{3,464} &   \tmblue{13,832} & \tmblue{55,303} &  \tmblue{221,191}\\\hline
  \multirow{2}{*}{DG} & \tmred{569} &  \tmred{1,436} & \tmred{4,891} &   \tmred{18,714} & \tmred{74,008} &  \tmred{295,192}\\
  					& \tmblue{478} &  \tmblue{1,135} & \tmblue{3,732} &   \tmblue{14,102} & \tmblue{55,574} &  \tmblue{221,463} \\\hline
  \multirow{2}{*}{DGr}  & \tmred{820} &  \tmred{1,723} & \tmred{5,197} &   \tmred{19,029} & \tmred{74,328} &  \tmred{295,513} \\
  						& \tmblue{698} &  \tmblue{1,361} & \tmblue{3,956} &   \tmblue{14,323} & \tmblue{55,793} &  \tmblue{221,680}\\\hline
  \multirow{2}{*}{DGw}     & \tmred{639} &  \tmred{1,524} & \tmred{4,989} &   \tmred{18,816} & \tmred{74,114} &  \tmred{295,298} \\
  & \tmblue{490} &  \tmblue{1,140} & \tmblue{3,733} &   \tmblue{14,101} & \tmblue{55,574} &  \tmblue{221,463} \\\hline
\end{tabular}
\caption{\textit{Hybrid quadrature on a NURBS sphere}: The mean number of quadrature points $\bar n_\mathrm{qp}= n_\mathrm{qp}/n_\mathrm{no}$ (rounded) for NURBS spheres of discretization level $\ell=1,\ldots,6$ (red font:~\tmred{single-patch sphere}, blue font:~\tmblue{six-patch sphere}) and quadrature density $n_0=3$.} \label{tab:n_qp}
\end{table}
\\\\The accuracy of the quadrature schemes is investigated by approximating the singular integrals of the BIE \eqref{eq:flow_BIE_surf}. The quadrature error for collocation point $\by_A$ is defined with respect to identity \eqref{eq:identity_SL} as
\begin{equation}\label{eq:err_SL}
 e_\mathrm{SL}^{\by_A}:=  \sqrt{ \sum_{i=1}^{3}\left[\left(\int_{\mcalS} G_{ij}(\br_A)\, n_j(\bx) \, \mathrm da_x\right)^2\right] } ~,
\end{equation}
and with respect to identity \eqref{eq:identity_DL} as
\begin{equation}\label{eq:err_DL}
 e_\mathrm{DL}^{\by_A}:=  \sqrt{ \sum_{i=1}^{3}\sum_{j=1}^{3}\left[\left(\int_{\mcalS} \frac{1}{4\pi}\, T_{ijk}(\br_A)\, n_k(\bx) \, \mathrm da_x + \delta_{ij}\right)^2\right] } ~,
\end{equation}
where $\br_A:=\bx-\by_A$.
The mean quadrature errors with respect to identities \eqref{eq:identity_SL} and \eqref{eq:identity_DL} then follow as 
\begin{equation}\label{eq:err_SL2}
 e_\mathrm{SL}:=  \frac{1}{n_A}\,\sum_{A=1}^{n_A} e_\mathrm{SL}^{\by_A} ~, \hh \mathrm{and}\hh e_\mathrm{DL} =  \frac{1}{n_A}\,\sum_{A=1}^{n_A} e_\mathrm{DL}^{\by_A} ~.
\end{equation}

\subsubsection{Single-patch NURBS sphere}\label{sec:hybrid_sphere_single}
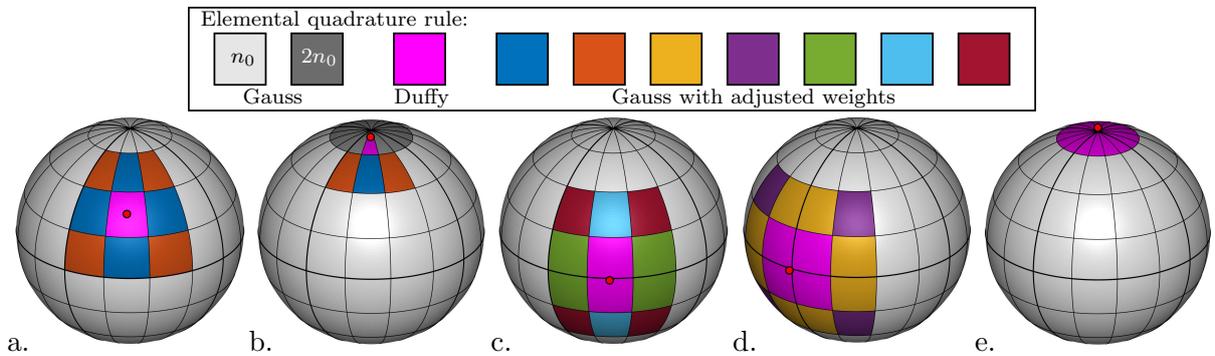
\begin{figure}[h]
\unitlength\linewidth
\begin{picture}(1,.29)
	\put(0,0){\includegraphics[trim = 167 100 133 47, clip, width=.2\linewidth ]{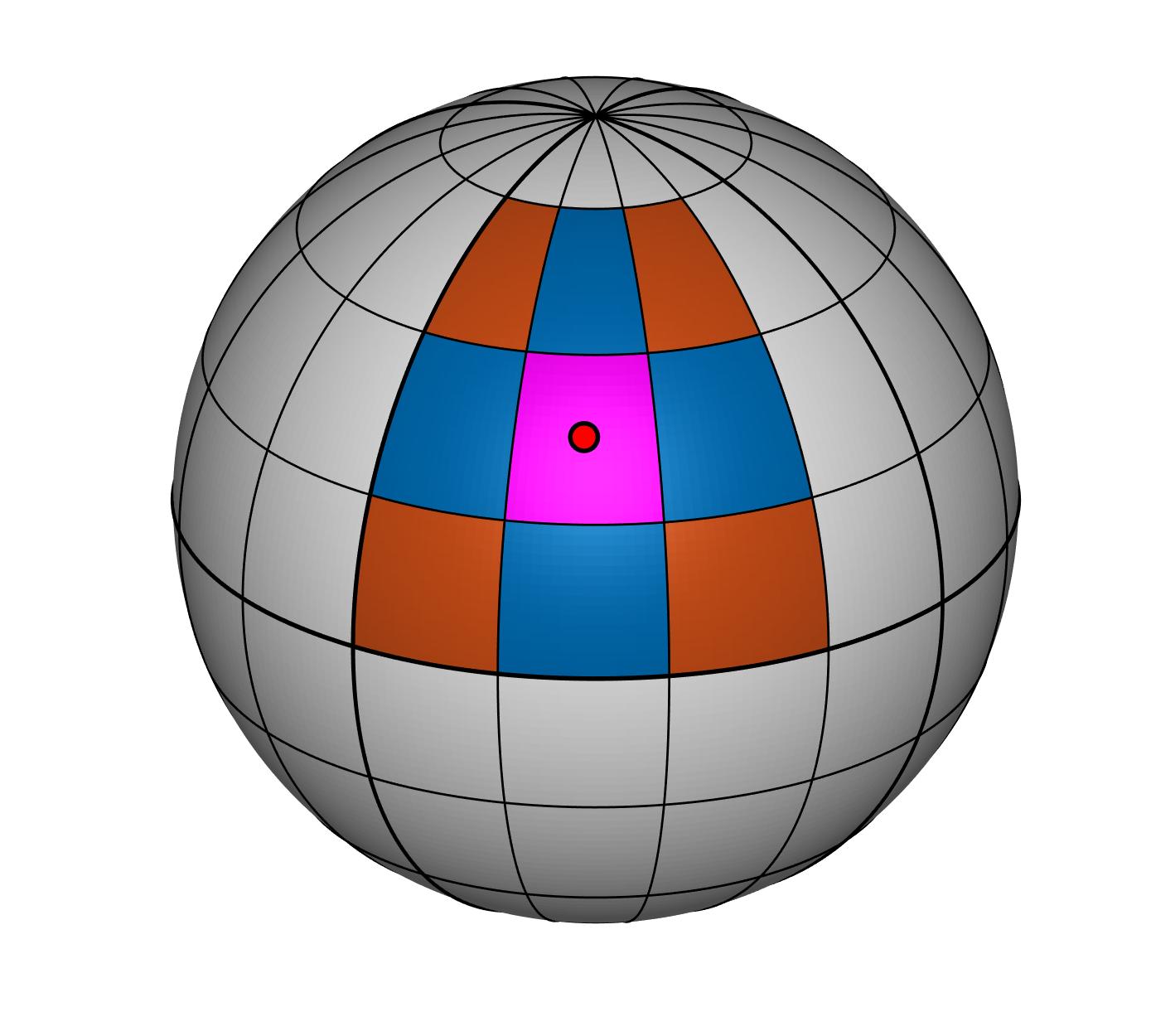}}	
	\put(0.2,0){\includegraphics[trim = 167 100 133 47, clip, width=.2\linewidth ]{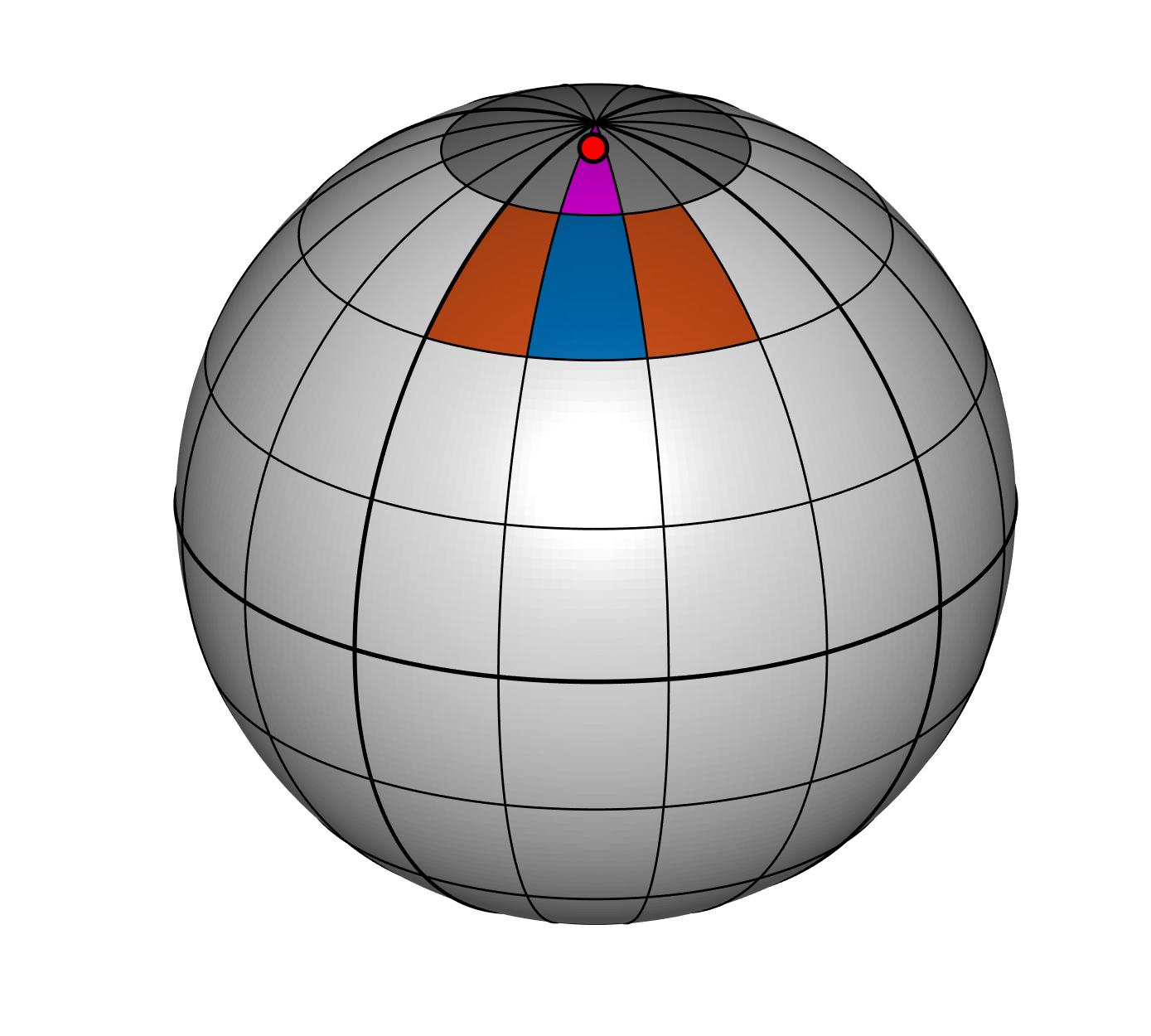}}	
	\put(0.4,0){\includegraphics[trim = 167 100 133 47, clip, width=.2\linewidth ]{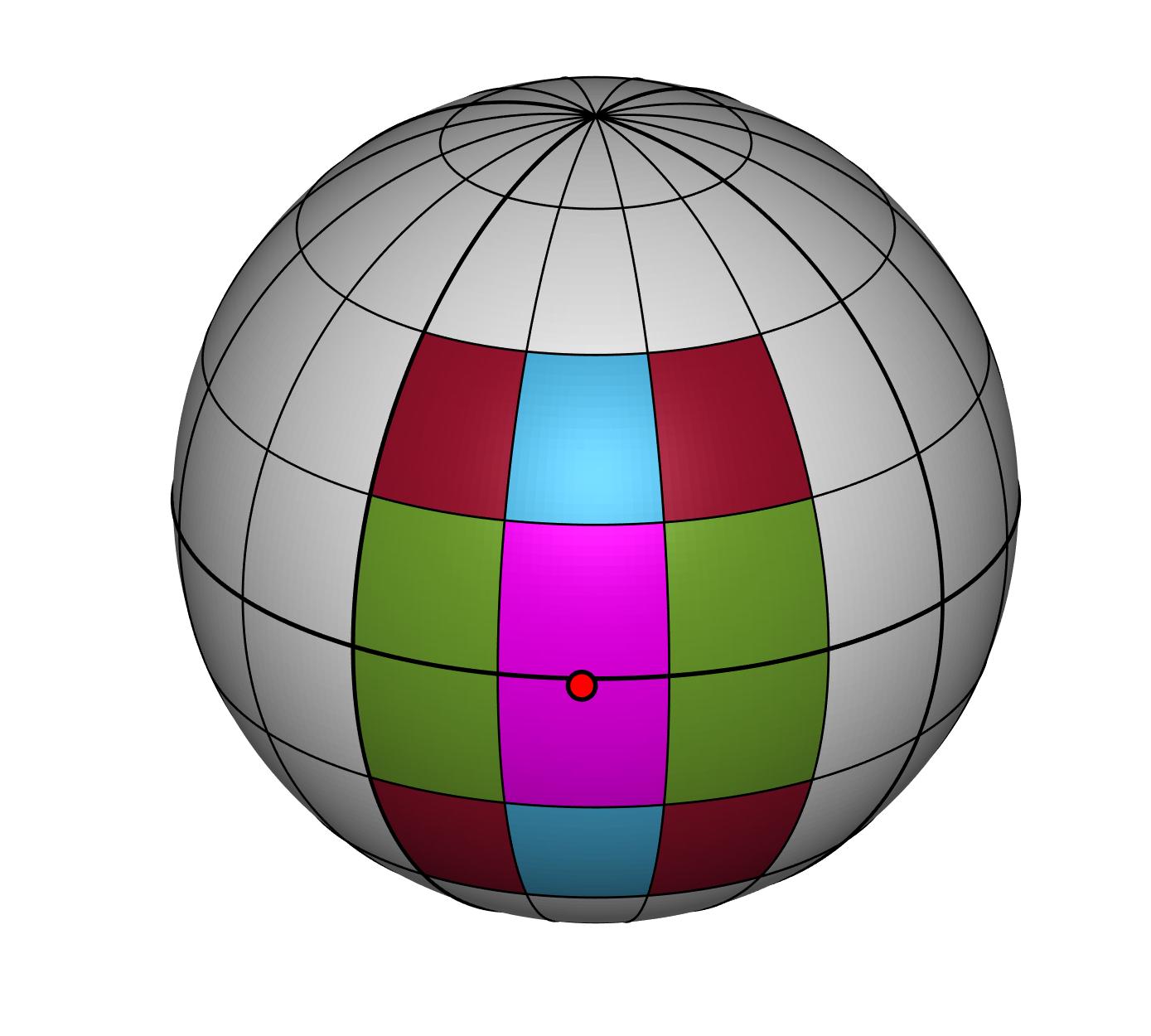}}	
	\put(0.6,0){\includegraphics[trim = 167 100 133 47, clip, width=.2\linewidth ]{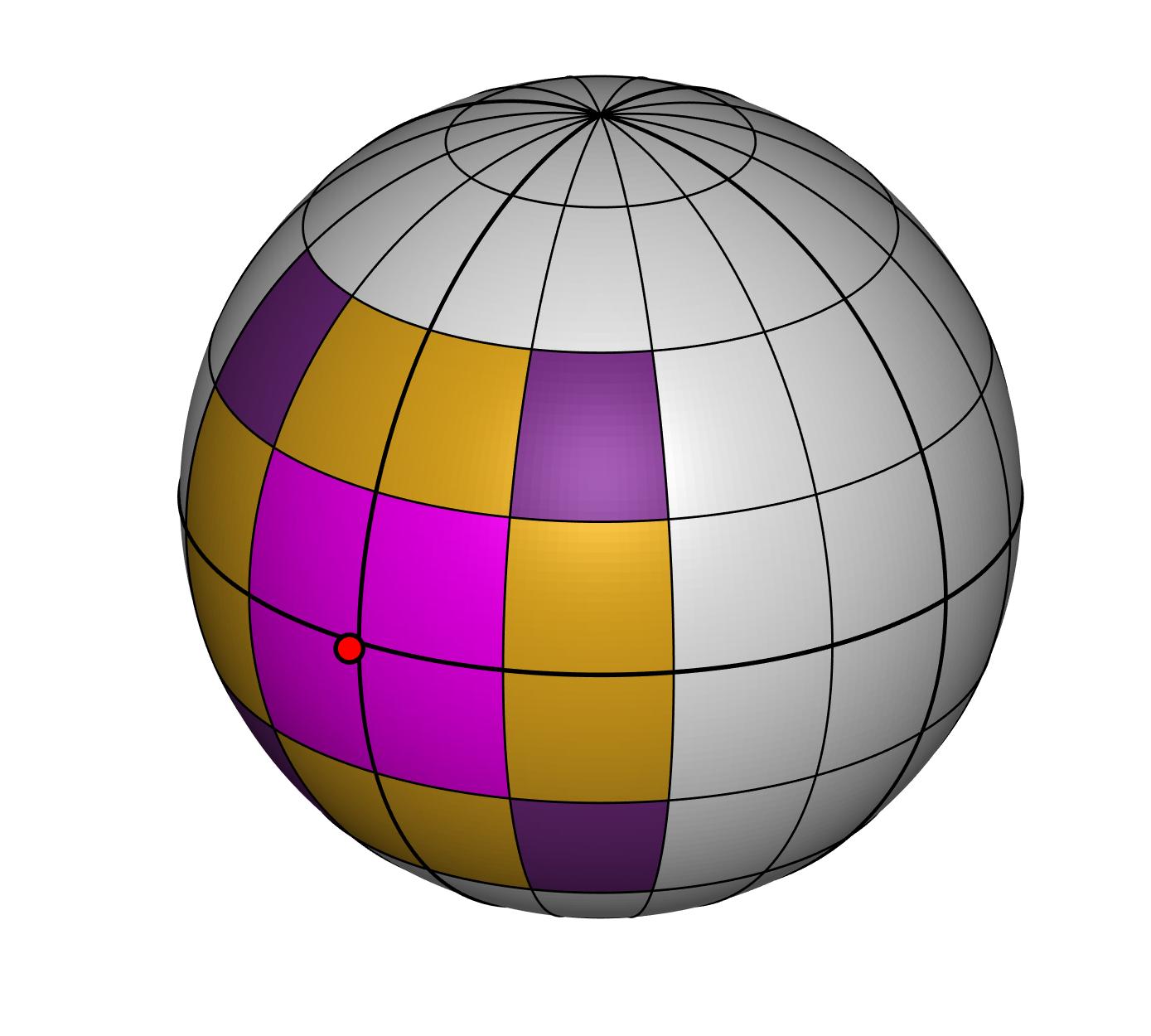}}	
	\put(0.8,0){\includegraphics[trim = 167 100 133 47, clip, width=.2\linewidth ]{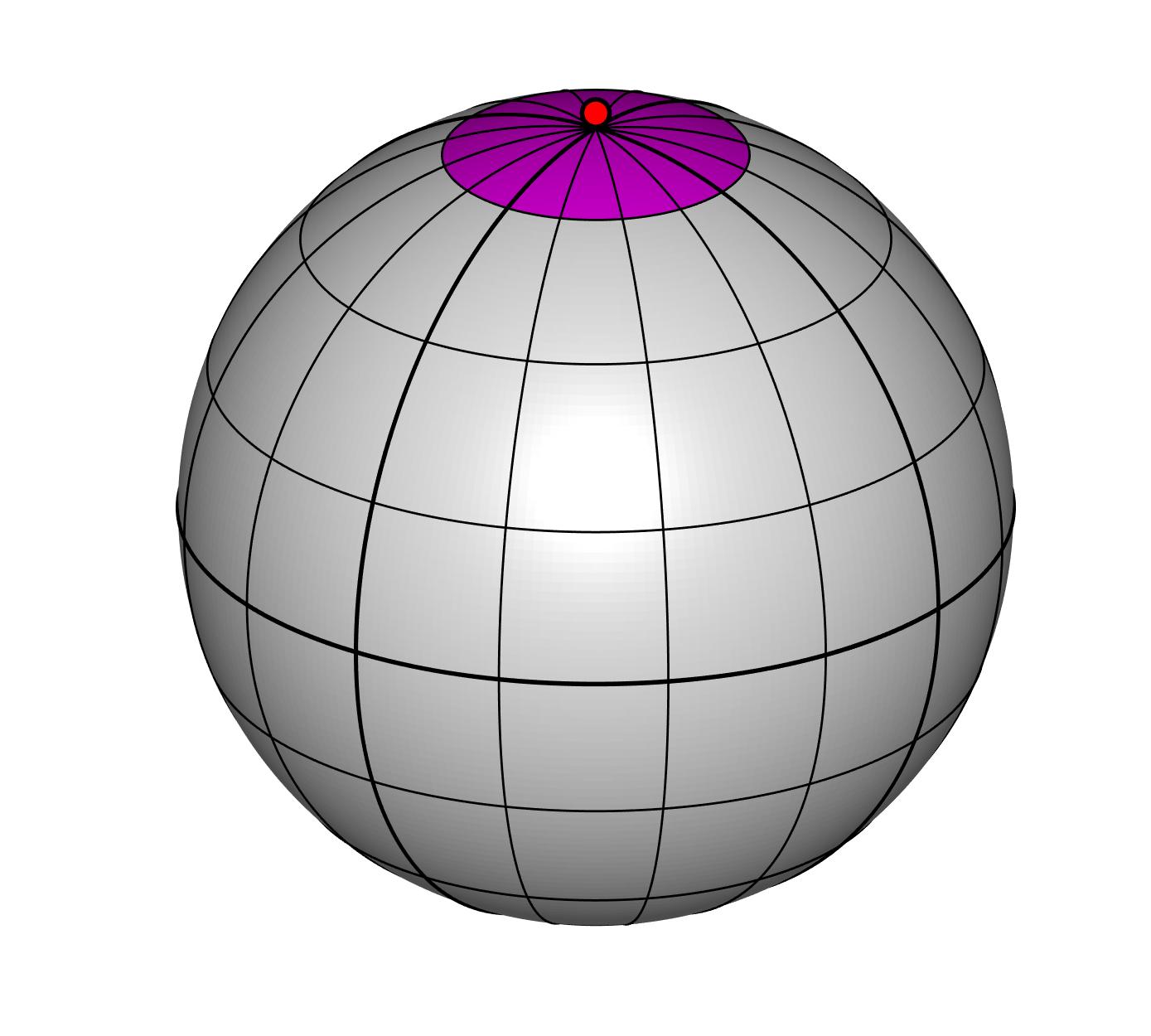}}	
 	\put(0,0){a.}\put(0.2,0){b.}\put(0.4,0){c.}\put(0.6,0){d.}\put(0.8,0){e.}
	\put(0.15,0.2){\includegraphics[trim = 0 0 0 0, clip, width=.7\linewidth ]{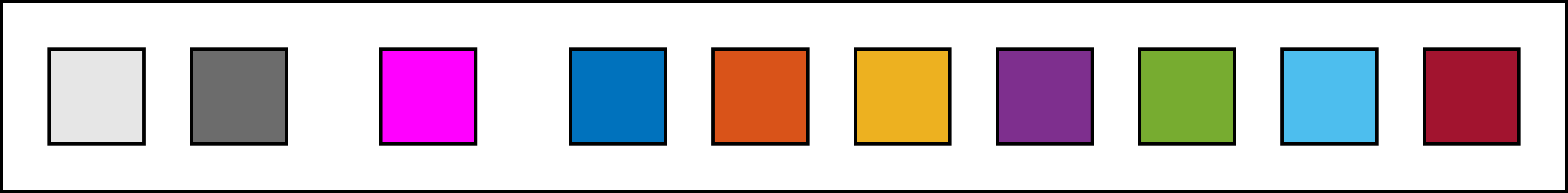}}
	\put(0.16,0.271){\scriptsize Elemental quadrature rule:}
  \put(0.195,0.206){\scriptsize Gauss}  \put(0.32,0.206){\scriptsize Duffy} \put(0.5,0.206){\scriptsize Gauss with adjusted weights}
	\put(0.185,0.24){\scriptsize $n_0$} \put(0.235,0.24){ \twhite{\scriptsize $2 n_0$}}
	\end{picture}
\caption{\textit{Hybrid quadrature on a single-patch NURBS sphere}: Quadrature rules used for hybrid Duffy-Gauss quadrature with adjusted weights~(DGw) considering collocation points that are located within an element~(a.~and b.), along the common edge of two elements~(c.), at the junction of four elements~(d.) and at the pole, where $2^{\ell+2}$ elements meet~(e.). Here $\ell = 2$.}\label{fig:quad_sphere_rules}
\end{figure}

\begin{figure}[h]
\unitlength\linewidth
\begin{picture}(1,.31)
	\put(0,0.06){\includegraphics[trim =  133 66 266 33, clip, width=.25\linewidth ]{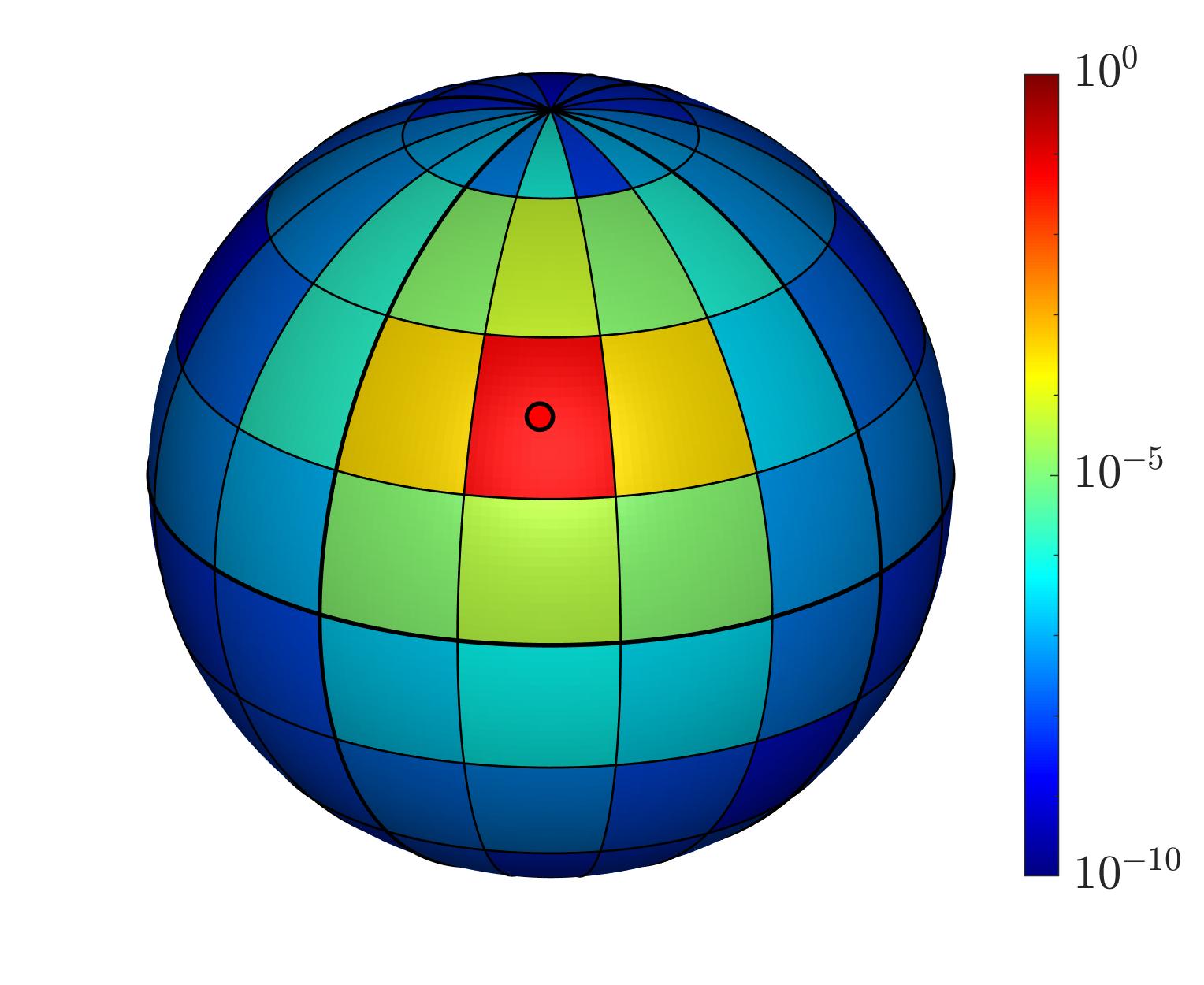}}	
	\put(0.25,0.06){\includegraphics[trim = 133 66 266 33, clip, width=.25\linewidth ]{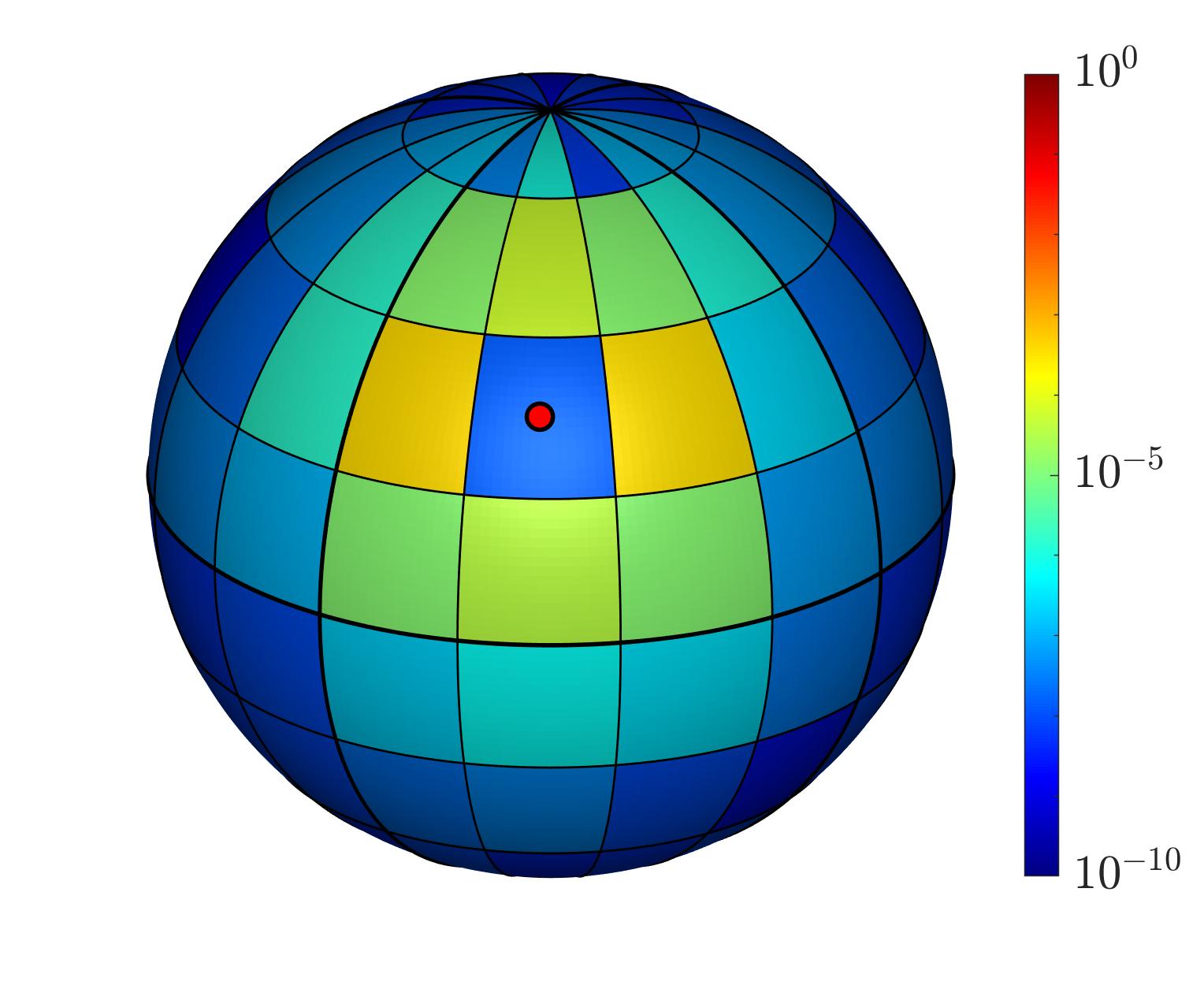}}	
	\put(0.5,0.06){\includegraphics[trim = 133 66 266 33, clip, width=.25\linewidth ]{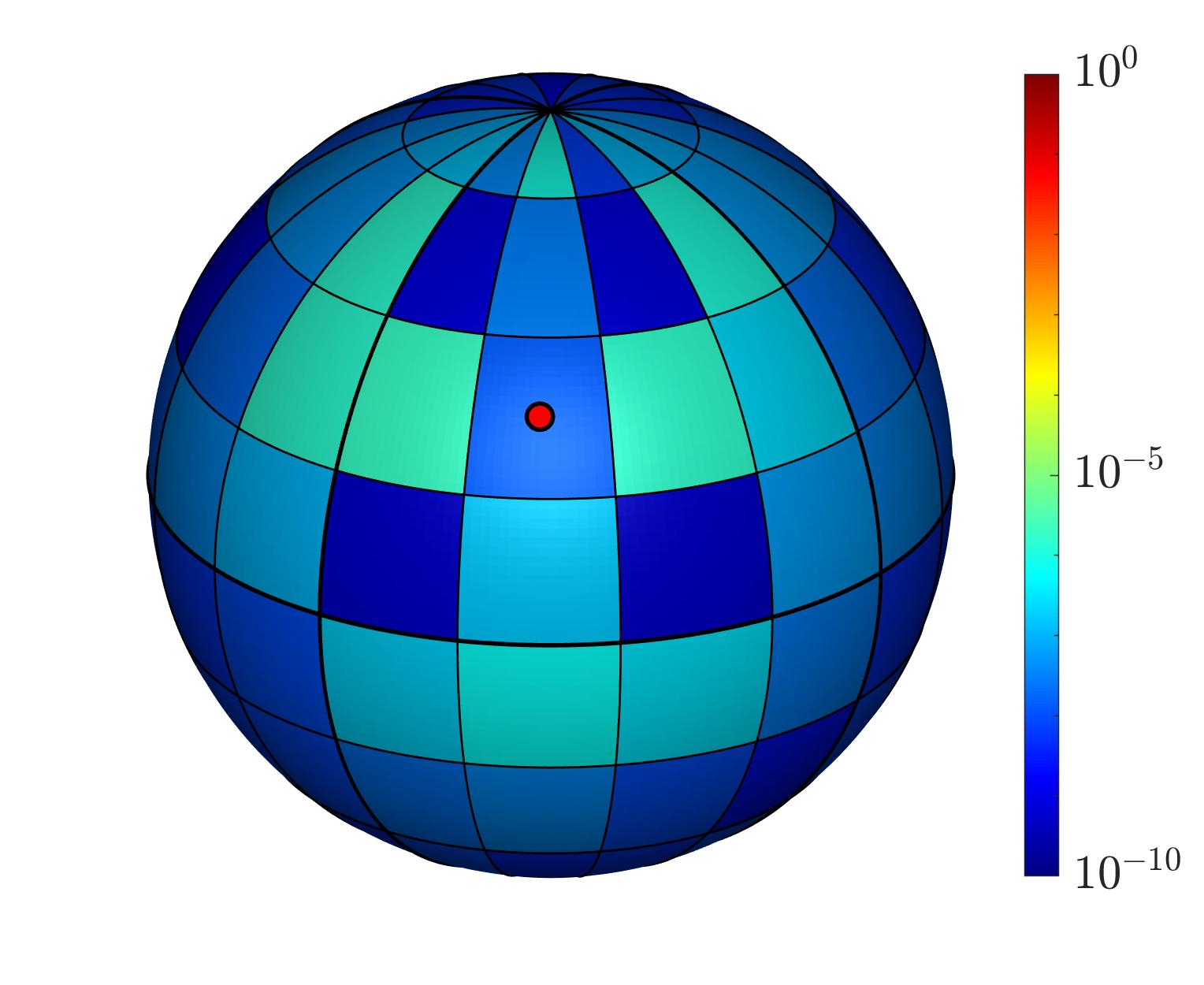}}	
	\put(0.75,0.06){\includegraphics[trim = 133 66 266 33, clip, width=.25\linewidth ]{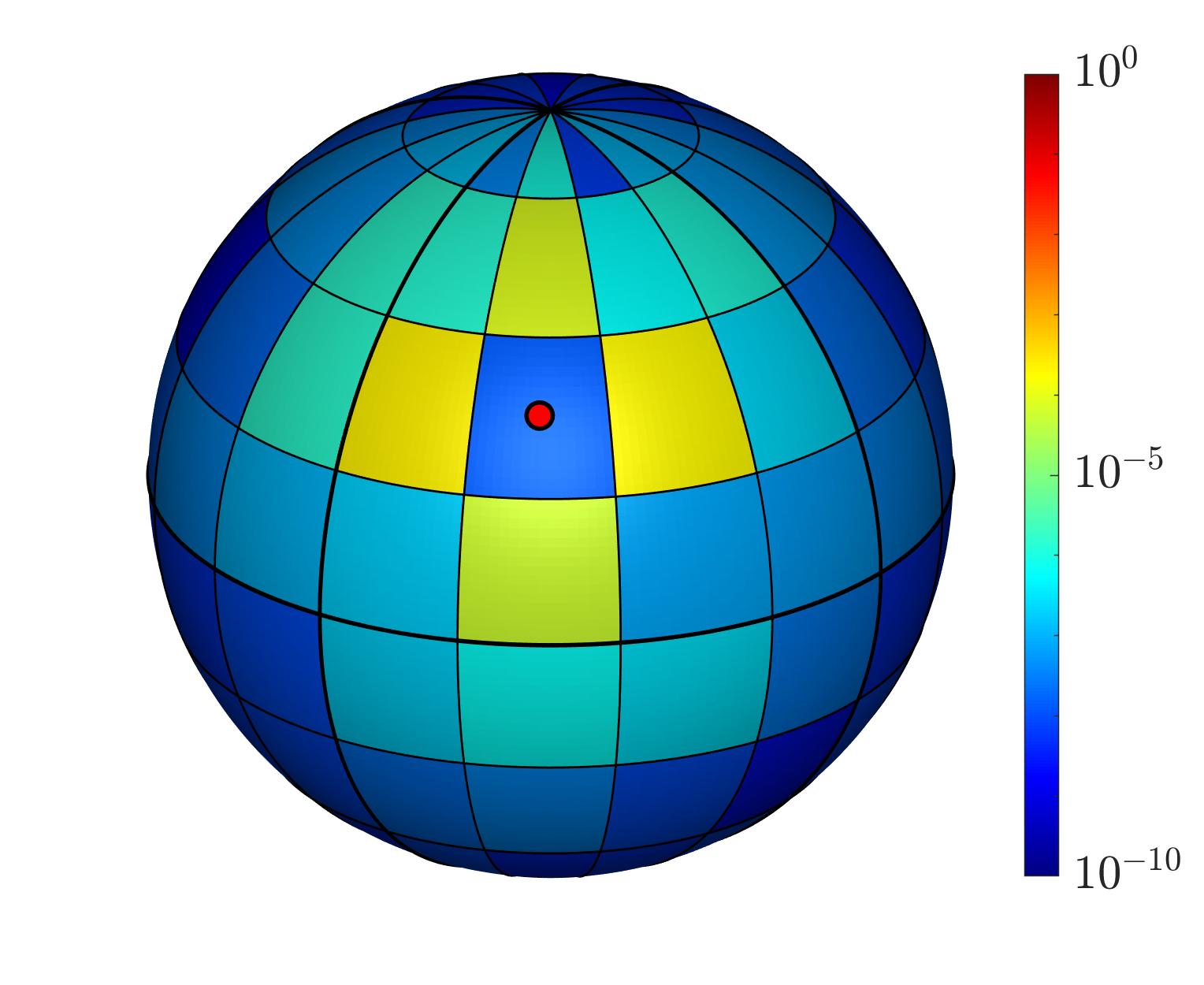}}	
	\put(0.1,0.01){\includegraphics[trim = 220 75 120 950, clip, width=0.8\linewidth ]{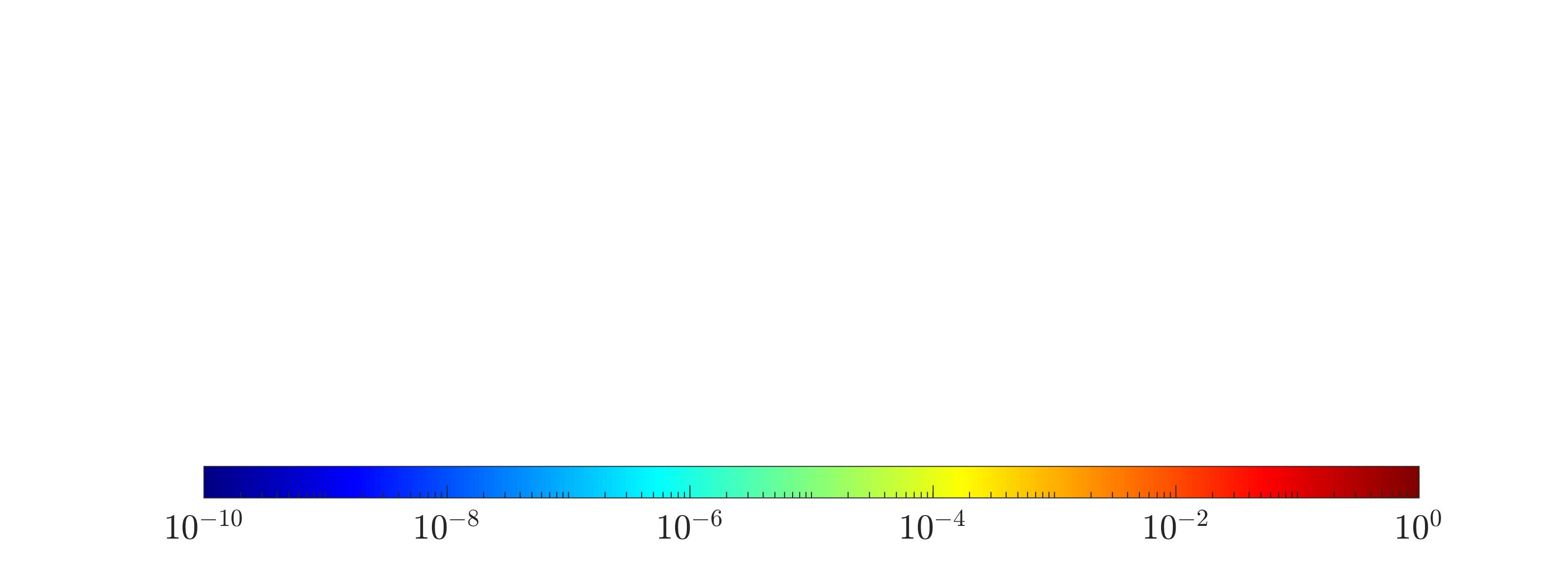}}
	\put(0.475,0){$e_\mathrm{abs}^{e}$}
 	\put(0,0.075){a.}\put(0.25,0.075){b.}\put(0.5,0.075){c.}\put(0.75,0.075){d.}
	\end{picture}
\caption{\textit{Hybrid quadrature on a single-patch NURBS sphere}: Absolute quadrature error \eqref{eq:quad_err0_elem} for Gauss-Legendre quadrature~(G,~a.), hybrid Duffy-Gauss quadrature (DG,~b.), Duffy-Gauss quadrature with progressive refinement~(DGr,~c.) and hybrid Duffy-Gauss quadrature with adjusted weights~(DGw,~d.), each considering collocation point $\by_A$ and $n_0=3$.}\label{fig:quad_sphere_error1}
\end{figure}
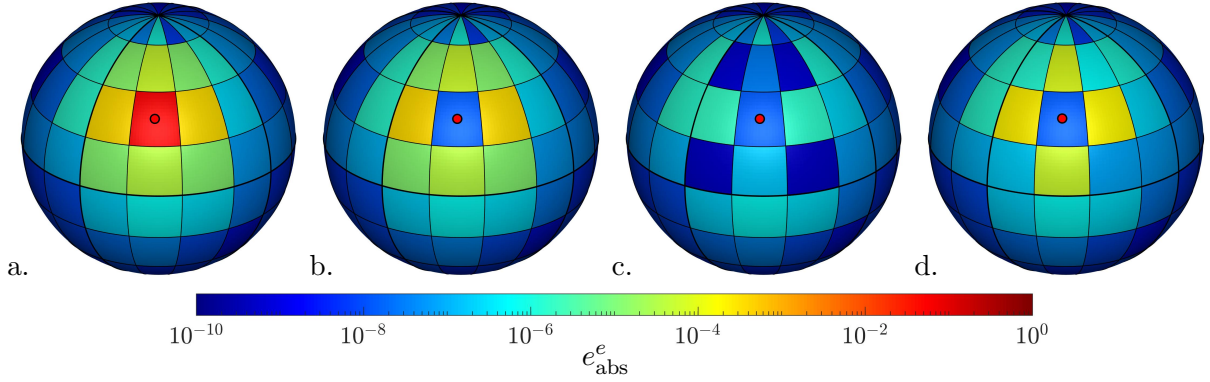

\begin{figure}[h]
\unitlength\linewidth
\begin{picture}(1,.42)
	\put(0,0){\includegraphics[trim = 10 0 30 20, clip, width=.5\linewidth]{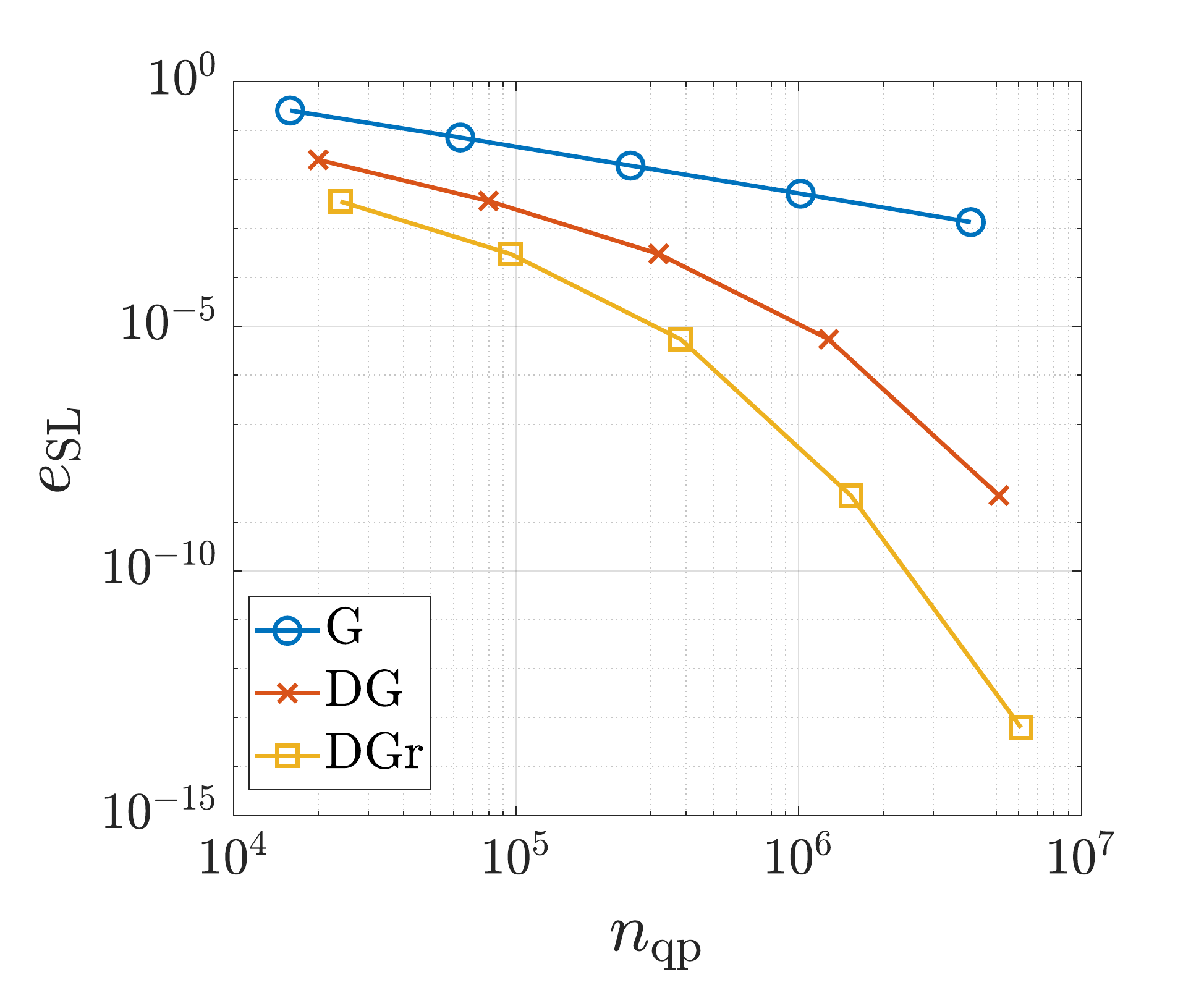}}
	\put(0.5,0){\includegraphics[trim = 10 0 30 20, clip, width=.5\linewidth]{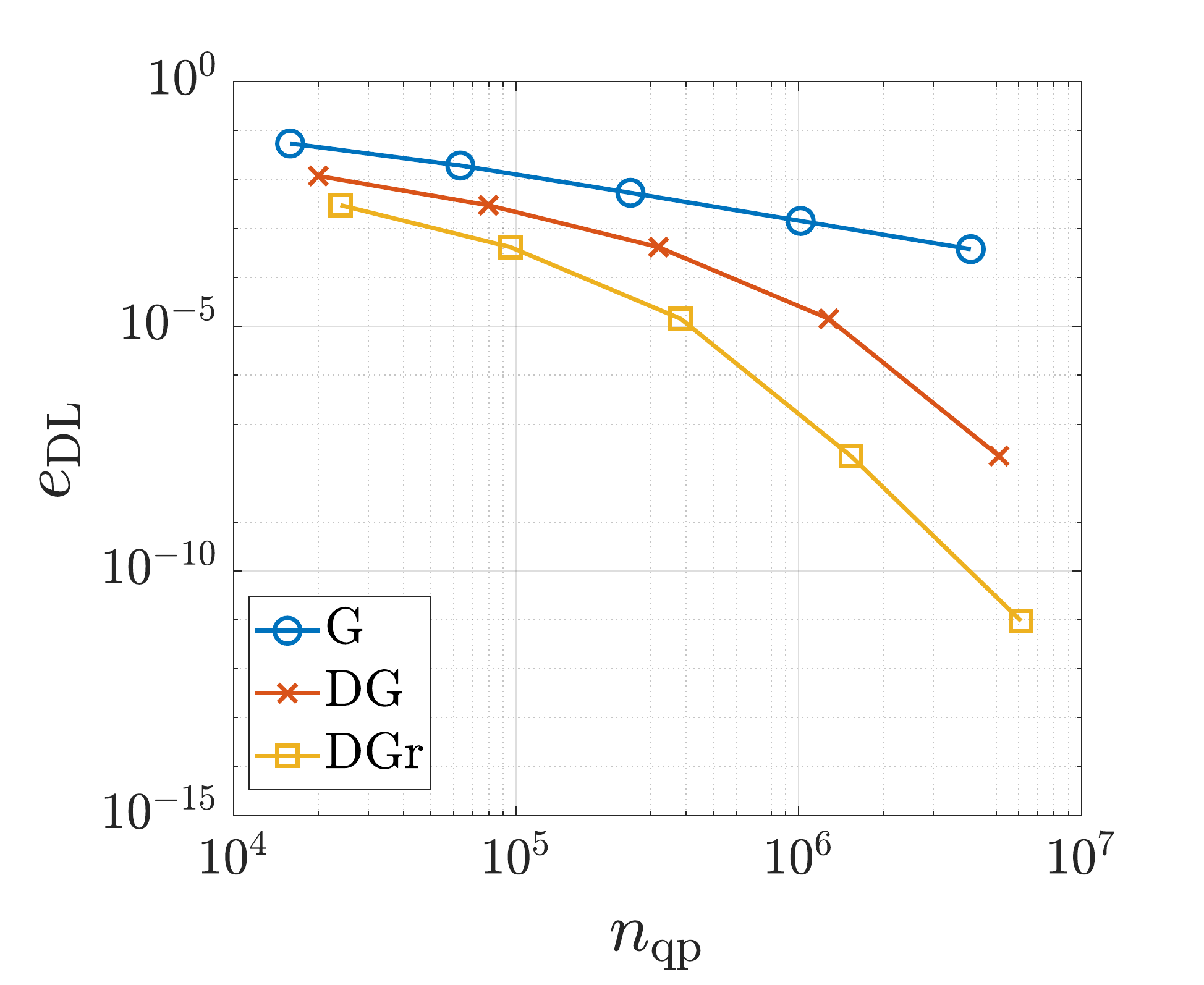}}
	 \put(0,0){a.}\put(0.5,0){b.}
\end{picture}
\caption{\textit{Hybrid quadrature on a single-patch NURBS sphere}: Mean quadrature error \eqref{eq:err_SL2} w.r.t.~BE identities (\ref{eq:identity_SL},~a.) and (\ref{eq:identity_DL},~b.) for $\ell=2$ and varying quadrature refinement $n_0=2,4,8,16,32$.}\label{fig:sphere_convergence}
\end{figure}

Single-patch NURBS surfaces are constructed by revolving a NURBS semicircle by a full rotation \citep{PieglTiller97_book} and they are thus exactly spherical. Fig.~\ref{fig:quad_sphere_discretization}a shows a coarse NURBS discretization of the surface (refinement level $\ell=1$) and the corresponding control points and collocation points. The sphere has constant radius ($C^0$-continuity), continuous tangent vectors ($C^1$-continuity) and constant curvature, making it $C^2$-continuous everywhere. However, these continuity properties are in general not maintained during deformation: The single-patch NURBS sphere is interpolatory across the patch boundary ($C^{-1}$-continuity) and between the other octants of the sphere ($C^0$-continuity), 
while it is $C^{p-1}$-continuous across the remaining element boundaries. Fig.~\ref{fig:quad_sphere_discretization}b shows a NURBS discretization of refinement level $\ell=2$ obtained by knot insertion \citep{Hughes05} maintaining the continuity properties from $\ell=1$.
\\\\Fig.~\ref{fig:quad_sphere_rules}a and b shows the quadrature rules used for hybrid Duffy-Gauss quadrature with adjusted weights (DGw) and exemplary collocation points located within elements. Particular attention has to be paid to collocation points within the degenerated elements located at the poles of the sphere (Fig.~\ref{fig:quad_sphere_rules}b). For these points, Gauss quadrature with adjusted weights cannot be applied to the near singular elements in a meaningful way. Therefore, refined Gauss-Legendre quadrature is applied instead to the elements adjacent to the pole (except the singular element where Duffy quadrature is applied). In contrast, the application of DGw to collocation points within other elements (Fig.~\ref{fig:quad_sphere_rules}a) follows directly from Sec. \ref{sec:hybrid_sheet} and is thus straightforward. Collocation points on element boundaries only occur along the lines of reduced continuity between octants of the sphere (see Fig.~\ref{fig:quad_sphere_discretization}a and b). DGw is applied to collocation points between two (Fig.~\ref{fig:quad_sphere_rules}c) and four octants (Fig.~\ref{fig:quad_sphere_rules}d) by treating each octant of the sphere separately as discussed in Sec. \ref{sec:hybrid_sheet}. For a collocation point at the pole (Fig.~\ref{fig:quad_sphere_rules}e), Duffy quadrature is applied to the $2^{\ell+2}$ elements adjacent to the collocation point.
\\Fig.~\ref{fig:quad_sphere_error1} shows the absolute element-wise quadrature error \eqref{eq:quad_err0_elem} for hybrid quadrature schemes G~(a.), DG~(b.), DGr~(c.) and DGw~(d.) considering collocation point $\by_A$ that is located within an element.\footnote{In absence of an analytical solution, integral $\mathcal I_e$~\eqref{eq:quad_int_elem} is approximated by DG using Duffy quadrature with $\tilde n_\mathrm{gp} = 16$ and Gauss-Legendre quadrature with $\tilde n_\mathrm{gp} = 60$ to obtain an reference solution.} It can be seen that scheme G is very inaccurate on the singular element, while the elemental error decreases with increasing distance to $\by_A$. Using scheme DG instead reduces the quadrature error on the singular element drastically so that the highest errors now occur on the near singular elements. Scheme DGr reduces the error on those elements so that the singular integral is approximated on all elements accurately. Gauss-Legendre quadrature with adjusted weights is designed in Sec.~\ref{sec:quad_nearly_adjusted} to determine singular integrals on near singular elements of plane and regular sheets exactly. Although integral \eqref{eq:quad_int_elem} is not approximated exactly on the near singular elements of curved surfaces, Fig.~\ref{fig:quad_sphere_error1}d shows that scheme DGw yields a lower error than scheme DG for the same number of quadrature points. 
\\\\The mean quadrature error \eqref{eq:err_SL2} with respect to the BE identities \eqref{eq:identity_SL} and \eqref{eq:identity_DL} is shown in Fig.~\ref{fig:sphere_convergence} vs.~the number of total quadrature points \eqref{eq:nqp} for mesh refinement level $\ell=2$ and varying quadrature refinement $n_0=2,4,8,16,32$. Both errors decrease with increasing quadrature refinement for quadrature schemes G, DG and DGr,\footnote{Scheme DGw is not considered for quadrature refinement studies, since it has only been formulated with $n_\mathrm{gp}^e=3\times 3$ for the near-singular elements.} where DG and especially DGr shows a much faster convergence than G: With $n_0=32$, the quadrature error with respect to the first identity (see Fig.~\ref{fig:sphere_convergence}a) is $e_\mathrm{SL}\approx 10^{-3}$ for scheme G, $e_\mathrm{SL}\approx 3\cdot 10^{-9}$ for scheme DG and $e_\mathrm{SL}\approx 6\cdot10^{-14}$ for scheme DGr. The error with respect to the second identity (see Fig.~\ref{fig:sphere_convergence}b) shows a similar behavior, albeit the difference between thee three schemes is slightly smaller. For both errors, schemes DG and DGr show a much better convergence behavior than scheme G whose convergence rate is linear.
\\The singular BE integrals are approximated with high efficiency using a coarse discretization and a moderate number of quadrature points as shown in the previous paragraph. However, refined discretizations are required to represent high order boundary conditions or more complex surface geometries. Another application that requires sufficiently small elements are coupled FE-BE simulations, where FE and BE analysis is conducted on a deforming surface. The accuracy of the quadrature schemes on refined discretization is thus of great interest for BE analysis and is investigated below.
\\\\Fig.~\ref{fig:quad_sphere_error_col} shows the quadrature error w.r.t.~identity \eqref{eq:identity_SL} for all collocation points on one octant of a sphere of refinement level $\ell=1$~(a.), $\ell=2$~(b.), $\ell=3$~(c.) and $\ell=4$~(d.) using quadrature scheme DGr with $n_0=3$. The quadrature error decreases for the majority of the collocation points with increasing quadrature refinement. However, the collocation points next to the pole yield similarly high errors independently of $\ell$. The mean quadrature error \eqref{eq:err_SL2} considering all collocation points is thus not monotonically decreasing as can be seen in Fig.~\ref{fig:sphere_convergence_m}: The quadrature error for DGr increases for the first refinement steps before it starts to decrease at $\ell=3$. Scheme G results in a much higher error that, however, shows linear convergence with an almost constant rate $\mu=0.5$. The curves for scheme DG and scheme DGw lie between those for G and DGr, with DGw being more accurate than DG. For higher $\ell$, all investigated quadrature scheme show linear convergence with the same rate, where DGr is more accurate than DGw, which is more accurate than DG, which is in turn more accurate than G.
\begin{figure}[h]
\unitlength\linewidth
\begin{picture}(1,.25)
	\put(0,0){\includegraphics[trim = 100 133 233 67, clip, width=.23\linewidth ]{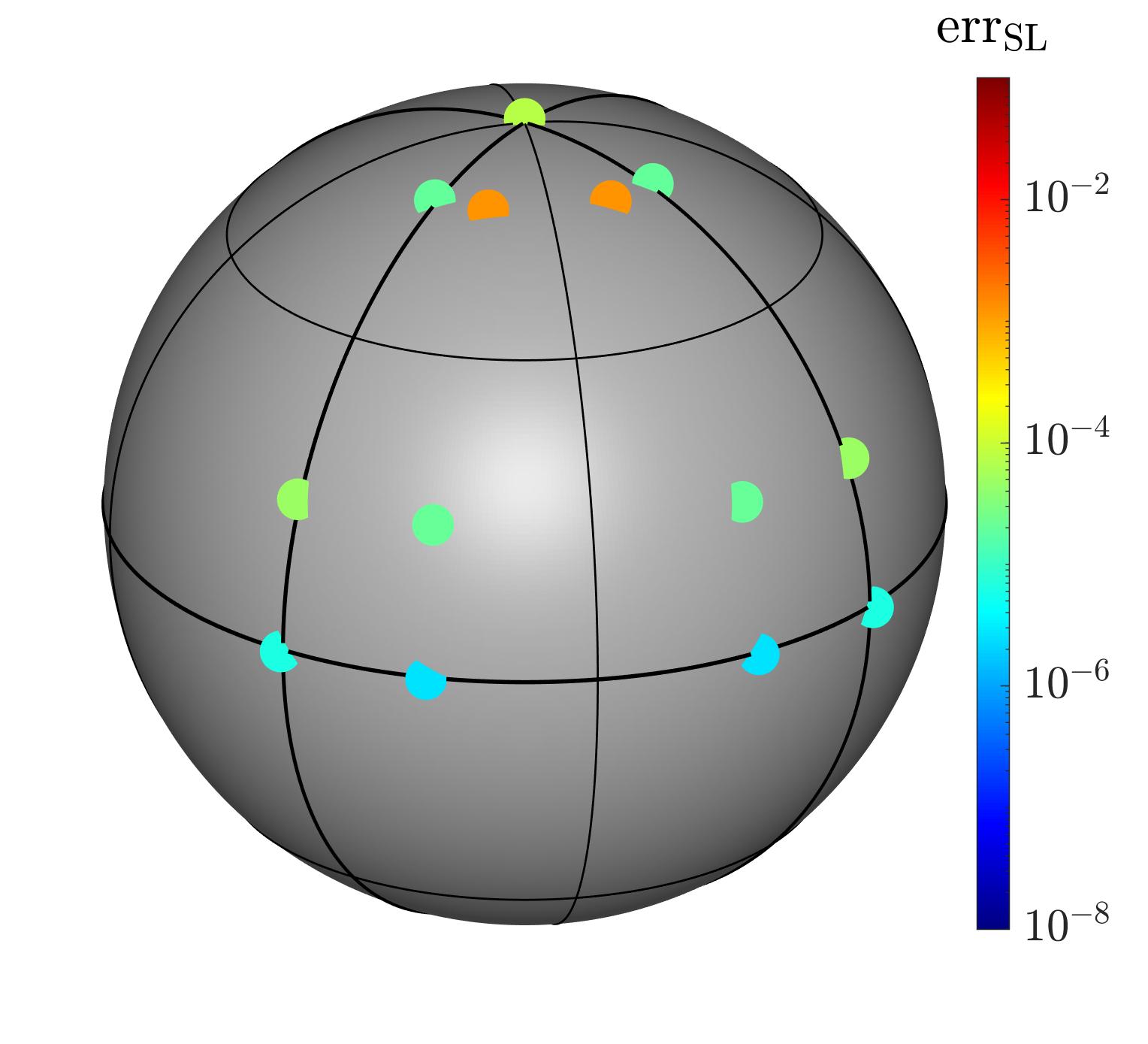}}	
	\put(0.23,0){\includegraphics[trim = 100 133 233 67, clip, width=.23\linewidth ]{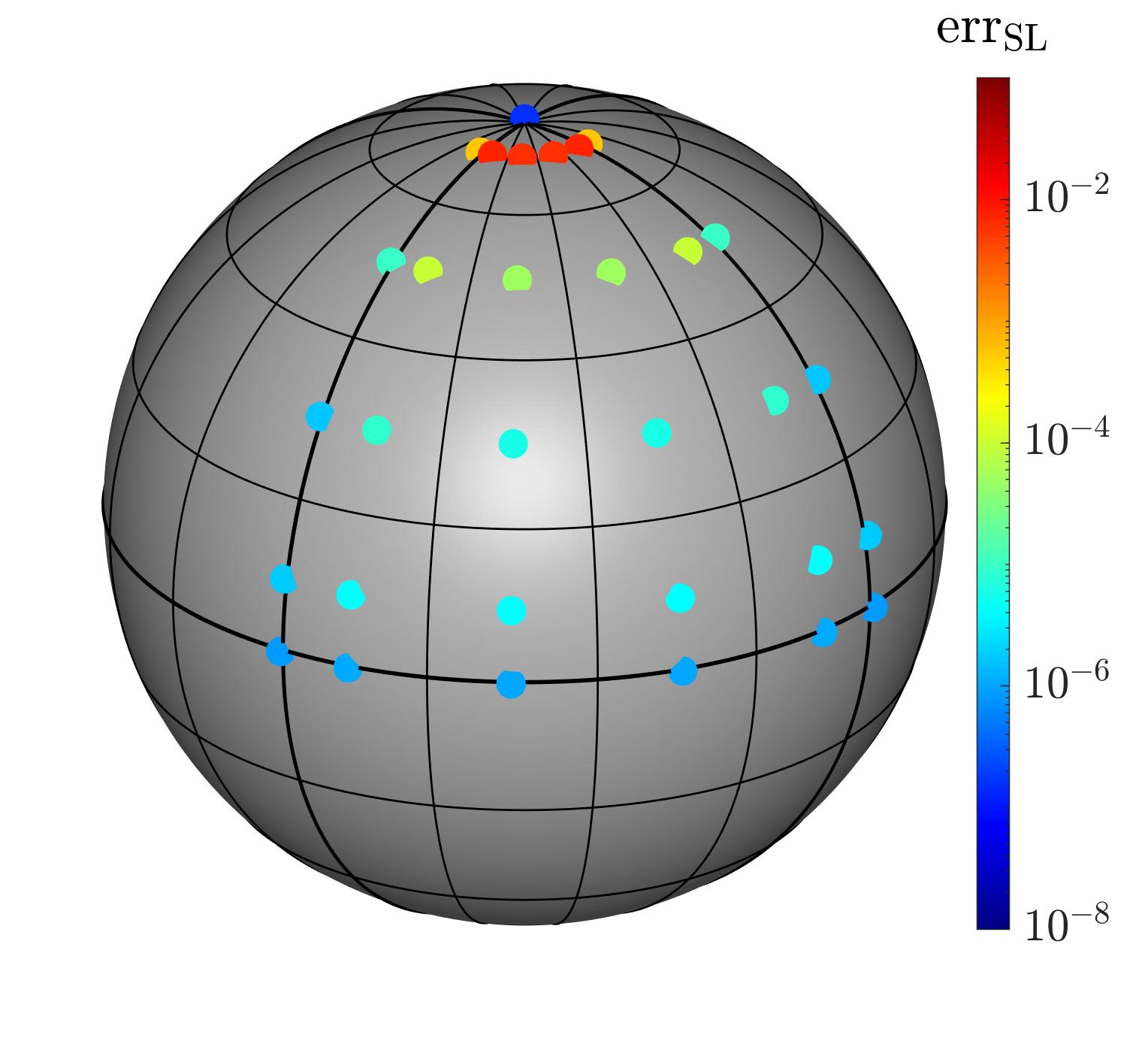}}	
	\put(0.46,0){\includegraphics[trim = 100 133 233 67, clip, width=.23\linewidth ]{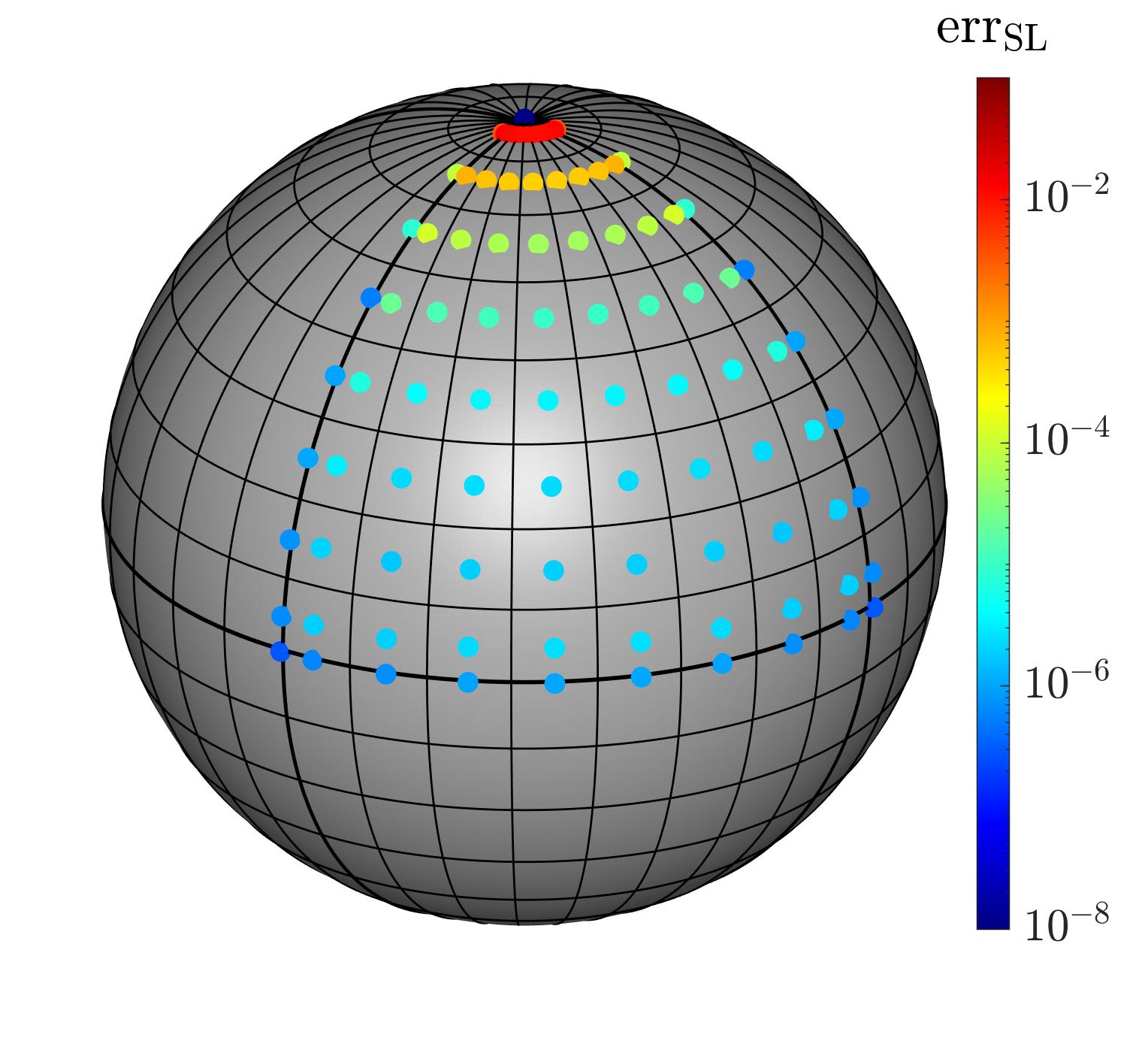}}	
	\put(0.69,0){\includegraphics[trim = 100 133 233 67, clip, width=.23\linewidth ]{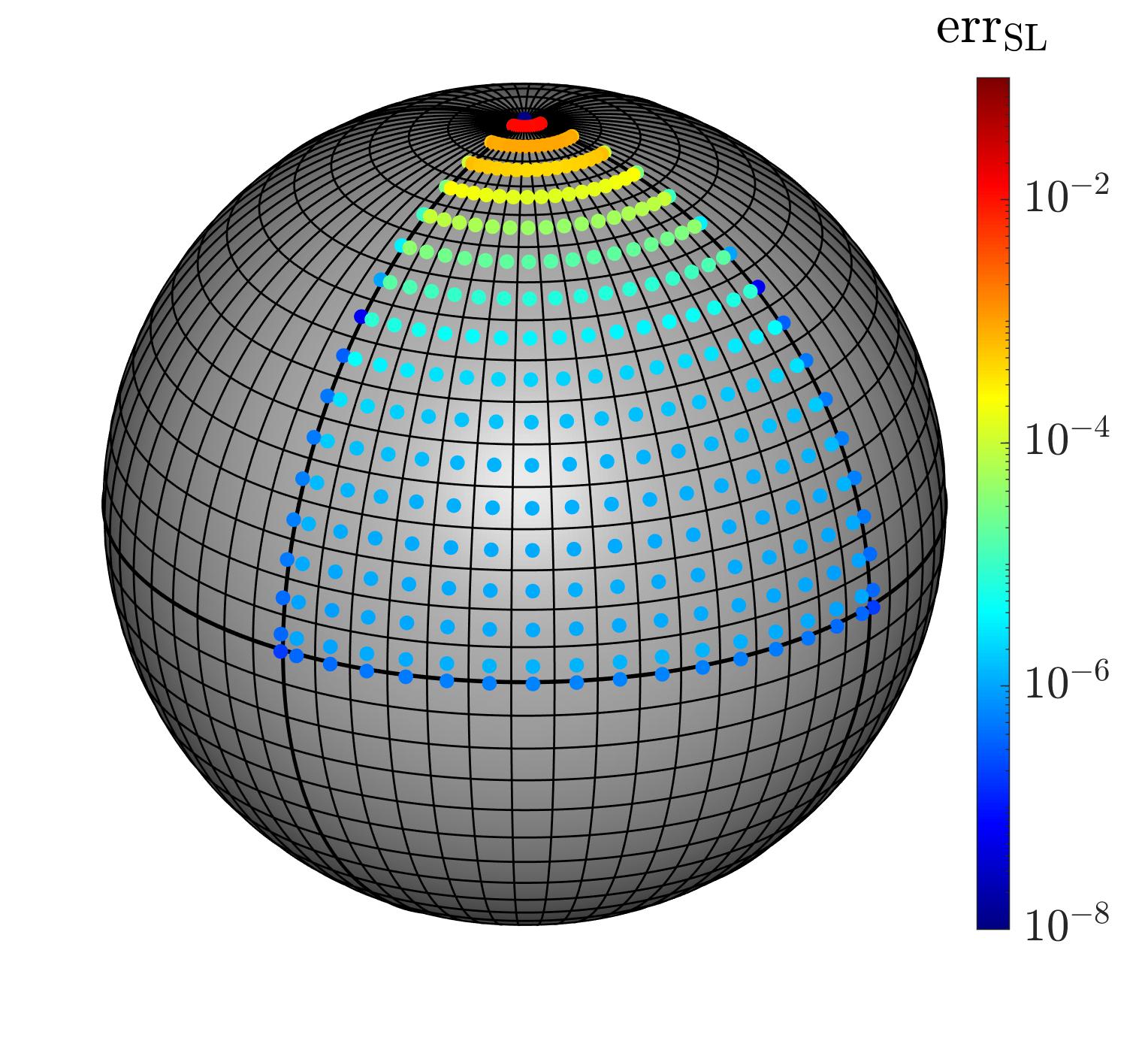}}	
	\put(0.92,0){\includegraphics[trim = 0 0 0 0, clip, width= 0.045\linewidth ]{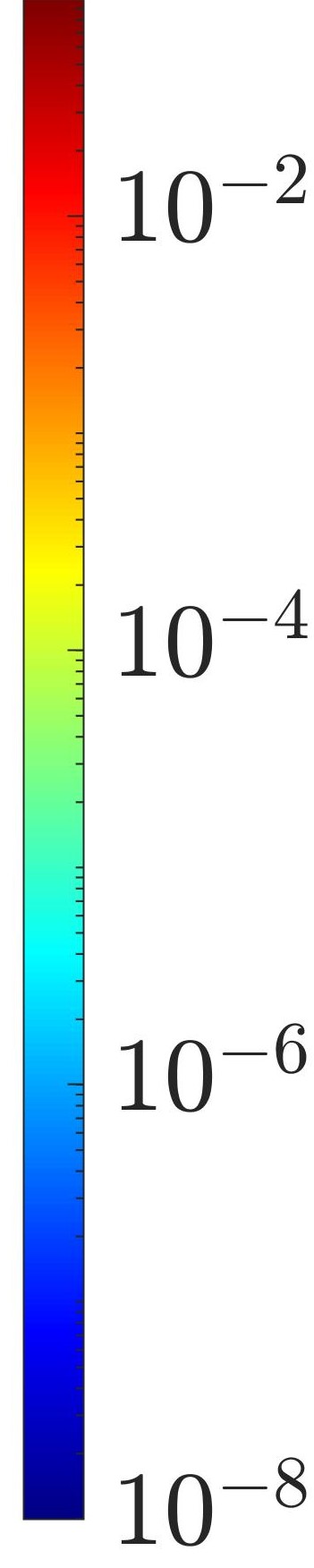}}	
	\vwput{0.96}{0.085}{\scriptsize $e_\mathrm{SL}^{\by_A}$}
 	\put(0,0){a.}\put(0.23,0){b.}\put(0.46,0){c.}\put(0.69,0){d.}
	\end{picture}
\caption{\textit{Hybrid quadrature on a single-patch NURBS sphere}: Quadrature error $e_\mathrm{SL}^{\by_A}$ \eqref{eq:err_SL} for DGr w.r.t.~identity \eqref{eq:identity_SL} considering all collocation points on one octant of the sphere for $n_0=3$ and discretization level $\ell=1$~(a.), $\ell=2$~(b.),  $\ell=3$~(c.) and  $\ell=4$~(d.).}\label{fig:quad_sphere_error_col}
\end{figure}
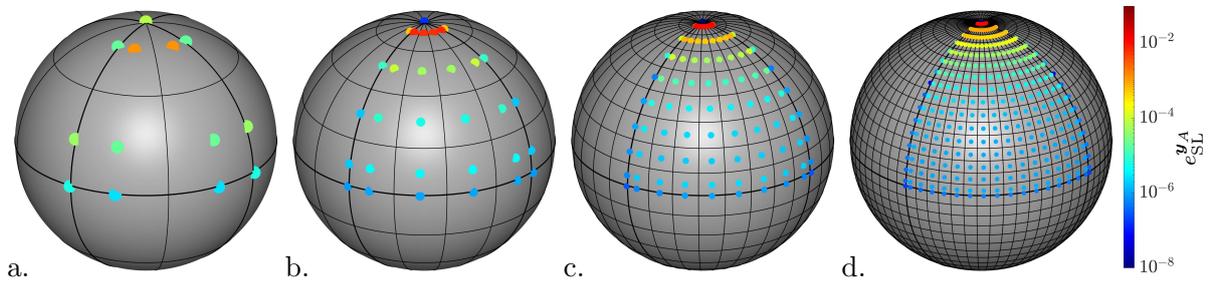
\begin{figure}[h]
\unitlength\linewidth
\begin{picture}(1,.42)
	\put(0,0){\includegraphics[trim = 10 0 30 20, clip, width=.5\linewidth]{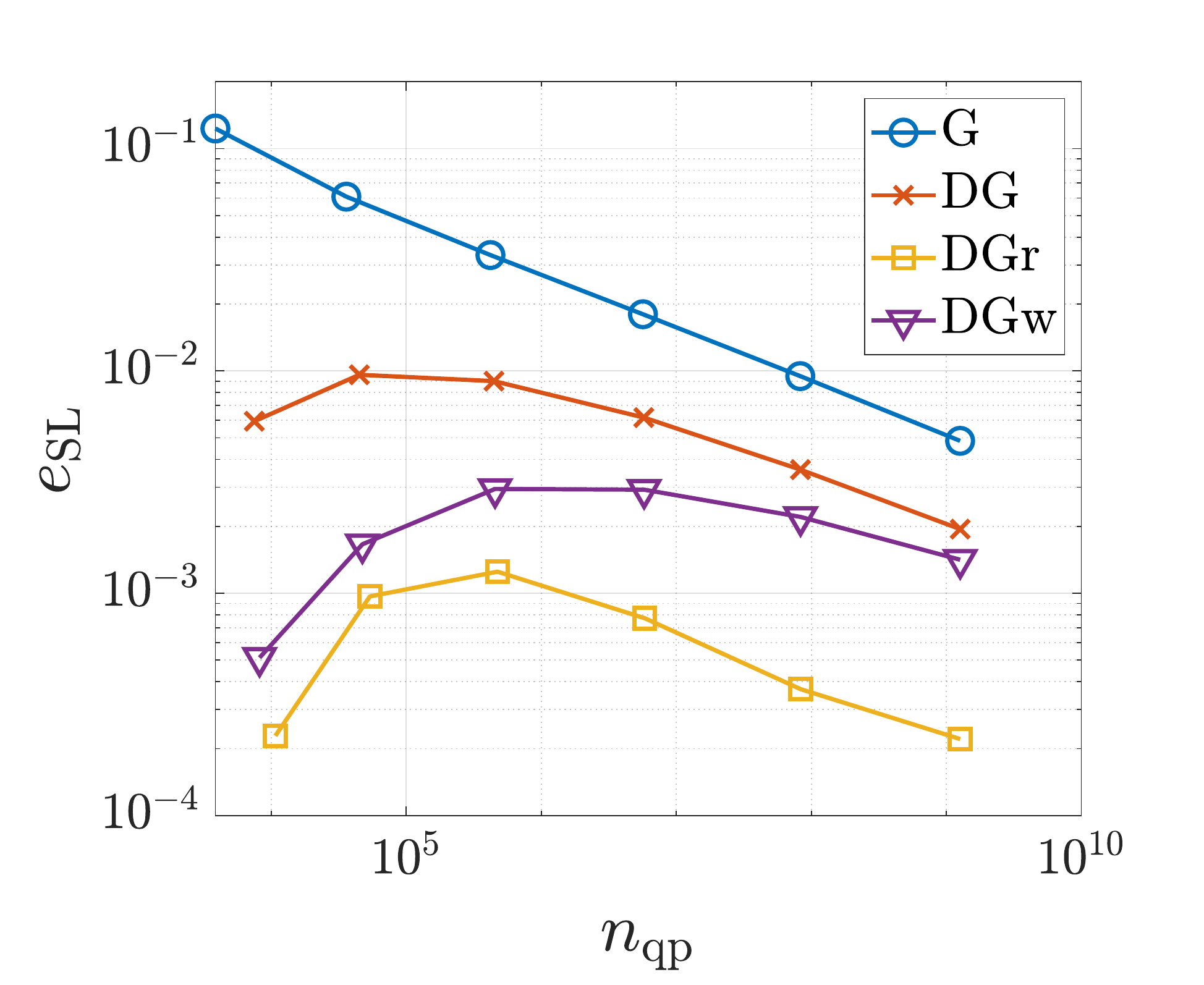}}
	\put(0.5,0){\includegraphics[trim = 10 0 30 20, clip, width=.5\linewidth]{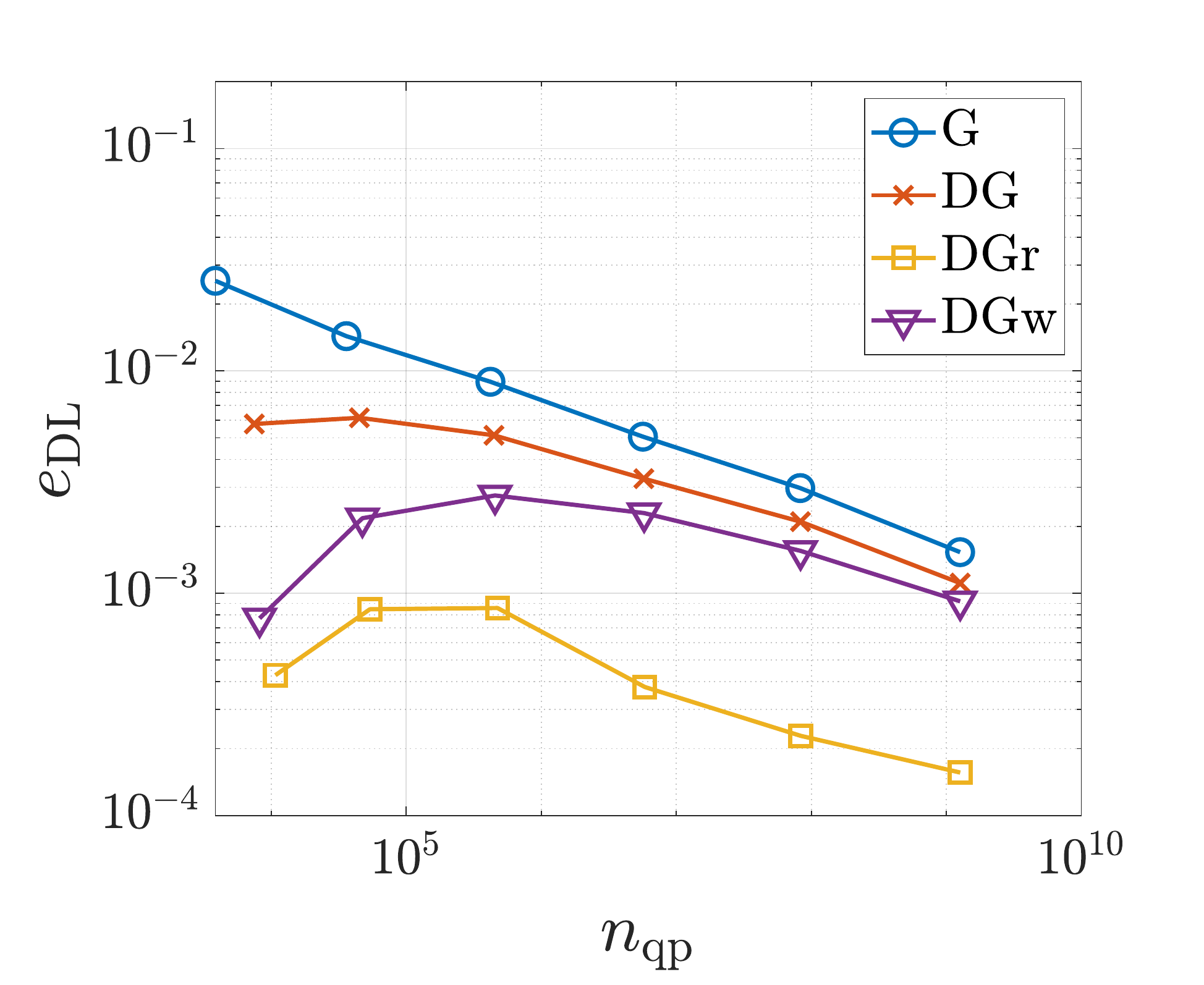}}
	 \put(0.01,0){a.}\put(0.51,0){b.}
\end{picture}
\caption{\textit{Hybrid quadrature on a single-patch NURBS sphere}: Mean quadrature error \eqref{eq:err_SL2} w.r.t.~BE identities (\ref{eq:identity_SL},~a.) and (\ref{eq:identity_DL},~b.) for $n_0=3$ and varying mesh refinement $\ell=1,\ldots,6$.}\label{fig:sphere_convergence_m}
\end{figure}

\subsubsection{Six-patch NURBS sphere}\label{sec:hybrid_sphere_six}
\begin{figure}[h]
\unitlength\linewidth
\begin{picture}(1,.29)
	\put(0,0){\includegraphics[trim = 150 100 150 70, clip, width=.2\linewidth ]{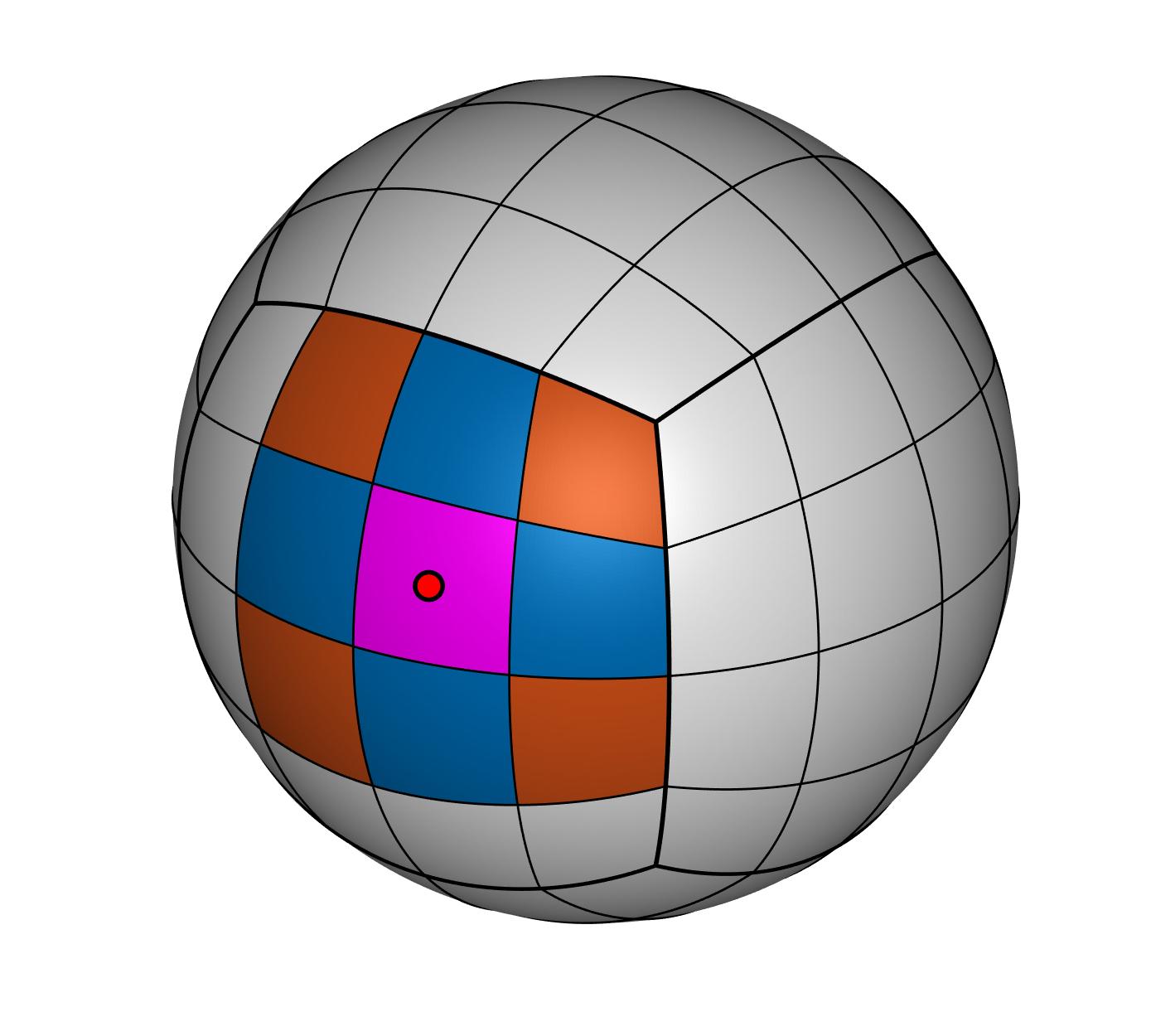}}	
	\put(0.2,0){\includegraphics[trim = 150 100 150 70, clip, width=.2\linewidth ]{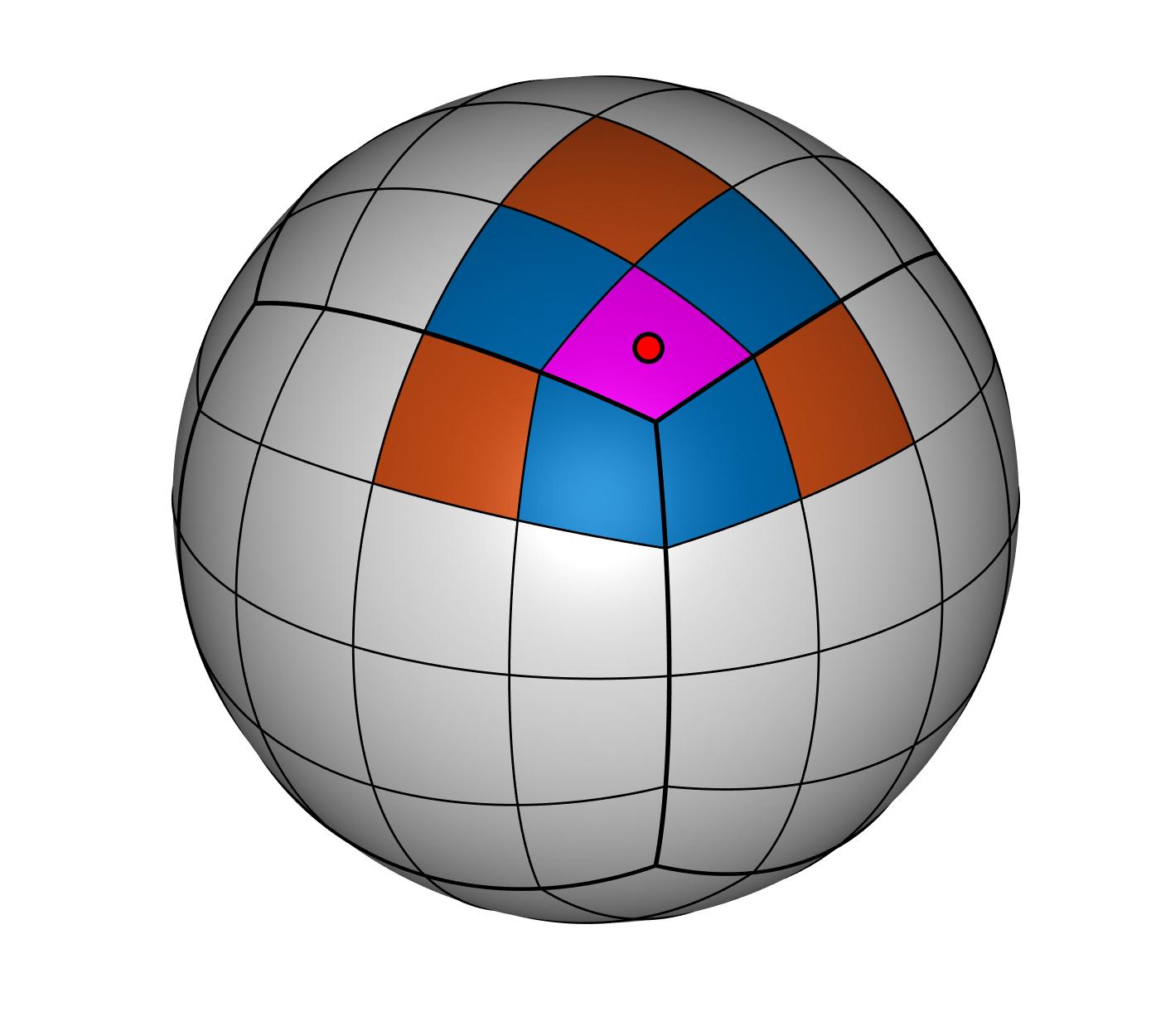}}	
	\put(0.4,0){\includegraphics[trim = 150 100 150 70, clip, width=.2\linewidth ]{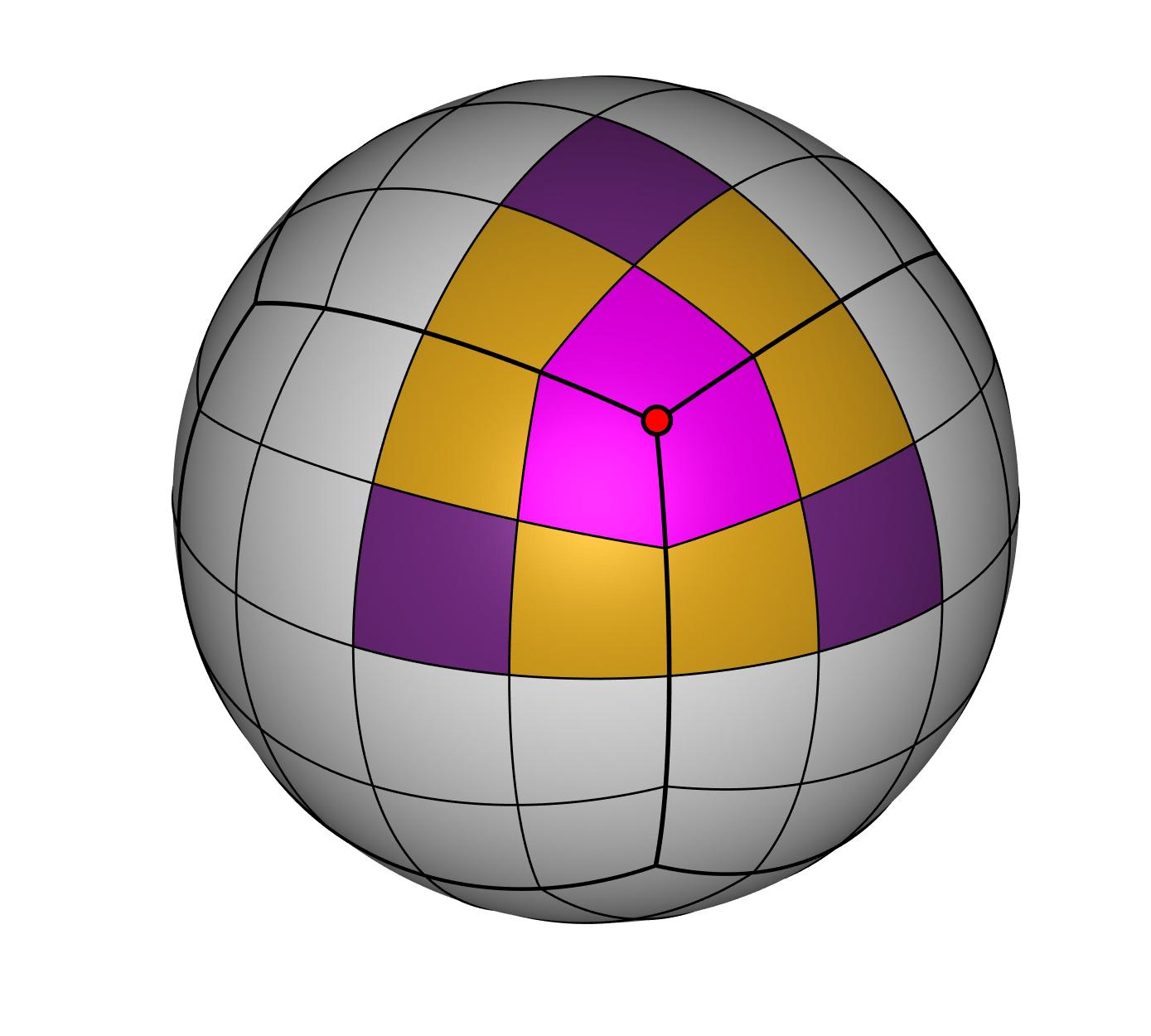}}	
	\put(0.6,0){\includegraphics[trim = 150 100 150 70, clip, width=.2\linewidth ]{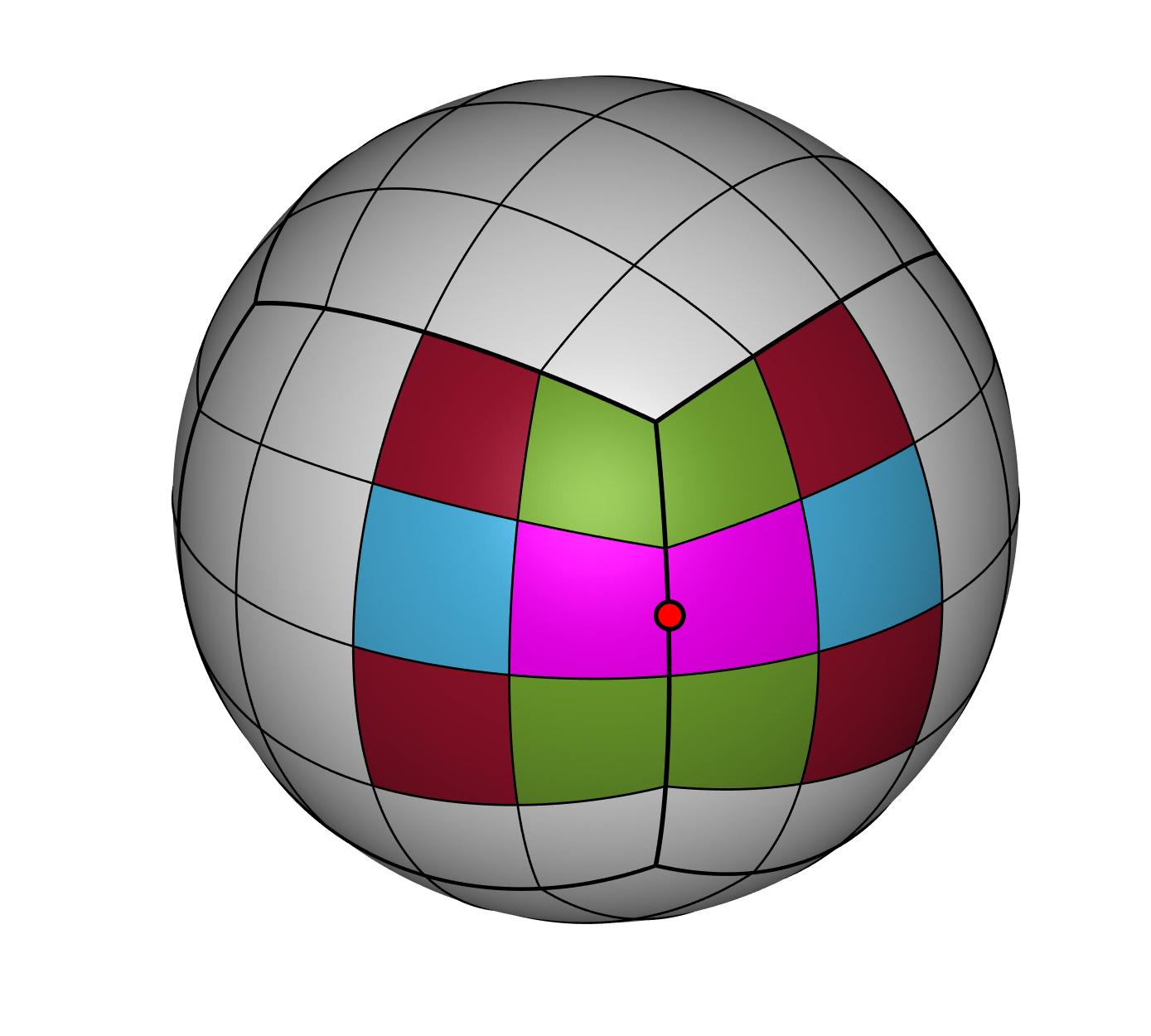}}	
	\put(0.8,0){\includegraphics[trim = 150 100 150 70, clip, width=.2\linewidth ]{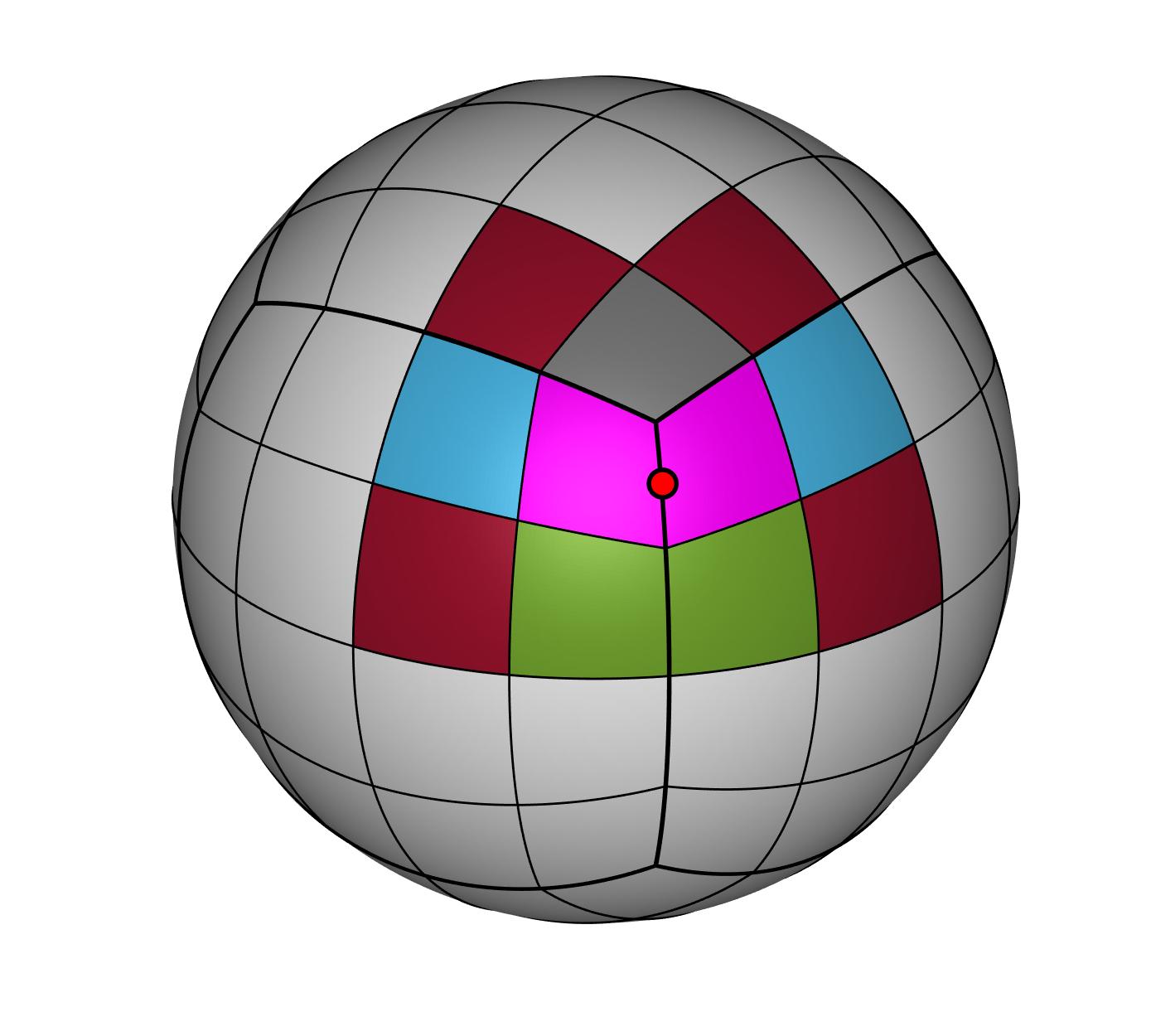}}	
 	\put(0,0){a.}\put(0.2,0){b.}\put(0.4,0){c.}\put(0.6,0){d.}\put(0.8,0){e.}
	\put(0.15,0.2){\includegraphics[trim = 0 0 0 0, clip, width=.7\linewidth ]{figures/sphere/legend_sphere2.png}}
	\put(0.16,0.271){\scriptsize Elemental quadrature rule:}
  \put(0.195,0.206){\scriptsize Gauss}  \put(0.32,0.206){\scriptsize Duffy} \put(0.5,0.206){\scriptsize Gauss with adjusted weights}
	\put(0.185,0.24){\scriptsize $n_0$} \put(0.235,0.24){ \twhite{\scriptsize $2 n_0$}} 	
	\end{picture}
\caption{\textit{Hybrid quadrature on a six-patch NURBS sphere}: Quadrature rules used for hybrid Duffy-Gauss quadrature with adjusted weights~(DGw) considering collocation points that are located within one patch (a.~and b.), at the junction of three patches (c.) and along the common edge of two patches (c.~and d.).}\label{fig:quad_sphere_rules_6p}
\end{figure}

\begin{figure}[h]
\unitlength\linewidth
\begin{picture}(1,.42)
	\put(0,0){\includegraphics[trim = 10 0 30 20, clip, width=.5\linewidth]{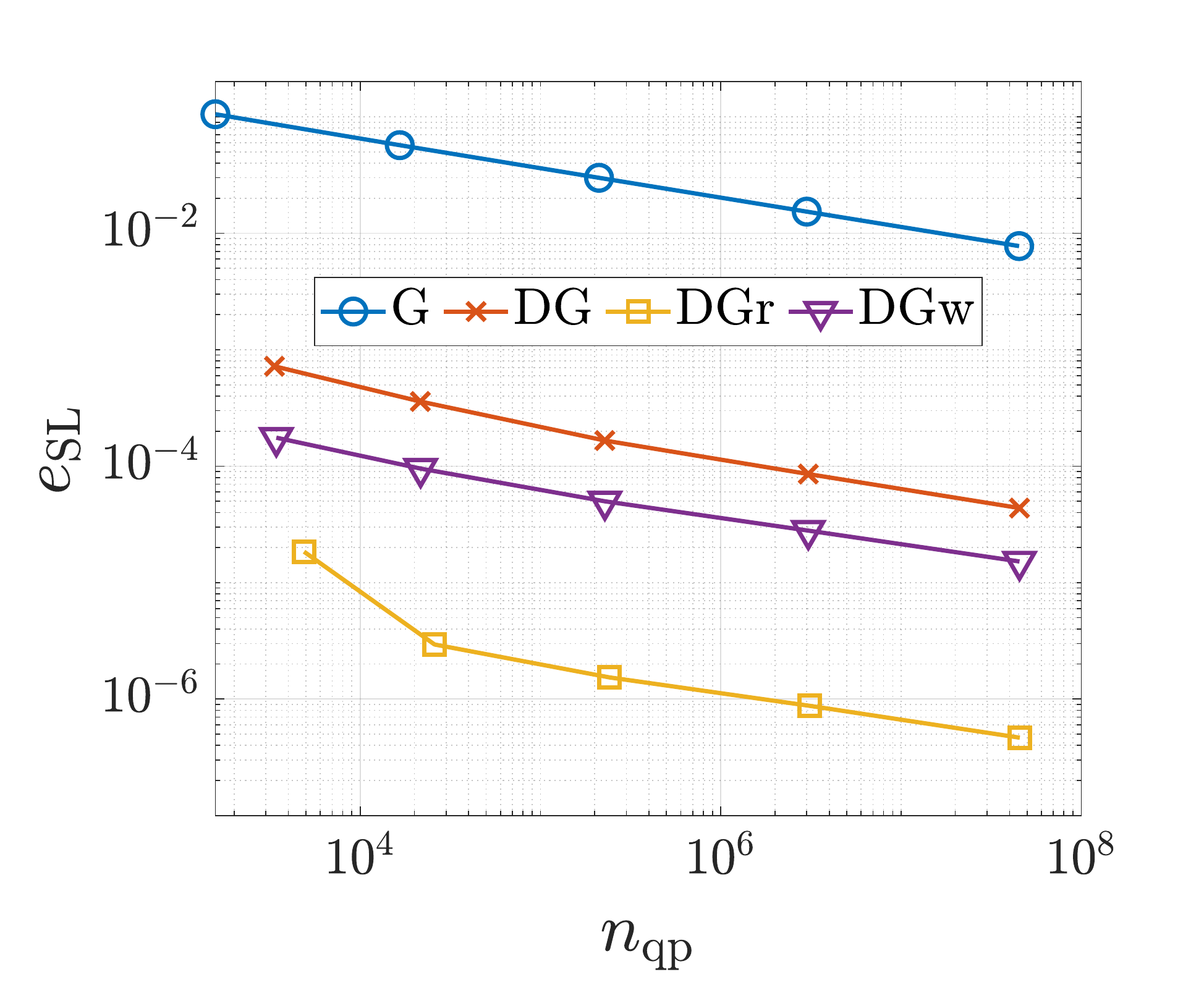}}
	\put(0.5,0){\includegraphics[trim = 10 0 30 20, clip, width=.5\linewidth]{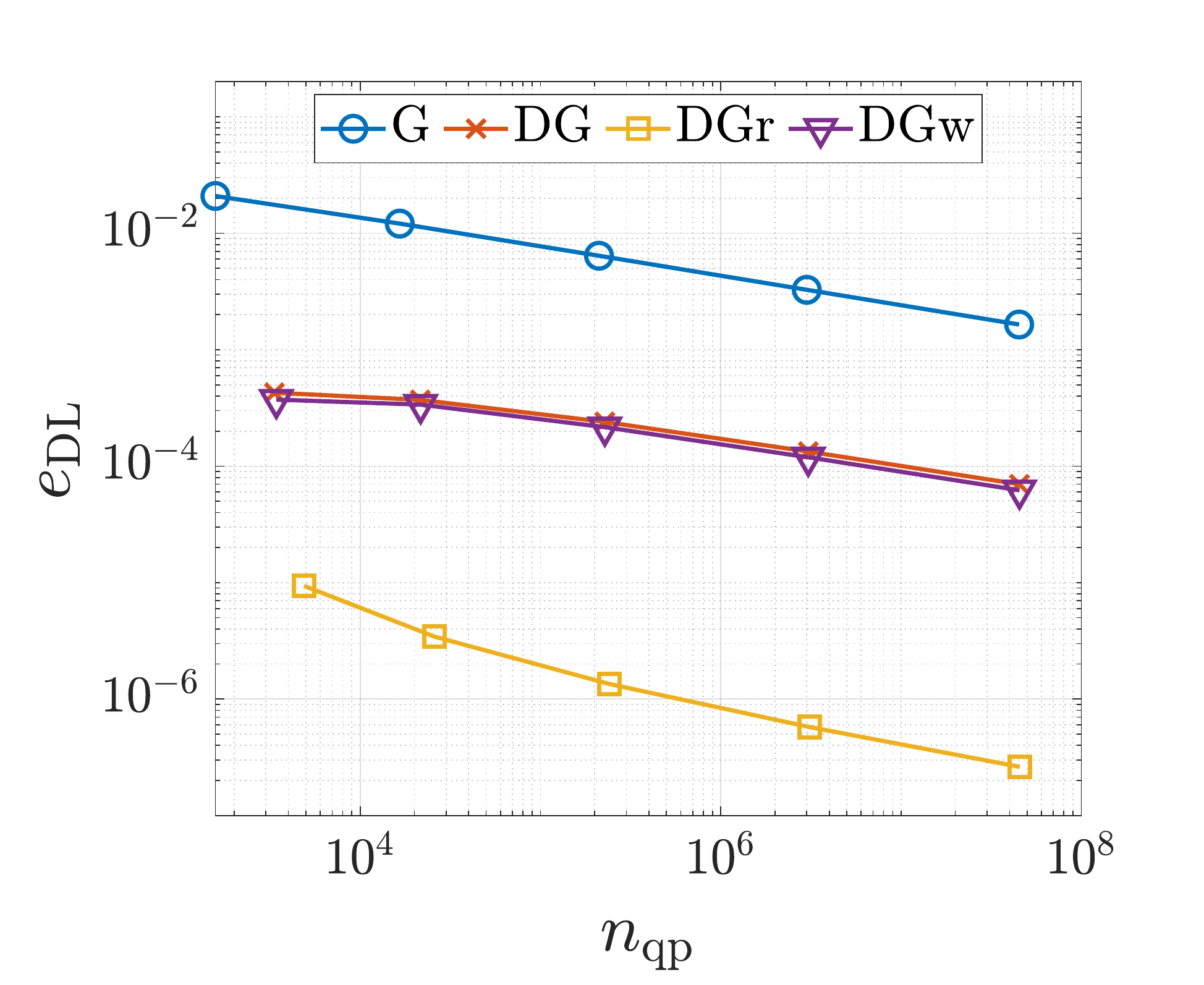}}
	 \put(0.01,0){a.}\put(0.51,0){b.}
\end{picture}
\caption{\textit{Hybrid quadrature on a six-patch NURBS sphere}: Mean quadrature error \eqref{eq:err_SL2} w.r.t.~BE identities (\ref{eq:identity_SL},~a.) and (\ref{eq:identity_DL},~b.) for quadrature density $n_0=3$ and varying mesh refinement $\ell=1,\ldots,5$.}\label{fig:sphere_convergence_6p_m}
\end{figure}

This section discusses the application of the four hybrid quadrature schemes to the discretizations shown in Fig.~\ref{fig:quad_sphere_discretization}c and d. These discretizations consist of six biquadratic NURBS patches and are thus referred to as six-patch NURBS spheres. In contrast to the single-patch discretziation from Sec.~\ref{sec:hybrid_sphere_six}, the six-patch NURBS sphere is only approximately spherical (see Fig.~\ref{fig:quad_sphere_discretization}e). The application of hybrid quadrature schemes to six-patch NURBS spheres is nevertheless promising, since there are no degenerated and small elements near the pole, which cause the largest quadrature errors on the single-patch sphere (cf.~Fig.~\ref{fig:quad_sphere_error_col}). It should be further noted that all collocation points lying on element boundaries also lie on patch boundaries as Fig.~\ref{fig:quad_sphere_discretization}c and d shows.
\\\\The quadrature rules and the element types for hybrid Duffy-Gauss quadrature with adjusted weights (DGw) are shown in Fig.~\ref{fig:quad_sphere_rules_6p}a and b for collocation points within patches, in Fig.~\ref{fig:quad_sphere_rules_6p}c for collocation points at junctions of three patches and in Fig.~\ref{fig:quad_sphere_rules_6p}d and e for collocation points at boundaries between two patches. Particular attention has to be paid to collocation points located near, but not directly at, junctions of three patches (see Fig.~\ref{fig:quad_sphere_rules_6p}b and e). For those collocation points, the number of near singular elements is reduced by one (cf.~Fig.~\ref{fig:quad_sphere_rules_6p}a and d).
\begin{figure}[h]
\unitlength\linewidth
\begin{picture}(1,.42)
	\put(0,0){\includegraphics[trim = 10 0 30 20, clip, width=.5\linewidth]{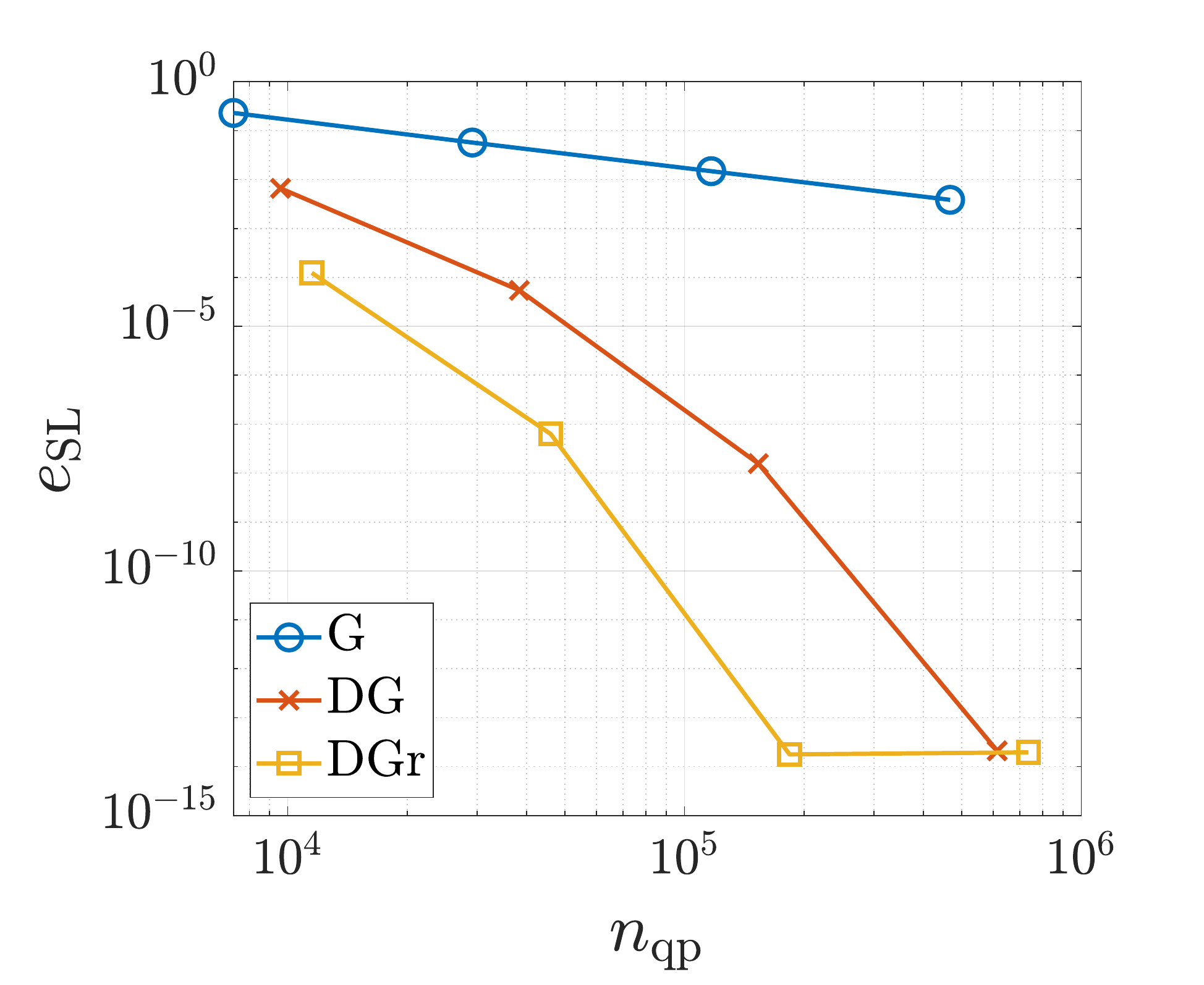}}
	\put(0.5,0){\includegraphics[trim = 10 0 30 20, clip, width=.5\linewidth]{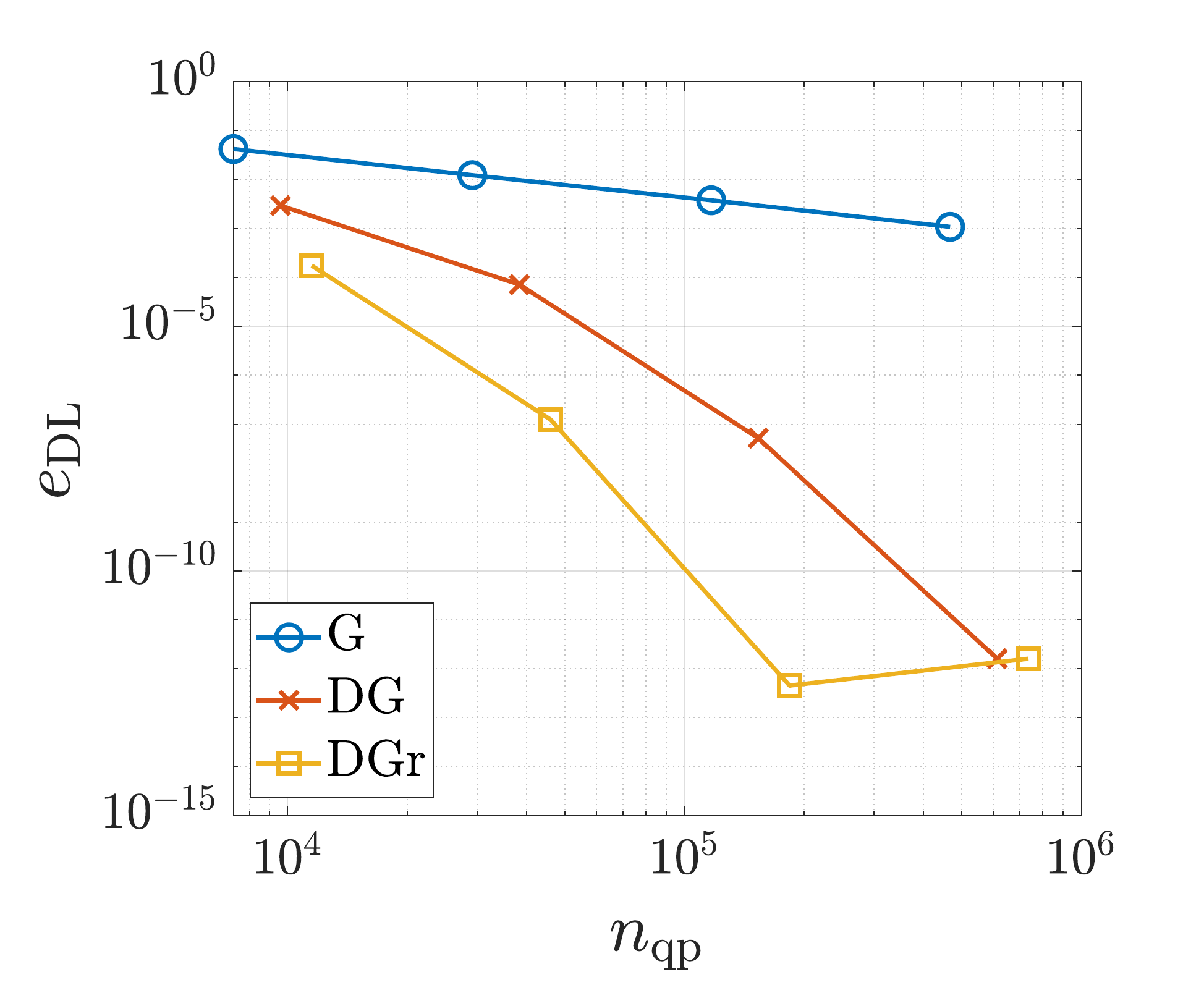}}
	 \put(0,0){a.}\put(0.5,0){b.}
\end{picture}
\caption{\textit{Hybrid quadrature on a six-patch NURBS sphere}: Mean quadrature error (\ref{eq:err_SL},~a.) and (\ref{eq:err_DL},~b.) for $\ell=2$ and varying quadrature refinement $n_0=2,4,8,16$.}\label{fig:sphere_convergence_6p}
\end{figure}
\\The mean quadrature error \eqref{eq:err_SL2} with respect to BE identities \eqref{eq:identity_SL} and \eqref{eq:identity_DL} is shown in Fig.~\ref{fig:sphere_convergence_6p} vs.~the number of total quadrature points \eqref{eq:nqp} for mesh refinement level $\ell=2$ and varying quadrature refinement $n_0=2,4,8,16$.  Similar to the single-patch sphere from Sec.~\ref{sec:hybrid_sphere_single}, the hybrid quadrature schemes G, DG, and DGr show decreasing quadrature errors as the number of quadrature points increases. The convergence behavior of quadrature schemes that use Duffy quadrature is even better for the six-patch sphere than for the single-patch sphere so that DG and DGr yield mean quadrature errors in the range of machine precision for $n_0=16$, respectively for $n_0=8$.
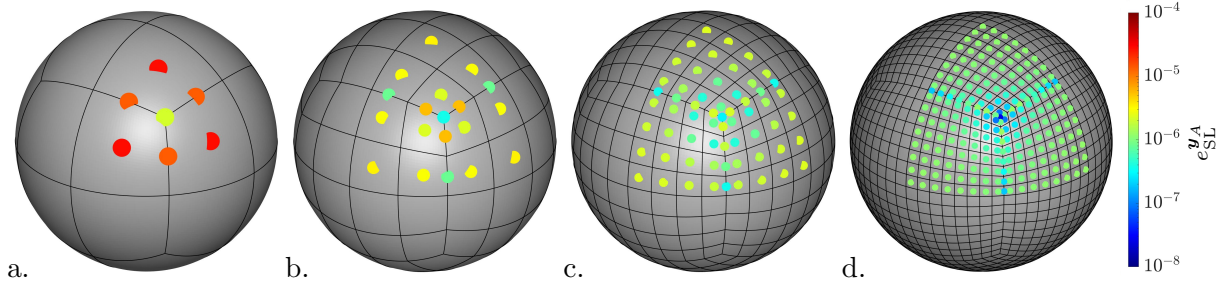
\begin{figure}[h]
\unitlength\linewidth
\begin{picture}(1,.23)
	\put(0,0){\includegraphics[trim = 100 133 233 67, clip, width=.23\linewidth ]{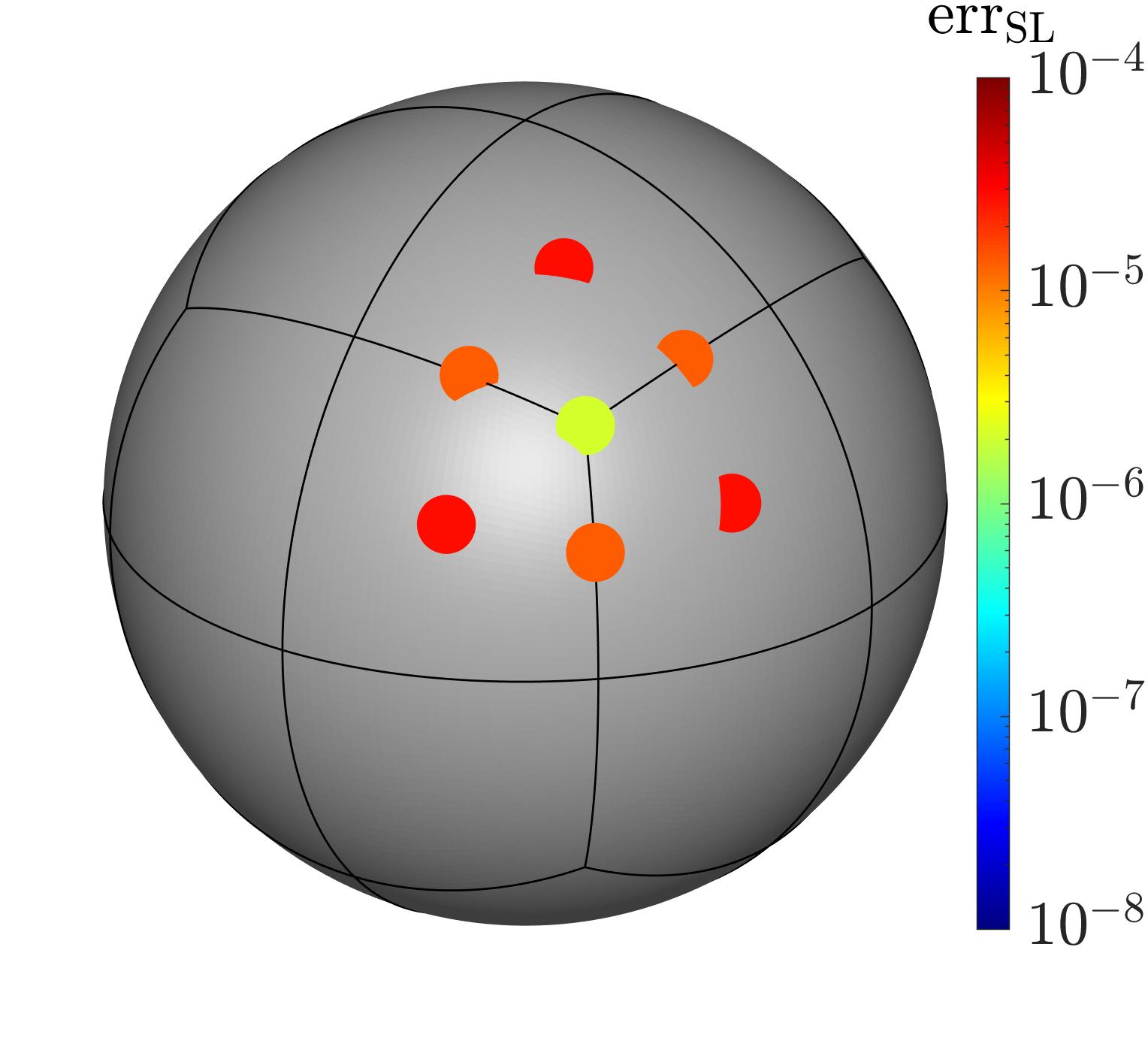}}	
	\put(0.23,0){\includegraphics[trim = 100 133 233 67, clip, width=.23\linewidth ]{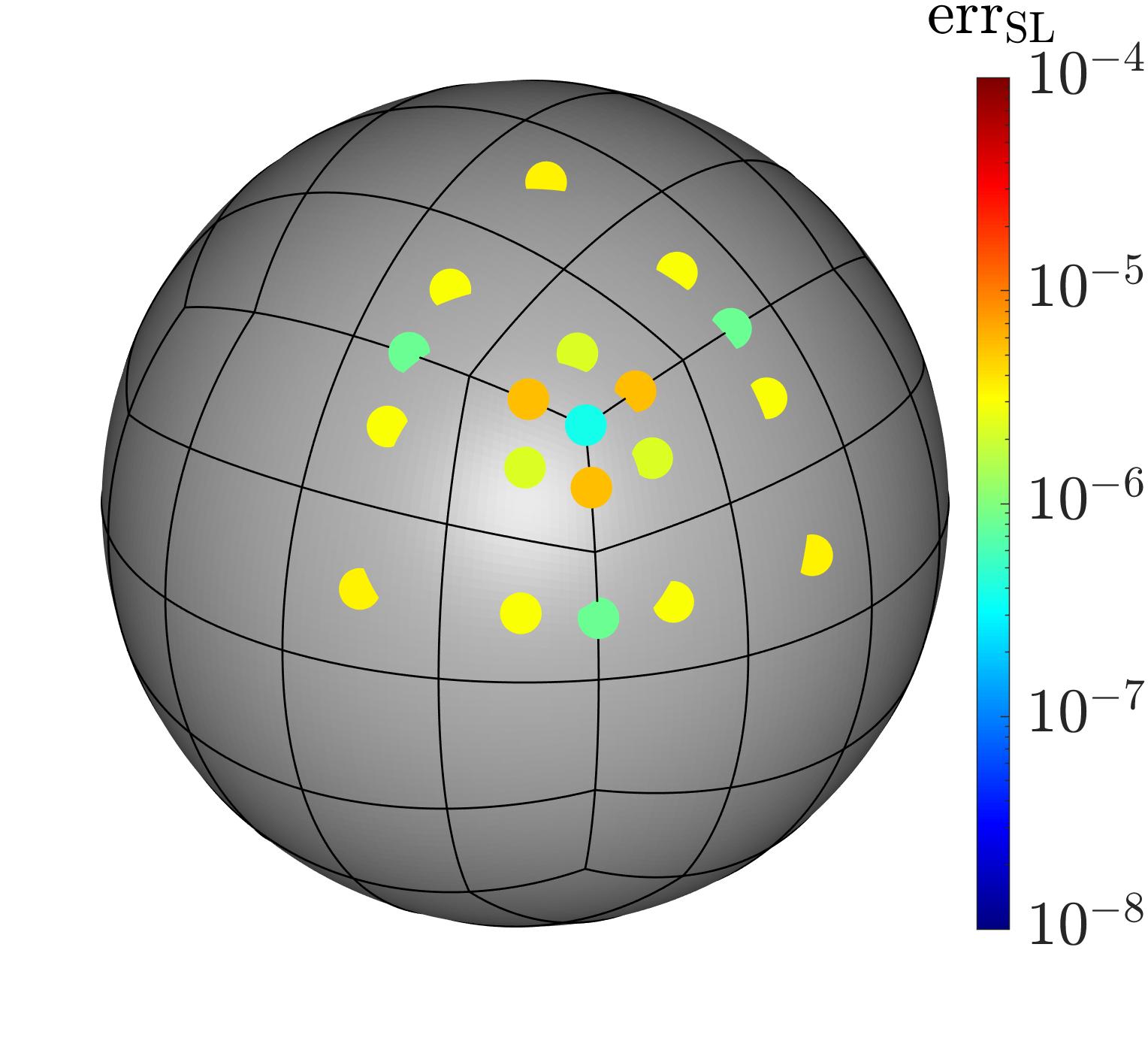}}	
	\put(0.46,0){\includegraphics[trim = 100 133 233 67, clip, width=.23\linewidth ]{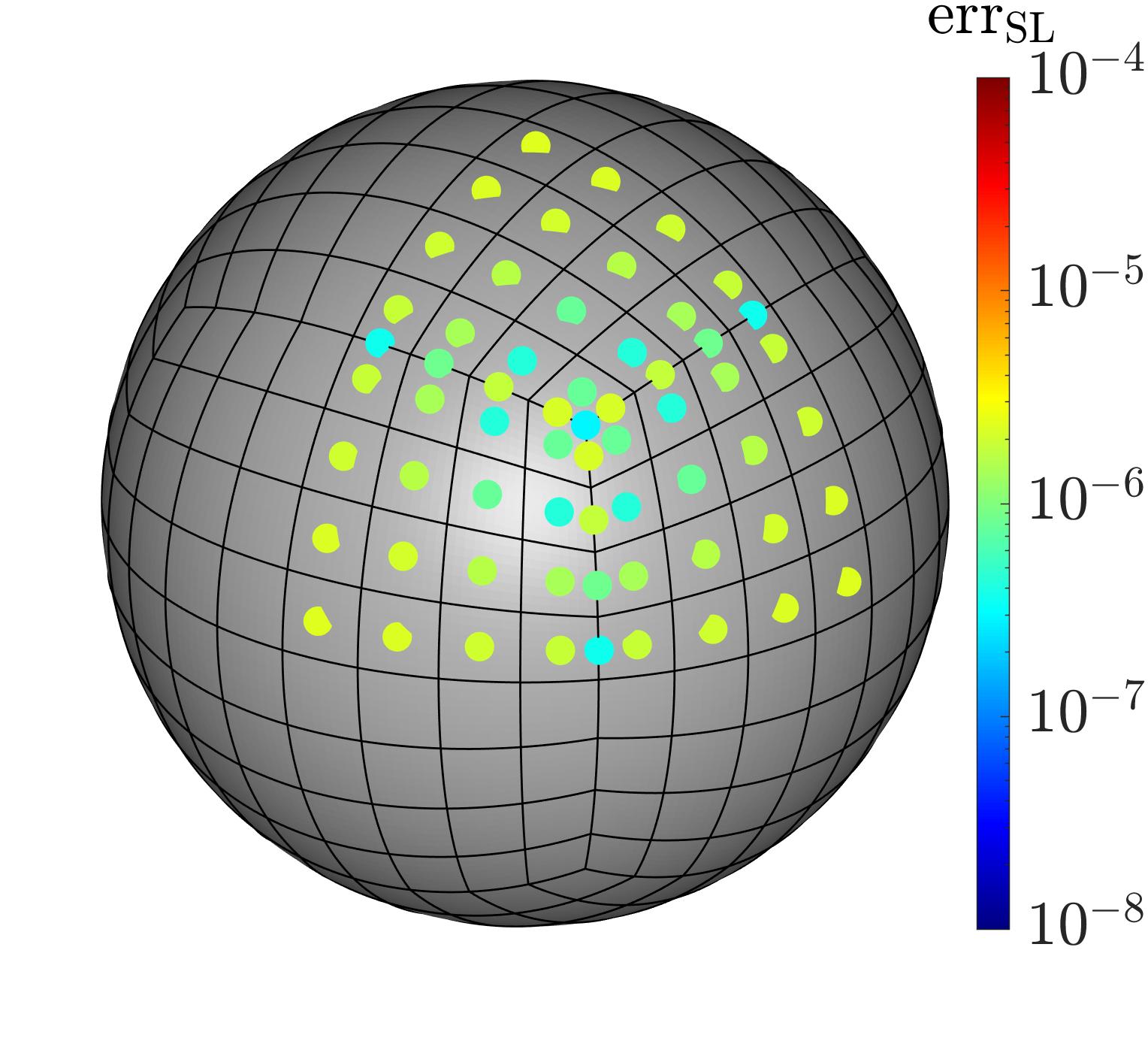}}	
	\put(0.69,0){\includegraphics[trim =100 133 233 67, clip, width=.23\linewidth ]{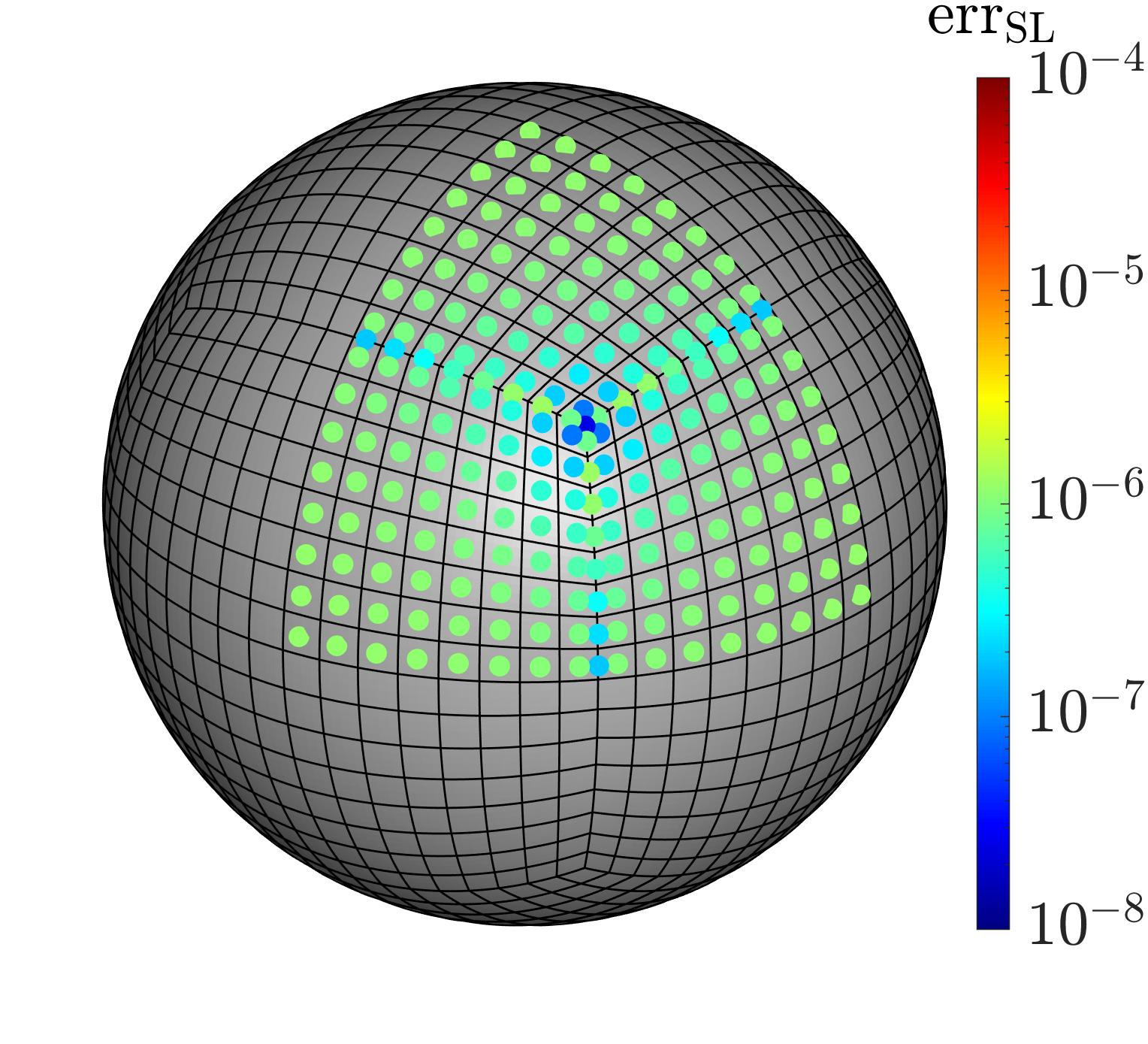}}	
	\put(0.925,0){\includegraphics[trim = 0 0 0 0, clip, height= 0.23\linewidth ]{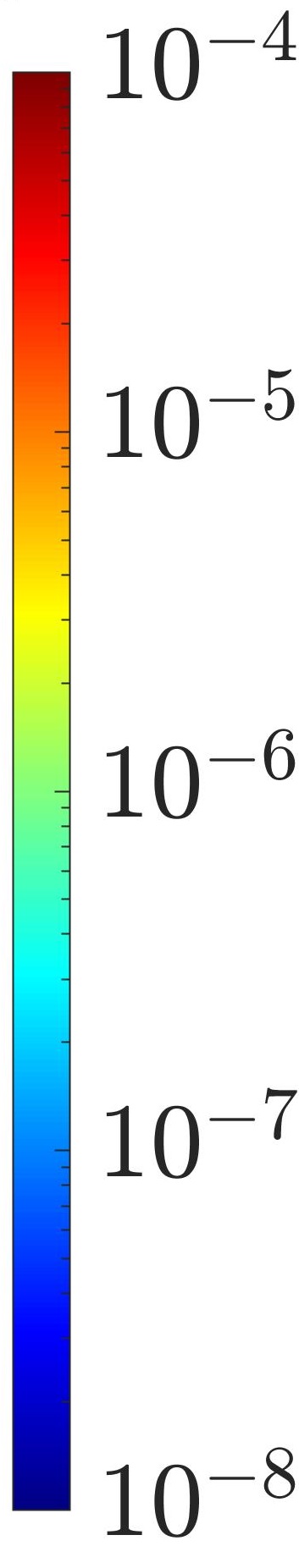}}	
	\vwput{0.97}{0.09}{\scriptsize $e_\mathrm{SL}^{\by_A}$} 
 	\put(0,0){a.}\put(0.23,0){b.}\put(0.46,0){c.}\put(0.69,0){d.}
	\end{picture}
\caption{\textit{Hybrid quadrature on a six-patch NURBS sphere}: Quadrature error $e_\mathrm{SL}^{\by_A}$ \eqref{eq:err_SL} for DGr w.r.t.~identity \eqref{eq:identity_SL} considering each collocation point on one octant of the sphere for $n_0=3$ and discretization level $\ell=1$~(a.), $\ell=2$~(b.),  $\ell=3$~(c.) and  $\ell=4$~(d.).}\label{fig:quad_sphere_error_col_6p}
\end{figure}
\\Fig.~\ref{fig:quad_sphere_error_col_6p} shows the quadrature error w.r.t.~BE identity \eqref{eq:identity_SL} for all collocation points on one quadrant of six-patch spheres of refinement level $\ell=1$~(a.), $\ell=2$~(b.), $\ell=3$~(c.) and $\ell=4$~(d.) using quadrature scheme DGr with $n_0=3$. It can be seen that the quadrature error decreases with increasing quadrature refinement for all collocation points, which was not the case for the single-patch sphere (cf.~Fig.~\ref{fig:quad_sphere_error_col}). The mean quadrature error \eqref{eq:err_SL2} with respect to identities \eqref{eq:identity_SL} and \eqref{eq:identity_DL} is shown in Fig.~\ref{fig:sphere_convergence_6p_m} for $n_0=3$ and varying mesh refinement $\ell=1,\ldots,5$. Both quadrature errors decrease monotonically with increasing mesh refinement for each of the four quadrature scheme. DGr provides by far the best result, followed by DGw and DG and eventually G, which is still quite inaccurate even for highly refined meshes. In contrast to $e_\mathrm{SL}$~(Fig.~\ref{fig:sphere_convergence_6p_m}a), DGw does not yield a significant improvement over DG for $e_\mathrm{DL}$~(Fig.~\ref{fig:sphere_convergence_6p_m}b).

\subsection{Conclusions on hybrid quadrature}\label{sec:hybrid_conclusion}
The suitability of the presented hybrid quadrature schemes for BE analysis has been demonstrated in the previous sections: Sec.~\ref{sec:hybrid_sheet} shows that hybrid Duffy-Gauss quadrature with adjusted weights (DGw) is by far the most accurate and efficient scheme for the approximation of singular integrals on flat surfaces, followed by Duffy-Gauss quadrature with progressive refinement~(DGr). Sec.~\ref{sec:hybrid_sphere} shows that DGr is by far the most efficient scheme on curved surfaces. However, the second most efficient scheme DGw still offers a slight improvement in accuracy over Duffy-Gauss quadrature without special treatment of the near singular elements (DG).
\\\\It has been further shown that the hybrid quadrature schemes approximate the singular BE integrals very efficiently on coarse single-patch NURBS spheres (Sec.~\ref{sec:hybrid_sphere_single}). These discretizations are exactly spherical and therefore recommended for BE analysis on rigid spheres or spheroids.
The application of hybrid quadrature schemes to the six-patch NURBS sphere (Sec.~\ref{sec:hybrid_sphere_six}) also yields accurate results for coarse discretizations, but moreover a significantly improved convergence behavior for mesh refinement. The six-patch discretization is therefore recommended for more complex surface geometries and deforming surfaces.

\section{Application to Stokes flow problems}\label{sec:examples}
The suitability of the introduced quadrature schemes for BE analysis is investigated here with three numerical examples: Sec.~\ref{sec:ex_rot} and Sec.~\ref{sec:ex_trans} consider spheres rotating in a viscous fluid and translating through a viscous fluid, respectively.\footnote{A fixed sphere in a steady rotational flow leads to a mathematically equivalent problem as a sphere rotating in a quiescent fluid. A fixed sphere in a steady transversal flow leads to a mathematically equivalent problem as a sphere translating through a quiescent fluid.} A rising bubble of non-spherical shape is considered as a third problem in Sec.~\ref{sec:ex_bubble}. In all three problems, the velocity on the surface is given by Dirichlet boundary conditions and the traction on the surface is determined by BE analysis. The numerical traction results are compared to analytical results based on \cite{Chwang75} and \cite{Kong12}. Although this paper focuses on Stokes flow, the presented quadrature schemes are also expected to be applicable to other BE formulations including those for linear elasticity and the Helmholtz equation. \cite{Venaas20} show that Duffy quadrature is capable to deal with the oscillatory nature of Helmholtz kernels. The adjusted quadrature weights incorporated into scheme DGw can be determined by replacing the Stokes kernel in \eqref{eq:quad_adjusted_condition} by any weakly singular kernel.

\subsection{Flow caused by a rotating sphere}\label{sec:ex_rot}
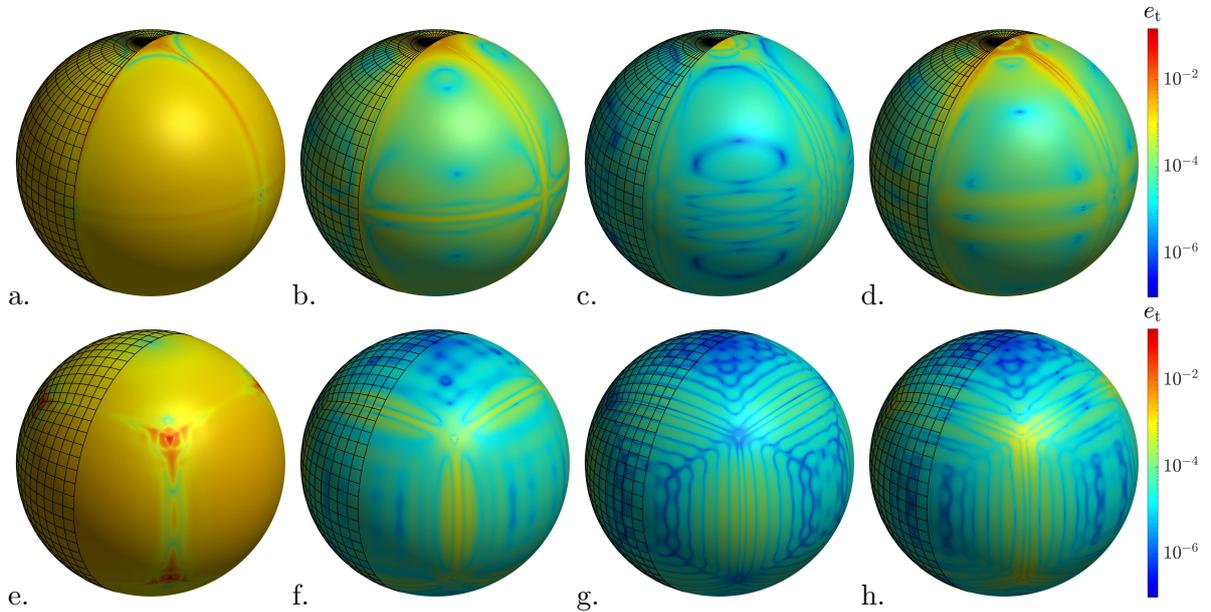
\begin{figure}[h]
\unitlength\linewidth
\begin{picture}(1,0.5)
\put(0,0.25){\includegraphics[trim = 75 120 240 0, clip, height=.25\linewidth]{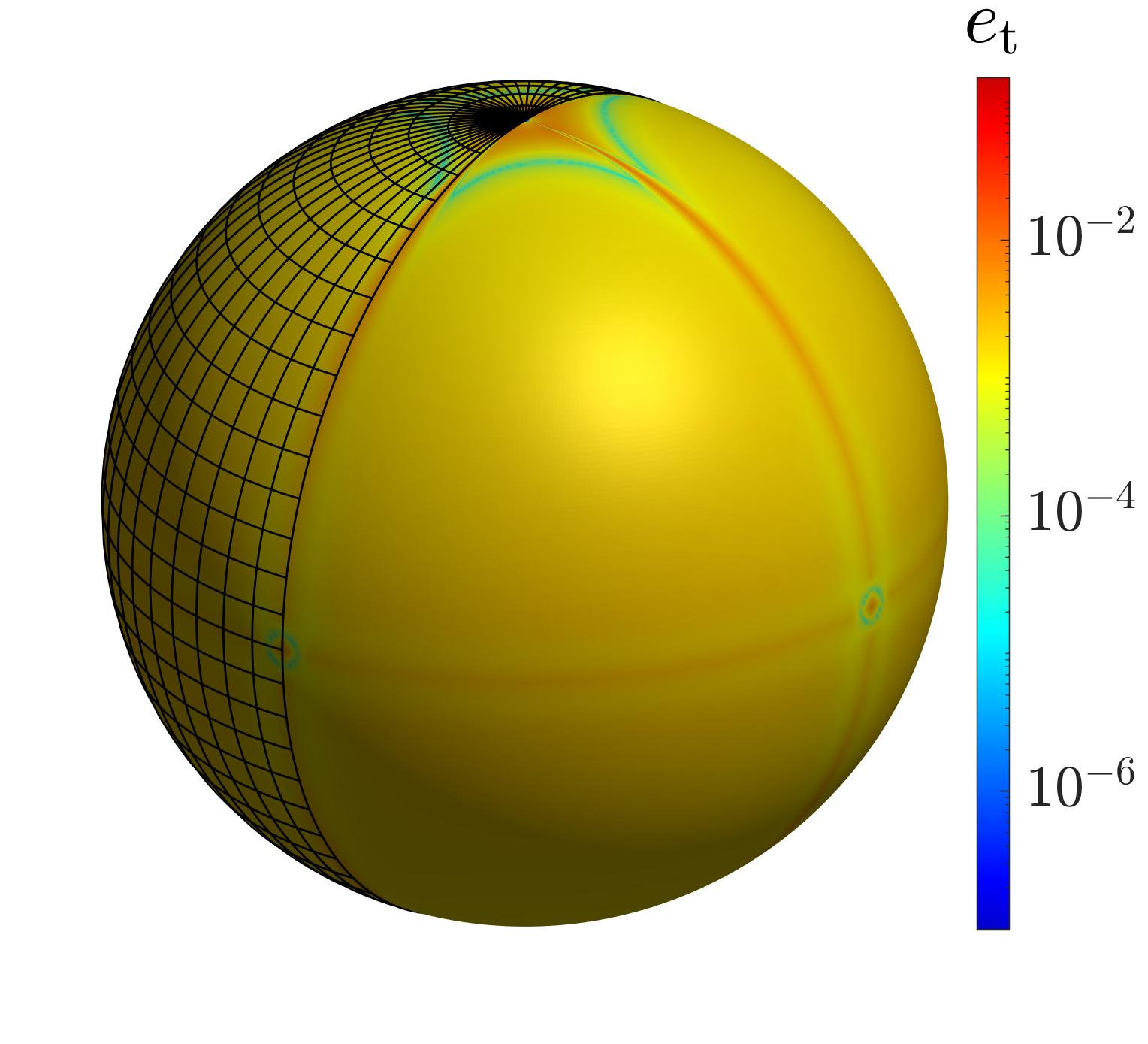}}
\put(0.235,0.25){\includegraphics[trim = 75 120 240 0, clip, height=.25\linewidth]{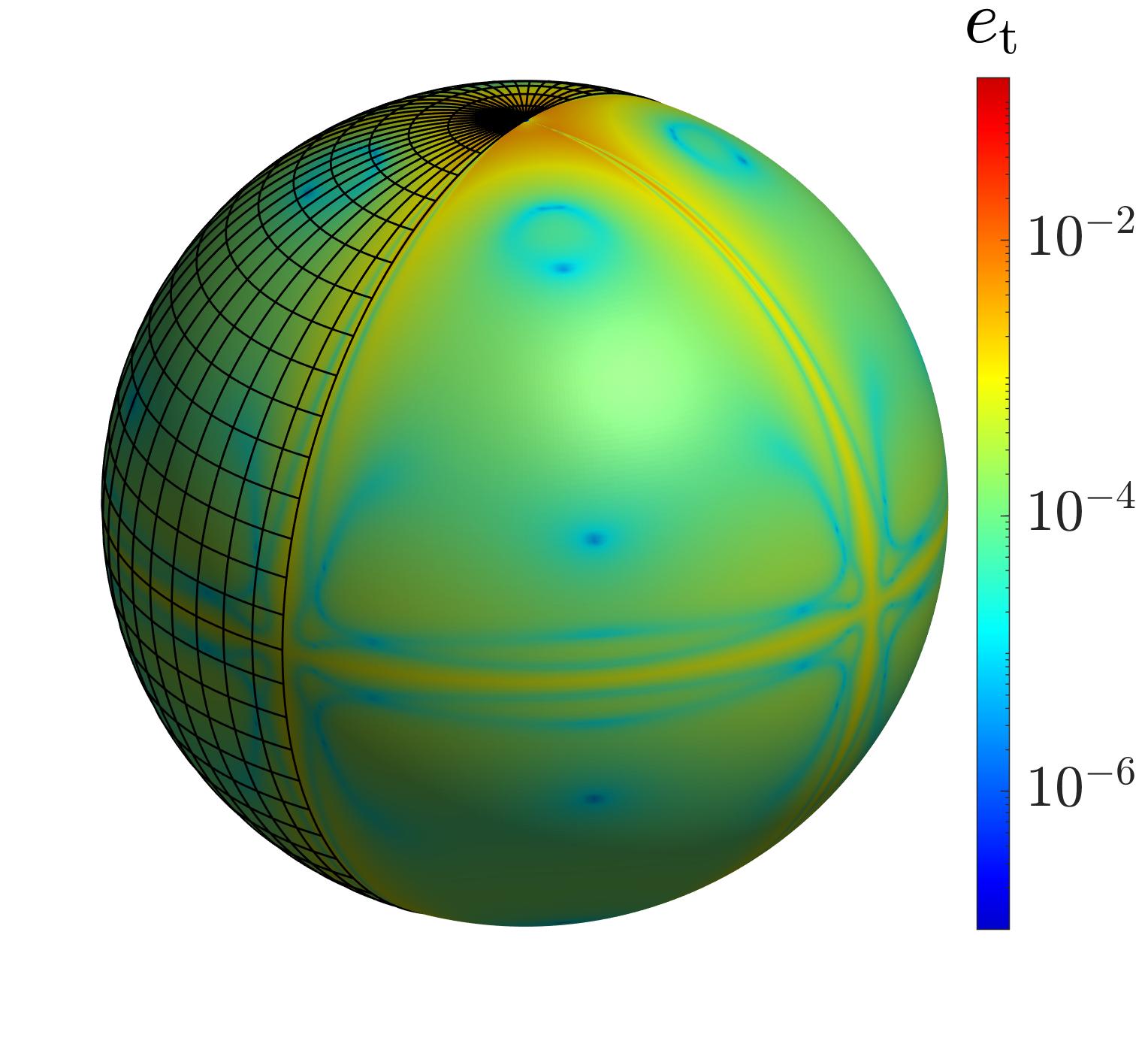}}
\put(0.47,0.25){\includegraphics[trim = 75 120 240 0, clip, height=.25\linewidth]{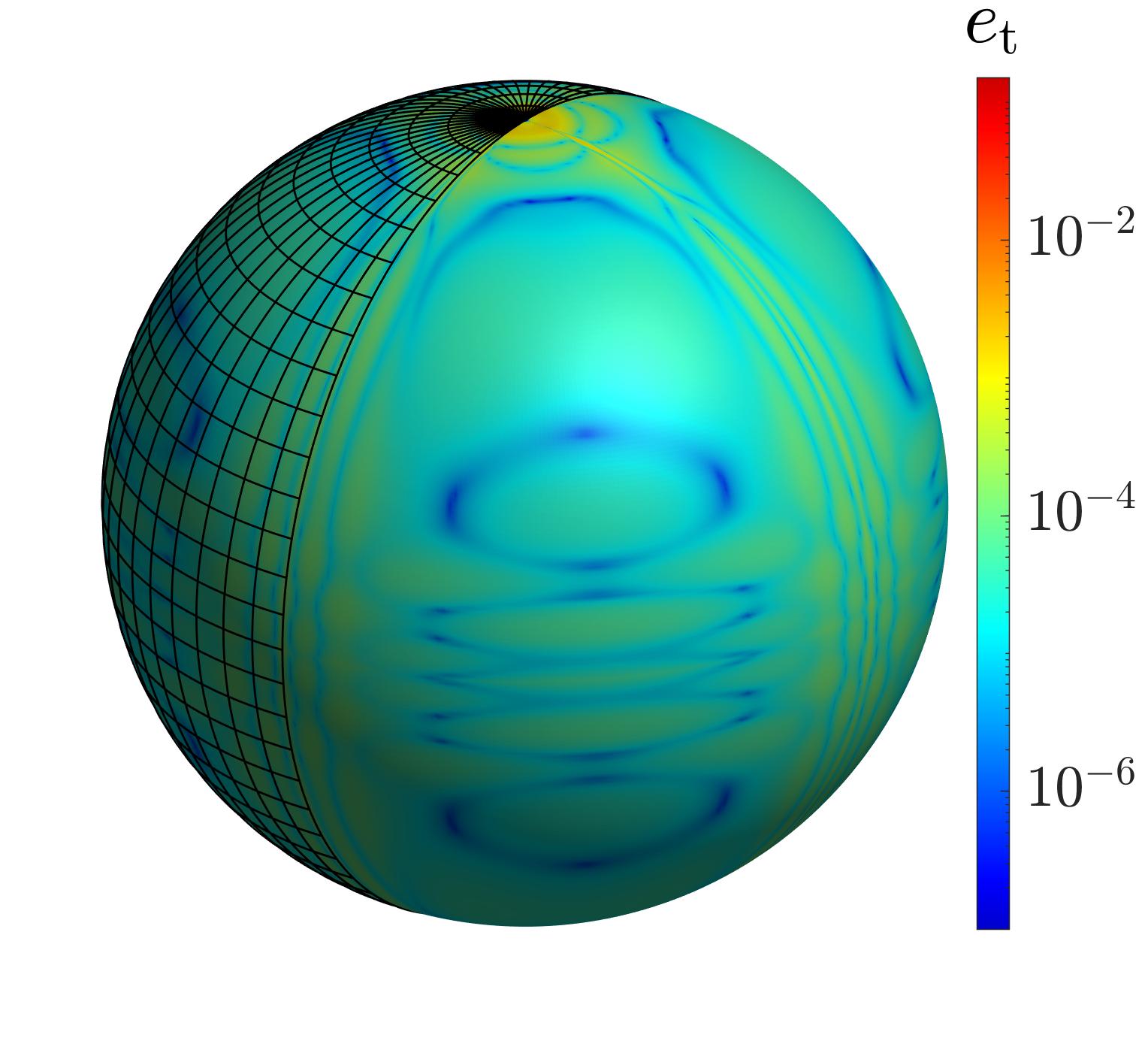}}
\put(0.705,0.25){\includegraphics[trim = 75 120 0 0, clip,height=.25\linewidth]{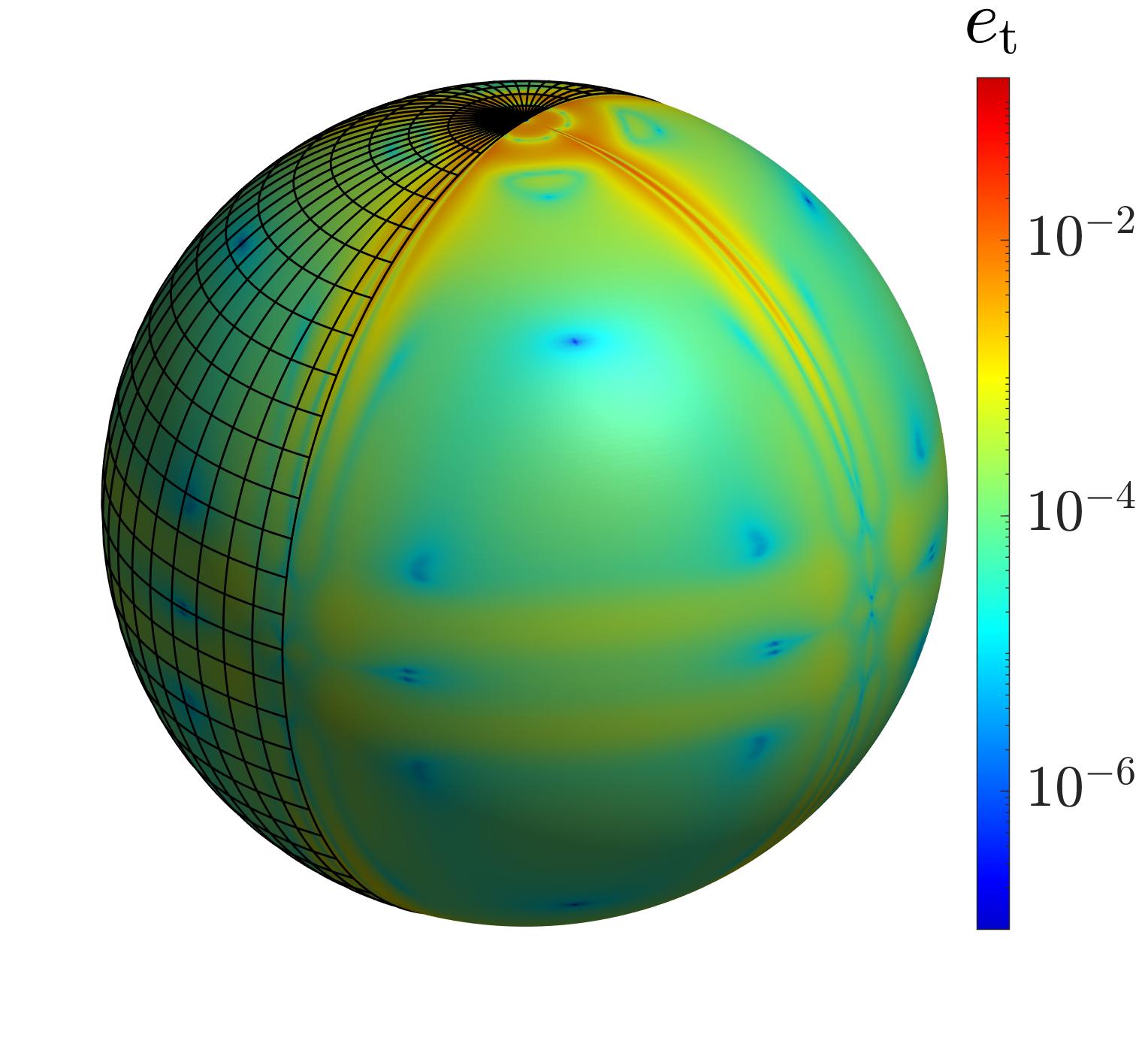}}
\put(0,0){\includegraphics[trim = 75 120 240 0, clip, height=.25\linewidth]{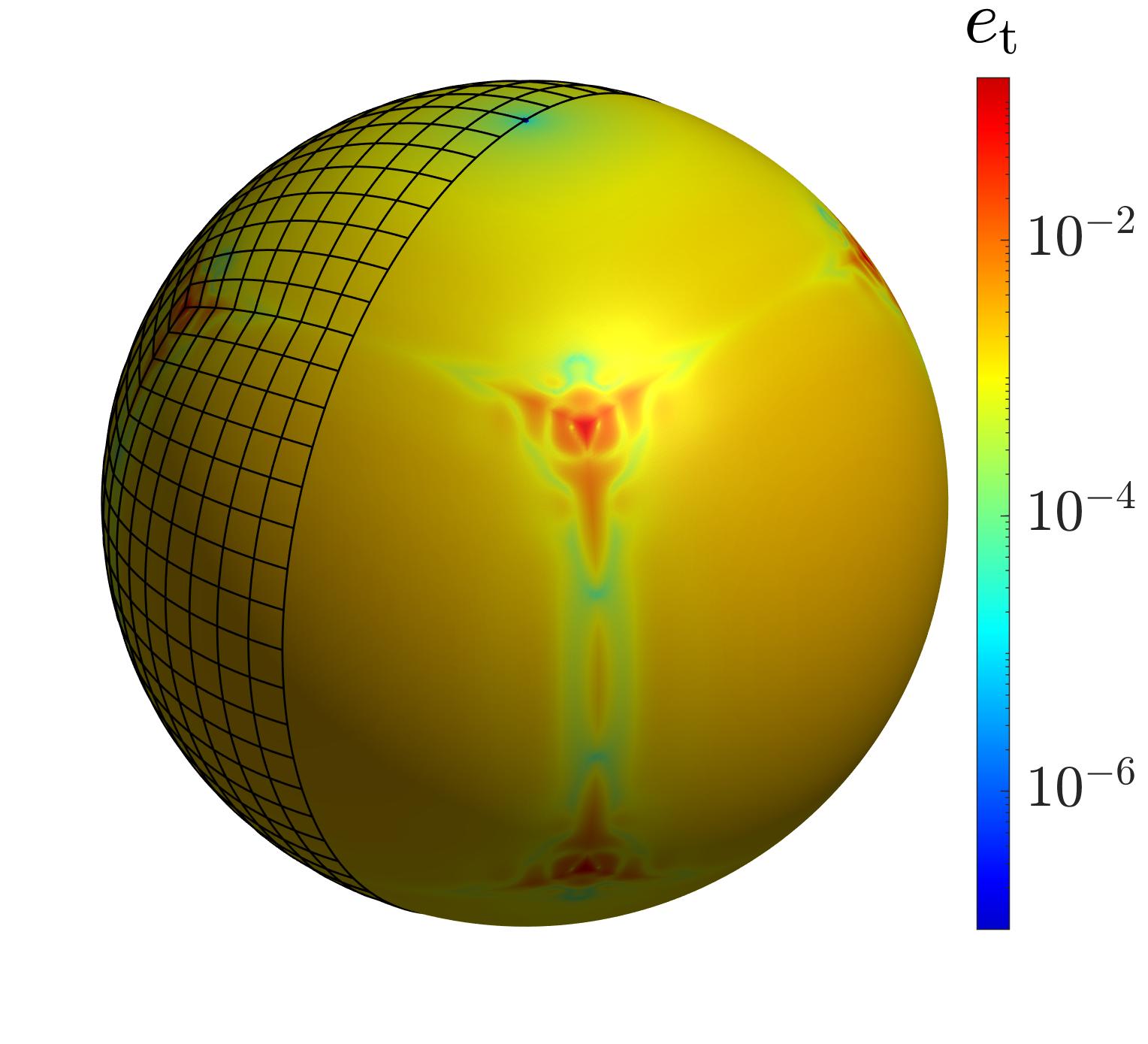}}
\put(0.235,0){\includegraphics[trim = 75 120 240 0, clip, height=.25\linewidth]{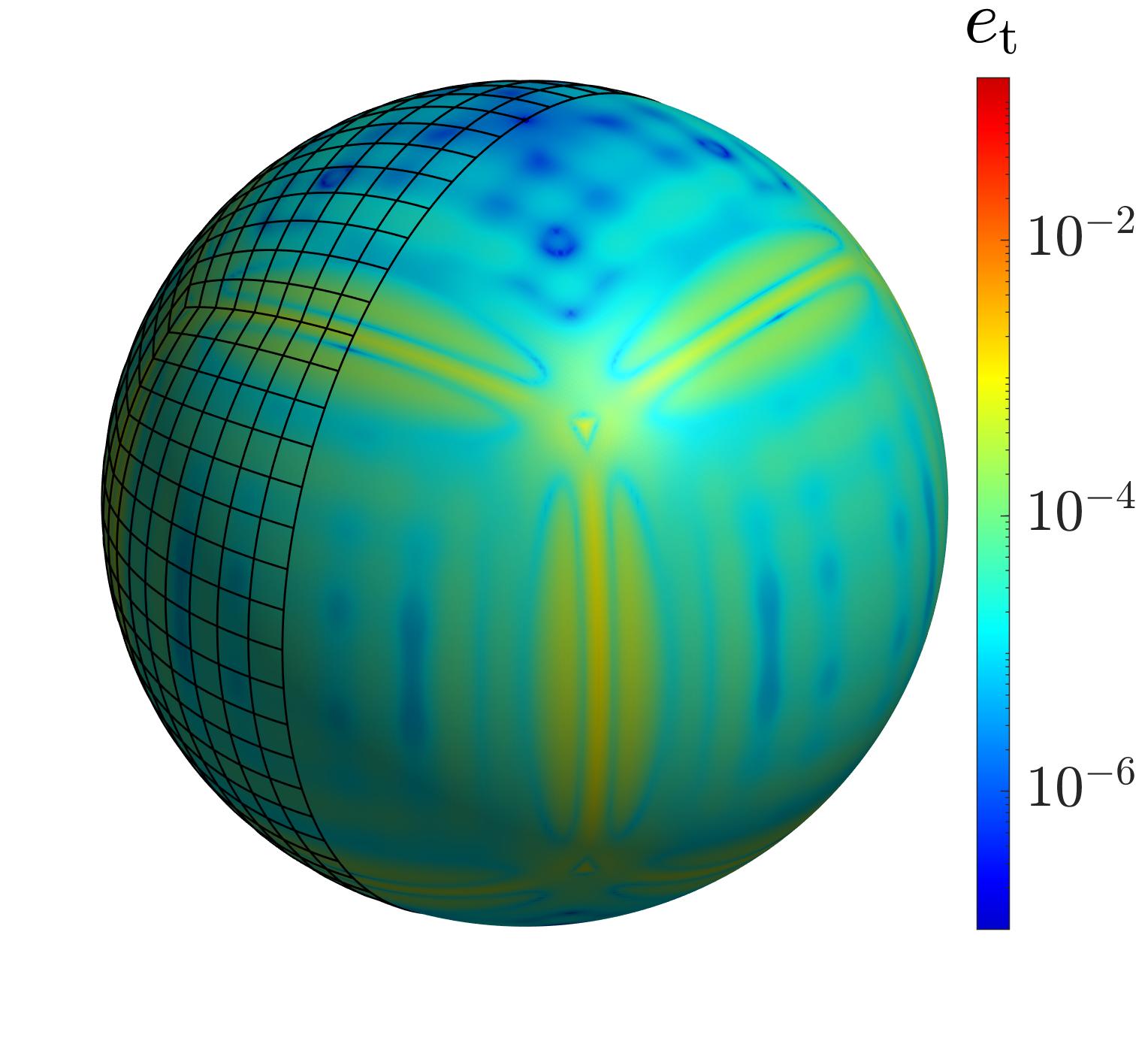}}
\put(0.47,0){\includegraphics[trim = 75 120 240 0, clip, height=.25\linewidth]{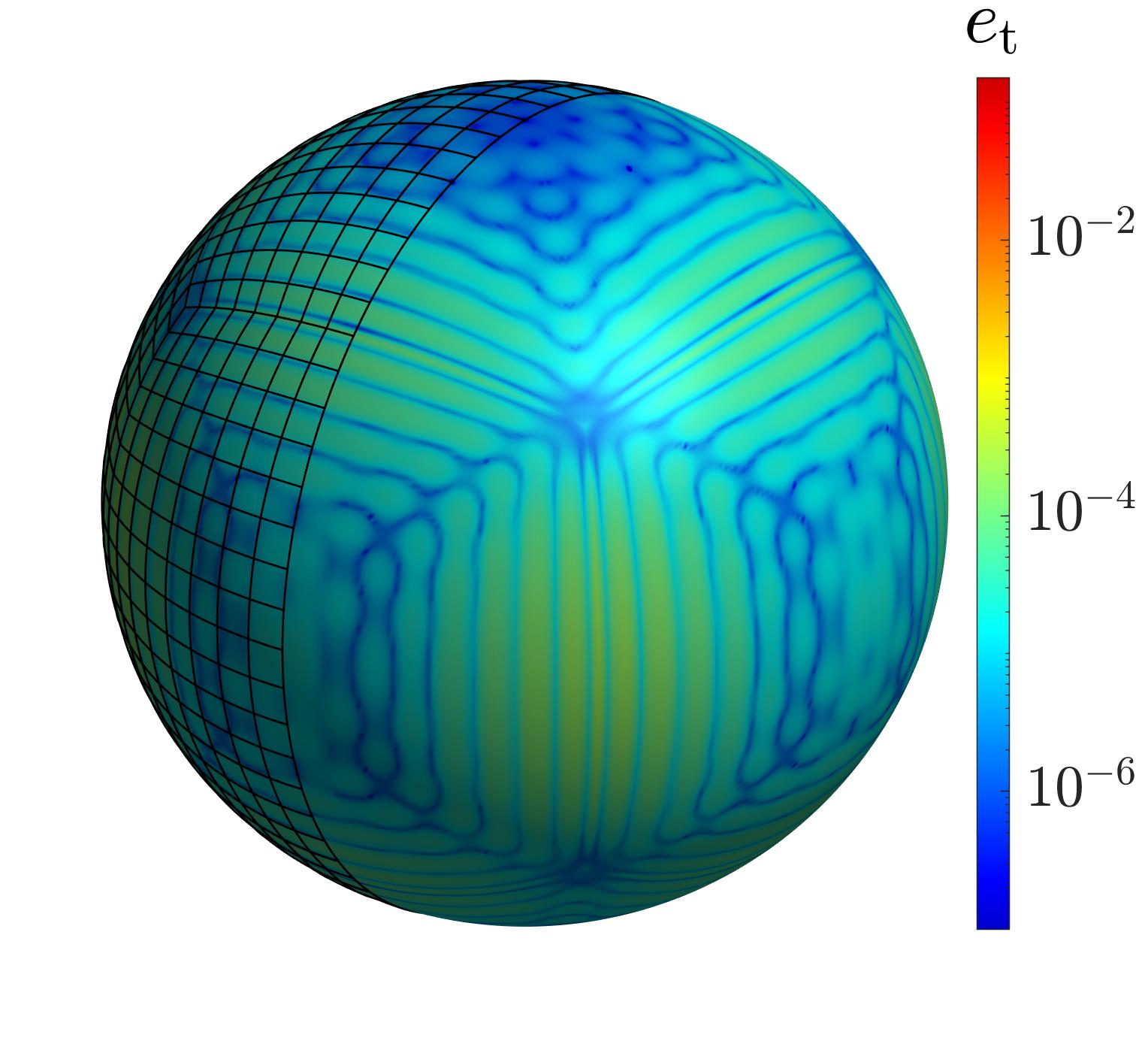}}
\put(0.705,0){\includegraphics[trim = 75 120 0 0, clip,height=.25\linewidth]{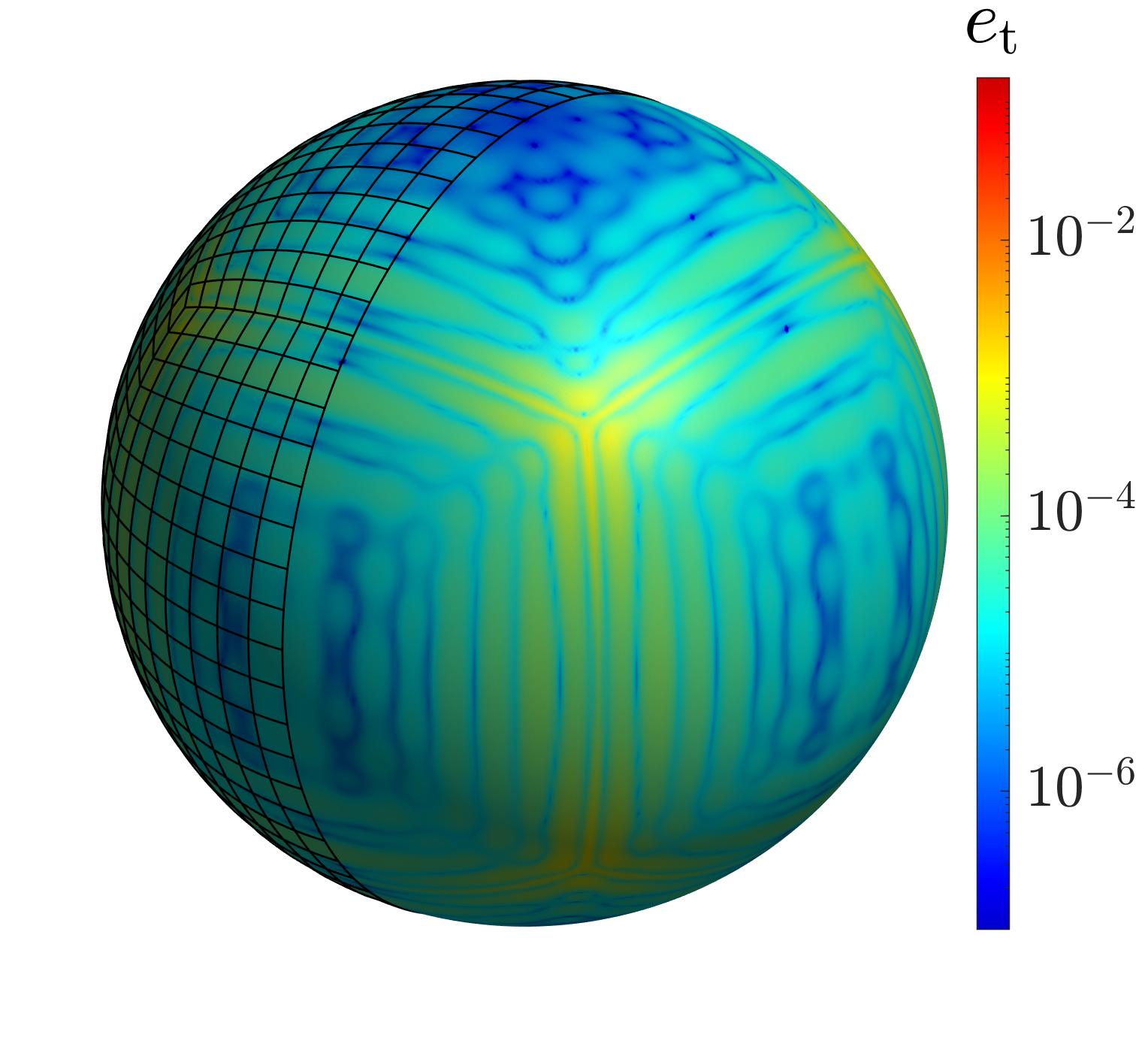}}
\put(0.005,0.25){a.} \put(0.24,0.25){b.} \put(0.475,0.25){c.} \put(0.71,0.25){d.}
\put(0.005,0){e.} \put(0.24,0){f.} \put(0.475,0){g.} \put(0.71,0){h.}
\end{picture}	
\caption{\textit{Rotating sphere:} Traction error for the four hybrid quadrature schemes with $n_0=3$. Single-patch NURBS sphere of refinement level $\ell=4$: G~(a.~$\bar n_\mathrm{qp}=18,440$), DG~(b.~$\bar n_\mathrm{qp}=18,714$), DGr~(c.~$\bar n_\mathrm{qp}=19,029$) and DGw~(d.~$\bar n_\mathrm{qp}=18,816$). Six-patch NURBS sphere of refinement level $\ell=4$: G~(e.~$\bar n_\mathrm{qp}=13,832$), DG~(f.~$\bar n_\mathrm{qp}=14,102$), DGr~(g.~$\bar n_\mathrm{qp}=14,323$) and DGw~(h.~$\bar n_\mathrm{qp}=14,101$).}\label{fig:ex_rot_surf}
\end{figure}

\begin{figure}[h]
\unitlength\linewidth
\begin{picture}(1,0.41)
\put(0,0){\includegraphics[trim = 0 0 30 20, clip, width=.5\linewidth] {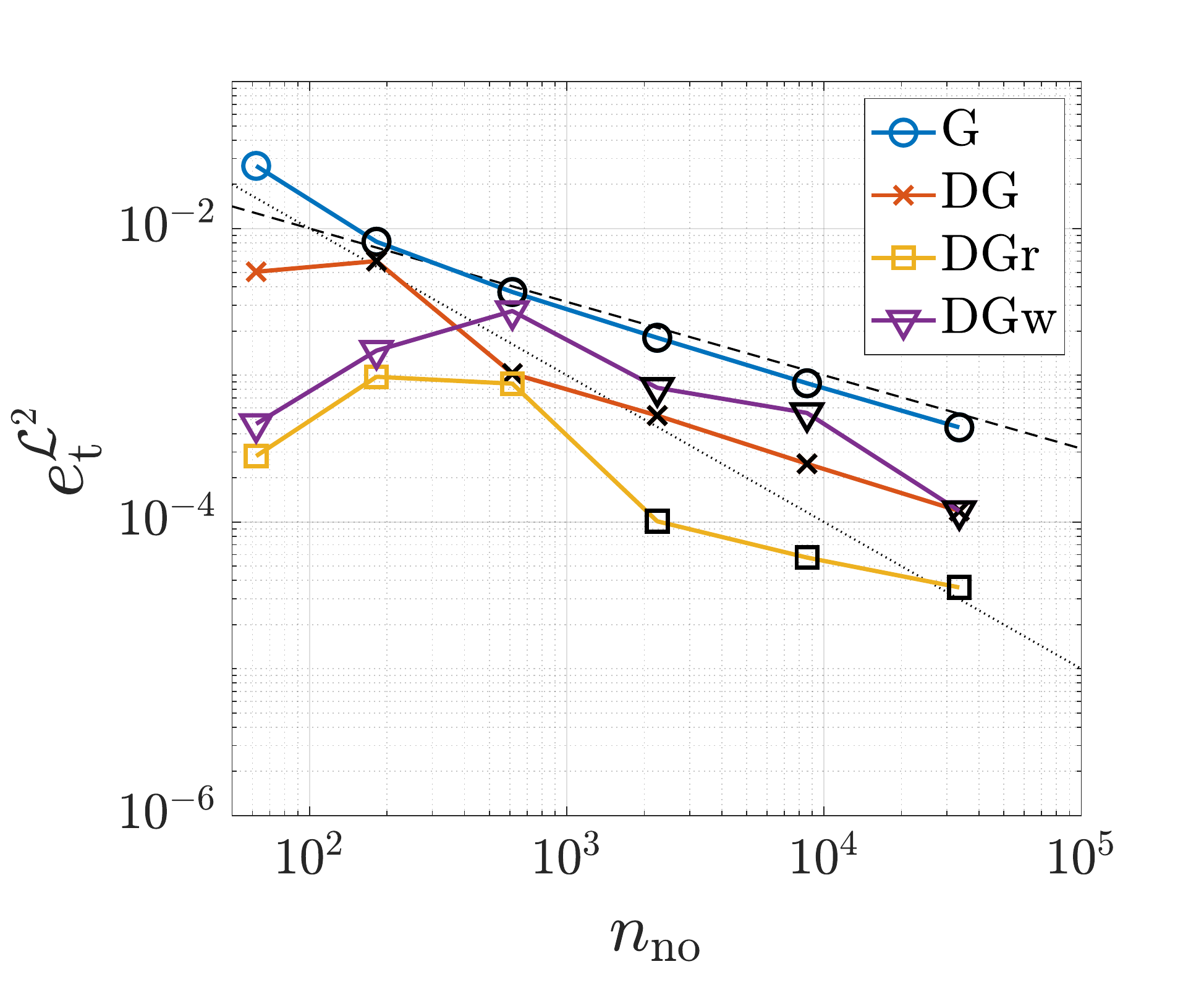}}
\put(0.25,0.31){\rotatebox{-15}{\footnotesize $\propto 1/\sqrt{n_\mathrm{no}}$}}
\put(0.395,0.16){\rotatebox{-30}{\footnotesize $\propto 1/n_\mathrm{no}$}}
\put(0.5,0){\includegraphics[trim = 0 0 30 20, clip, width=.5\linewidth]{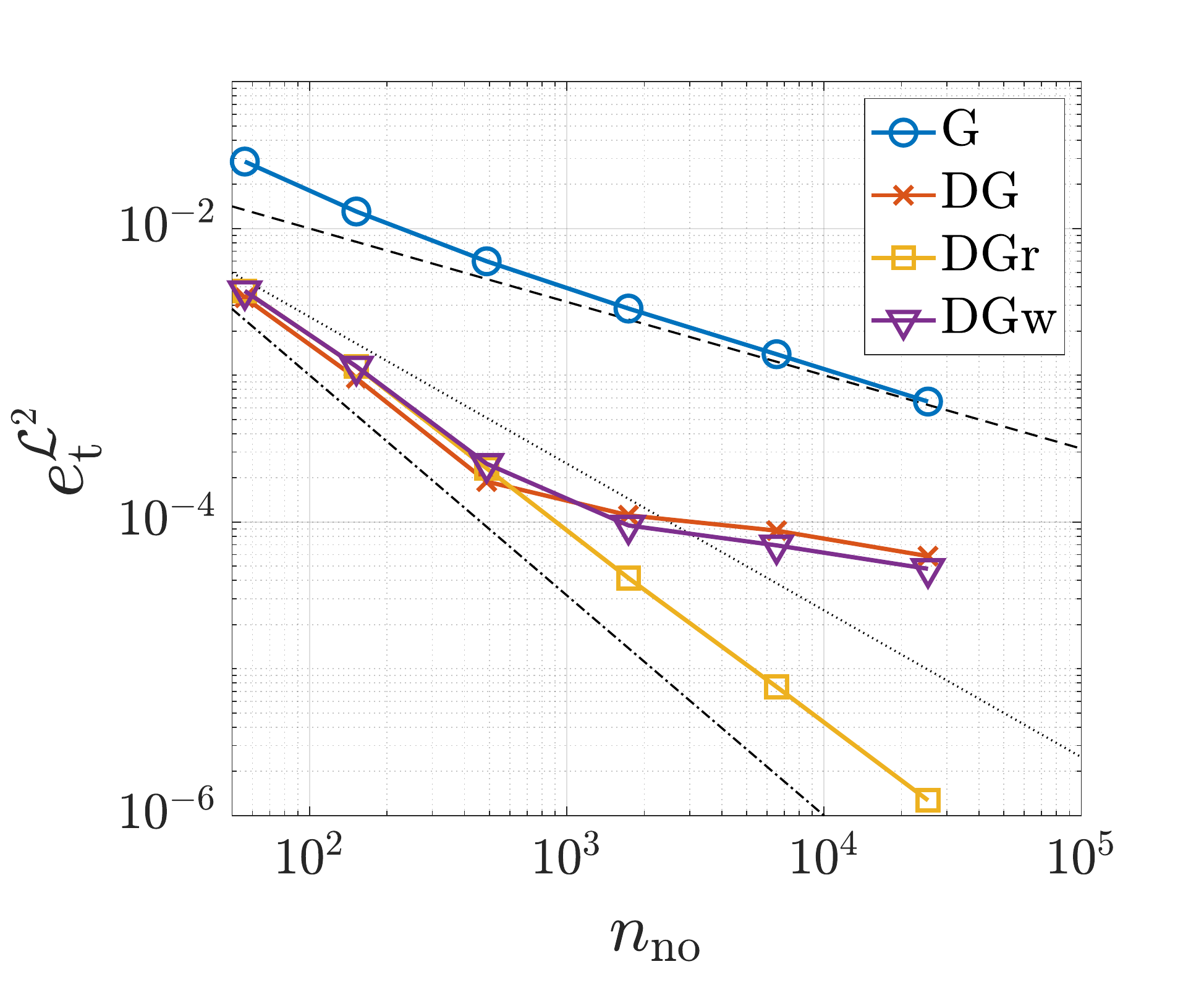}}
\put(0.875,0.24){\rotatebox{-15}{\footnotesize $\propto 1/\sqrt{n_\mathrm{no}}$}}
\put(0.86,0.14){\rotatebox{-30}{\footnotesize $\propto 1/n_\mathrm{no}$}}
\put(0.7,0.18){\rotatebox{-43}{\footnotesize $\propto 1/n_\mathrm{no}^{1.5}$}}
\put(0.01,0.01){a.} \put(0.51,0.01){b.}
\end{picture}
\caption{\textit{Rotating sphere:} Traction error $e_t^{\mathcal L^2}$ for quadrature density $n_0=3$ on the single-patch NURBS sphere~(a) and on the six-patch NURBS sphere~(b), both of varying refinement level ${\ell=1,\ldots,6}$. The black markers depict results from an iterative solver.}\label{fig:sphere_rot_mesh_res2}
\end{figure}
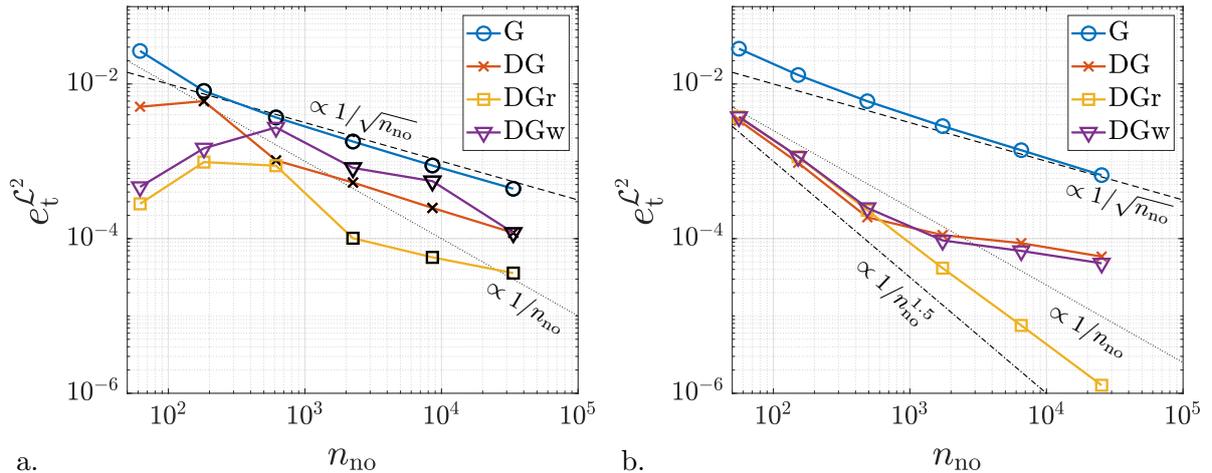

The first example considers a rigid sphere (surface $\mcalS$, Radius $R$) that is surrounded by a viscous fluid $\mcalF$ (dynamic viscosity $\eta$). The sphere rotates around its center with the prescribed angular velocity $\bar\bome \in \mathbb{R}^3$, so that the surface velocity is given by
\begin{equation}\label{eq:ex_rot_BC}
\bv(\bx)= \bar\bome\times \bx, \hh \forall\; \bx \in \mcalS~.
\end{equation}
A rotation around the vertical axis is considered here, i.e.~$\bar\bome=\bar \omega\, \be_3$. The BIE can be solved for the surface traction $\bt:=\bsig\bn$ considering \eqref{eq:ex_rot_BC} as a Dirichlet boundary condition. Due to the pure rotational velocity, the traction field has only tangential components, while the normal traction (i.e.~the pressure) is zero.
\\\\The BE traction error
\begin{equation}\label{eq:ex_rot_err}
e_\mrt(\bx)=\frac{\| \bt^h(\bx) - \bt(\bx)\|}{t_\mathrm{max}}~,
\end{equation}
with the maximum traction value
\begin{equation}\label{eq:ex_rot_tmax}
t_\mathrm{max} := \max_{\bx \in\mcalS} \|\bt \|= 3\,\eta\,\bar\omega
\end{equation}
located along the equator of the sphere, is introduced to compare the numerical results systematically to the analytical solution. Fig.~\ref{fig:ex_rot_surf} shows the BE traction error on a single-patch NURBS sphere and on a six-patch NURBS sphere, both of mesh refinement level $\ell=4$, considering the introduced quadrature schemes from Sec.~\ref{sec:hybrid}. Half of the symmetric BE meshes is hidden to improve the visibility of the results.
\\\\The accuracy of the introduced hybrid quadrature schemes is investigated systematically in two convergence studies: First, the influence of quadrature refinement is investigated by varying the quadrature density $n_0=2,4,8,16$ on a fixed mesh ($\ell=2$). Second, the influence of mesh refinement is investigated by varying the mesh refinement level $\ell=1,\ldots, 6$, while $n_0=3$. Both convergence studies consider single-patch and six-patch NURBS spheres. The traction error on the surface is characterized by the $\mcalL^2$ norm of the traction error defined as
\begin{equation}\label{eq:ex_rot_L2_dof}
e_\mrt^{\mathcal L^2} = \frac{1}{\sqrt{A_\mcalS}} \sqrt{\int_\mcalS  e_\mrt^2 \,\mrd a}~,
\end{equation}
where $A_\mcalS= 4\pi R^2$ denotes the surface area of $\mcalS$.
\begin{figure}[h]
\unitlength\linewidth
\begin{picture}(1,0.82)
\put(0,0.41){\includegraphics[trim = 0 0 30 20, clip, width=.5\linewidth]{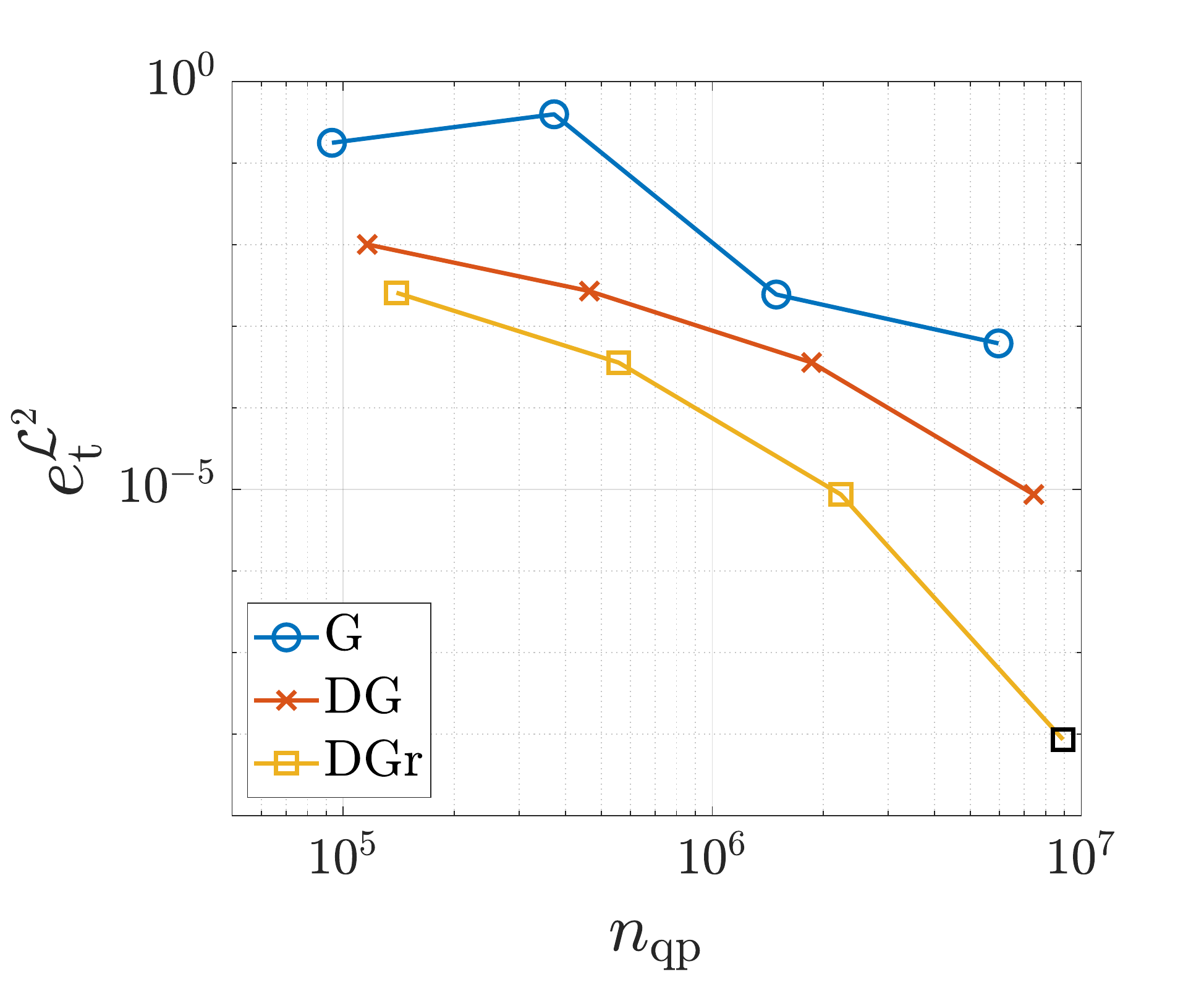}}
\put(0.5,.41){\includegraphics[trim = 0 0 30 20, clip, width=.5\linewidth]{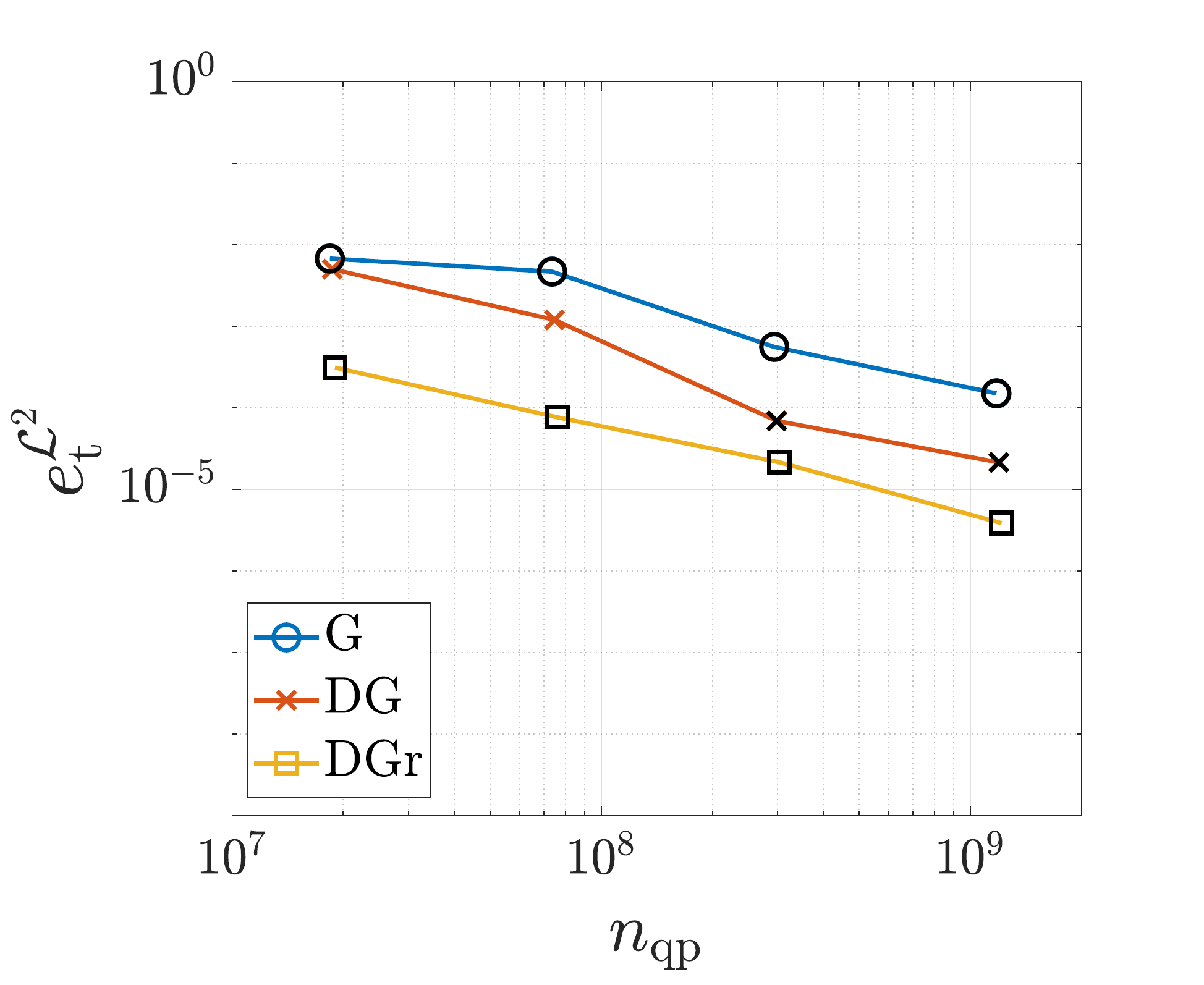}}
\put(0,0){\includegraphics[trim = 0 0 30 20, clip, width=.5\linewidth]{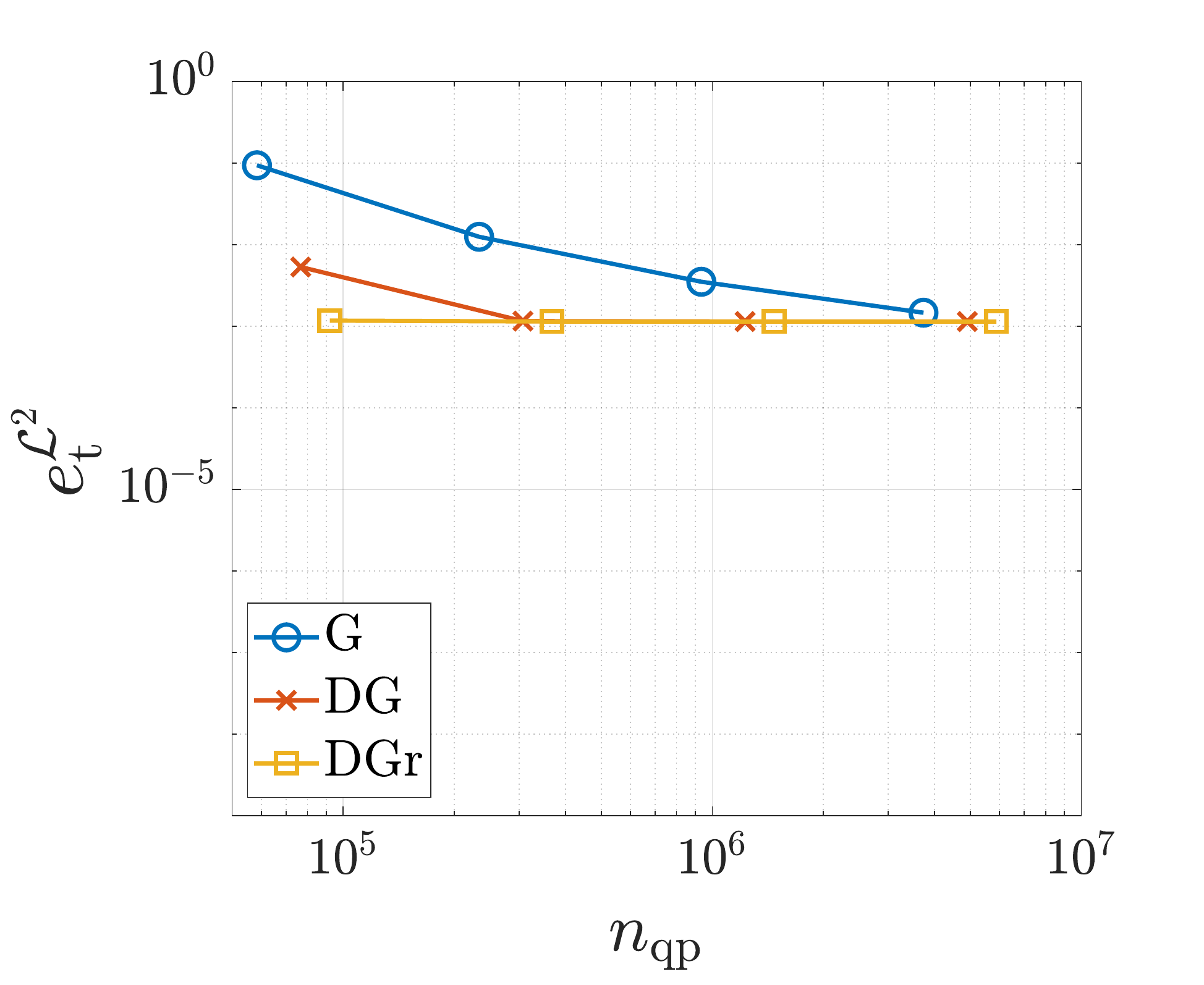}}
\put(0.5,0){\includegraphics[trim = 0 0 30 20, clip, width=.5\linewidth]{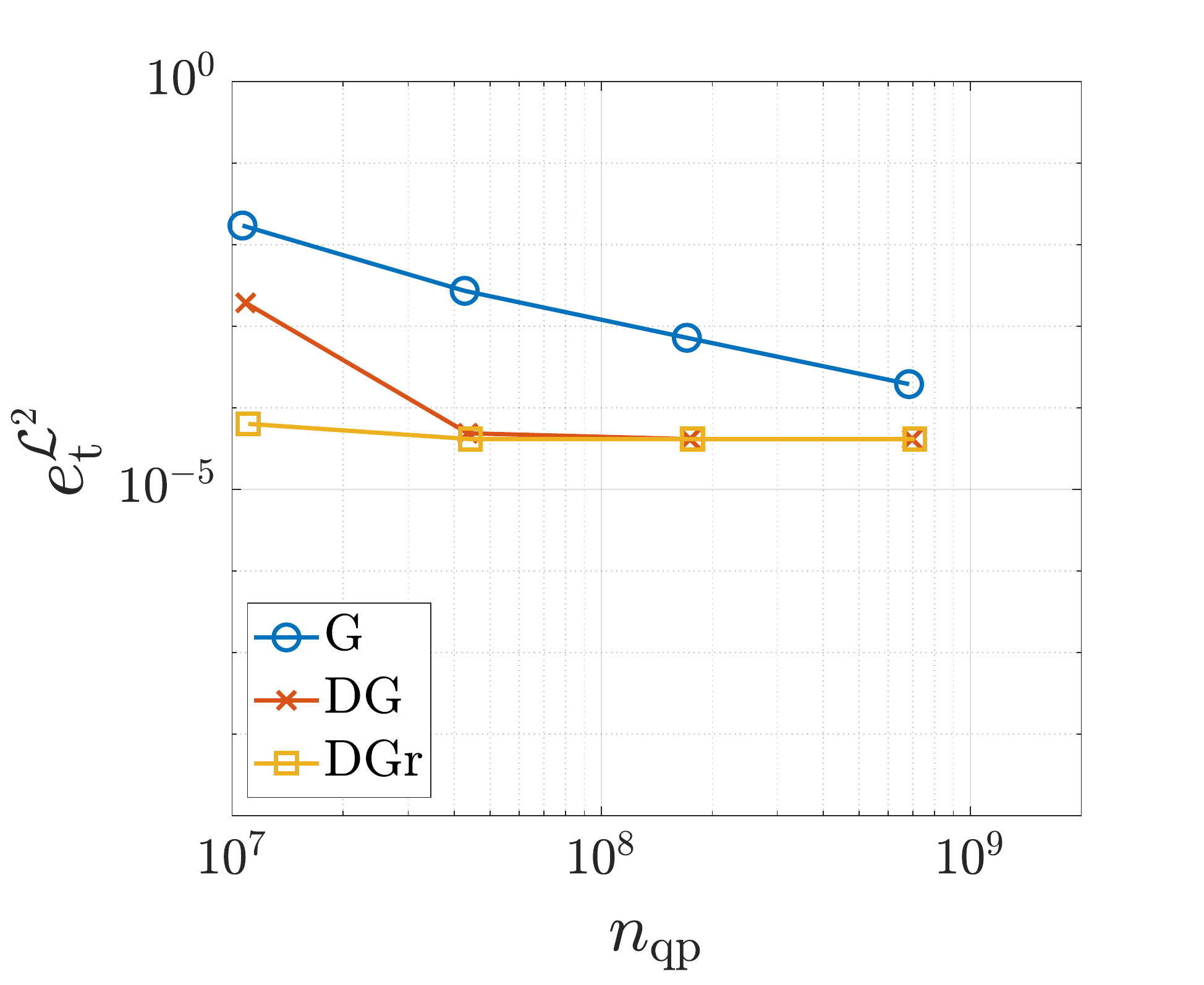}}
\put(0.01,0.42){a.}\put(0.51,0.42){b.}\put(0.01,0.01){c.}\put(0.51,0.01){d.}
\end{picture}
\caption{\textit{Rotating sphere:} Traction error $e_t^{\mathcal L^2}$ for varying quadrature density $n_0=2,4,8,16$ as introduced in Sec.~\ref{sec:hybrid}. The error on the single-patch sphere is shown in (a) for $\ell=2$ and in (b) for $\ell=4$. Black markers indicate where an iterative solver is used. The error on the six-patch sphere is shown in (c) for $\ell=2$ and in (d) for $\ell=4$.}\label{fig:sphere_rot_quad_res}
\end{figure}
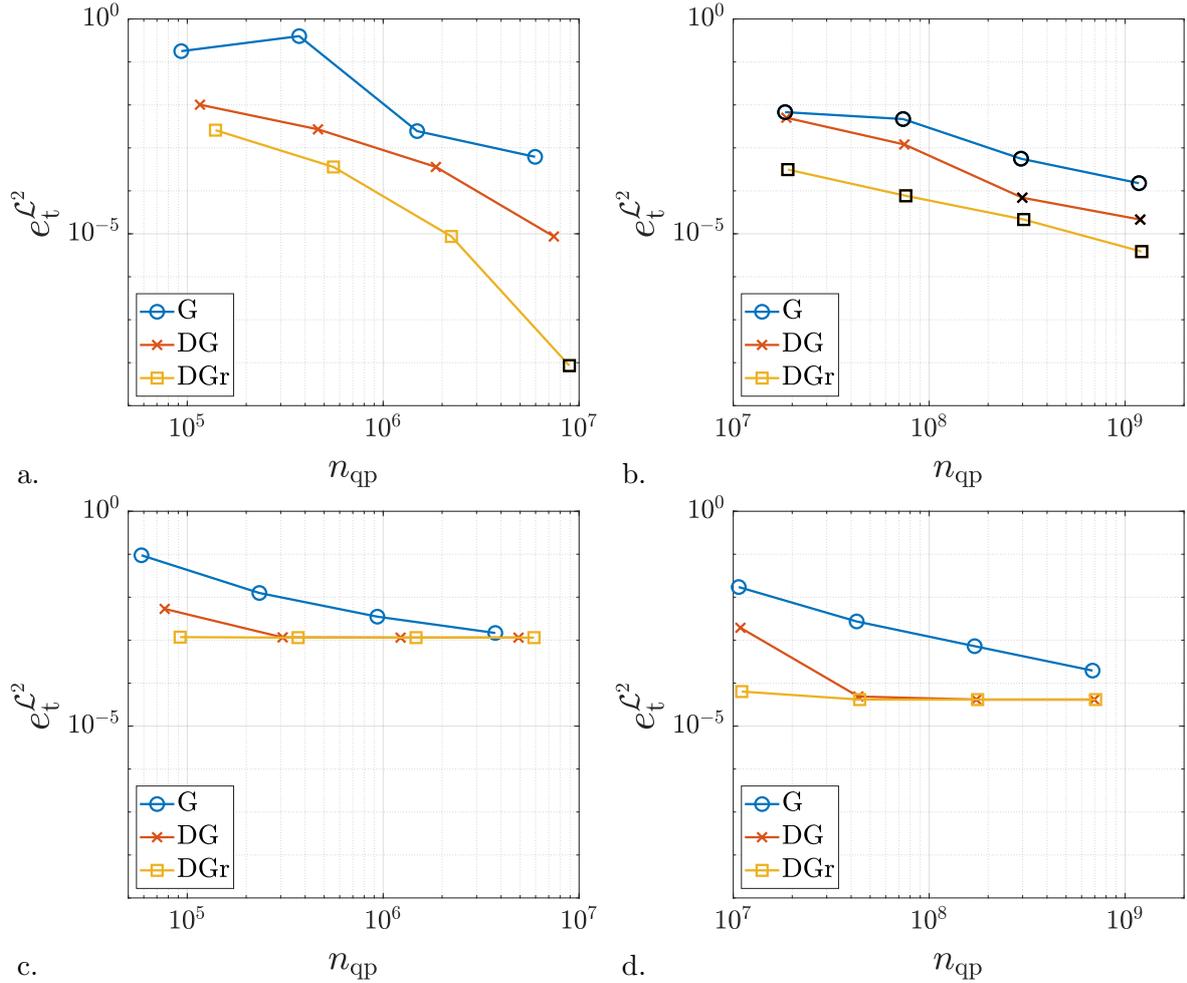
\\\\Fig.~\ref{fig:sphere_rot_quad_res} depicts $e_t^{\mathcal L^2}$ for the first convergence study. On the single-patch sphere (Fig.~\ref{fig:sphere_rot_quad_res}a and b), the traction error decreases for each of the hybrid quadrature schemes, where hybrid Duffy-Gauss quadrature with progressive refinement~(DGr) provides the best result, followed by hybrid Duffy-Gauss quadrature~(DG) and eventually by Gauss-Legendre quadrature~(G). Black markers denote results obtained from an iterative solver.\footnote{The iterative solution method used is the preconditioned conjugate gradients method. For each individual combination of quadrature scheme, mesh refinement and quadrature density, the depicted error is chosen based on
$$e_\mrt^{\mathcal L^2} = 
\begin{cases}
e_\mathrm{t,iter}^{\mathcal L^2}, & e_\mathrm{t,iter}^{\mathcal L^2} < 0.9 \,e_\mathrm{t,direct}^{\mathcal L^2}\\
e_\mathrm{t,direct}^{\mathcal L^2}, & \, \mathrm{otherwise}
\end{cases}~,$$
where $e_\mathrm{t,iter}^{\mathcal L^2}$ and $e_\mathrm{t,direct}^{\mathcal L^2}$ denote the $\mcalL^2$ norm of the traction error using an iterative or a direct solver, respectively. The iterative solver is more robust w.r.t.~ill-conditioned BE matrices and yields thus better results on highly refined single-patch spheres with coinciding control points at the poles.} This convention is used for all examples in Sec.~\ref{sec:examples}.
Similar results are obtained for the six-patch sphere (Fig.~\ref{fig:sphere_rot_quad_res}c and d), but the $\mathcal L^2$ norm of the traction error converges only to $e_t^{\mathcal L^2}\approx 10^{-3}$ for $\ell=2$ and to $e_t^{\mathcal L^2}\approx 5\cdot 10^{-5}$ for $\ell=4$. Better results are prevented by the inaccuracy of the six-patch sphere discretization (see Fig.~\ref{fig:quad_sphere_discretization}e). The six-patch sphere does not require an iterative solver and is thus preferable in terms of robustness and computational effort. The quadrature weights for adjusted weight quadrature are given in Table \ref{tab:Gauss} and \ref{tab:weights} only for $n_\mathrm{qp}=3\times 3$. The results for DGw are thus not investigated for a varying quadrature density.
\\\\Fig.~\ref{fig:sphere_rot_mesh_res2} depicts the $\mcalL^2$ norm of the traction error for the second convergence study. It can be seen that the traction error decreases for each of the hybrid quadrature schemes and on both of the considered discretizations. As expected, scheme DGr and G provide the best and the worst results, respectively. The results of the remaining two schemes are in between, where DGw seems to offer no significant advantage in accuracy over DG here. Applying G or DG to a single-patch sphere with $\ell>1$ requires an iterative solver, as indicated by the black markers in Fig.~\ref{fig:sphere_rot_mesh_res2}a. Schemes DGw and DGr, which also treat the near singular elements properly, allow to use a direct solver for $\ell=2,3$ also, and are therefore preferable in regard to the computational time.
\\Unlike the single-patch sphere, the six-patch sphere (Fig.~\ref{fig:sphere_rot_mesh_res2}b) does not require an iterative solver and is therefore very robust, both with respect to mesh and quadrature refinement. The six-patch sphere further leads to a smaller error, where DGr provides by far the best result with a convergence rate of almost 1.5~($e_t^{\mathcal L^2} \propto 1/n_\mathrm{no}^{1.5}$). The kink in the curves for DG and DGw shows that applying Duffy quadrature to the singular elements without increasing the quadrature density on the near singular elements is insufficient to yield accurate results for $\ell>3$.

\subsection{Flow caused by a rising sphere}\label{sec:ex_trans}

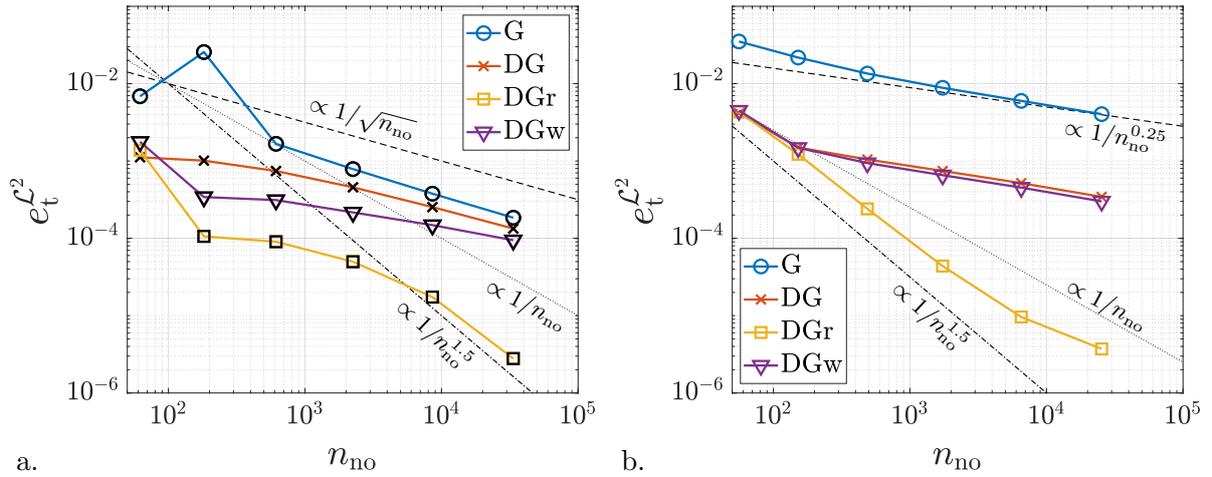
\begin{figure}[h]
\unitlength\linewidth
\begin{picture}(1,0.41)
\put(0,0){\includegraphics[trim = 0 0 30 20, clip, width=.5\linewidth]{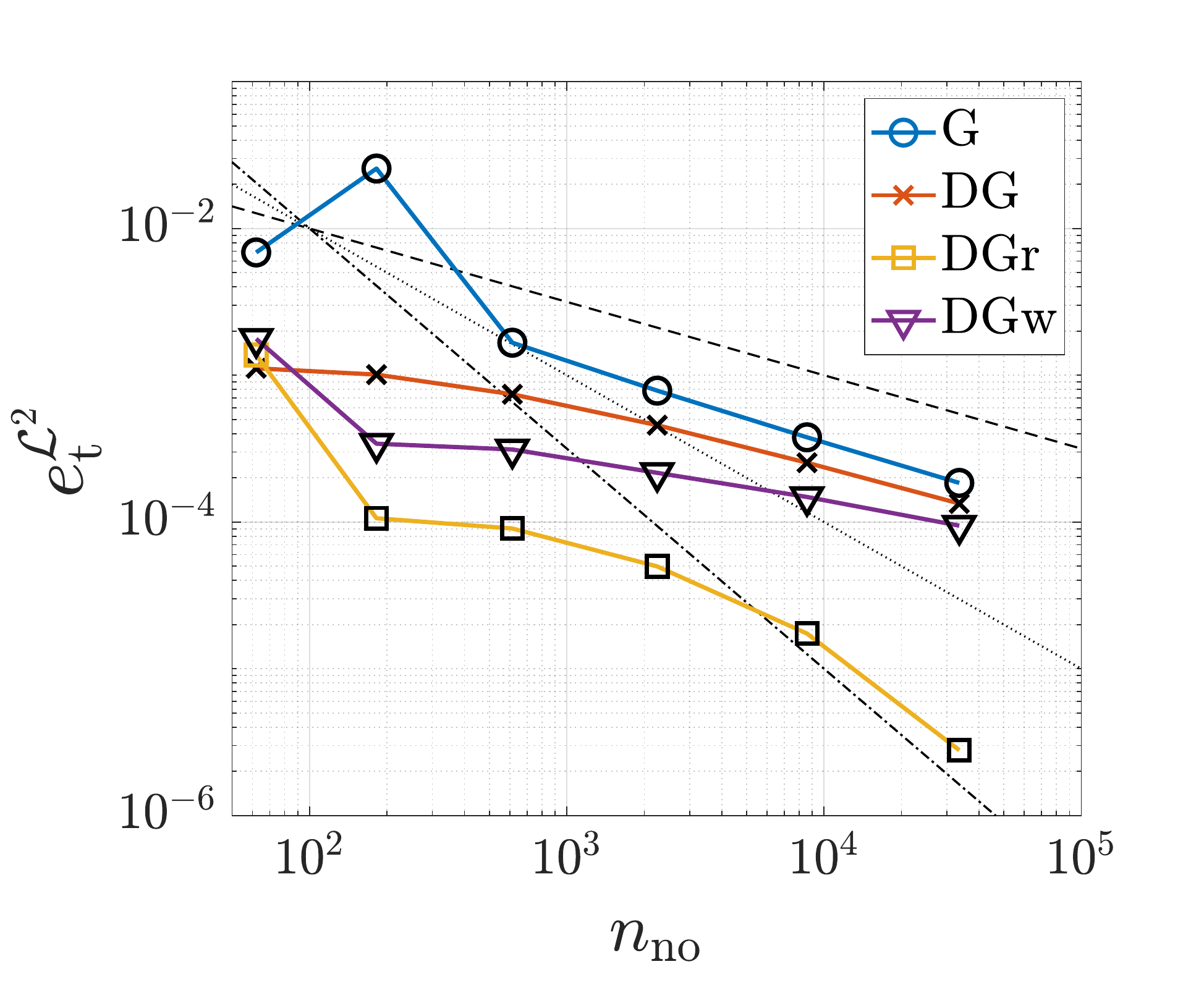}}
\put(0.25,0.31){\rotatebox{-15}{\footnotesize $\propto 1/\sqrt{n_\mathrm{no}}$}}
\put(0.395,0.16){\rotatebox{-28}{\footnotesize $\propto 1/n_\mathrm{no}$}}
\put(0.32,0.145){\rotatebox{-43}{\footnotesize $\propto 1/n_\mathrm{no}^{1.5}$}}
\put(0.5,0){\includegraphics[trim = 0 0 30 20, clip, width=.5\linewidth]{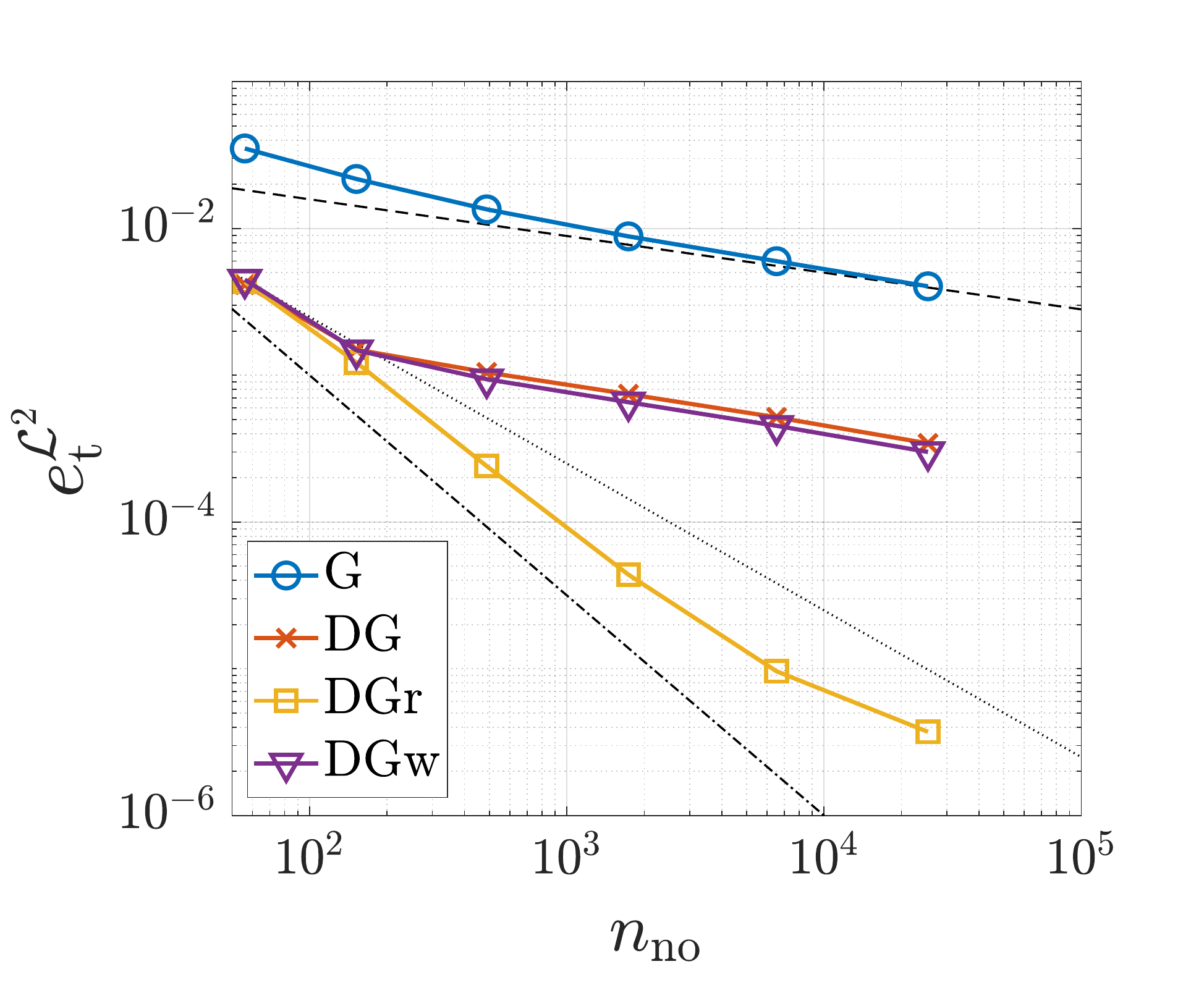}}
\put(0.875,0.29){\rotatebox{-10}{\footnotesize $\propto 1/n_\mathrm{no}^{0.25}$}}
\put(0.875,0.16){\rotatebox{-28}{\footnotesize $\propto 1/n_\mathrm{no}$}}
\put(0.73,0.16){\rotatebox{-43}{\footnotesize $\propto 1/n_\mathrm{no}^{1.5}$}} 
\put(0.01,0.01){a.} \put(0.51,0.01){b.}
\end{picture}
\caption{\textit{Rising sphere:} Traction error $e_t^{\mathcal L^2}$ for quadrature density $n_0=3$ on the single-patch NURBS sphere~(a) and on the six-patch NURBS sphere~(b), both of varying refinement level ${\ell=1,\ldots,6}$. The black markers depict results from an iterative solver.}\label{fig:sphere_trans}
\end{figure}

The suitability of the introduced quadrature schemes for a pure rotational problem has been shown in Sec.~\ref{sec:ex_rot}. However, this problem yields only tangential tractions, while the normal tractions are zero. The suitability of the quadrature schemes for problems with non-zero velocity and traction components in both normal and tangential direction is investigated in the second example: A sphere with radius $R$ is transtaled with constant velocity $\bar\bv = \bar v\, \be_3$\footnote{The prescribed velocity can be chosen arbitrarily, i.e.~$\bar \bv= \bar v_i \,\be_i$. Here, a pure vertical velocity is chosen for the sake of simplicity.} through the fluid. The Dirichlet boundary condition on the surface $\mcalS$ is therefore given by
\begin{equation} \label{eq:ex_trans_BC}
\bv(\bx)= \bar\bv, \hh \forall\; \bx \in \mcalS~,
\end{equation}
while the surface traction is unknown. 
\\\\Fig.~\ref{fig:sphere_trans} shows the $\mcalL^2$ norm of the relative BE traction error
\begin{equation}\label{eq:ex_trans_err}
	e_\mrt(\bx)=\frac{\| \bt^h(\bx) - \bt(\bx)\|}{\|\bt(\bx)\|}
\end{equation}
for translating spheres of refinement level ${\ell=1,\ldots,6}$. The results are generally similar to those from Sec.~\ref{sec:ex_rot}. On the single-patch sphere, the translating sphere problem can be solved with higher accuracy than the rotating sphere problem. On the six-patch sphere, it is the other way around. It is noteworthy that the convergence rate for DGr and $\ell>4$ decreases slightly. Treating additional rings with refined quadrature would allow to maintain a convergence rate of 1.5. However, a single refinement ring already gives an excellent gain in accuracy.

\subsection{Flow caused by a rising ellipsoid}\label{sec:ex_bubble}
\begin{figure}[h]
\unitlength\linewidth
\begin{picture}(1,0.36)
\put(0.02,0){\includegraphics[trim = 100 133 13 0, clip, width=0.3\linewidth]{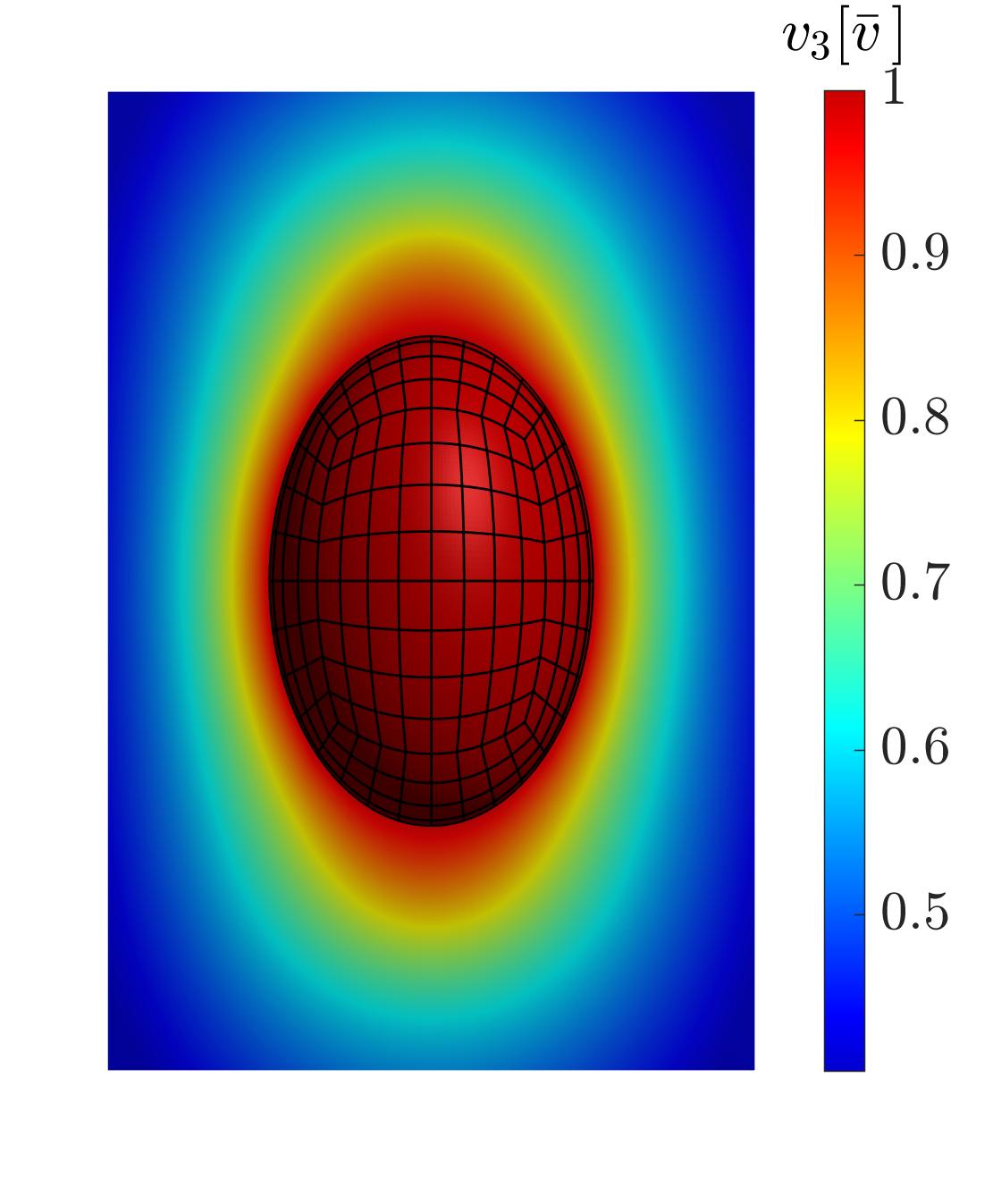}}
\wbox{0.25}{0.3495}{0.035}{0.008} \put(0.241,.347){\footnotesize $v_3 \,[\bar v]$}
\put(0.36,0){\includegraphics[trim = 100 133 0 0, clip, width=0.304\linewidth]{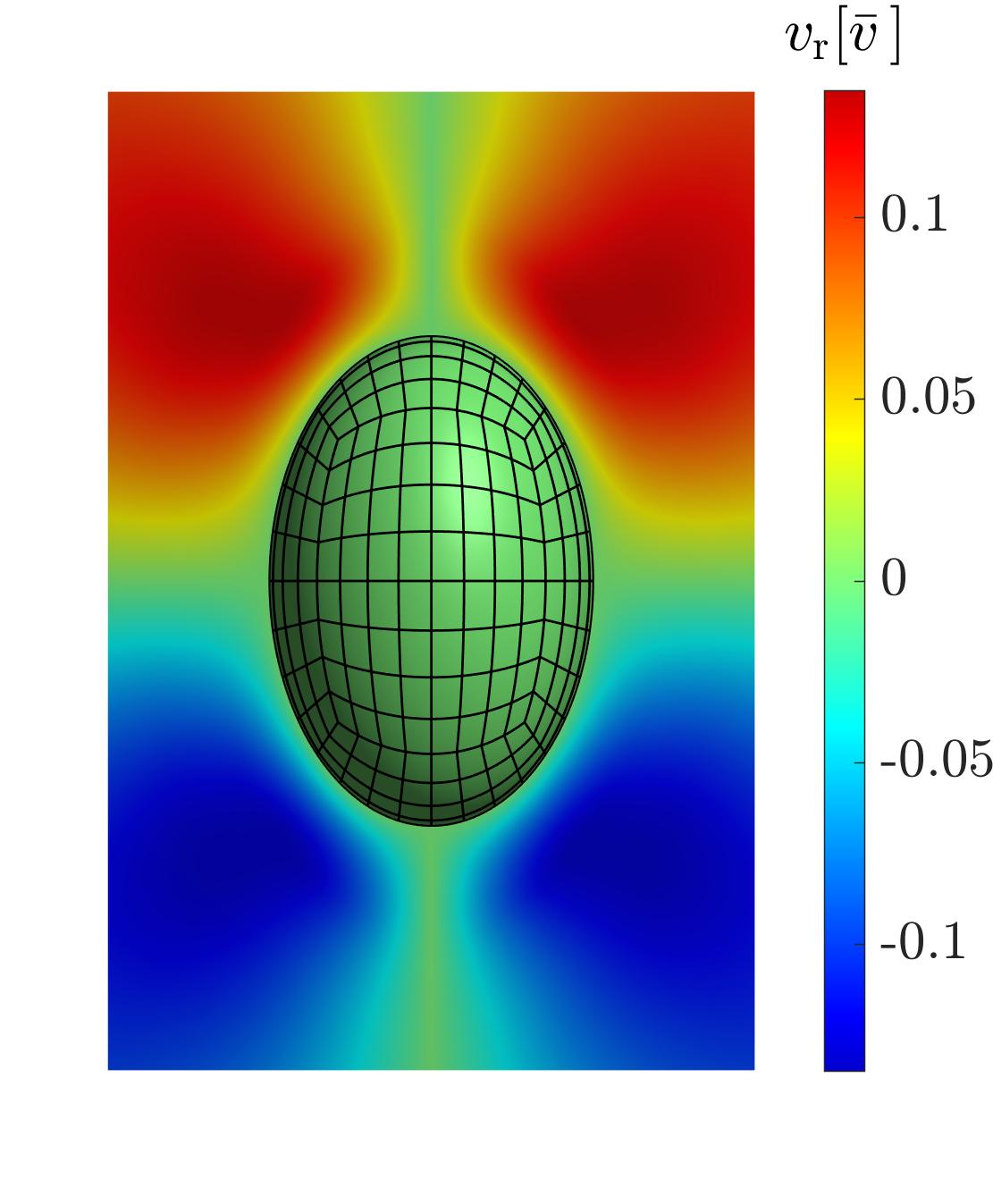}}
\wput{0.58}{0.347}{\footnotesize $v_\mrr \,[\bar v]$} 
\put(0.70,0){\includegraphics[trim =  100 133 13 0, clip, width=0.3\linewidth]{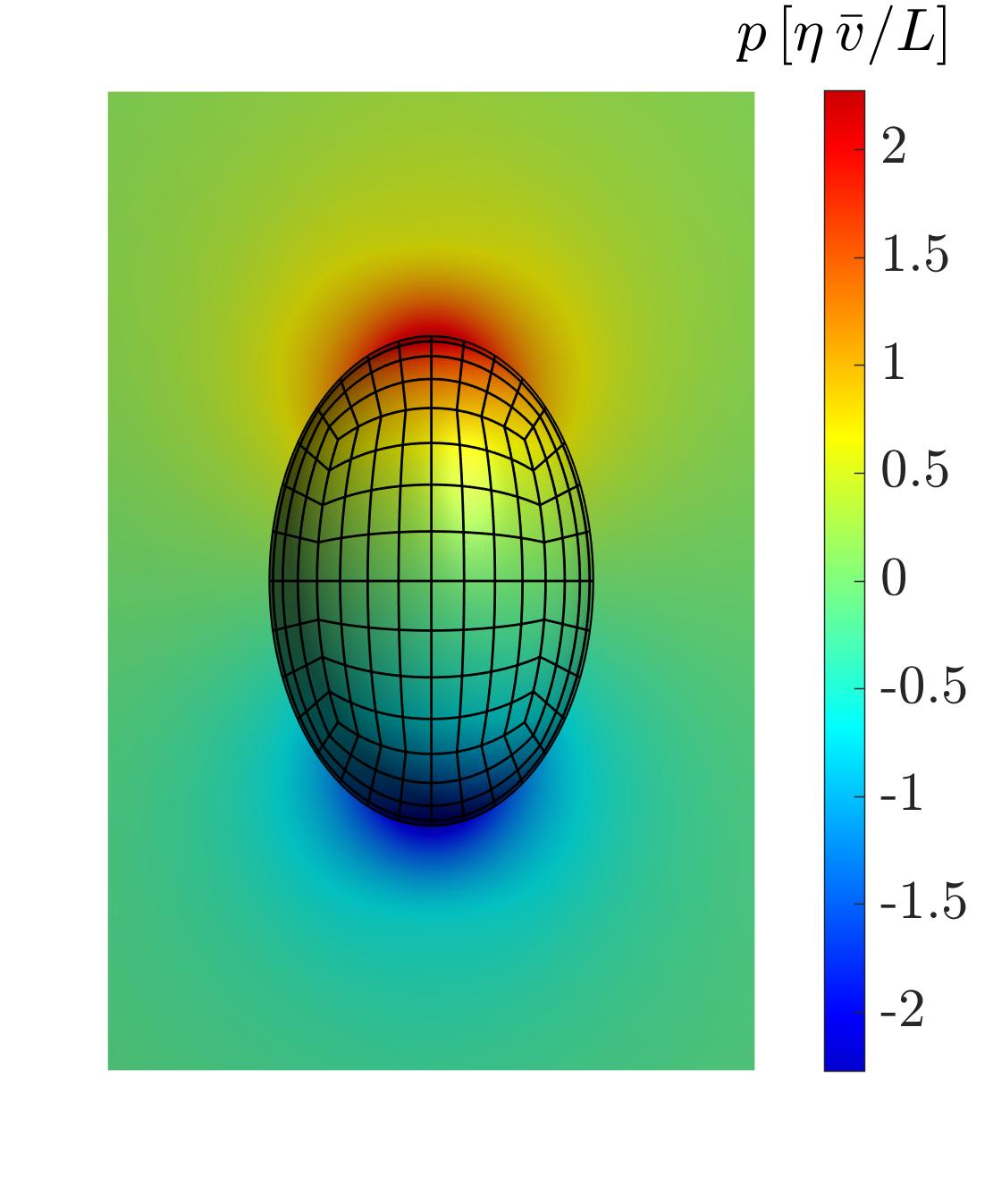}}
\wput{0.91}{0.347}{\footnotesize $p \,[\eta \bar v/L]$} 
\put(0,0){a.}\put(0.34,0){b.}\put(0.68,0){c.}
\end{picture}
\caption{\textit{Rising ellipsoid}: Velocity and pressure field on the ellipsoid surface and in the surrounding fluid: a.~vertical velocity $v_3=\bv \cdot \be_3$; b.~radial velocity $v_\mrr = \bv \cdot \be_\mrr$, where $\be_\mrr := \bx_\mrr /\| \bx_\mrr \|$ with $\bx_\mrr= x_1 \be_1 +x_2 \be_2$ ; c.~pressure $p=-\mathrm{tr}\, \bsig/3$.
}\label{fig:ex_bubble_post}
\end{figure}

\begin{figure}[h]
\unitlength\linewidth
\begin{picture}(1,0.41)
\put(0,0){\includegraphics[trim = 0 0 30 20, clip, width=.5\linewidth]{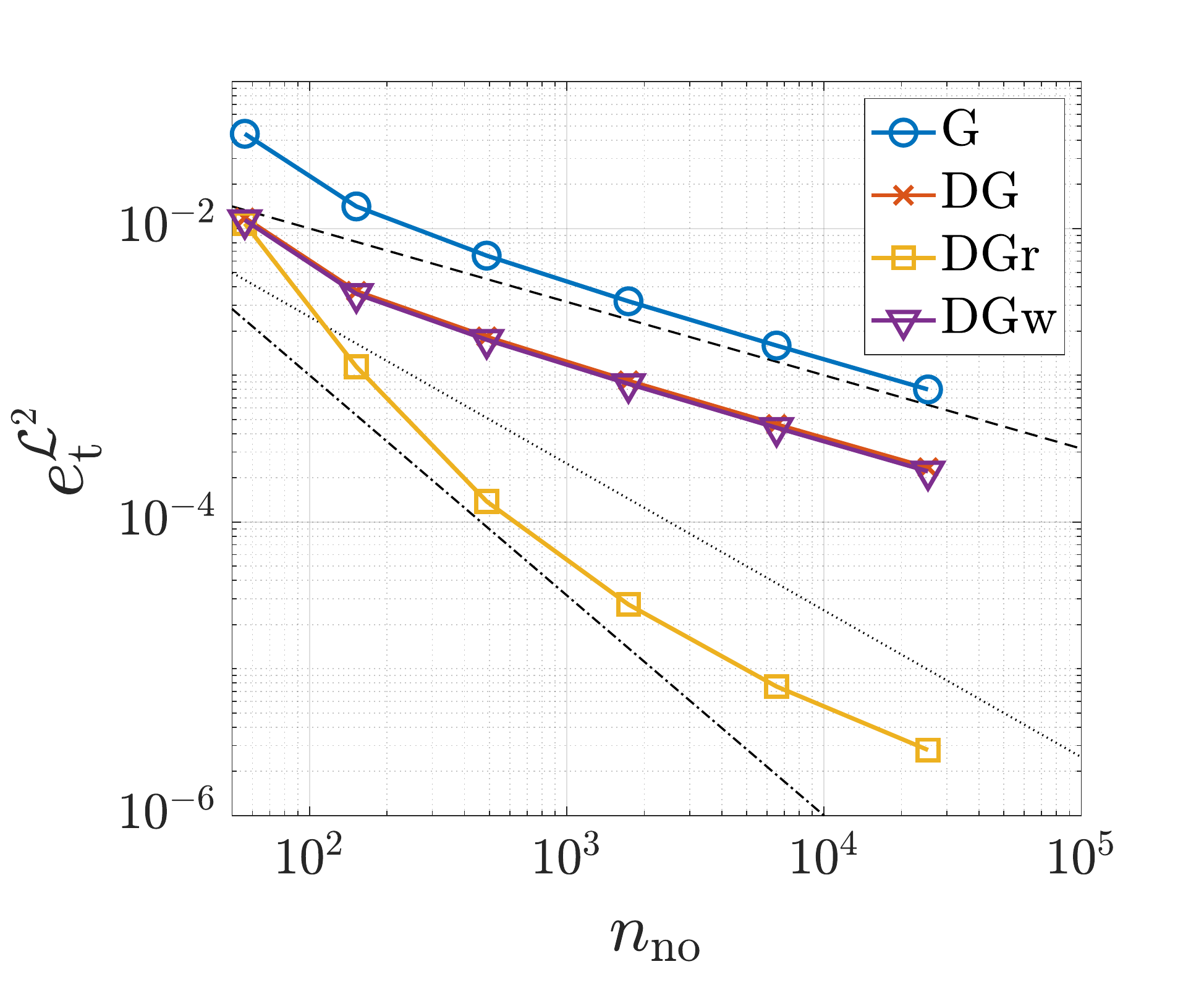}}
\put(0.5,0){\includegraphics[trim = 0 0 30 20, clip, width=.5\linewidth]{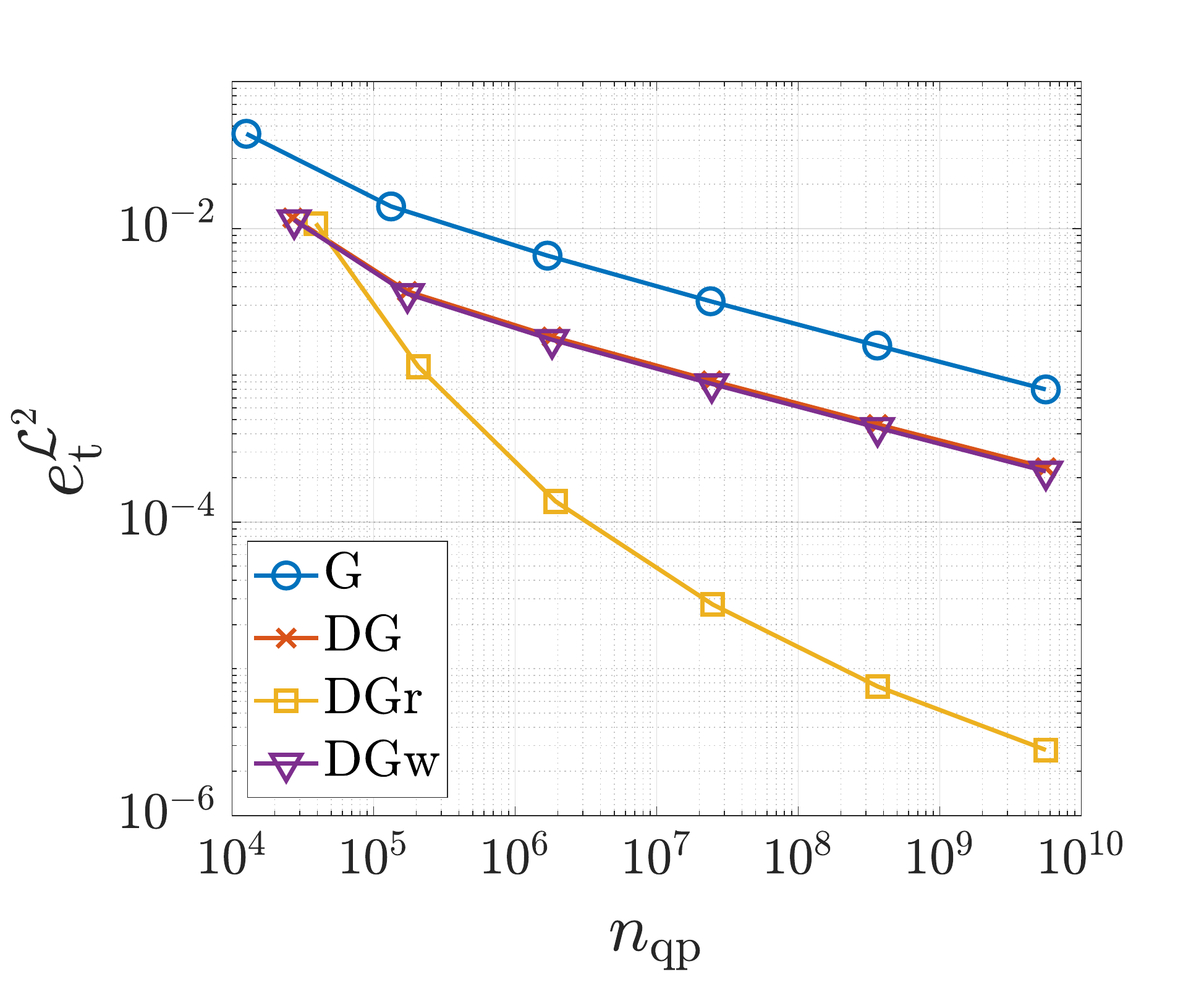}}
\put(0.34,0.253){\rotatebox{-15}{\footnotesize $\propto 1/n_\mathrm{no}^{0.25}$}}
\put(0.375,0.16){\rotatebox{-28}{\footnotesize $\propto 1/n_\mathrm{no}$}}
\put(0.23,0.16){\rotatebox{-43}{\footnotesize $\propto 1/n_\mathrm{no}^{1.5}$}} 
\put(0.01,0.01){a.} \put(0.51,0.01){b.}
\end{picture}
\caption{\textit{Rising ellipsoid:} The $\mcalL^2$ norm of the relative pressure error $e_\mrp$ for quadrature density $n_0=3$ and six-patch NURBS ellipsoids of varying refinement level ${\ell=1,\ldots,6}$.
}\label{fig:ex_bubble_study}
\end{figure}
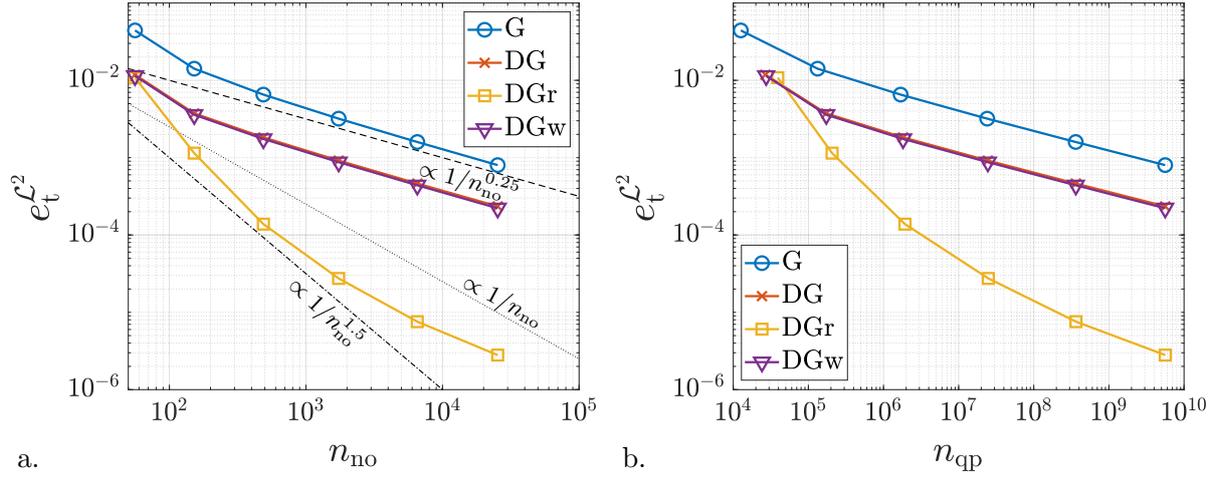

The suitability of the quadrature schemes for problems with non-spherical surface geometries is investigated in the third example:
%
A ellipsoid with surface $\mcalS$ rises with constant velocity $\bar v$ through a fluid of dynamic viscosity $\eta$. The surface of the ellipsoid is described by 
\begin{equation}\label{eq:ex_ellipsoid_surf}
\frac{x_1^2 +x_2^2}{a^2\,(1-e^2)} + \frac{x_3^2}{a^2} =1~, \hh \mathrm{for}\; \bx  \in \mcalS
\end{equation}
where $e$ denotes the eccentricity of the ellipsoid ($0\leq e<1$) and $a$ is the length of its semi-major axis.\footnote{The eccentricity of a ellipsoid is defined as $e:=\sqrt{1-b^2/a^2}$, where $b$ denotes the length of the semi-minor axis. Two parameters out of the triplet ($a$,$b$,$e$) needs to be known to define the ellipsoidal surface. The example from Sec.~\ref{sec:ex_bubble} considers a ellipsoid with $e=0.75$ and $a=L/e^{1/3}$.} The velocity on $\mcalS$ is given by the Dirichlet boundary condition \eqref{eq:ex_trans_BC}.
\begin{figure}[h]
\unitlength\linewidth
\begin{picture}(1,0.255)
	\put(0,0){\includegraphics[trim = 250 180 390 120, clip, width=.21\linewidth]{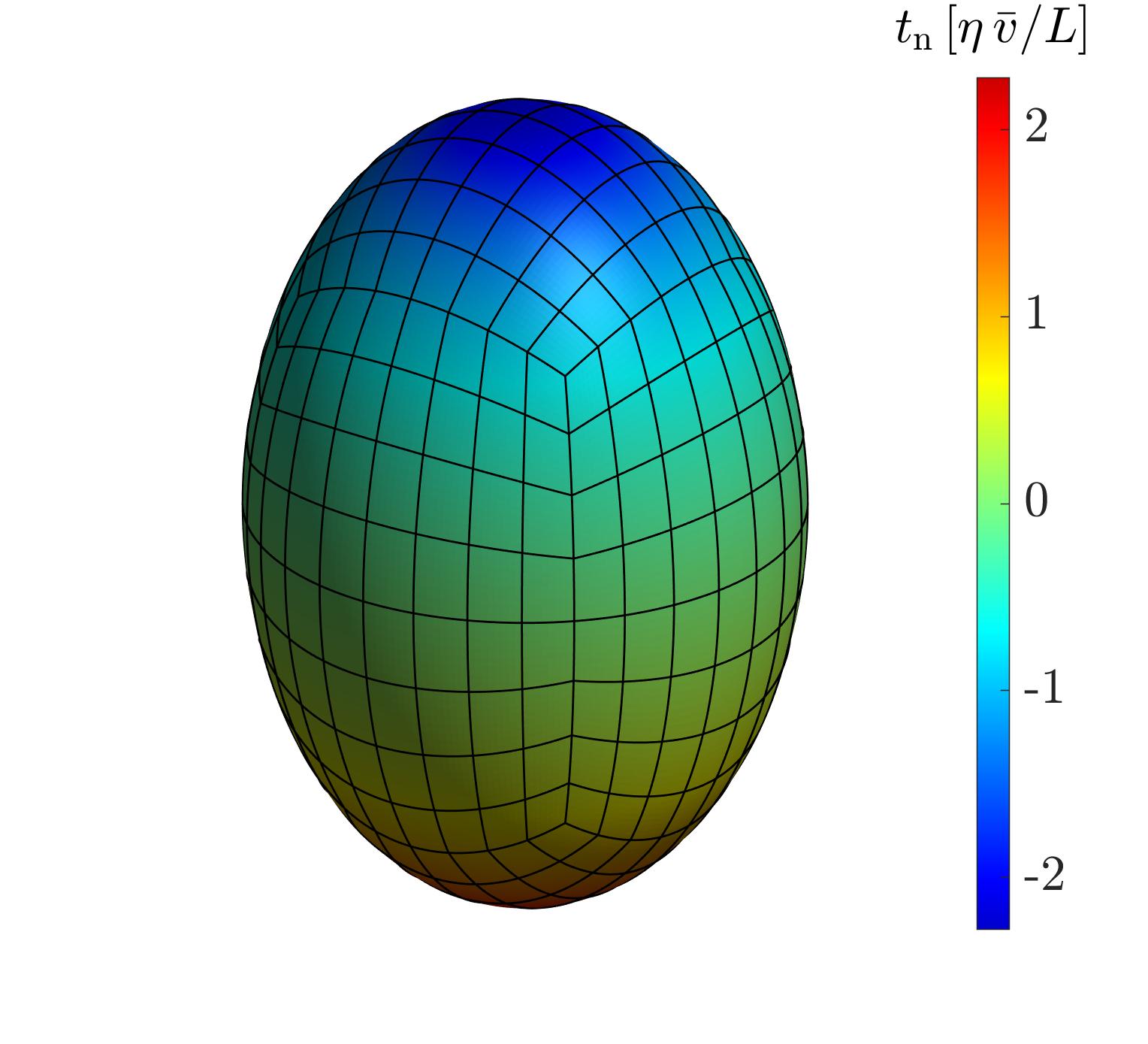}}
	\put(.205,0){\includegraphics[trim = 1250 130 40 76, clip, height=.24\linewidth]{figures/examples/bubble/small/t_n.jpg}}
	\put(.25,0){\includegraphics[trim = 250 180 390 120, clip, width=.21\linewidth]{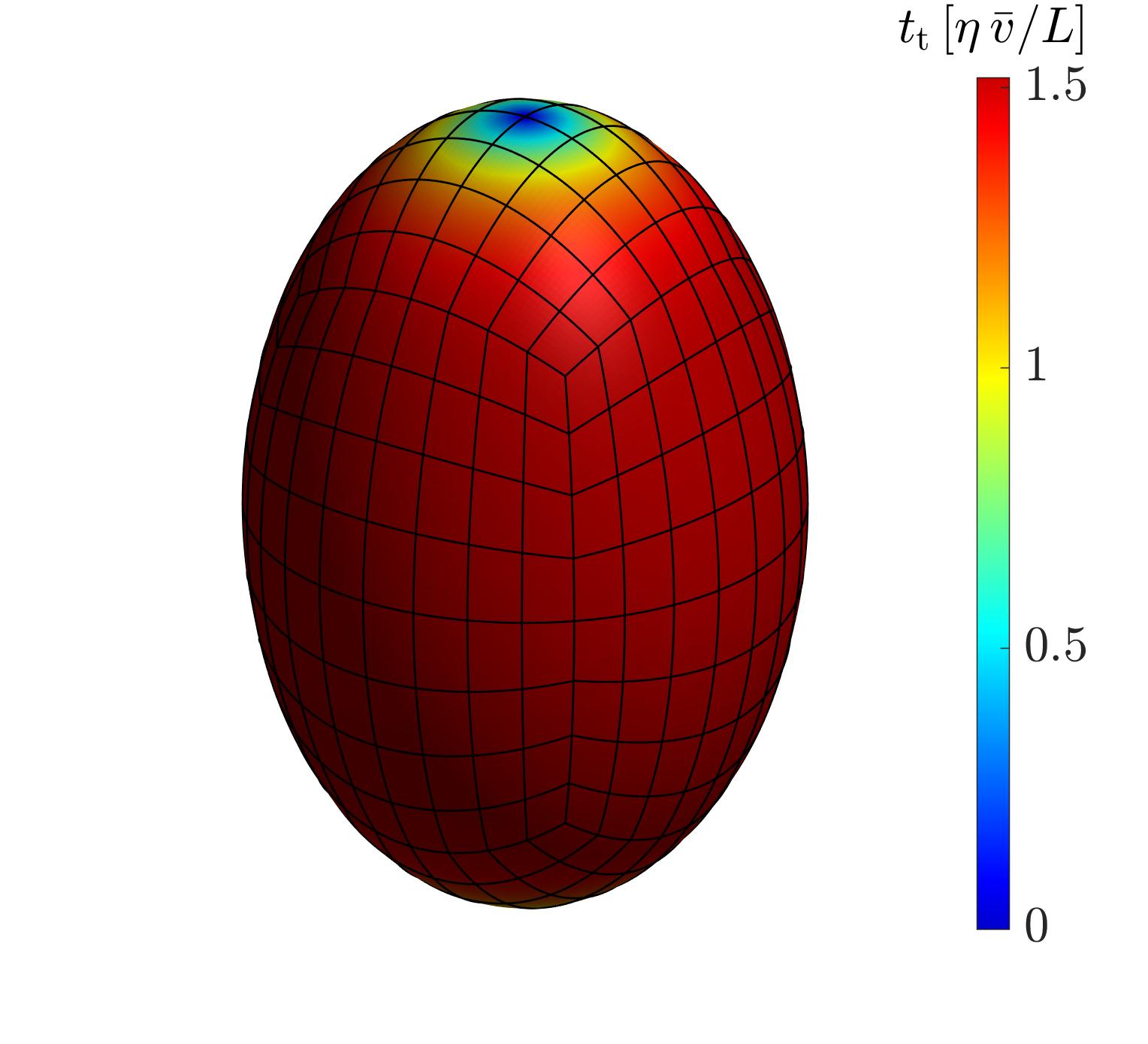}}
	\put(.455,0){\includegraphics[trim = 1250 130 40 76, clip, height=.24\linewidth]{figures/examples/bubble/small/t_t.jpg}}
	\put(0.5,0){\includegraphics[trim = 250 180 390 120, clip, width=.21\linewidth]{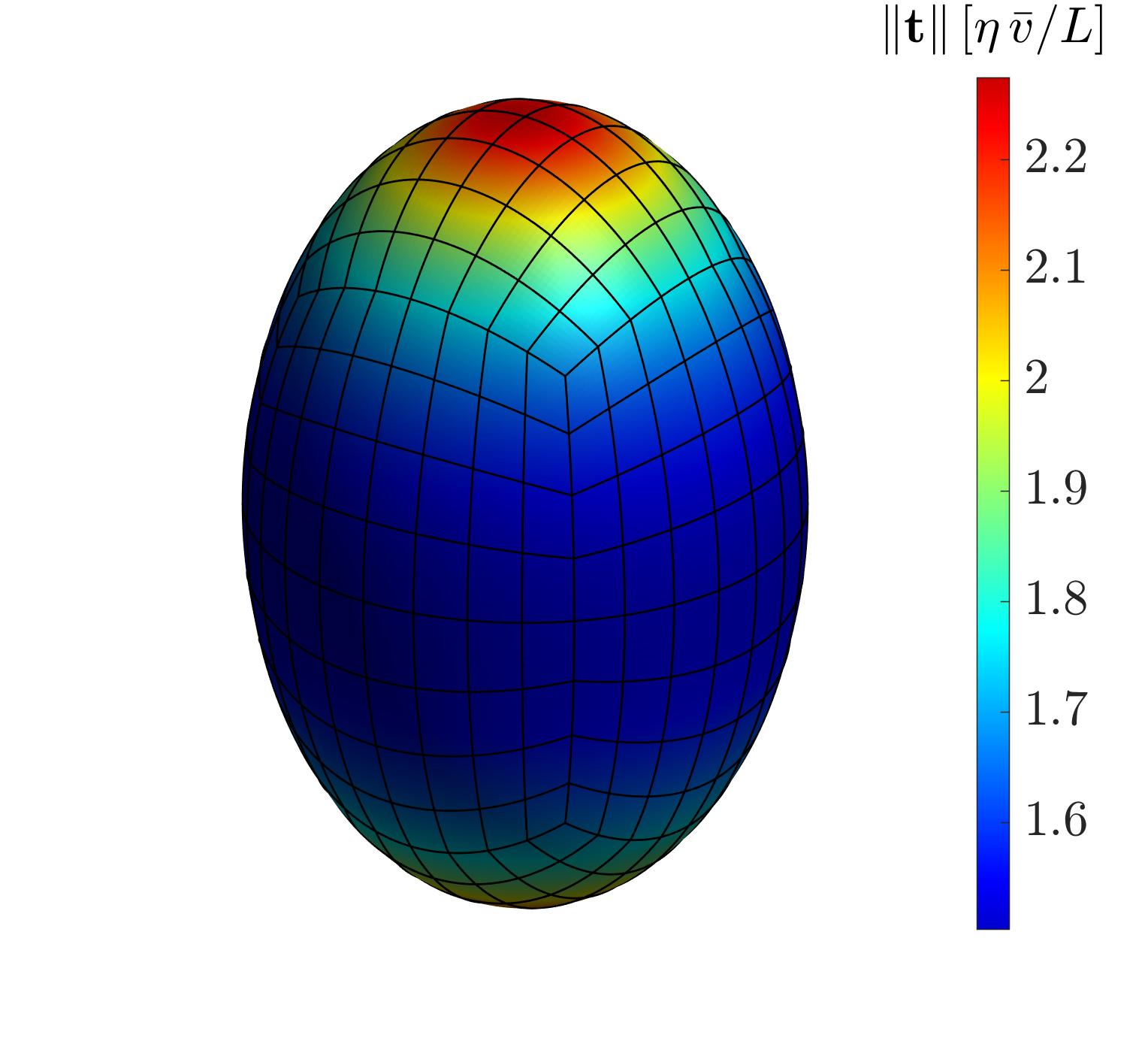}}
	\put(0.705,0){\includegraphics[trim =  1250 130 40 76, clip, height=.24\linewidth]{figures/examples/bubble/small/t.jpg}}
	\put(0.75,0){\includegraphics[trim = 250 180 390 120, clip, width=.21\linewidth]{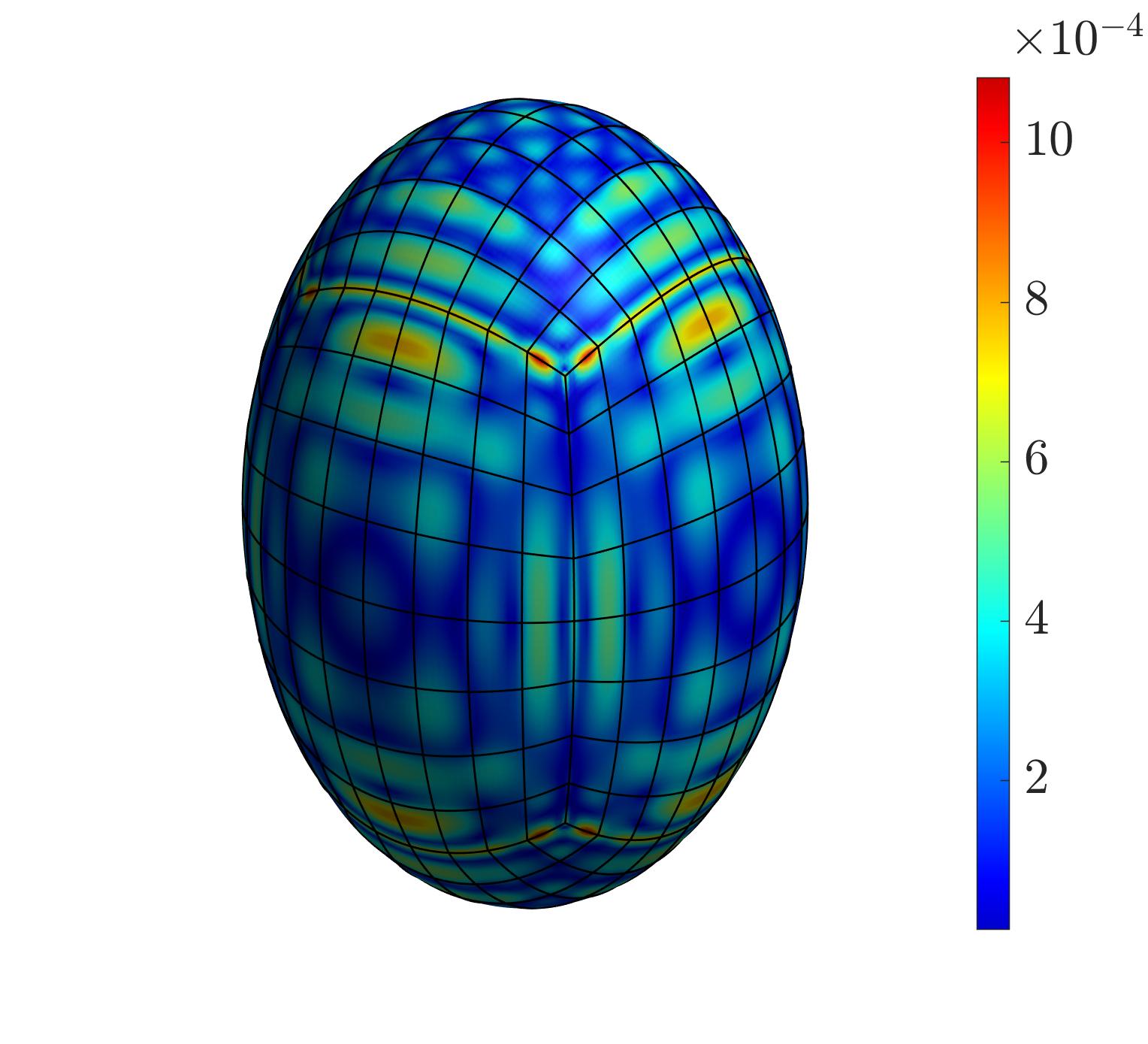}}
	\put(0.955,0){\includegraphics[trim =  1250 130 50 76, clip, height=.24\linewidth]{figures/examples/bubble/small/err_t.jpg}}
	\put(.18,0.245){\scriptsize $t_\mrn[\eta \,\bar v /L ]$} \put(.43,0.245){\scriptsize $t_\mrt[\eta \,\bar v /L ]$}
	\put(.68,0.245){\scriptsize $\|\bt\|[\eta \,\bar v /L ]$} \put(.93,0.245){\scriptsize $e_\mrt\,[10^{-4}]$} 
	\put(0.01,0){\small a.} \put(0.26,0){\small b.} \put(0.51,0){\small c.} \put(0.76,0){\small d.}
\end{picture}	
\caption{\textit{Rising ellipsoid:} BE traction for a ellipsoid of refinement level $\ell=3$ and quadrature scheme DGr with $n_0=3$. a.~normal traction $t_\mrn = -p$; b.~tangential traction $t_\mrt$; c.~norm of the traction $\|\bt\|=t_3$ here; d.~relative traction error $e_\mrt$.}\label{fig:ex_bubble_traction}
\end{figure}
\\\\Fig.~\ref{fig:ex_bubble_traction} shows the BE traction on a six-patch ellipsoid of refinement level $\ell=3$ using quadrature scheme DGr with $n_0=3$. The normal and tangential components of the traction vector and its magnitude are shown in Fig.~\ref{fig:ex_bubble_traction}a, b and c, respectively. Fig.~\ref{fig:ex_bubble_traction}d shows that the relative BE traction error \eqref{eq:ex_trans_err} yields $e_\mrt \approx 10^{-3}$ along patch boundaries and an even much smaller error away from patch boundaries. The velocity and pressure field on the ellipsoid and in the surrounding fluid are shown in Fig.~\ref{fig:ex_bubble_post}. The fluid velocity is determined by the BIE for points within the domain \eqref{eq:flow_BIE_domain}, while the surface velocity on $\mcalS$ is given by BC \eqref{eq:ex_trans_BC}. The pressure is determined on the ellipsoid surface by $p = - \bt \cdot \bn$ and within the fluid domain by the pressure BIE \eqref{eq:flow_BIE_pressure}.
%
%
\\The $\mcalL^2$ norm \eqref{eq:ex_rot_L2_dof} of the relative BE traction error \eqref{eq:ex_trans_err} is shown in Fig.~\ref{fig:ex_bubble_study}a vs.~the number of nodes, and in Fig.~\ref{fig:ex_bubble_study}b vs.~the total number of quadrature points \eqref{eq:nqp}. The results are very similar to those of the rising sphere in Sec.~\ref{sec:ex_trans}, since the sphere is a special case of the ellipsoid (for $e = 0$): Quadrature schemes DGr and G provide the best and the worst of the BE results, respectively. Considering DGr with additional rings of refined quadrature would allow to maintain a convergence rate of 1.5. Quadrature schemes DG and DGw provide almost the same results, which are in between the results of the former. Thus, all hybrid quadrature schemes, in particular DGr, prove to be robust to changes in surface shape and are therefore suitable for various applications, including coupled FE-BE simulations.
\section{Conclusion and outlook}\label{sec:conclusion}
This work presents new quadrature schemes for the efficient approximation of weakly singular integrals ($1/r$ kernel) and provides important findings on singular quadrature in BE~(Boundary element) analysis. A new quadrature rule that approximates weakly singular integrals on plane and regular near singular elements exactly using adjusted weights is presented in Sec.~\ref{sec:quad_nearly_adjusted}. Numerical examples have shown that this quadrature rule can also be advantageous on curved surfaces, as it improves the accuracy compared to classical Gaussian quadrature. The presented quadrature rule can be easily extended to determine singular integrals exactly on other surfaces like spheres and cylinders.
\\\\Sec.~\ref{sec:hybrid} presents four new hybrid quadrature schemes for BE analysis that combine different quadrature rules for singular, near singular, and regular elements. Numerical investigations on isogeometric surfaces show that all presented schemes converge robustly for both quadrature and mesh refinement.
Hybrid Duffy-Gauss quadrature with adjusted weights (DGw) is by far the most efficient scheme on flat surfaces, whereas hybrid Duffy-Gauss quadrature with progressive refinement~(DGr) is the most efficient scheme on curved surfaces. The accuracy of DGw on curved surfaces could be further increased by adapting the quadrature rule with adjusted weights as it is briefly described in the last paragraph of Sec.~\ref{sec:quad_nearly_adjusted}. This adaption is not considered here and thus left for future work.
The present paper considers a single ring of elements with refined quadrature, which gives already an excellent gain in accuracy. However, considering additional rings would yield even better approximations and is thus worth to investigate in future works. The quadrature resolution for each ring could then be determined with respect to numerical criteria such as those of \cite{Bu95}.
\\\\The application of the presented hybrid quadrature schemes to coupled FE-BE analysis of fluid structure interaction problems is planned in a forthcoming paper. The presented schemes also offer a wide range of other applications, since they are suitable for arbitrary weakly singular integrals and are not limited to Stokes flow or to BE analysis at all.
\section*{Acknowledgements}
The Deutsche Forschungsgemeinschaft (DFG) is gratefully acknowledged for financial support through SFB 1120. The authors also wish to thank Micha{\l} P.~Rajski for his valuable comments and help.


\appendix
\numberwithin{equation}{section}
\numberwithin{table}{section}
\numberwithin{figure}{section}
\section{Boundary element analysis}\label{sec:BE}
Solving the boundary integral equation~(BIE) for Dirichlet, Neumann or mixed problems by use of boundary elements~(BE) is discussed in this section. The spatial discretization of the surface geometry and the BIE~\eqref{eq:flow_BIE_surf} is introduced in Sec.~\ref{sec:BE_discr}, while the subsequent collocation into a solvable system is discussed in Sec.~\ref{sec:BE_collo}. The elemental mapping from the physical domain to a parameter domain is described in Sec.~\ref{sec:BE_mapping} and the conceptual approximation of surface integrals by numerical quadrature is presented in Sec.~\ref{sec:BE_quad}.

\subsection{Boundary element discretization}\label{sec:BE_discr}
The surface geometry $\mcalS$ and the BIE~\eqref{eq:flow_BIE_surf} are discretized into $\nel$ finite boundary elements, numbered ${e = 1,\ldots,n_\mathrm{el}}$, and $\nno$ nodal points. Element $e$ occupies the surface domain $\Omega^e\subset \mcalS^h$ such that the surface geometry is approximated by
\begin{equation}
 \mcalS  \approx \mcalS^h = \bigcup_{e=1}^{n_\mathrm{el}} \Omega^e~,
\end{equation}
where superscript $h$ denotes approximated quantities. Point $\bx\in\Omega^e$ is approximated by the nodal interpolation 
\begin{equation}\label{eq:BE_discr_x}
\bx \approx \bx^h= \sum_{I=1}^{n_e} N_I \, \mx_I~, \hh \bx^h \in \Omega^e \subset \mcalS^h,
\end{equation}
where $\mx_I$ denote one of the $n_e$ nodal points that define element $e$ and $N_I$ denotes the corresponding shape functions \footnote{Non-italic discrete arrays $\mx_I$ and $\mx^e$ should not be confused with the italic field variable $\bx$.} (e.g.~Lagrange functions, B-splines or NURBS).
Discrete elemental arrays are defined to simplify the notation: The nodal positions for element $e$ are assembled in the discrete column vector of length $\dof^e:=3\,n_e$
\begin{equation}\label{eq:BE_discr_XI}
\mx^e =
\begin{bmatrix}
\mx_1\\
\mx_2\\
\vdots\\
\mx_{n_e}
\end{bmatrix}~,
\end{equation}
whereas the nodal shape functions are assembled in the ($3\times n_\mathrm{dof}^e$) array 
\begin{equation}\label{eq:BE_discr_elem}
\mN^e =
\begin{bmatrix}
N_1\, \bone,\, N_2 \,\bone,\,\cdots,\, N_{n_e}\, \bone~,
\end{bmatrix}
\end{equation}
where $\bone$ is the usual identity tensor in $\mathbb{R}^3$. Approximation \eqref{eq:BE_discr_x} can be re-written in matrix form by use of the elemental position vectors as
\begin{equation}\label{eq:BE_discr_x2}
	\bx \approx \mN^e \, \mx^e~.
\end{equation}
To discretize the BIE, velocity and traction on the surface are approximated in the same fashion. Velocity and traction are approximated at field point $\bx\in\Omega^e$ by
\begin{equation}\label{eq:BE_discr_v}
\bv(\bx)\approx \mN^e \, \mv^e~,\hh \mathrm{and}\hh  \bt(\bx)\approx \mN^e \, \mt^e~,
\end{equation}
where $\mv^e$ and $\mt^e$ denote the nodal velocity and traction vector of length $n_\mathrm{dof}^e$ for element $e$. At source point $\by\in\Omega^{\bar e}$, the velocity is analogously approximated by
\begin{equation}\label{eq:BE_discr_v2}
	\bv(\by)\approx \mN^{\bar e} \, \mv^{\bar e}~,
\end{equation}
where $\mN^{\bar e}$ and $\mv^{\bar e}$ are the shape function array and the nodal velocity vector for element $\bar e$. The components of the discrete vectors $\mv^e$, $\bar \mv^e$ and $\mt^e$ are denoted by $\mrv_A^e$,  $\mrv_A^{\bar e}$ and $\mrt_A^e$, respectively, while the components of the shape function arrays $\mN^e$ and $\mN^{\bar e}$ are denoted by $\mrN_{iA}^e$ and $\mrN_{iA}^{\bar e}$. The lowercase indices $i$,~$j$ and $k$ run from $1$ to $3$, while the uppercase index $A$ runs from $1$ to $n_\mathrm{dof}^e$.
\\\\The BIE \eqref{eq:flow_BIE_surf} is then approximated by use of \eqref{eq:BE_discr_v} and \eqref{eq:BE_discr_v2} as
\begin{equation}\label{eq:BE_BIE}
	\mrN_{iA}^{\bar e}\, \mrv_A^{\bar e} = \sum_{e=1}^{n_e} \left[ \mcalG_{iA}^e(\by)\,\mrt^e_A +  \mcalT_{iA}^e(\by)\, \mrv^e_A \right]~,
\end{equation}
where summation over $A$ is implied. The components of the discrete BE matrices $\boldsymbol{\mcalG}^e$ and $\boldsymbol{\mcalT}^e$, both of size ($3\times n_\mathrm{dof}^e$), are given by
 \begin{equation}
 \begin{aligned}\label{eq:BE_GreensInt}
	 \mcalG_{iA}^e(\by) &= 	-\frac{1}{4\pi\,\eta}\, \int_{\Omega^e} G_{ij}(\bx-\by) \,\mrN^e_{jA} \,\mathrm da_x~,\\
	\mcalT_{iA}^e(\by) &= \frac{1}{4\pi}\ \int_{\Omega^e} T_{ijk}(\bx-\by)\, n_k(\bx) \, \mrN^e_{jA}\,\mathrm da_x ~.
 \end{aligned}
\end{equation}
The elemental vectors $\mv^e$ and $\mt^e$ are assembled into global knot vectors $\mv$ and $\mt$ of length ${n_\mathrm{dof}:=3\,n_n}$, where $n_n$ denotes the total number of nodal discretization points, i.e.~nodes for Lagrange and control points for isogeometric discretizations. Similarly, the elemental arrays $\boldsymbol{\mcalG}^e$, $\boldsymbol{\mcalT}^e$ and $\mN^{\bar e}$ are assembled into global global matrices $\boldsymbol{\mcalG}$, $\boldsymbol{\mcalT}$ and $\boldsymbol{\mathcal N}$ of size ($3\times n_\mathrm{dof}$).
The discretized boundary integral equation for a source point $\by$ on the interface $\mcalS$ can thus be given in matrix-vector form by
\begin{equation}\label{eq:BE_BIE2}
	\bcalN(\by)\, \mv = \boldsymbol\mcalG(\by)\,\mt +  \boldsymbol\mcalT(\by)\,\mv~.
\end{equation}
A continuous boundary integral representation for the pressure is given in \cite{Pozrikidis92_book} by
\begin{equation}\label{eq:flow_BIE_pressure}
	p(\by)= - \frac{1}{8\pi} \int_{\mcalS} P_i(\bx-\by)\, {t_i}(\bx) \,\mrd a_x  + \frac{\eta}{8\pi}\int_{\mcalS} v_i(\bx) \, \Pi_{ij}(\bx-\by)\, n_j(\bx)  \, \mrd a_x
\end{equation}
with pressure Green's function
\begin{equation}\label{eq:flow_green_p}
P_i (\br) =  2 \,\frac{\bar r_i}{r^2}~,
\end{equation}
and the corresponding tensor for pressure field associated with the stresslet
\begin{equation}\label{eq:flow_green_pi}
\Pi_{ij}(\br) =  \frac{4}{r^3}(3\,\bar r_i \bar r_j - \delta_{ij})  ~.
\end{equation}
More details on theory, discretization, assembly and implementation of BIEs can be found in \cite{Harmel21_thesis}.

\subsection{Boundary element collocation}\label{sec:BE_collo}

The discretized BIE is then collocated at $\nno$ source points $\by_A\in\mcalS$, for $A=1,\ldots,\nno$, to obtain a square system. The source points $\by_A$ are thus referred to as collocation points in the following. 
\\\\For Lagrange discretizations on the one hand, the collocation points are simply selected to be the same as the nodes located on the surface. For isogeometric discretizations on the other hand, the control  points are not necessarily located on the surface, which makes then unsuitable for collocation. Therefore, the locations of the collocation points are determined by the Greville abscissae (see e.g.~\cite{Greville64}, \cite{Johnson05} and \cite{Aurrichio10}). Fig.~\ref{fig:quad_sphere_discretization}a and b depict the resulting collocation points on spheres with biquadratic NURBS discretizations.
\\\\The discretized BIE \eqref{eq:BE_BIE2} is evaluated for each collocation point to obtain the discrete system
\begin{equation}\label{eq:BE_system}
\mT^\mathrm{BE}\,	 \mv
+\mN^\mathrm{BE}\, \mv
+\mG^\mathrm{BE}\, \mt 
= \bzero,
\end{equation}
where $\mN^\mathrm{BE}$, $\mT^\mathrm{BE}$ and $\mG^\mathrm{BE}$ are square matrices of size $(\dof \times \dof)$. Evaluating the BIE for collocation point $\by_A$ fills three rows of the matrices, indicated by vector $\md_A:= {[ 3\,A-2,\,3\,A-1,\,3\,A]^\mrT}$ of length 3, with
\begin{equation}\label{eq:green_int}
	\begin{aligned}
	\mN^\mathrm{BE}_{(\md_A,\,:)}&= \bcalN(\by_A),\\
	\mG^\mathrm{BE}_{(\md_A,\,:)} &= \bcalG(\by_A),\\
	\mT^\mathrm{BE}_{(\md_A,\,:)}&= \bcalT(\by_A)~,
	\end{aligned}
\end{equation}
where now $A=1,\ldots,\nno$ and $\bcalN$, $\bcalG$ and $\bcalT$ as defined in \eqref{eq:BE_discr_v2} and \eqref{eq:BE_GreensInt}. The subscripts on the left hand side of \eqref{eq:green_int} are given in MATLAB-like notation so that $\mN^\mathrm{BE}_{(\md_1,\,:)}$, for example, refers to the first three rows of $\mN^\mathrm{BE}$.
The BE system \eqref{eq:BE_system} is solvable for Dirichlet, Neumann and mixed problems. However, only Dirichlet problems are considered in the numerical examples in Sec.~\ref{sec:examples}.
\\\\Remark: The probably most intuitive way to create the BE system is to evaluate the discretized BIE \eqref{eq:BE_BIE2} for all collocation points successively. However, changing the order of the loops, i.e.~approximating the integrals on element $e$ for all collocation points $\by_A$ and repeating this step for the remaining elements subsequently, is computationally much more efficient.

%

\subsection{Mapping to a parameter domain}\label{sec:BE_mapping}
The surface $\mcalS$ is fully characterized by the parametric description
\begin{equation}\label{eq:surf_map}
\bx=\bx (\xi^\alpha)~,\hh\bx \in \mcalS
\end{equation} 
where $\xi^\alpha$ with $\alpha=1,2$ are curvilinear coordinates associated with a parameter domain $\mcalP$. The mapping  \eqref{eq:surf_map} reflects the property that the surface is a two-dimensional object embedded within three-dimensional space. Each elemental surface $\Omega^e$ is mapped to a quadrilateral master element in the parameter domain ${\xi^\alpha \in [-1,1]}$ of side length 2. Comparing approximations \eqref{eq:BE_discr_x} and \eqref{eq:BE_discr_x2} with the parametric mapping \eqref{eq:surf_map} shows that the nodal shape functions are defined on the master element, i.e.
\begin{equation}\label{eq:BE_discr_shape}
	N_I = N_I (\xi^1,\xi^2)~.
\end{equation}
This paper considers quadrilateral elements since these can be conveniently related to the master element introduced above. The reader is referred to \cite{Sauer18_CISM} for more information on mapping \eqref{eq:surf_map} and the corresponding surface description in curvilinear coordinates.


\subsection{Boundary quadrature}\label{sec:BE_quad}
The boundary integrals are first mapped to the parameter space with mapping \eqref{eq:surf_map}. Each elemental integral is defined on the master element as
\begin{equation}\label{eq:BE_int_map}
\int_{\Omega^e} k(\bx) \, \mrd a 
= \int_{-1}^1 \int_{-1}^1 k(\bxi)  \,  J_a(\bxi) \,\mrd \xi^1 \, \mrd\xi^2~,
\end{equation}
where $k(\bx)$ denotes an arbitrary integral kernel and $J_a$ denotes the local surface stretch between surface $\mcalS$ and parameter domain $\mcalP$ \citep{Sauer14_CMAME}. The right hand side of \eqref{eq:BE_int_map} is then approximated with numerical quadrature rules \citep{Gauss1815, Golub69, Laurie01} that are stated as weighted sums of function values at specified positions in the parameter space. Bivariate quadrature rules with $n_\mathrm{qp}^e$ quadrature points for the approximation of surface integrals on the master element are thus defined as 
\begin{equation}\label{eq:quad_general}
\int_{-1}^1 \int_{-1}^1 f(\bxi) \, \mrd \xi^1 \, \mrd\xi^2 
\approx \sum_{i=1}^{n_\mathrm{qp}} f(\bxi_i) \, \mrw (\bxi_i) ~, 
\end{equation}
where $f(\bxi)$ denotes the integral kernel on the master element, while $\bxi_i$ and $\mrw (\bxi_i)$ denote the position and weight of quadrature point $i$. Various bivariate quadrature rules are investigated in Sec.~\ref{sec:quad} with respect to approximation of singular integrals.

\section{Existing boundary quadrature rules}\label{app:quad}

The accurate approximation of the singular boundary integrals from \eqref{eq:BE_system} is crucial in BE analysis. Various quadrature approaches are therefore investigated in this section with respect to their suitability for singular integral approximation: The classical Gauss-Legendre quadrature rule is considered in Sec.~\ref{app:quad_Gauss1}, a modified Gauss-Legendre quadrature rule is introduced in Sec.~\ref{app:quad_Gauss2}, while Sec.~\ref{app:quad_Duffy} considers a Duffy transformation-based quadrature rule for singular integrals.

\subsection{Classical Gauss-Legendre quadrature}\label{app:quad_Gauss1}
\begin{figure}[h]
\unitlength\linewidth
\begin{picture}(1,0.27)
\put(0,0){\includegraphics[trim = 5 40 10 10, clip, width=.33\linewidth]{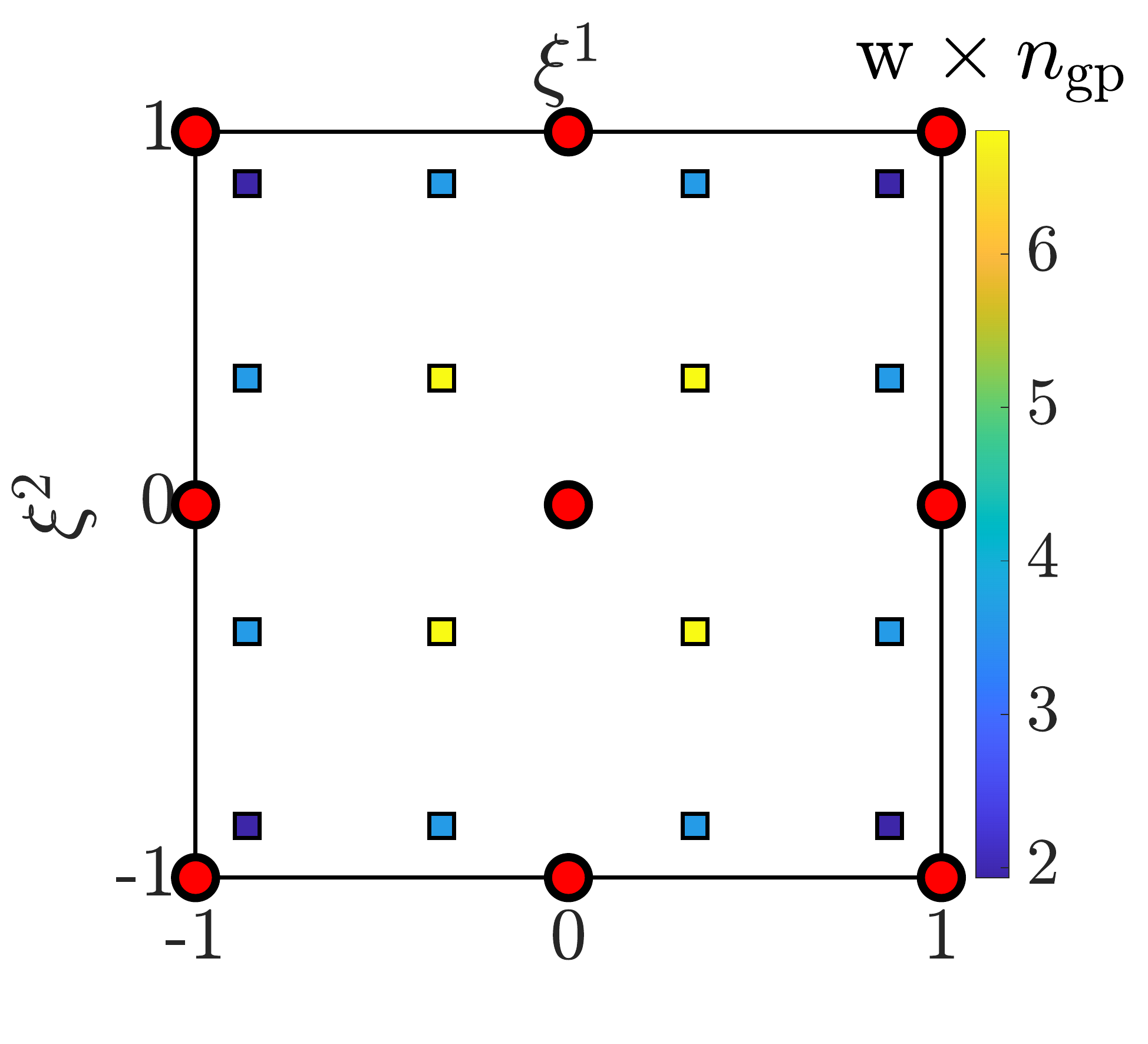}}
\put(0.335,0){\includegraphics[trim = 5 40 10 10, clip, width=.33\linewidth]{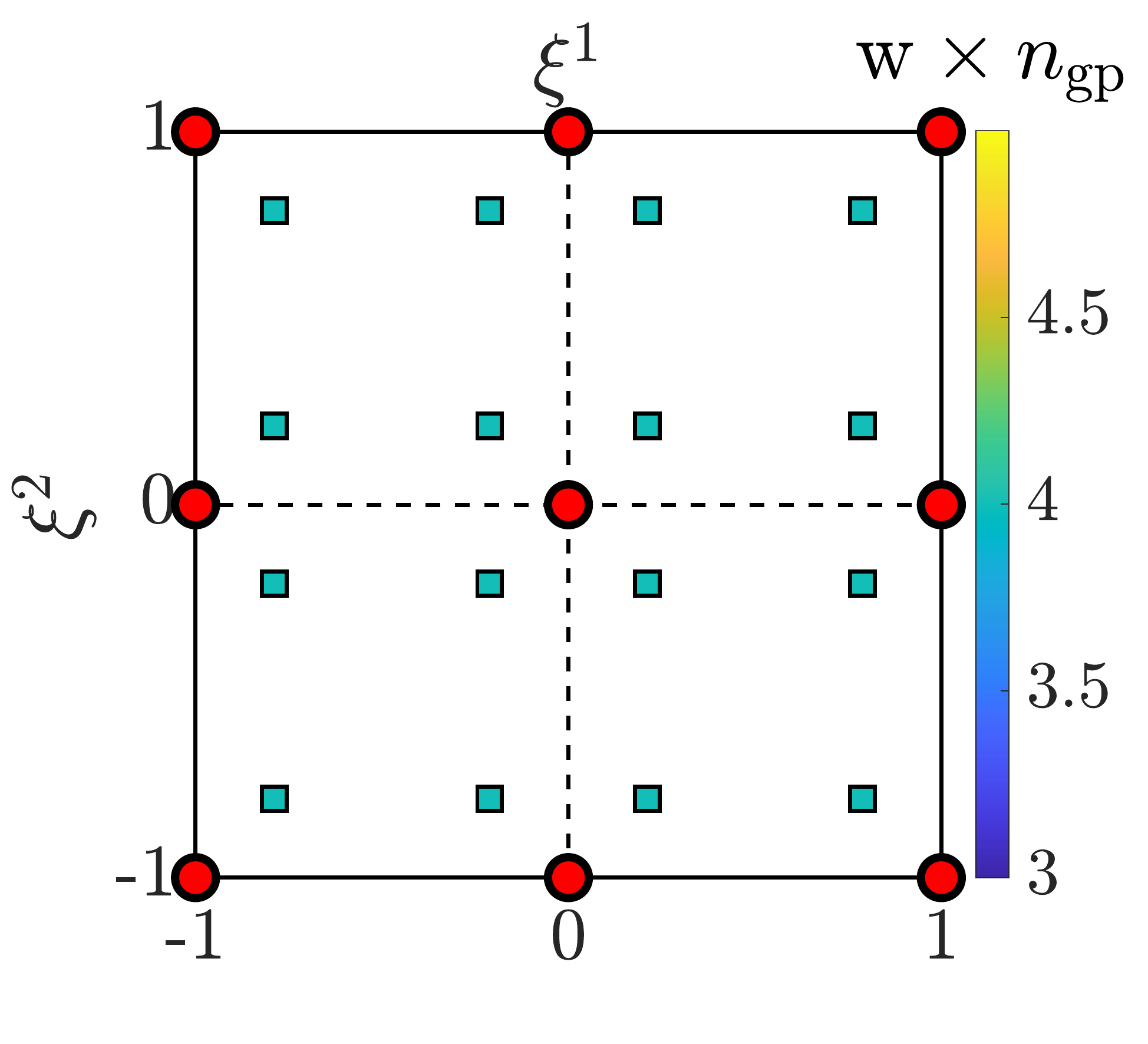}}
\put(0.67,0){\includegraphics[trim = 5 40 10 10, clip, width=.33\linewidth]{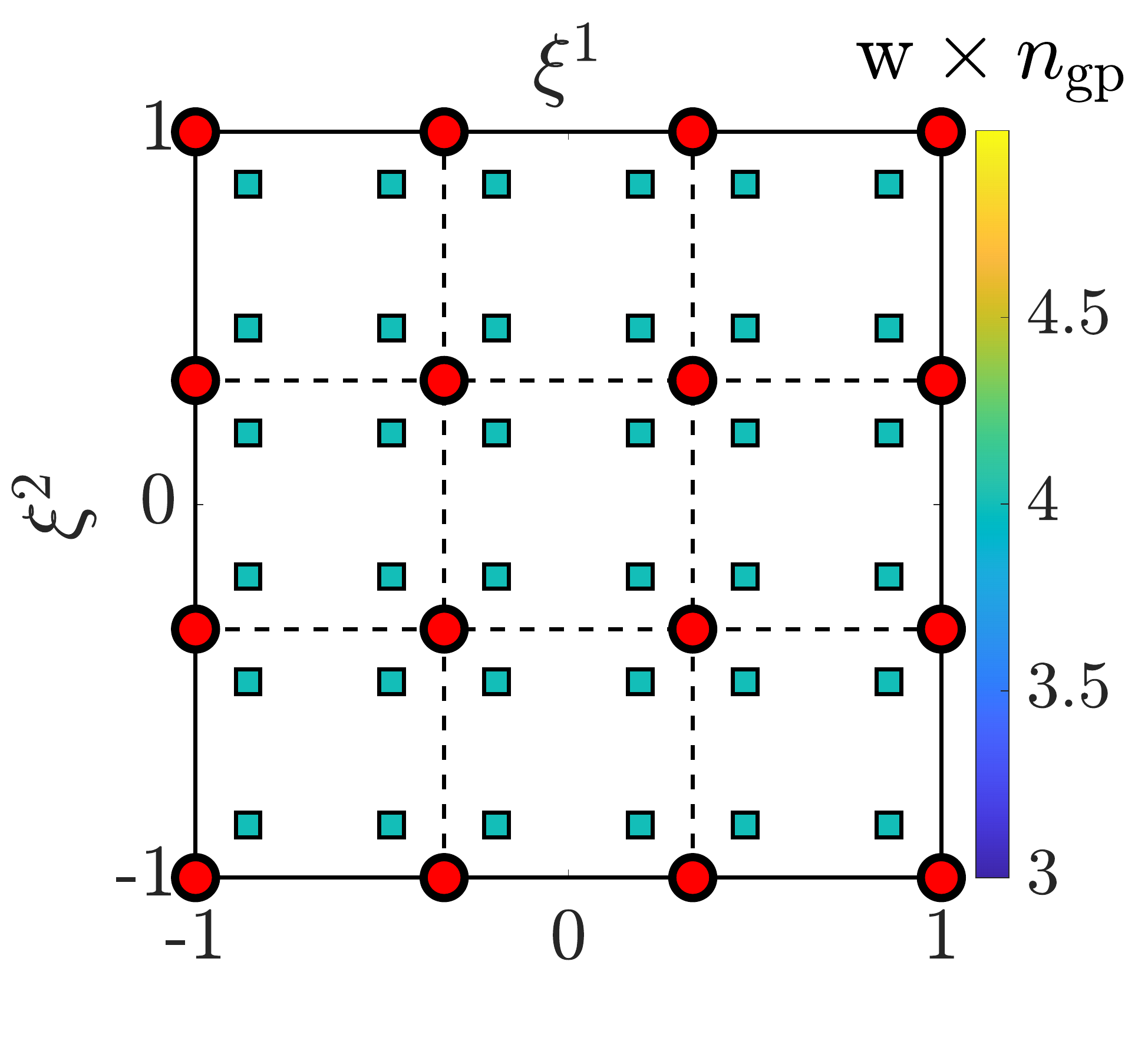}}
\put(0,0){a.} \put(0.335,0){b.} \put(0.67,0){c.}
\put(0.025,0.05){$\by_1$} \put(0.15,0.045){$\by_2$} \put(0.15,0.155){$\by_3$}
\end{picture}
\caption{\textit{Gauss-Legendre quadrature}: Quadrature points and weights for Gauss-Legendre quadrature~(GL) with $\tilde n_\mathrm{qp}=4$. Classical GL with $n_\mathrm{qp}^e = 4 \times 4$~(a.), modified GL on a biquadratic discretization with $n_\mathrm{qp}^e = 4 \times 4$~(b.) and on a bicubic discretization with $n_\mathrm{qp}^e = 6 \times 6$ according to \eqref{app_eq:quad_GL2_nsub}~(c.).}\label{fig:quad_GL}
\end{figure}
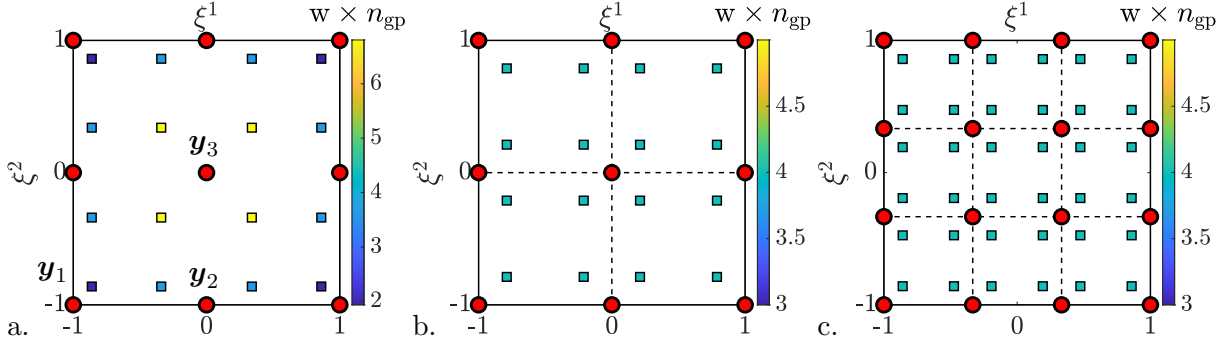

\begin{figure}[h]
\unitlength\linewidth
\begin{picture}(1,0.44)
\put(0,0){\includegraphics[trim = 0 0 40 20, clip, width=.5\linewidth]{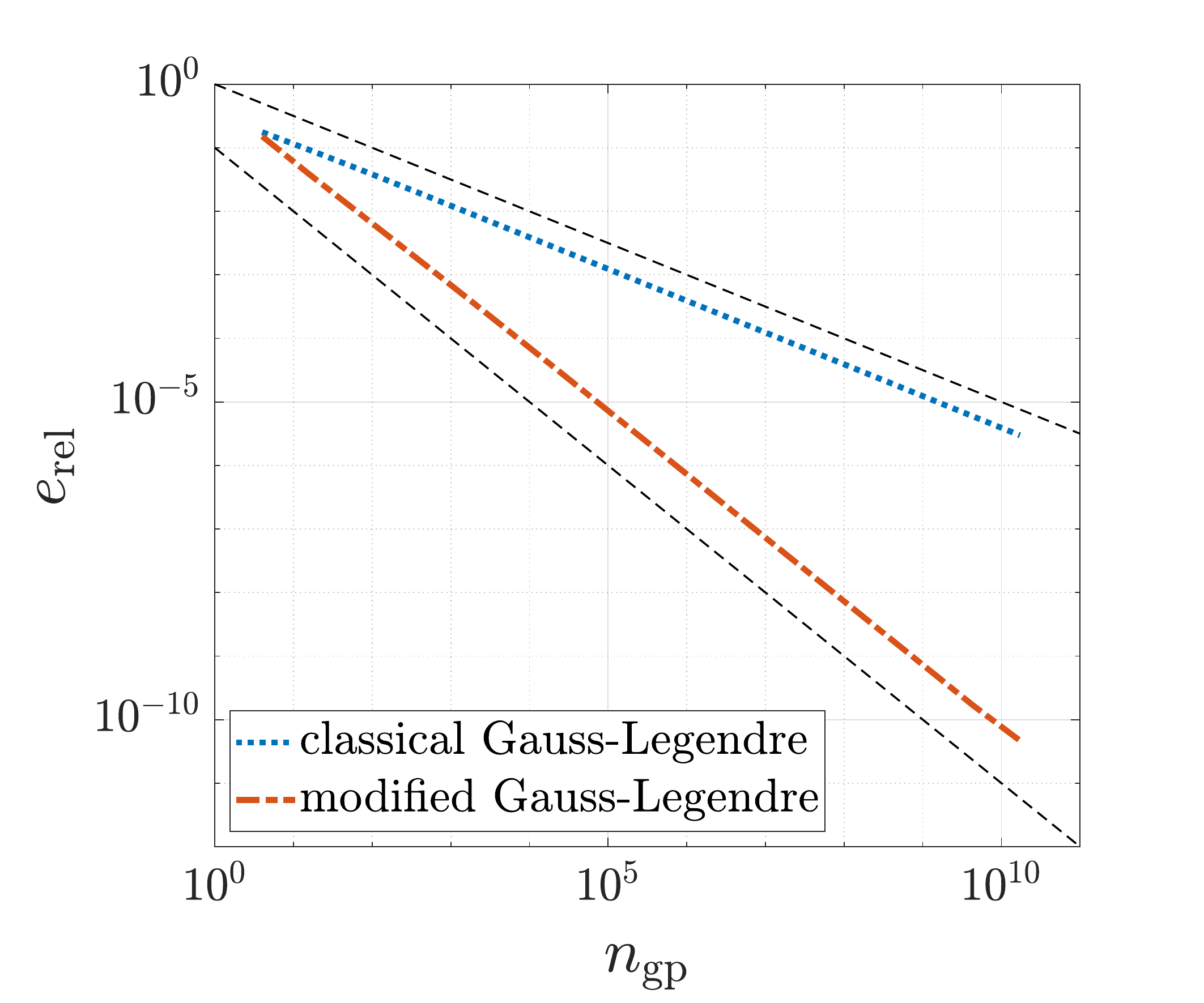}}
\put(0.22,0.2){$\propto 1/n_\mathrm{qp}^e$} \put(0.3,0.34){$\propto 1/\sqrt{n_\mathrm{qp}^e}$}
\put(0.52,0){\includegraphics[trim =40 0 40 20, clip,  width=.465\linewidth]{figures/quad/single/err_dist_GL_nqp_2_131072.pdf}}
\put(.045,0.02){a.}
\put(.53,0.02){b.} \put(.94,0.14){\rotatebox{90}{\large $\varepsilon_\mathrm{machine}$}}
\put(.585,0.255){$e_\mathrm{class}$}
\put(.585,0.099){$e_\mathrm{mod}$}
\wput{0.26}{0.018}{$n_\mathrm{qp}^e$}
\end{picture}
\caption{\textit{Gauss-Legendre quadrature}: a.~Mean quadrature error $e_\mathrm{rel}$ for kernel $1/r$ vs.~the number of quadrature points. b.~$e_\mathrm{rel}$ vs.~the minimum distance $r_\mathrm{min}$ (see Fig. \ref{fig:quad_res}a).}\label{fig:quad_GL_res}
\end{figure}

Gauss-Legendre quadrature has been introduced by \cite{Gauss1815} and is nowadays the most common and established quadrature
rule for numerical integration. Bipolynomial kernels of orders $p$ and $q$ are integrated exactly by Gauss-Legendre quadrature with a minimum number of 
\begin{equation}\label{app_eq:quad_Gauss_nqp}
n_\mathrm{qp}^e= {\lceil (p+1)/2\rceil \times \lceil (q+1)/2\rceil}
\end{equation}
quadrature points per element, while non-polynomials kernels are only integrated approximately. This work focus on bivariate Gauss-Legendre quadrature rules with $n_\mathrm{qp}^e=\tilde n_\mathrm{qp}\times \tilde n_\mathrm{qp}$ without loss of generality. The quadrature point locations are thus denoted by
\begin{equation}\label{app_eq:quad_Gauss_xi}
	\bxi_i= [\tilde\xi_j,\, \tilde\xi_k]^\mrT~,
\end{equation}
for $i:=j+(k-1)\,\tilde n_\mathrm{qp}$, where $j$ and $k$ run from $1$ to $\tilde n_\mathrm{qp}$. The scalar values $\tilde\xi_j$ and $\tilde\xi_k$ denote the $j$-th and the $k$-th root of the Legendre polynomial $P_{\tilde n_\mathrm{qp}}(\xi)$, respectively. They can be determined with numerical methods like Newton's method or by exploiting explicit expressions \citep{Golub69} or tables \citep{Laurie01}. 
The corresponding weights are given by
\begin{equation}\label{app_eq:quad_Gauss_w}
\mrw_i= \tilde \mrw_j\, \tilde \mrw_k~, 
\end{equation}
with
\begin{equation}\label{app_eq:quad_Gauss_w1D}
\tilde \mrw_j = \frac{2}{\left(1-\tilde \xi_j^2\right) \big[P_{\tilde n_\mathrm{qp}}'(\tilde \xi_j)\big]^2}
\end{equation}
and analogously for $\tilde \mrw_k$.
Fig.~\ref{fig:quad_GL}a shows the quadrature point locations on the master element and the corresponding weight values for classical Gauss-Legendre quadrature with $n_\mathrm{qp}^e= 4\times 4$.
\\The application of classical Gauss-Legendre quadrature to regular, i.e.~non-singular, kernels is accurate and robust \citep{Chawla68, Kambo70, Leone79}. However, the numerical integration of singular kernels is much more challenging. The \textit{weakly singular} integral
\begin{equation}\label{app_eq:quad_int}
\mcalI_A:=\int_\mcalS \frac{1}{r_A} \,\mrd a~,
\end{equation}
where $r_A:=\|\bx-\by_A\|$,\footnote{Integral \eqref{app_eq:quad_int} is representative of the behavior of the two singular integrals in the BIE (\eqref{eq:flow_BIE_surf} and \eqref{eq:BE_system}).} is investigated in the following considering the collocation points depicted in Fig.~\ref{fig:quad_GL}a:
\begin{packeditemize}
	\item $\by_1=[-1,-1]^\mrT$ at the corner of the element,
	\item $\by_2=[0,-1]^\mrT$ at the midpoint of an edge,
	\item and $\by_3=[0,0]^\mrT$ at the elemental midpoint.
\end{packeditemize}
Fig.~\ref{fig:quad_GL_res}a shows that the mean relative quadrature error defined by
\begin{equation}\label{app_eq:quad_err}
e_\mathrm{rel}:= \frac{1}{n_A}\,\sum_{A=1}^{n_A}{ \frac{\|\mcalI_A^h - \mcalI_A\|}{\mcalI_A}}
\end{equation}
decreases for classical Gauss-Legendre quadrature (blue line) with increasing number of quadrature points, i.e.~$\tilde n_\mathrm{qp}=2^1,2^2,\ldots, 2^{17}$, where $\mcalI^h_A$ denotes the numerical approximation of integral \eqref{app_eq:quad_int}. However, the convergence rate is only $1/2$, so the highest number of quadrature points $n_\mathrm{qp}^e=2^{17}\times 2^{17}\approx {1.72 \times10^{10}}$ still results in $e_\mathrm{rel} \approx {3\times 10^{-6}}$.
\\Fig.~\ref{fig:quad_GL_res}b depicts the mean quadrature error
vs.~the minimum distance between quadrature points and collocation points that is defined by
\begin{equation}\label{app_eq:quad_dist}
 r_\mathrm{min} :=\min_{\substack{ i= {1,\ldots, n_\mathrm{qp}} \\A={1,2,3}}} {\left\| \bx(\bxi_i)-\by_A\right\|}~,
\end{equation}
where $\bx(\bxi_i)$ denote the position on $\mcalS$ for quadrature point $i$. The vertical dashed line shows $\varepsilon_\mathrm{machine} \approx {2.2\times 10^{-16}}$, which is the smallest computationally admissible value for $r_\mathrm{min}$. Extending the curve by extrapolating the simulation data shows that the mean error of classical Gauss-Legendre quadrature is bounded by $e_\mathrm{rel}>e_\mathrm{class}\approx3\times10^{-9}$, even with unlimited memory capacity. A huge number of $n_\mathrm{qp}>10^{16}$ quadrature points would be required to obtain such a small error. It should be further noted that classical Gauss-Legendre quadrature does in general not prevent the coincidence of collocation points and quadrature points. For instance, using an odd number for $\tilde n_\mathrm{qp}$ results in a quadrature point located at the middle of the element, where in many cases a collocation point is also located. The approximation of the integral $I_3$ \eqref{app_eq:quad_int} is infinite for these cases and thus useless.\footnote{The coincidence of collocation points and quadrature points potentially also occurs for even $\tilde n_\mathrm{qp}$. The collocation point locations depend on the order of the shape functions used for discretization.}
\\\\In conclusion, three main drawbacks for singular BE integral approximation with classical Gauss-Legendre quadrature can be identified: The quadrature rule
\begin{packeditemize}
	\item[--] is very inefficient with respect to the number of quadrature points,	
	\item[--] has a lower error bound of $e_\mathrm{class}\approx3\times10^{-9}$,
	\item[--] does not prevent the coincidence of collocation points and quadrature points.
\end{packeditemize}
On the other hand, classical Gauss-Legendre quadrature is
\begin{packeditemize}
	\item[+] simple to implement and computationally efficient.\footnote{{The shape function values and the corresponding weights are determined once at the quadrature points on the master element. Isogeometric shape function values are determined from universal Bernstein polynomials and element specific B\'ezier extraction operators \citep{Borden11}.}}
\end{packeditemize}
The disadvantages in accuracy and robustness clearly outweigh the advantageous implementation, making the classical Gauss-Legendre quadrature unsuitable for the approximation of singular BE integrals.

\subsection{Modified Gauss-Legendre quadrature}\label{app:quad_Gauss2}
The classical Gauss-Legendre rule can be modified to overcome the discussed deficiencies, while maintaining the efficient equal treatment of all elements. This modification is briefly introduced by \cite{Heltai14} for biquadratic NURBS and is generalized here for discretizations with Lagrange or isogeometric basis functions of arbitrary order. The key idea of the modification is to split the surface elements into virtual sub-elements such that collocation points are exclusively located at the corners of sub-elements. Classical Gauss-Legendre quadrature is then applied to the virtual sub-elements.
\\\\Choosing the collocation points on a discretization of order $p$ and $q$ as shown in Fig.~\ref{fig:quad_sphere_discretization} yields collocation point locations on the master element at
\begin{equation}
\xi_\mathrm{col}^i=-1+\frac{2\,(i-1)}{p}~,\hh i=1,\ldots,p+1
\end{equation}
and analogous for $q$ and $\eta_\mathrm{col}^j$. The surface elements are split into virtual sub-elements at $\xi_\mathrm{col}^i$ for $i=1,\ldots,p+1$ and at $\eta_\mathrm{col}^j$ for $j=1,\ldots,q+1$. The sub-elements are illustrated by dashed lines in Fig.~\ref{fig:quad_GL}b for biquadratic elements ($p\!=\!q=\!2$) and in Fig.~\ref{fig:quad_GL}c for bicubic elements ($p\!=\!q=\!3$).
Classical Gauss-Legendre quadrature rules are then applied to the virtual sub-elements without coincidence of quadrature points and collocation points: The collocation points are located at the corner of the sub-elements, while Gauss-Legendre quadrature points are located within the element. The use of sub-elements furthermore provides an equal minimum distance \eqref{app_eq:quad_dist} for all collocation points. A surface element of order $p$ and $q$ consists of $n_\mathrm{sub} = p\times q$ sub-elements that are of equal size in the parametric space.\footnote{This statement holds for isogeometric as well as for Lagrange discretizations. Lagrange discretizations requires a higher number of nodal points than the isogeometric counterpart of the same order and thus also a higher number of collocation points. For Lagrange discretizations a collocation point is located at every corner point of a virtual sub-element, while several possible positions are empty for isogeometric discretizations. However, in both cases there is no collocation point at any other position than $\bxi_\mathrm{col}$ and $\bet_\mathrm{col}$.}
\\\\
Note, that the sub-elements are only introduced to illustrate the quadrature strategy. The discretization and thus the number of elements and nodes remains unchanged. Code-wise, the division into $n_\mathrm{sub}$ sub-elements does not appear at all. Instead, the quadrature point locations and the corresponding weights are defined by so-called modified quadrature rules. The classical quadrature rules of the virtual sub-elements are combined into one modified quadrature rule for the whole element. The quadrature weights of the modified Gauss-Legendre rule are obtained as
\begin{equation}\label{app_eq:quad_GL2_w}
\mw =\frac{1}{p\, q}\, \big[\underbrace{\mw_\mathrm{GL}^\mrT}_{i=1},\underbrace{\mw_\mathrm{GL}^\mrT}_{i=2}, \ldots, \underbrace{\mw_\mathrm{GL}^\mrT}_{i=p\,q}\, ]^\mrT~,
\end{equation}
where $\mw_\mathrm{GL}$ denotes a column vector that contains the weights for classical Gauss-Legendre quadrature \eqref{app_eq:quad_Gauss_w}. With the master element defined on the interval $I_0={[-1,1] \times [-1,1]}$ in the parametric domain, the sub-elements are defined on the intervals ${I_1=[\xi^1_\mathrm{col},\xi^2_\mathrm{col}] \times [\eta^1_\mathrm{col},\eta^2_\mathrm{col}]}$, ${I_2=[\xi^2_\mathrm{col},\xi^3_\mathrm{col}] \times [\eta^1_\mathrm{col},\eta^2_\mathrm{col}]}$, \ldots, $I_{p\,q}=[\xi^{p-1}_\mathrm{col},\xi^p_\mathrm{col}] \times [\eta^{q-1}_\mathrm{col},\eta^q_\mathrm{col}]$.
The quadrature point positions of the modified Gauss-Legendre quadrature rule then are
\begin{equation}\label{app_eq:quad_GL2_xi}
\bXi = \left[\bXi_1^\mrT,\, \bXi_2^\mrT, \ldots , \, \bXi_{pq}^\mrT\, \right]^\mrT ~,
\end{equation}
where $\bXi_i$ denotes the classical quadrature points $\bXi_0:=[\bxi_1^\mrT, \bxi_2^\mrT, \ldots, \bxi_{pq}^\mrT]$ from \eqref{app_eq:quad_Gauss_xi}, shifted to interval $I_i$.\footnote{For a biquadratic element, the original quadrature points are shifted to the intervals $I_1$, $I_2$, $I_3$ and $I_4$ such that
$\bXi_1= (\bXi-\bi\,[1, 1] )/2$ and $\bXi_2= (\bXi+\bi\,[1, -1] )/2$ and $\bXi_3= (\bXi+\bi\,[-1, 1] )/2$ and $\bXi_4= (\bXi+\bi\,[1, 1] )/2$,where $\bi$ denotes a ($n_\mathrm{qp}\times 1$) array with only `1' entries.}
The total number of elemental quadrature points for a modified bivariate Gauss-Legendre quadrature rule is accordingly given by
\begin{equation}\label{app_eq:quad_GL2_n}
n_\mathrm{qp}^e= p\,q\,(\tilde n_\mathrm{qp}^\mathrm{sub})^2~,
\end{equation}
where $\tilde n_\mathrm{qp}^\mathrm{sub}$ denotes the number of quadrature points per dimension on each sub-element.
In order to obtain a comparable, but not a smaller $n_\mathrm{qp}^e$ as for the equivalent classical Gauss-Legendre rule \eqref{app_eq:quad_Gauss_nqp}, the sub-elemental number of quadrature points is chosen as 
\begin{equation}\label{app_eq:quad_GL2_nsub}
\tilde n_\mathrm{qp}^\mathrm{sub}:= \left\lceil \frac{\tilde n_\mathrm{qp}}{p}\right\rceil~.
\end{equation}
Fig.~\ref{fig:quad_GL} shows the quadrature point positions and the corresponding weights for classical Gauss-Legendre quadrature (a), for modified Gauss-Legendre quadrature on a biquadratic element (b) and for modified Gauss-Legendre quadrature on a bicubic element (c), where $\tilde n_\mathrm{qp}=4$ for all three examples. The total number of elemental quadrature points yields accordingly $n_\mathrm{qp}^e=\! 4\!\times \!4\!=\!16$, $n_\mathrm{qp}^e\!= \!4 \lceil 4/2 \rceil^2\!=\!16$ and $n_\mathrm{qp}^e= 9 \lceil 4/3 \rceil^2=36$, respectively.
\\\\The introduced modified Gauss-Legendre quadrature is used to approximate the singular integral \eqref{app_eq:quad_int} for the same collocation points $\by_1^0$, $\by_2^0$ and $\by_3^0$ as in Sec.~\ref{app:quad_Gauss1} to compare the accuracy of both quadrature rules. The mean quadrature error for modified Gauss-Legendre quadrature is illustrated by red lines in Fig.~\ref{fig:quad_GL_res}. It yields convergence rates of $1$ with respect to the number of quadrature points~(a) and with respect to the minimum distance between collocation points and quadrature points~(b). It can be seen that modified Gauss-Legendre quadrature shows a twice higher convergence rate than the classical counterpart. The mean quadrature error for modified Gauss-Legendre quadrature can thus be reduced to the range of machine precision by increasing the number of quadrature points to $n_\mathrm{qp}^e>10^{15}$.
\\\\Summarizing, the modified Gauss-Legendre quadrature rule eliminates the first and second drawback of the classical one. Also the third drawback has been addressed, as the required number of quadrature points to achieve a certain precision has been reduced due to the better convergence rate. However, a convergence rate of $1$ still has additional potential for improvement. Even the modified quadrature rule still requires at least $n_\mathrm{qp}^e=1,\!000$ to obtain $e_\mathrm{rel}\leq 10^{-3}$.
Besides addressing the drawbacks, the modified Gauss-Legendre quadrature also maintains the desirable property from the classical Gauss-Legendre quadrature: It is sufficient to evaluate the shape functions once on the master element and reuse the values for all elements independently of the collocation point position.
\\\\In conclusion, the presented modified Gauss-Legendre quadrature
\begin{packeditemize}
	\item[$\circ$] shows moderate efficiency with respect to the number of quadrature points,	
	\item[$+$] leads to a vanishing error with progressing quadrature refinement,
	\item[$+$] prevents the coincidence of collocation points and quadrature points,
	\item[$+$] is simple to implement and computationally efficient.
\end{packeditemize}
The numerical investigation of the modified Gauss-Legendre quadrature shows a strong robustness and a moderate accuracy making it suitable for the approximation of singular integrals. A special quadrature rule for singular functions is presented in Sec.~\ref{app:quad_Duffy} to achieve even more accurate integral approximations.

\subsection{Duffy transformation-based quadrature}\label{app:quad_Duffy}
\begin{figure}[h]
\unitlength\linewidth
\begin{picture}(1,0.3)
\put(0,0){\includegraphics[trim = 0 35 40 0, clip, height=.3\linewidth]{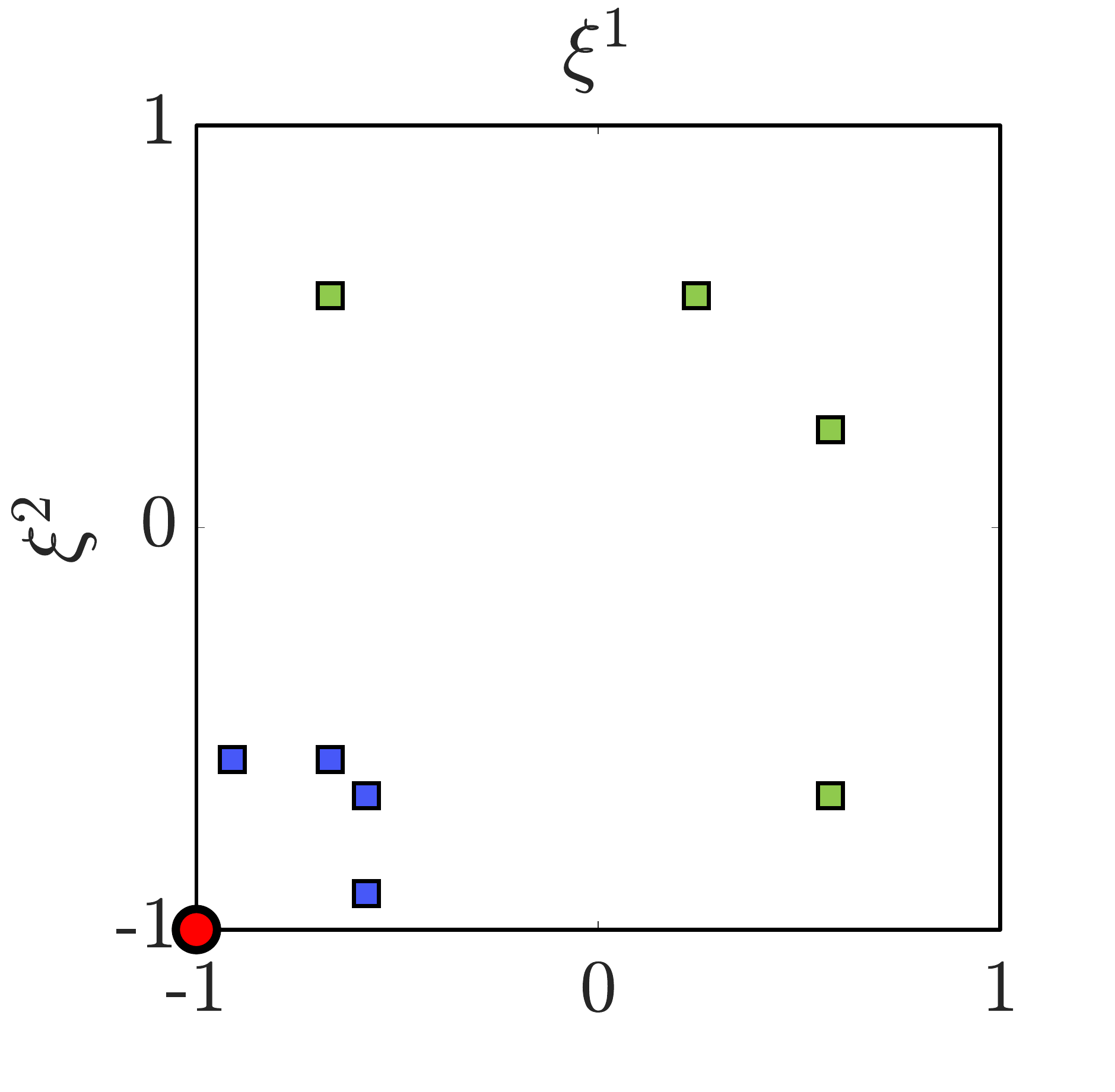}}
\put(0.3425,0){\includegraphics[trim = 25 40 40 0, clip, height=.3\linewidth]{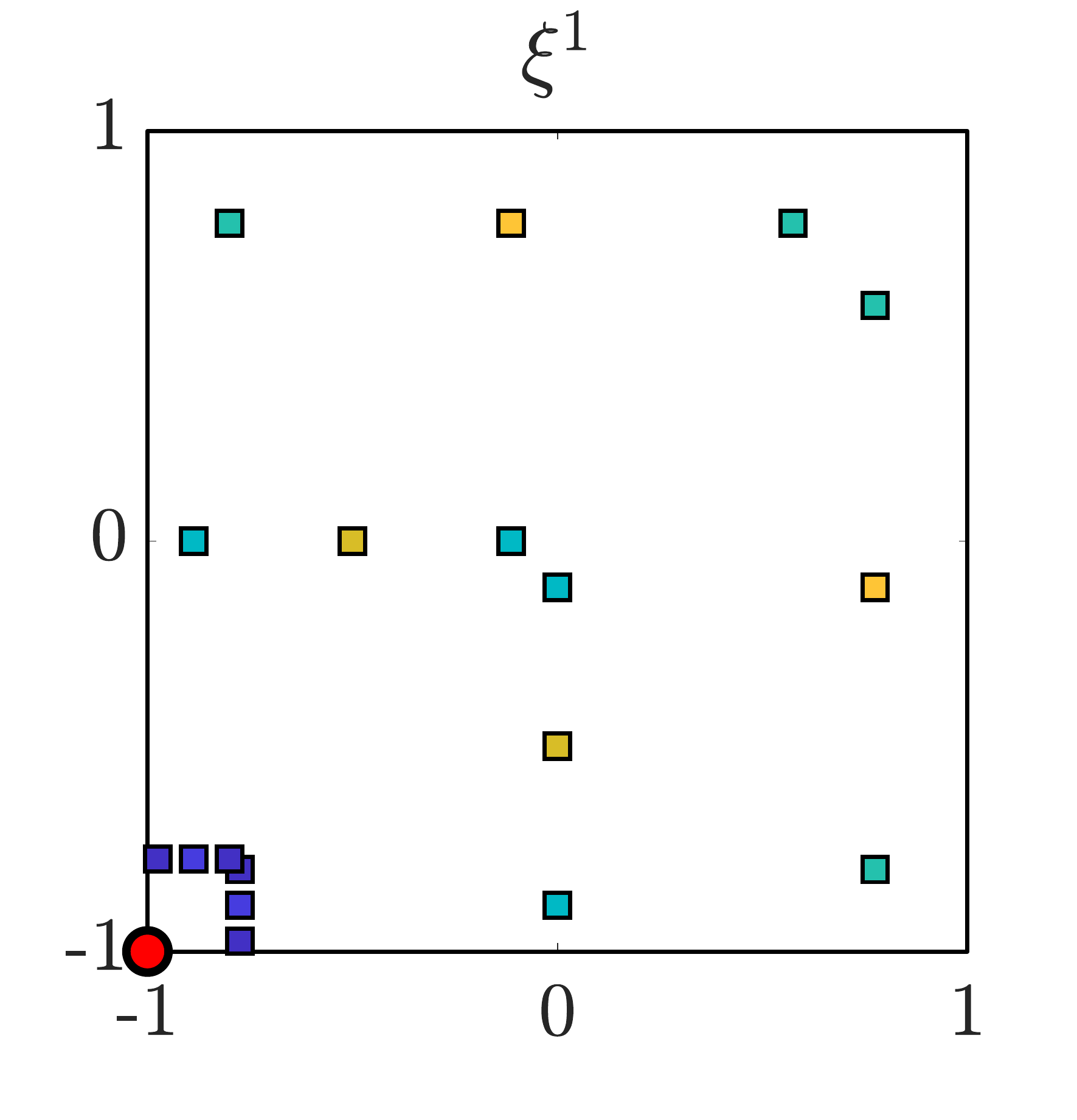}}
\put(0.65,0){\includegraphics[trim = 30 40 10 0, clip, height=.3\linewidth]{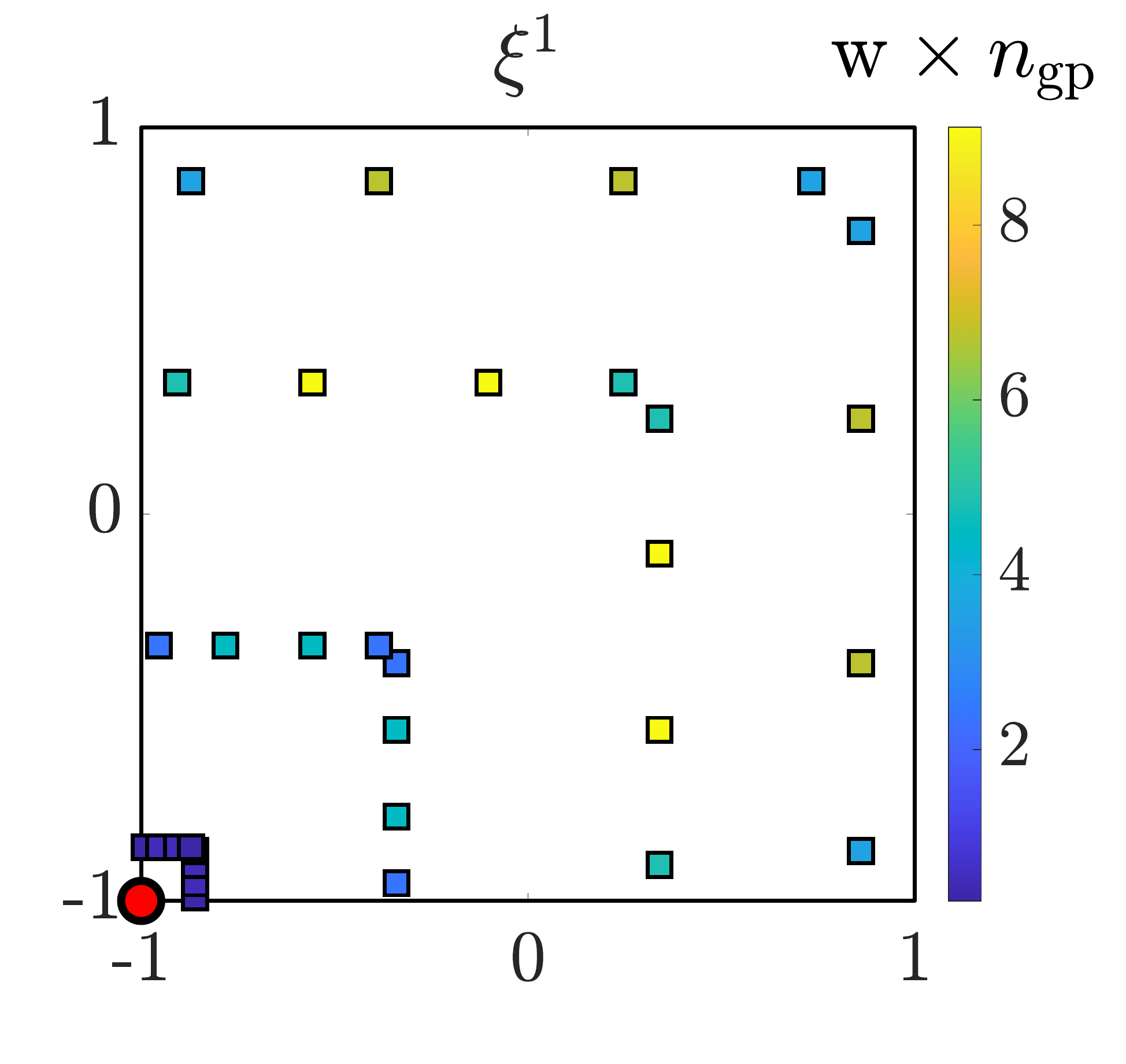}}
\put(0.02,-0.015){a.} \put(0.3425,-0.015){b.} \put(0.65,-0.015){c.}  \wput{0.926}{0.276}{\large $n_\mathrm{qp}$}
\end{picture}
\caption{\textit{Duffy quadrature}: Duffy quadrature points and their weights on the master element for collocation point $\by_1^0$ with ${\tilde n_\mathrm{qp}=2}$~(a.), ${\tilde n_\mathrm{qp}=3}$~(b.) and ${\tilde n_\mathrm{qp}=4}$~(c.).}\label{fig:quad_Duff}
\end{figure}
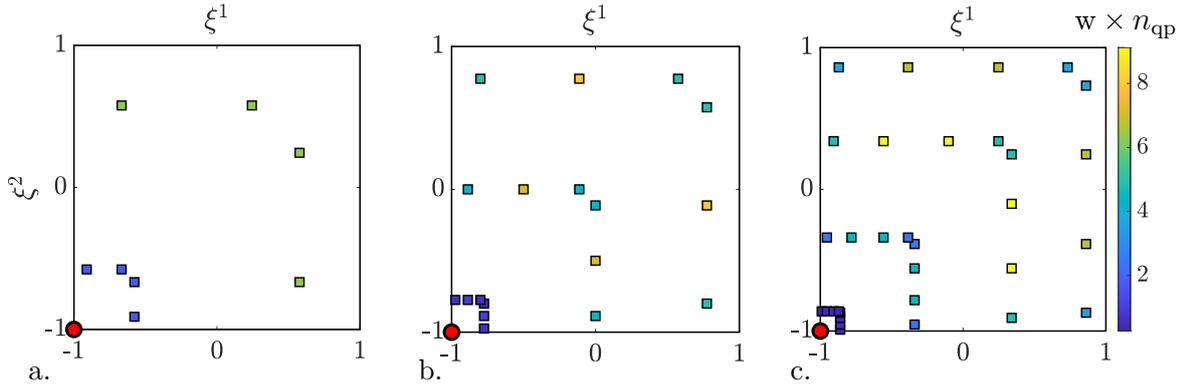

\begin{figure}[h]
\unitlength\linewidth
\begin{picture}(1,.28)
    \put(0.025,-0.01){\includegraphics[trim = 50 50 50 50, clip, width=.3\linewidth]{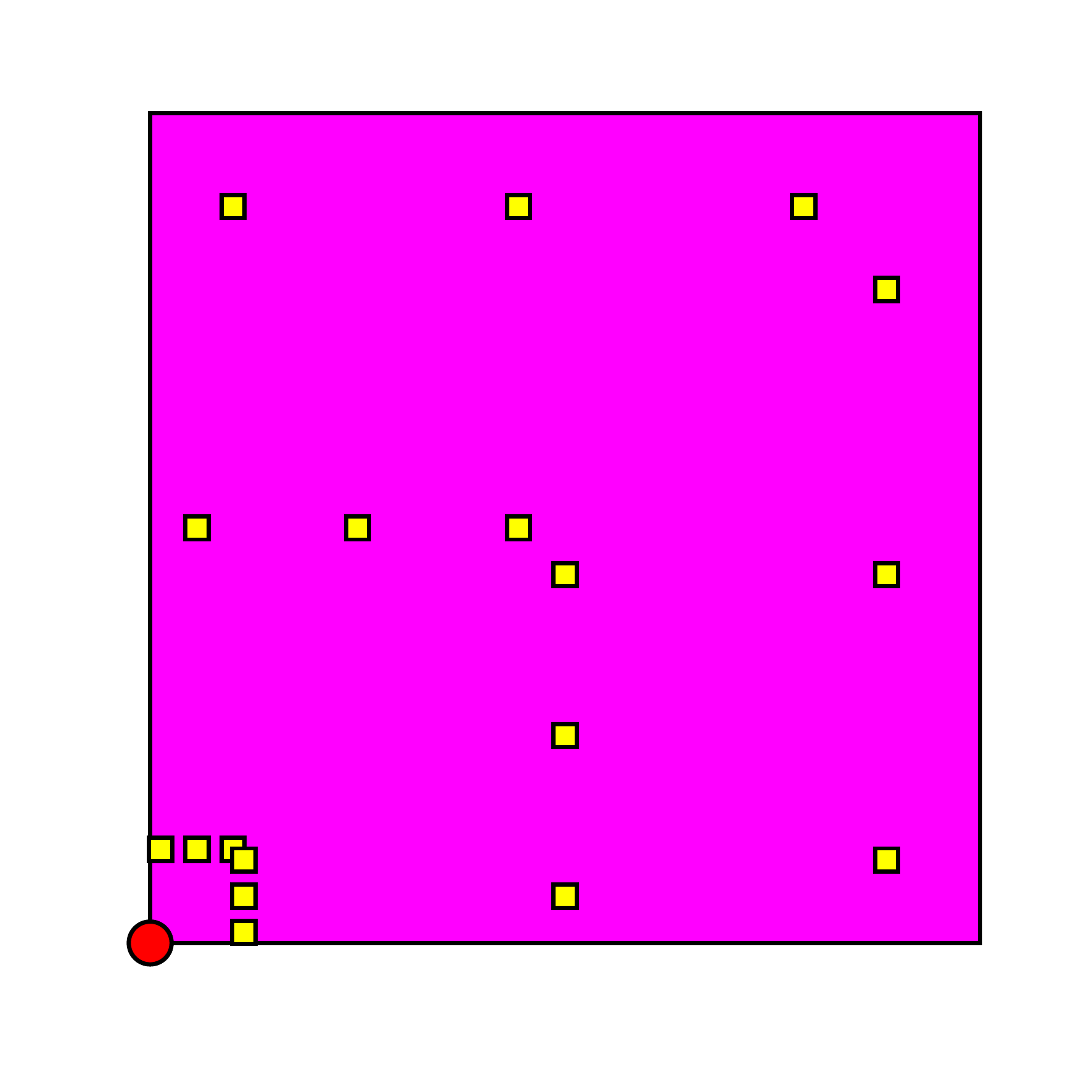}}
    \put(0.35,-0.01){\includegraphics[trim = 50 50 50 50, clip, width=.3\linewidth]{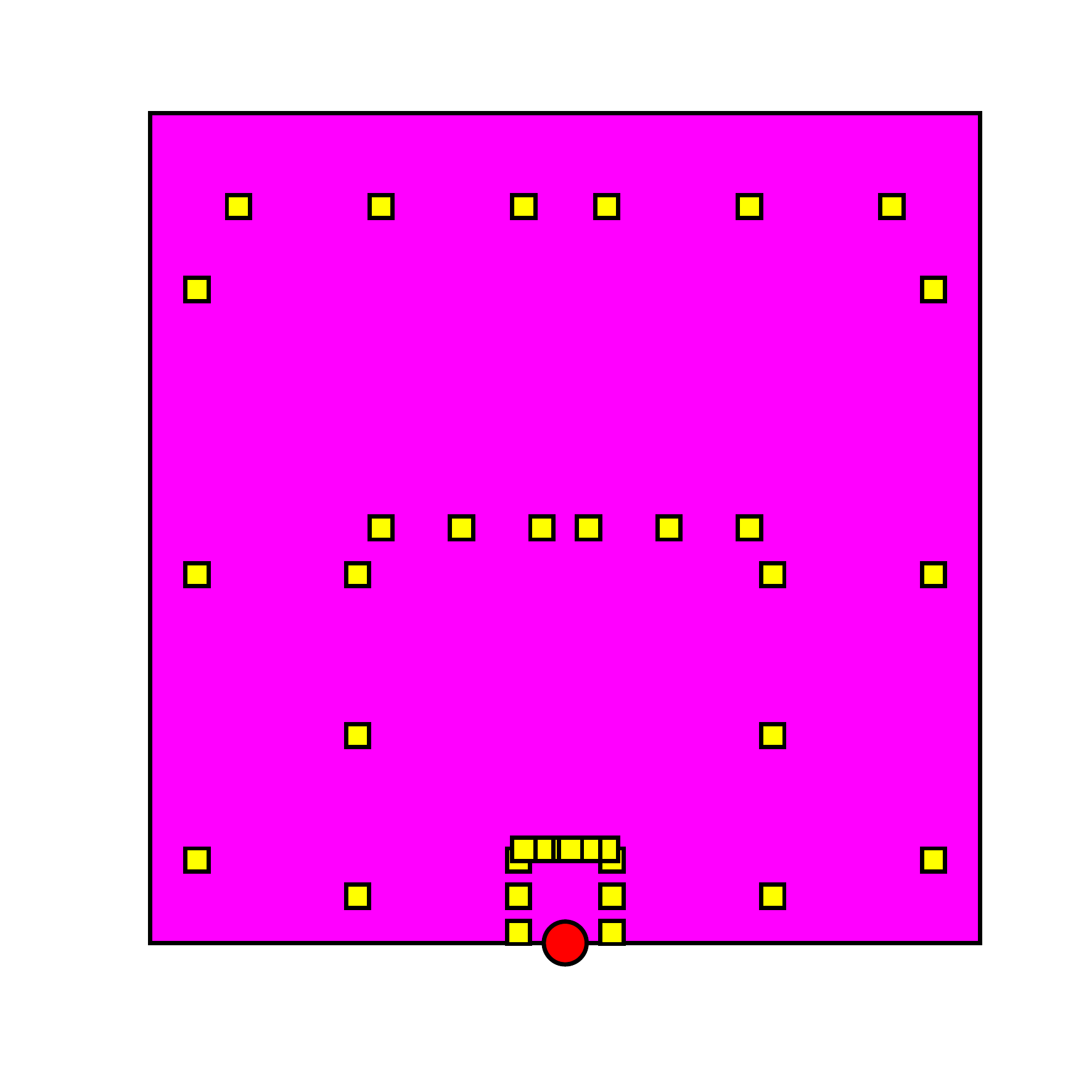}}
    \put(0.675,-0.01){\includegraphics[trim = 50 50 50 50, clip, width=.3\linewidth]{figures/quad/single/duffy/n_3_y5.pdf}}
	\put(0.01,0){a.}\put(0.34,0){b.}\put(0.665,0){c.}
	\end{picture}
\caption{\textit{Duffy quadrature}: Duffy quadrature points for $\tilde n_\mathrm{qp}=3$ considering collocation points at a corner point of an element (a.), at the midpoint of an elemental edge (b.) and at the midpoint of an element (c.).}\label{fig:quad_Duff2}
\end{figure}

The introduced Gauss-Legendre quadrature rule comes along with significant accuracy shortcomings as discussed in the previous two sections. The quadrature approach from \cite{Fairweather79} and \cite{Duffy82} is considered here to overcome this shortcomings. This approach exploits Duffy transformation from a triangle to a square to remove corner singularities of the type $1/r$. The quadrature rule is referred to as {Duffy transformation-based quadrature} or as {Duffy quadrature} in short form. While \cite{Duffy82} and \cite{Mousavi10} use Duffy quadrature for integral approximation on triangular elements, two Duffy triangles are combined here to define the quadrature rule on the master element. The Duffy quadrature is therefore applicable to arbitrary quadrilateral elements.
\\\\Fig.~\ref{fig:quad_Duff} shows the quadrature point locations on the master element and the corresponding weight values 
for collocation point $\by^0_1 = [-1,-1]^\mrT$ and $\tilde n_\mathrm{qp}=2,3,4$. For Duffy quadrature $\tilde n_\mathrm{qp}$ refers to the number of quadrature points per dimension on each of the two triangles. The total number of Duffy quadrature points per element thus yields 
\begin{equation}\label{app_eq:quad_Duffy_nqp}
n_\mathrm{qp}^e= 2 \tilde n_\mathrm{qp}^2~.
\end{equation}
The quadrature rules from Fig.~\ref{fig:quad_Duff} accordingly consist of $n_\mathrm{qp}^e=8$ (a), $n_\mathrm{qp}^e=18$ (b) and $n_\mathrm{qp}^e=32$ (c) quadrature points, respectively.
\\The singular integral \eqref{app_eq:quad_int} is approximated using Duffy quadrature with an increasing number of quadrature points to compare its accuracy to Gauss-Legendre quadrature. Fig.~\ref{fig:quad_res}b shows that a coarse Duffy quadrature with $n_\mathrm{qp}^e=2$ yields a similar accuracy as a refined modified Gauss-Legendre quadrature with $n_\mathrm{qp}^e=25$. The Duffy quadrature also shows a much better convergence behavior:  The relative error converges to the analytical solution with $e_\mathrm{rel} \propto 1/{n_\mathrm{qp}^e}^{2.5}$ for small $n_\mathrm{qp}^e\leq 18$ and with $e_\mathrm{rel} \propto 1/{n_\mathrm{qp}^e}^{12}$ for higher $n_\mathrm{qp}^e$ such that is in the range of machine precision for $n_\mathrm{qp}^e>200$.
%
\\\\In conclusion, the introduced Duffy transformation-based quadrature 
\begin{packeditemize}
	\item[+] shows outstanding efficiency with respect to the number of quadrature points,	
	\item[+] provides a perfect approximation already for moderate quadrature refinement,
	\item[+] prevents the coincidence of collocation points and quadrature points.
\end{packeditemize}
Duffy quadrature is several orders of magnitude more accurate than modified Gauss-Legendre quadrature, not to mention classical Gauss-Legendre quadrature, and is thus highly suitable for the approximation of strongly singular integrals. However, the Duffy quadrature presented is only directly applicable to corner collocation points and needs to be modified for other locations (see Fig.~\ref{fig:quad_Duff2}).
Since quadrature points and weights depend on the particular collocation point, Duffy quadrature 
\begin{packeditemize}
	\item[--] is more complex to implement than Gauss-Legendre quadrature.
\end{packeditemize}
Three different quadrature approaches are investigated in Sec.~\ref{sec:quad} with respect to their suitability for singular integral approximation. First, classical Gauss-Legendre quadrature is not robust and very inefficient and is thus not recommend for the quadrature of singular kernels. Second, modified Gauss-Legendre quadrature is very robust, simple to implement and yields moderate accuracy. As for classical Gauss-Legendre quadrature, it is sufficient to evaluate the shape functions only once on the master element for all collocation points and elements. Third, Duffy quadrature approximates singular integrals very efficiently, but the quadrature points and weights depend on the location of the collocation point. However, this disadvantage is clearly outweighed by the enormous gain in accuracy. In conclusion, two investigated methods are recommend for singular integral approximation: Highly accurate and efficient Duffy quadrature and modified Gauss-Legendre quadrature that is robust and simple to implement.

\section{Quadrature point locations and weights}\label{app:table}
\begin{figure}[h]
\unitlength\linewidth
\begin{picture}(1,0.57)
	\put(0.1,.46){\includegraphics[trim = 0 0 0 0, clip, width=.8\linewidth]{figures/quad/legend.png}} 
	\put(0.11,.555){\footnotesize Elemental quadrature rule:}	
	\put(0.53,.47){\footnotesize Duffy} \put(0.67,.47){\footnotesize Gauss with adj.~weights}	\put(0.15,.47){\footnotesize classical~Gauss}
	\put(0.325,.47){\footnotesize modified~Gauss}
	\put(0.145,.525){\footnotesize $\tilde n_\mathrm{qp}$}	\put(0.145,.5){\footnotesize $=3$}
	\put(0.26,.525){\footnotesize \twhite{$\tilde n_\mathrm{qp}$}}	\put(0.26,.5){\footnotesize \twhite{$=6$}}
	\put(0.37,.525){\footnotesize $\tilde n_\mathrm{qp}$}	\put(0.37,.5){\footnotesize $=3$}
	\put(0.54,.525){\footnotesize $\tilde n_\mathrm{qp}$}	\put(0.54,.5){\footnotesize $=6$}
	\put(0.71,.525){\footnotesize \twhite{$\tilde n_\mathrm{qp}$}}	\put(0.71,.5){\footnotesize \twhite{$=3$}}
	\put(0.83,.525){\footnotesize \twhite{$\tilde n_\mathrm{qp}$}}	\put(0.83,.5){\footnotesize \twhite{$=3$}}
	\put(0.02,0.23){\includegraphics[trim = 100 50 80 30, clip, width=.225\linewidth]{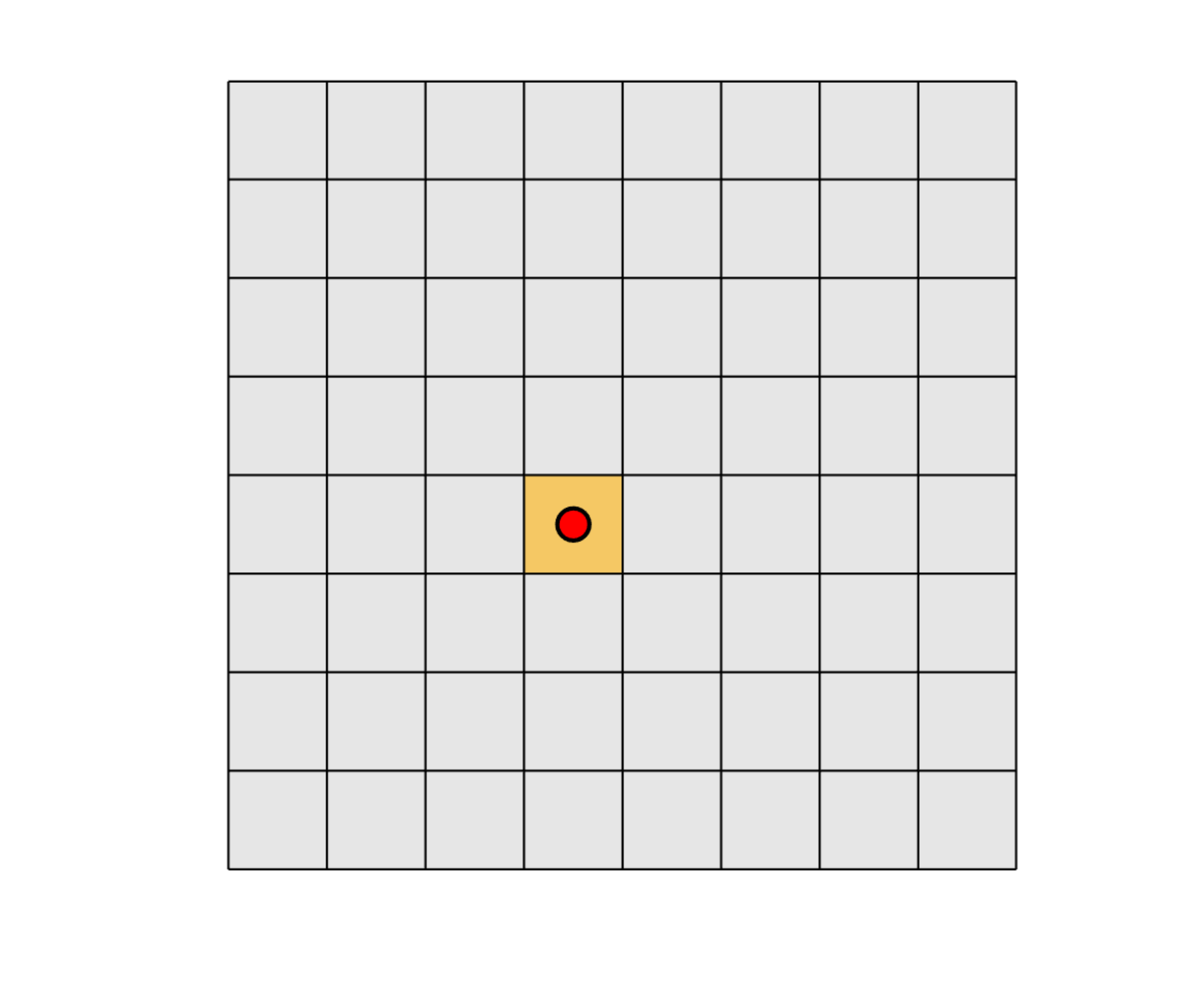}}
	\put(0.27,0.23){\includegraphics[trim = 100 50 80 30, clip, width=.225\linewidth]{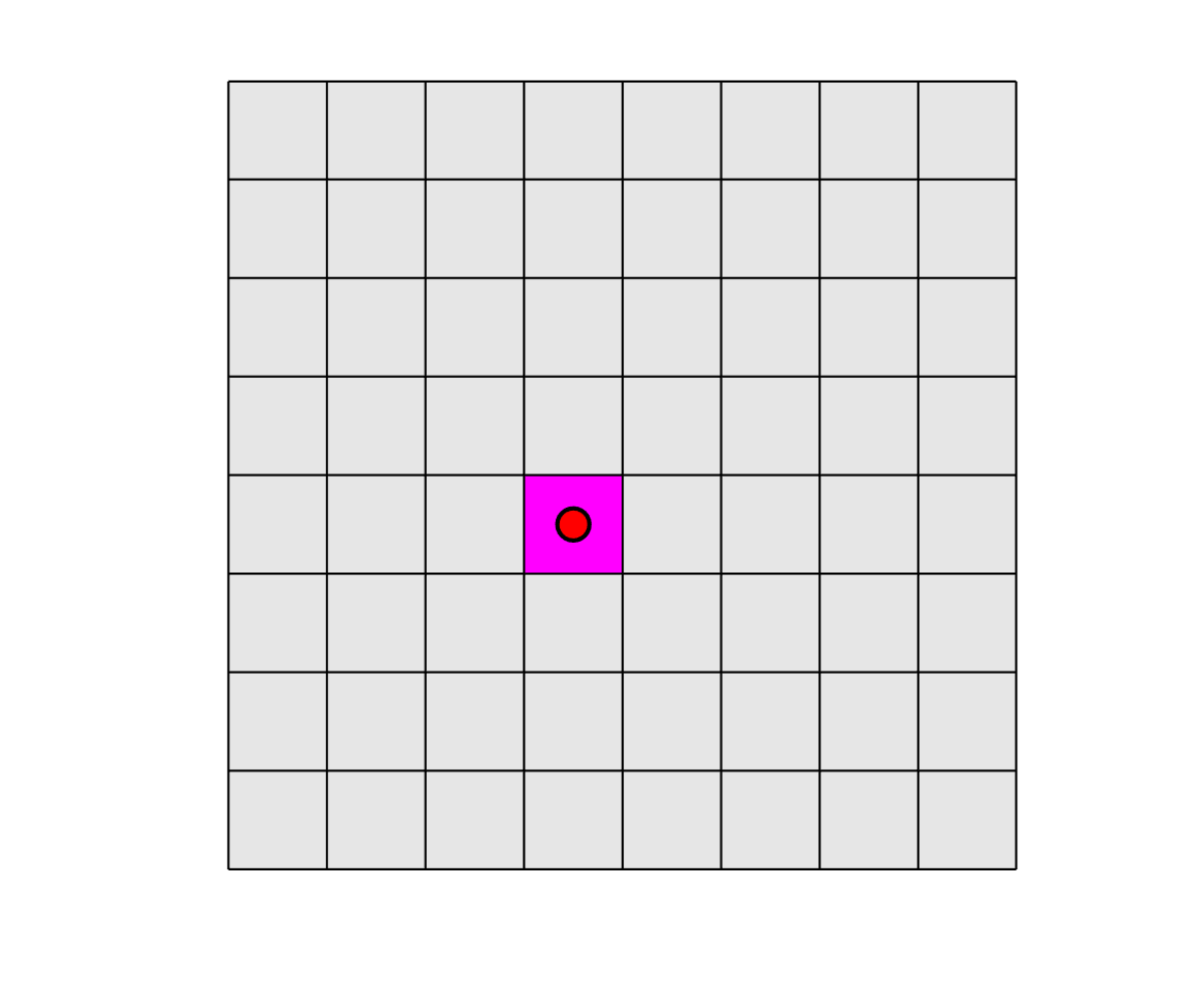}}
	\put(0.52,0.23){\includegraphics[trim = 100 50 80 30, clip, width=.225\linewidth]{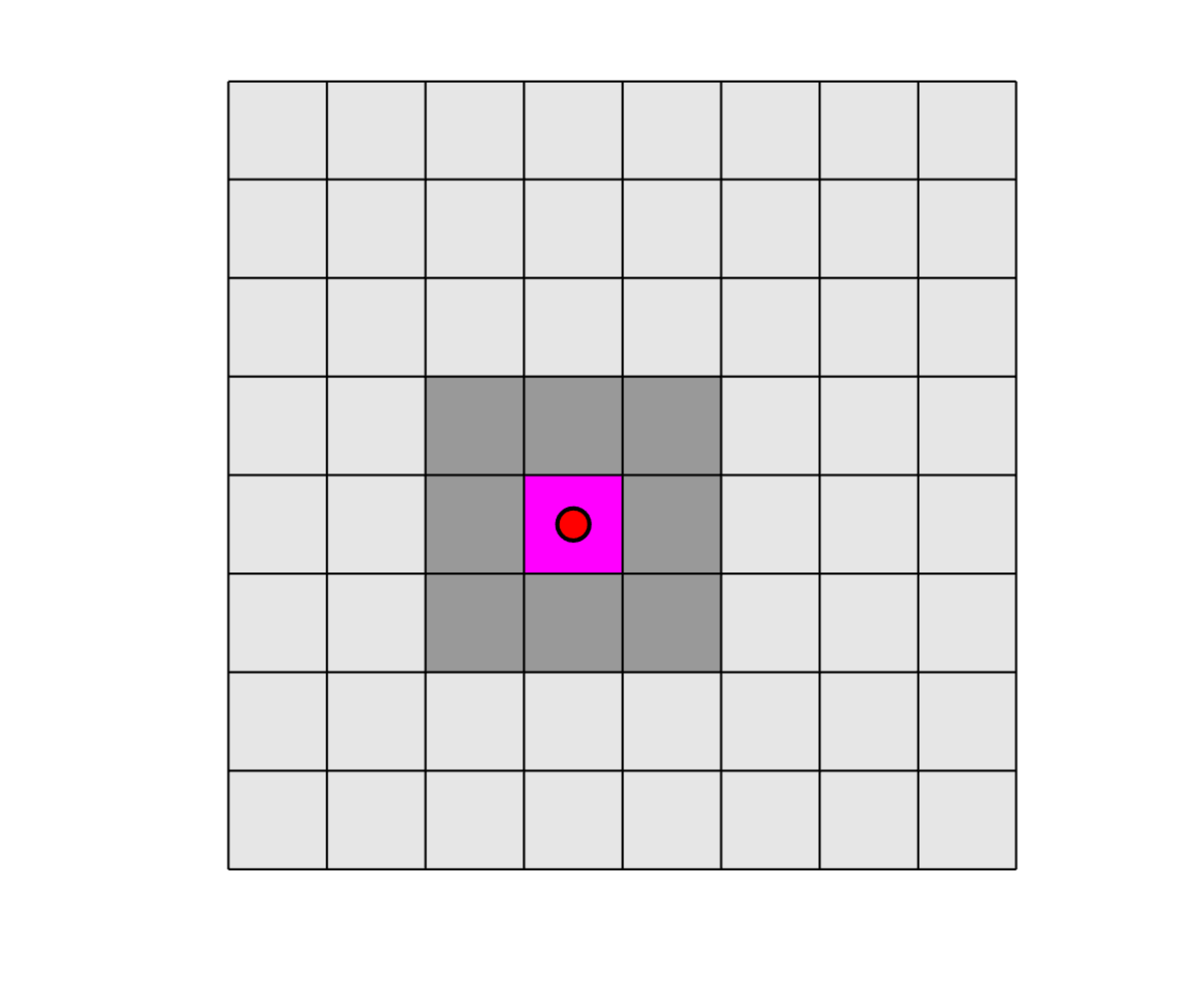}}
	\put(0.77,0.23){\includegraphics[trim = 100 50 80 30, clip, width=.225\linewidth]{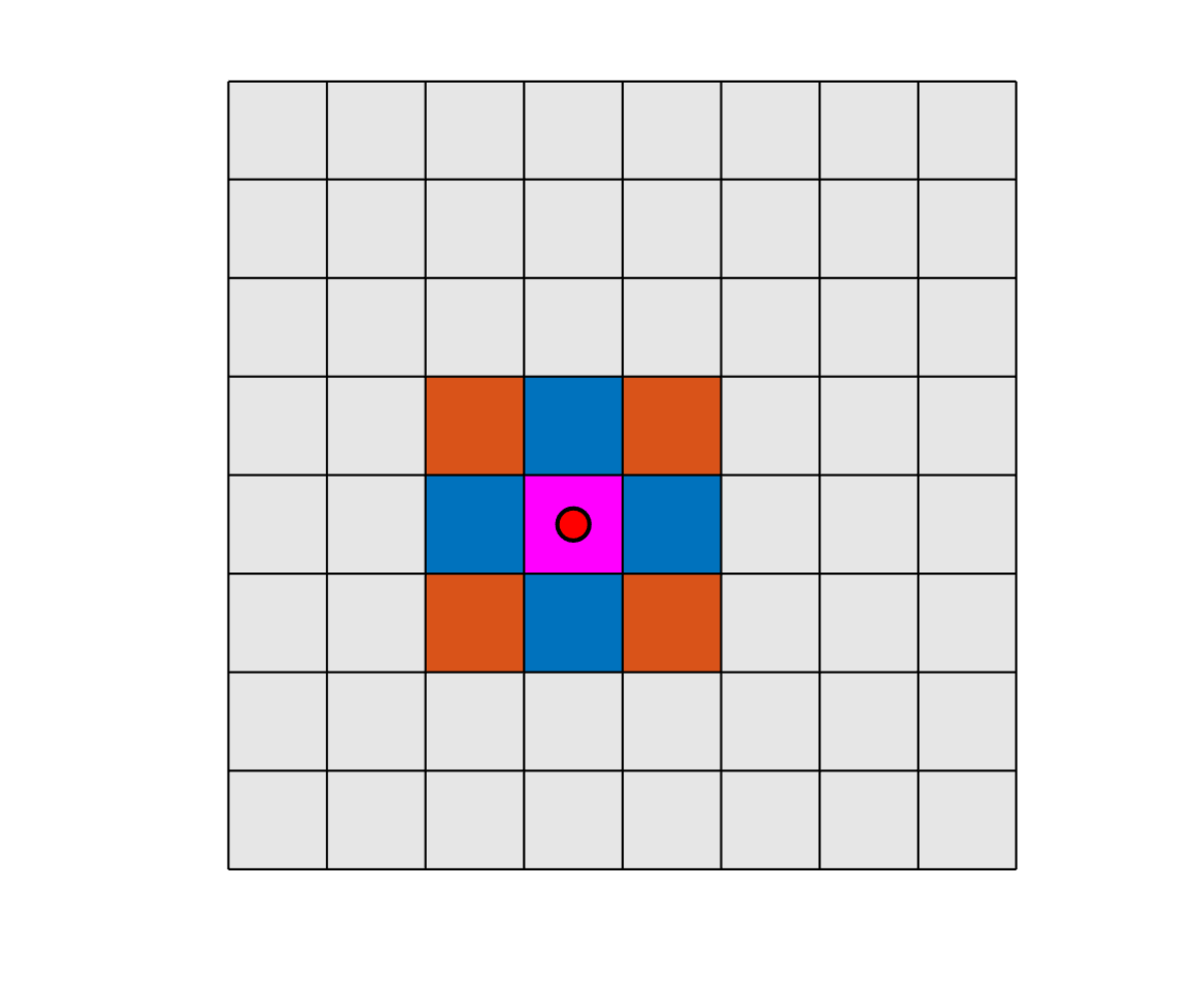}}
	\put(0.02,0){\includegraphics[trim = 100 50 80 30, clip, width=.225\linewidth]{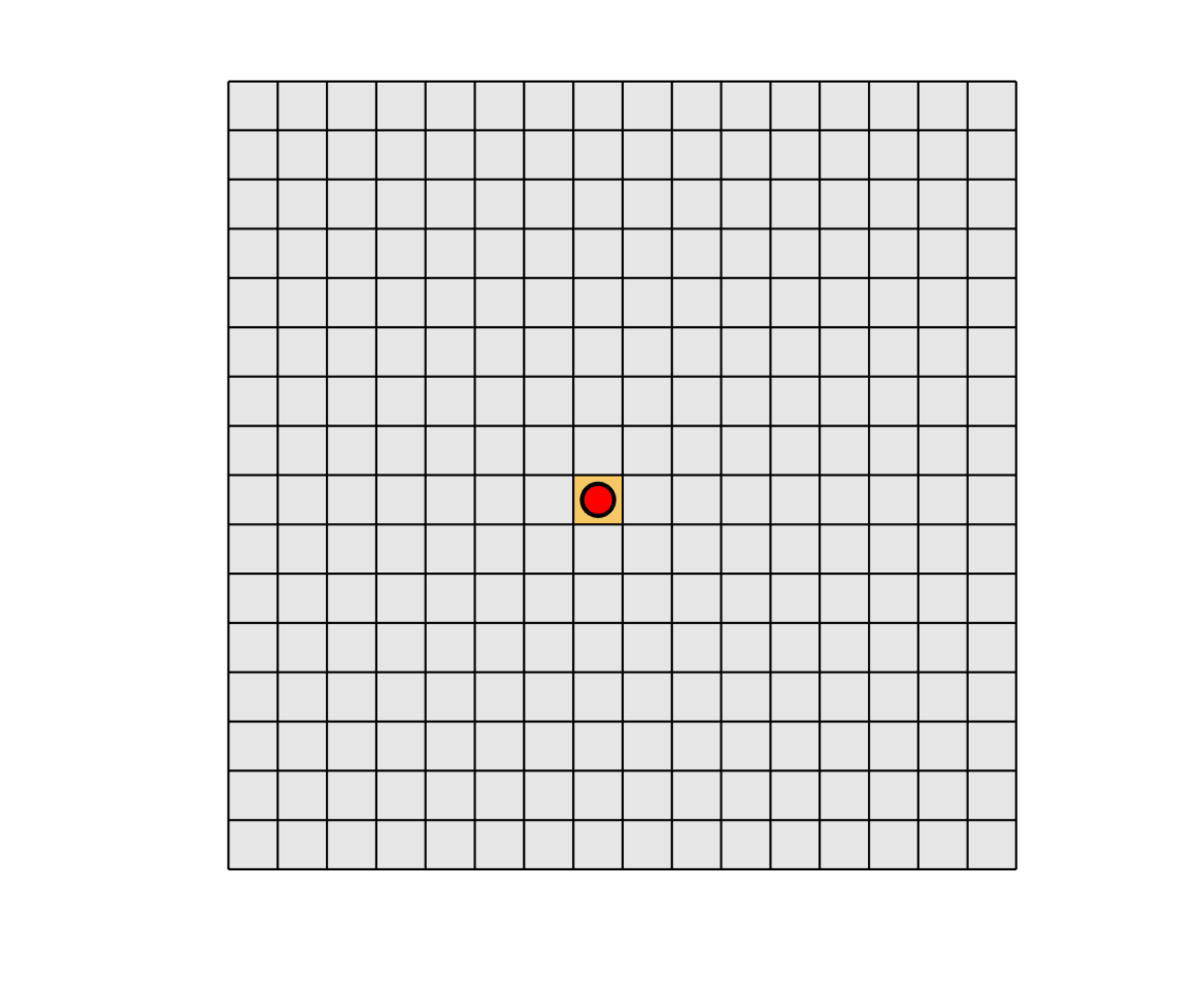}}
	\put(0.27,0){\includegraphics[trim = 100 50 80 30, clip, width=.225\linewidth]{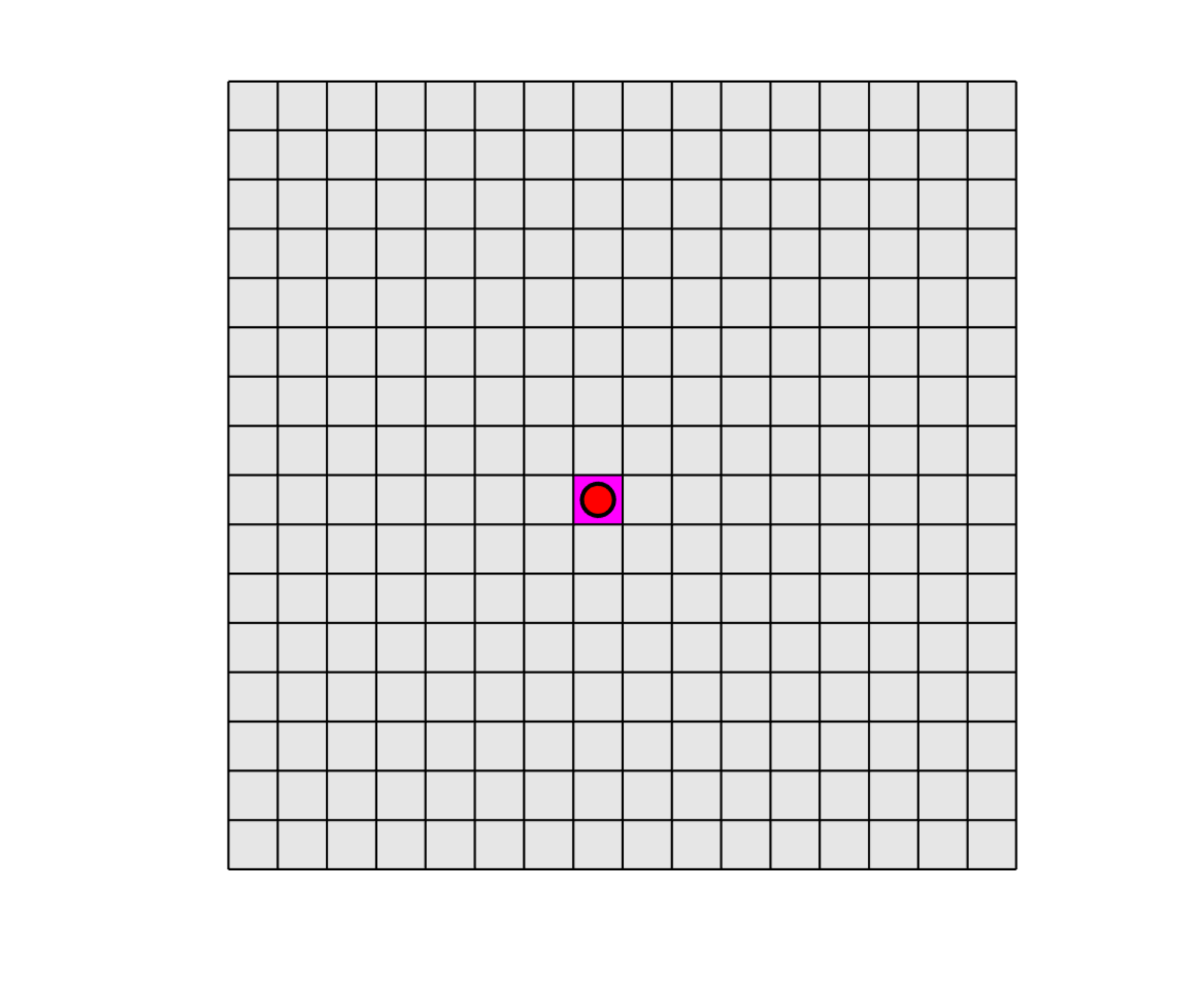}}
	\put(0.52,0){\includegraphics[trim = 100 50 80 30, clip, width=.225\linewidth]{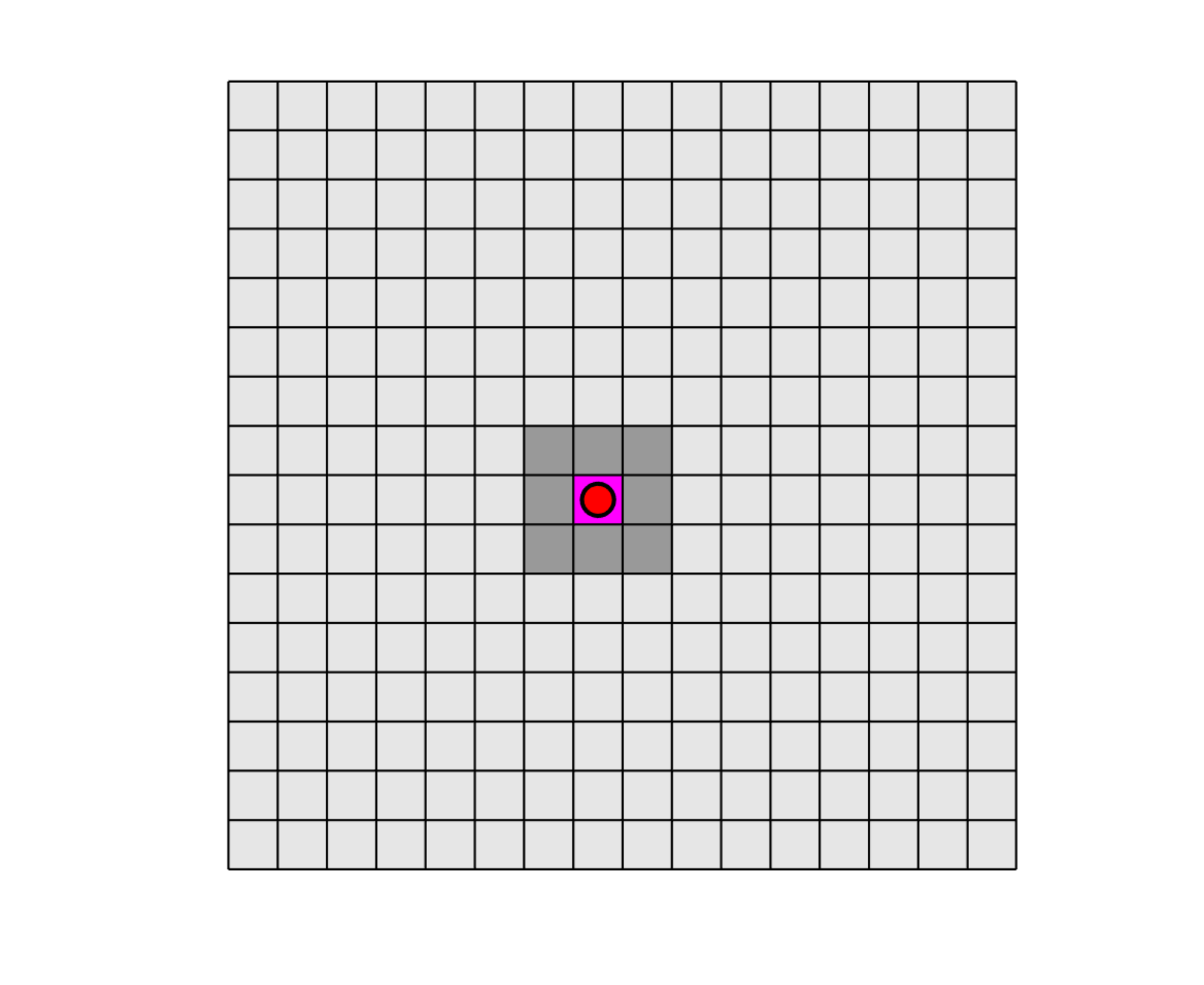}}
	\put(0.77,0){\includegraphics[trim = 100 50 80 30, clip, width=.225\linewidth]{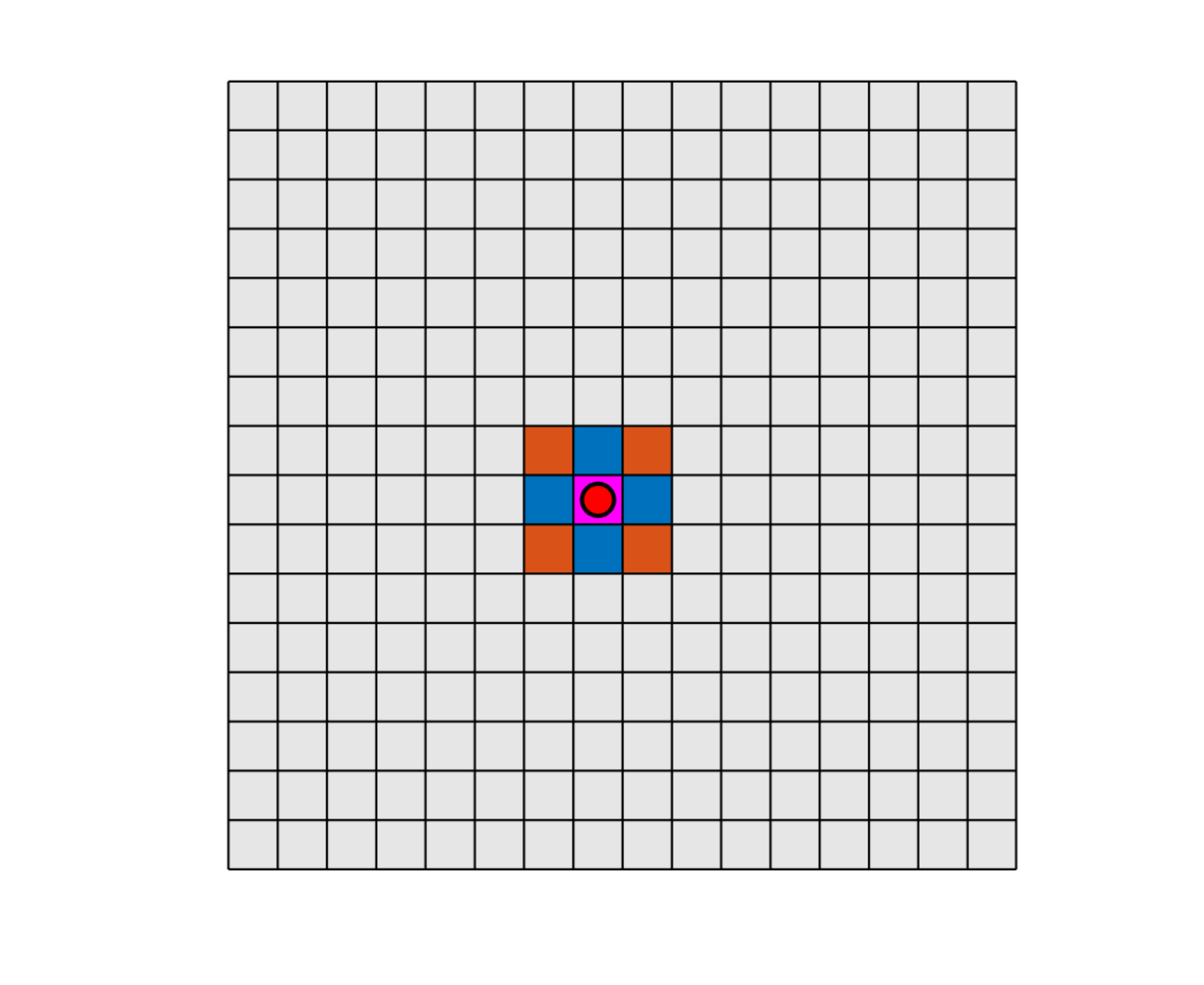}}
	\put(0,0){a.}\put(.25,0){b.}\put(.5,0){c.}\put(0.75,0){d.}
	\end{picture}
\caption{\textit{Quadrature point locations and weights}: Hybrid quadrature on a biquadratic B-spline sheet of refinement level $\ell=2$ (upper row) and $\ell=3$ (lower row). Quadrature rules used for Gauss-Legendre quadrature~(G,~a.), for hybrid Duffy-Gauss quadrature (DG,~b.), for Duffy-Gauss quadrature with progressive refinement~(DGr,~c.) and for hybrid Duffy-Gauss quadrature with adjusted weights~(DGw,~d.) considering collocation point $\by_0$ and $n_0=3$.}\label{fig:app_rules}
\end{figure}

\begin{figure}[h]
\unitlength\linewidth
\begin{picture}(1,1.22)		
	\put(0.01,0.9){\includegraphics[trim = 90 50 130 20, clip, width=.3\linewidth]{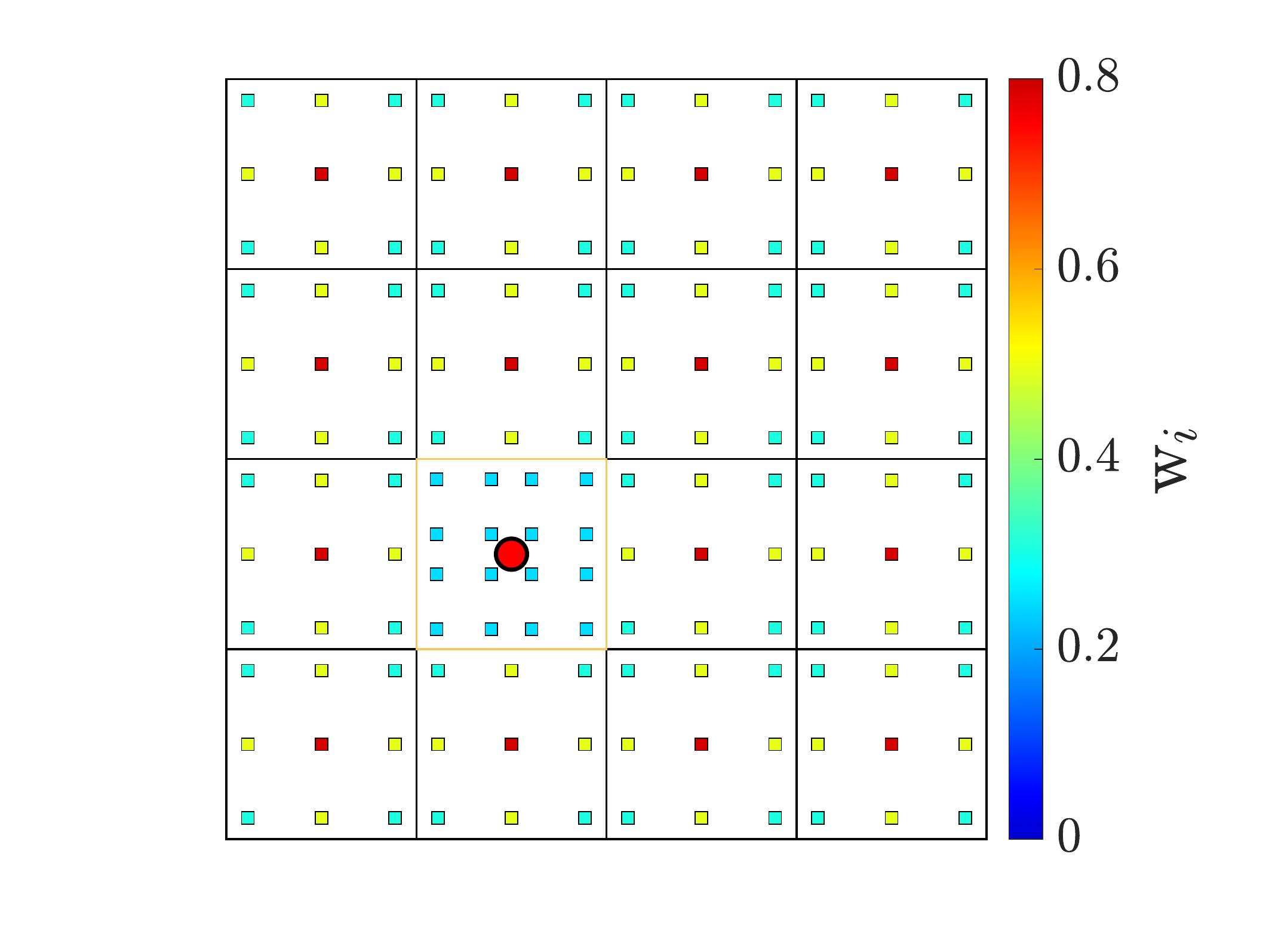}}
	\put(0.31,0.9){\includegraphics[trim = 90 50 130 20, clip, width=.3\linewidth ]{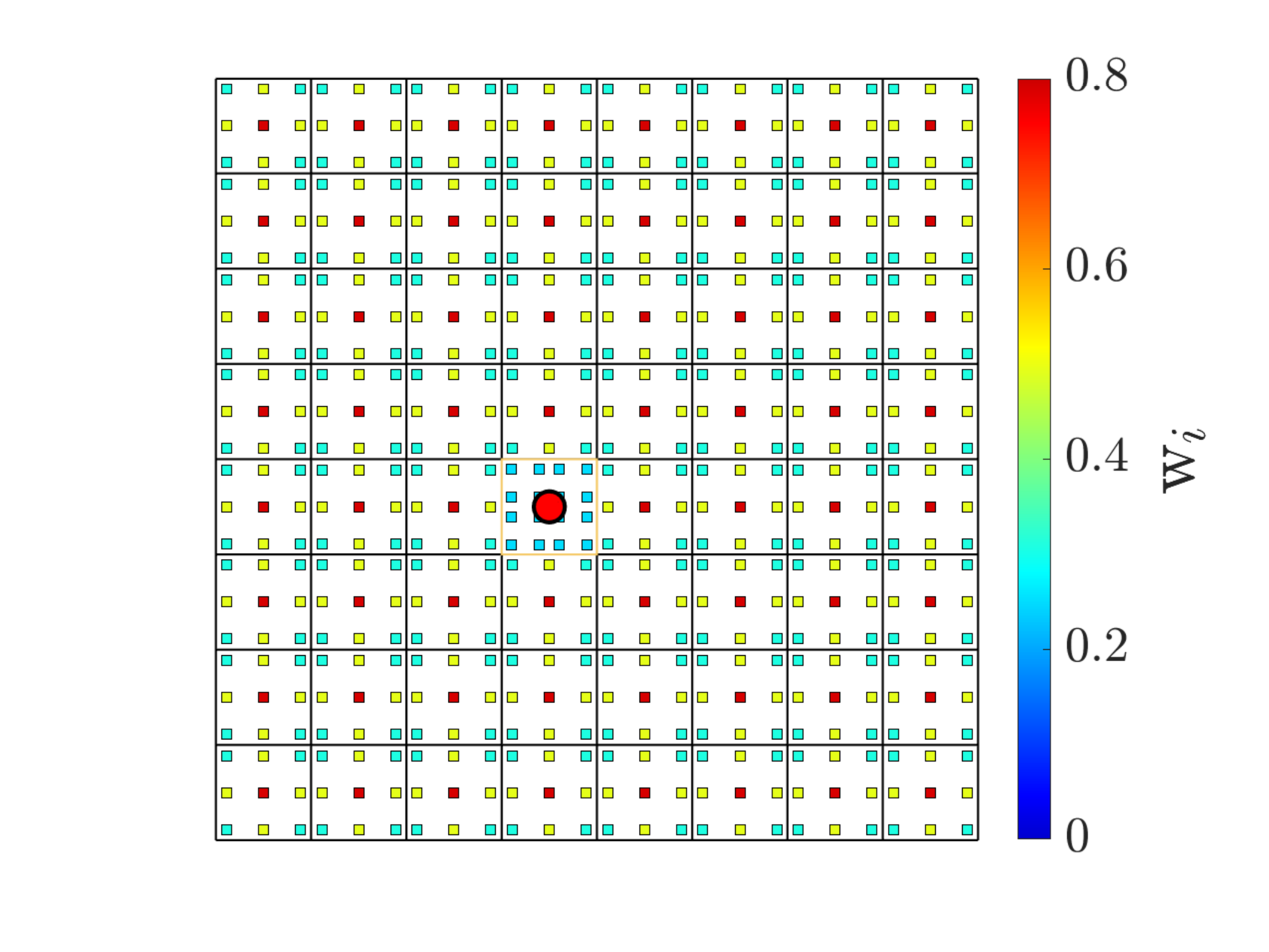}}
	\put(0.61,0.9){\includegraphics[trim = 90 50 130 20, clip, width=.3\linewidth ]{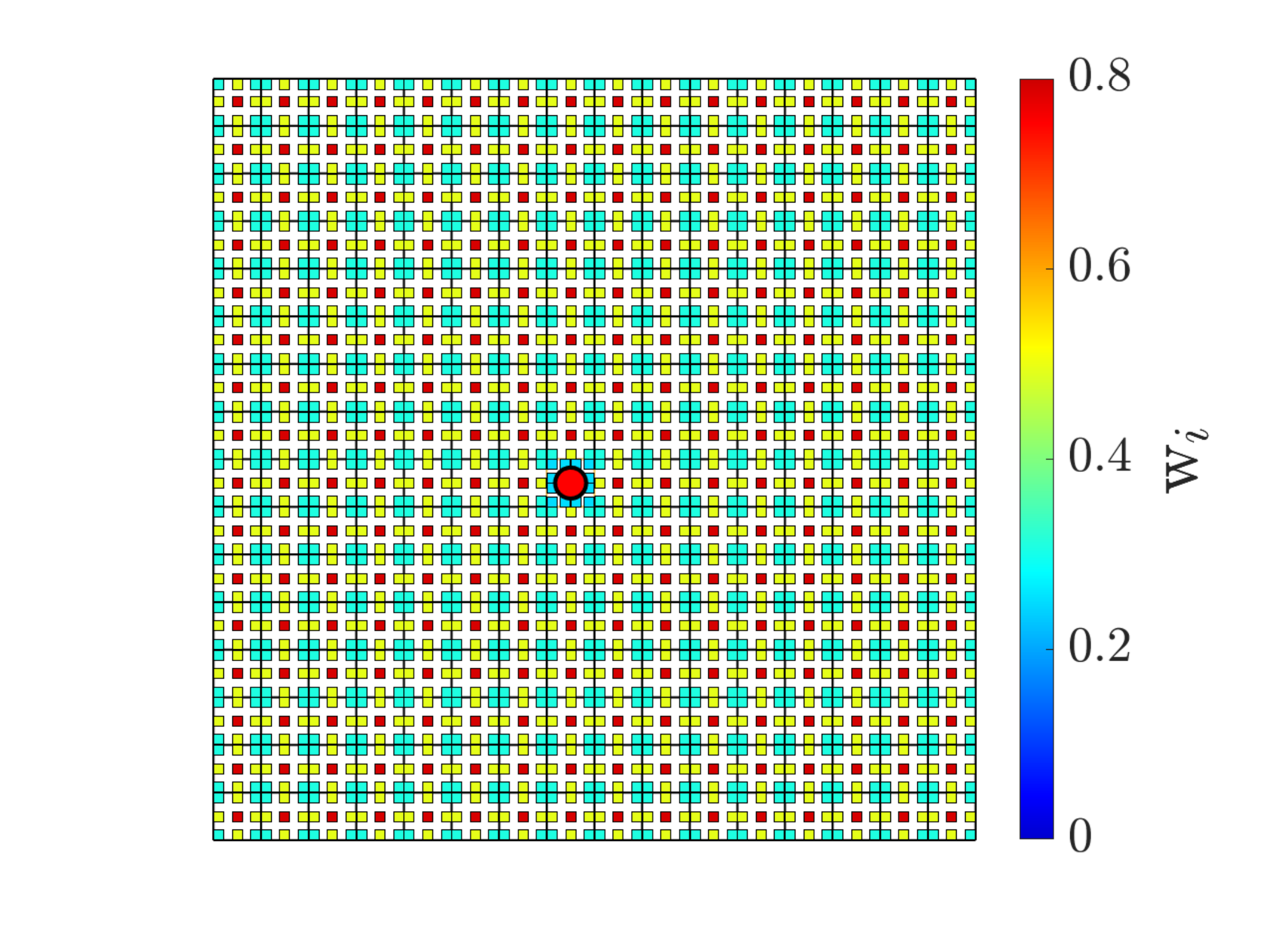}}
	\put(0.01,0.6){\includegraphics[trim = 90 50 130 20, clip, width=.3\linewidth]{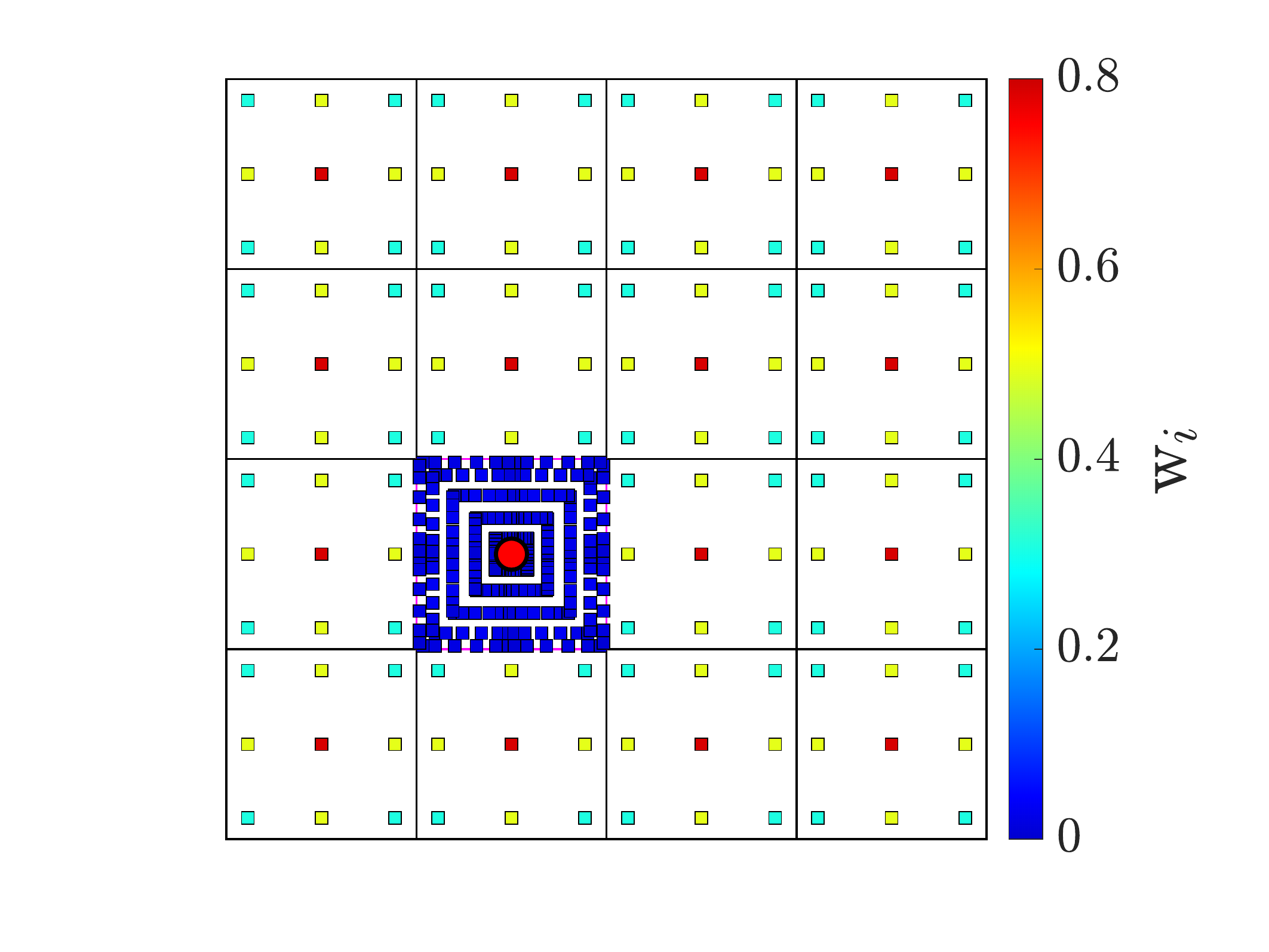}}
	\put(0.31,0.6){\includegraphics[trim = 90 50 130 20, clip, width=.3\linewidth ]{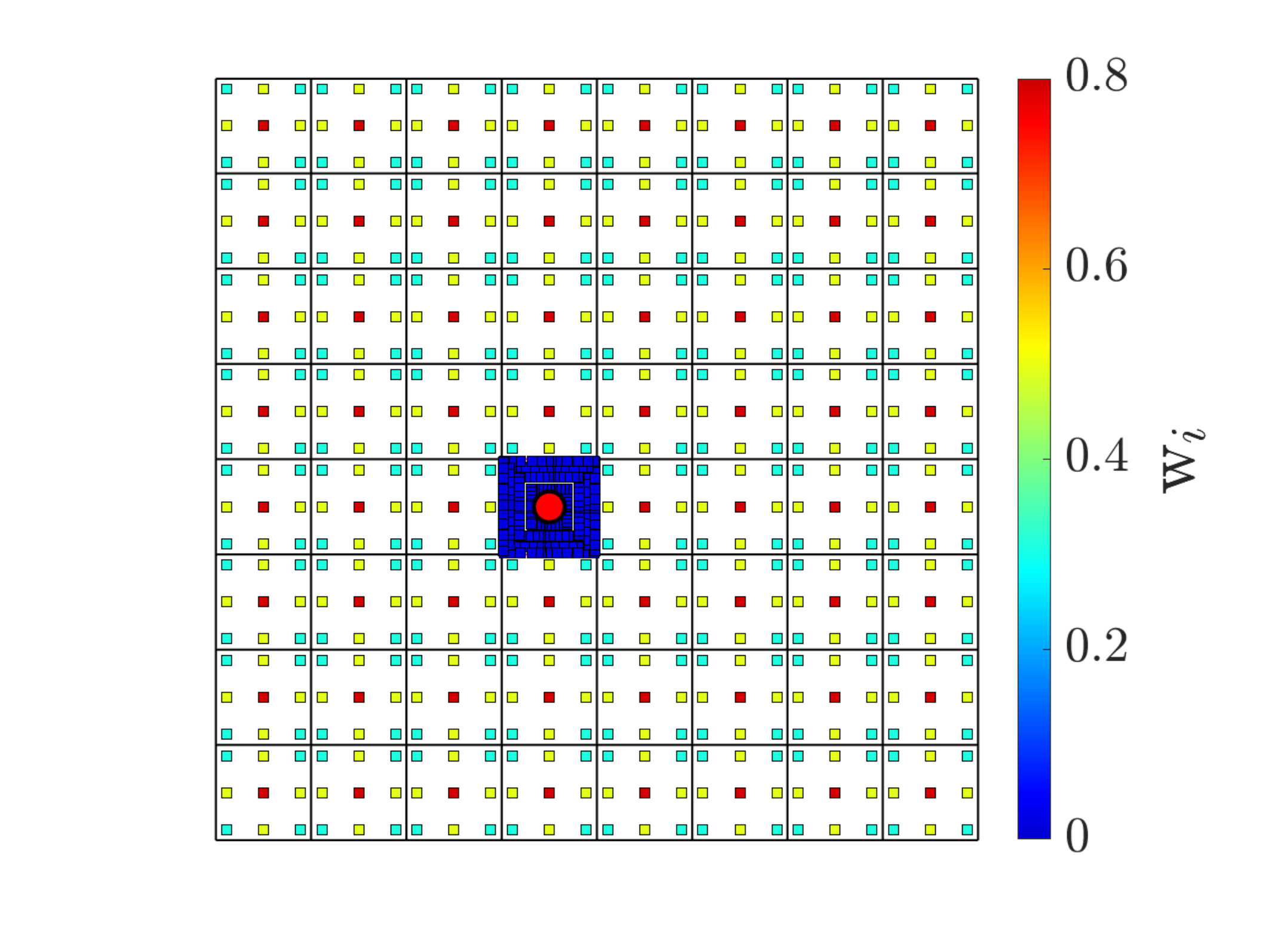}}
	\put(0.61,0.6){\includegraphics[trim = 90 50 130 20, clip, width=.3\linewidth ]{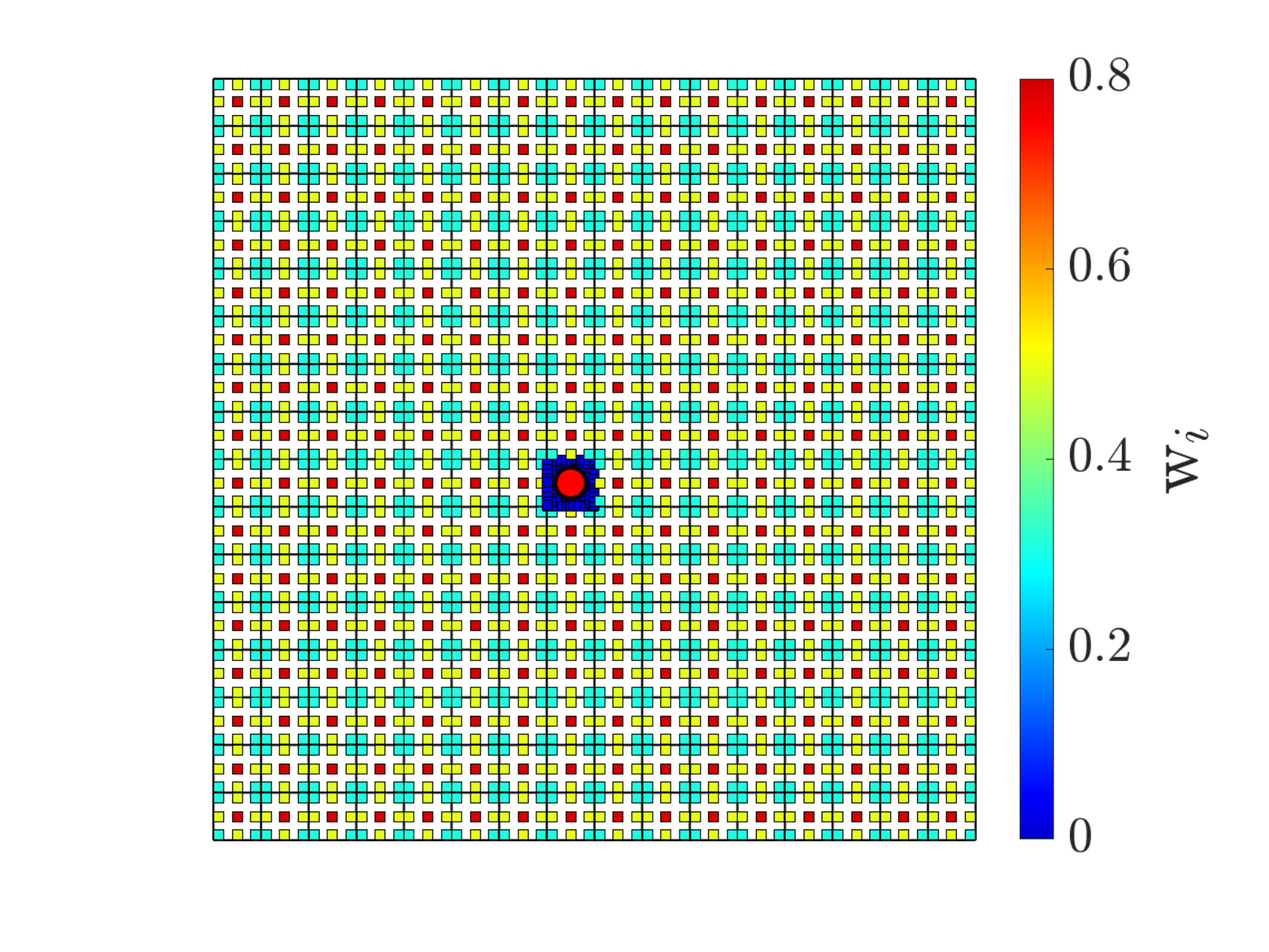}}
	\put(0.01,0.3){\includegraphics[trim = 90 50 130 20, clip, width=.3\linewidth]{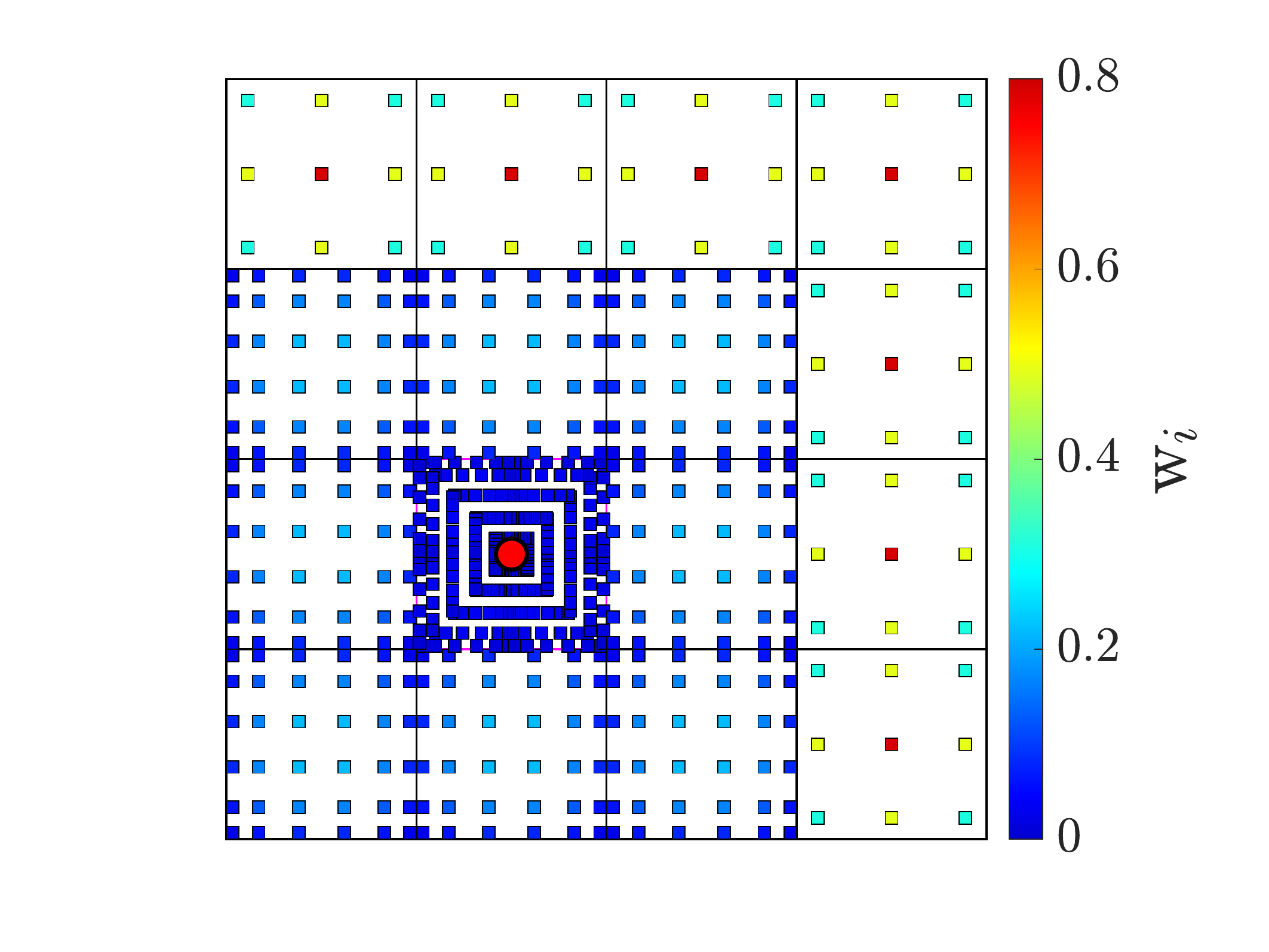}}
	\put(0.31,0.3){\includegraphics[trim = 90 50 130 20, clip, width=.3\linewidth ]{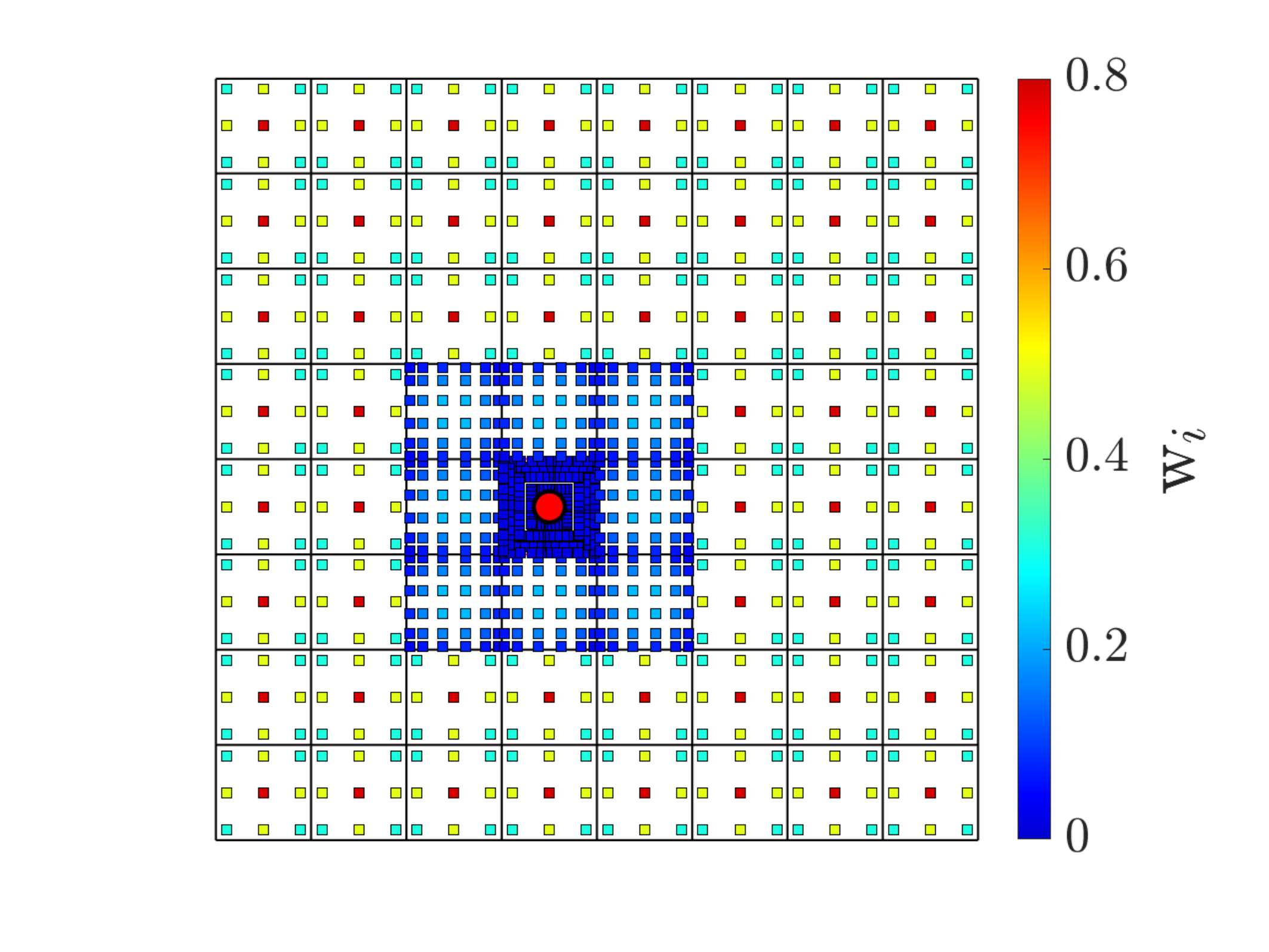}}
	\put(0.61,0.3){\includegraphics[trim = 90 50 130 20, clip, width=.3\linewidth ]{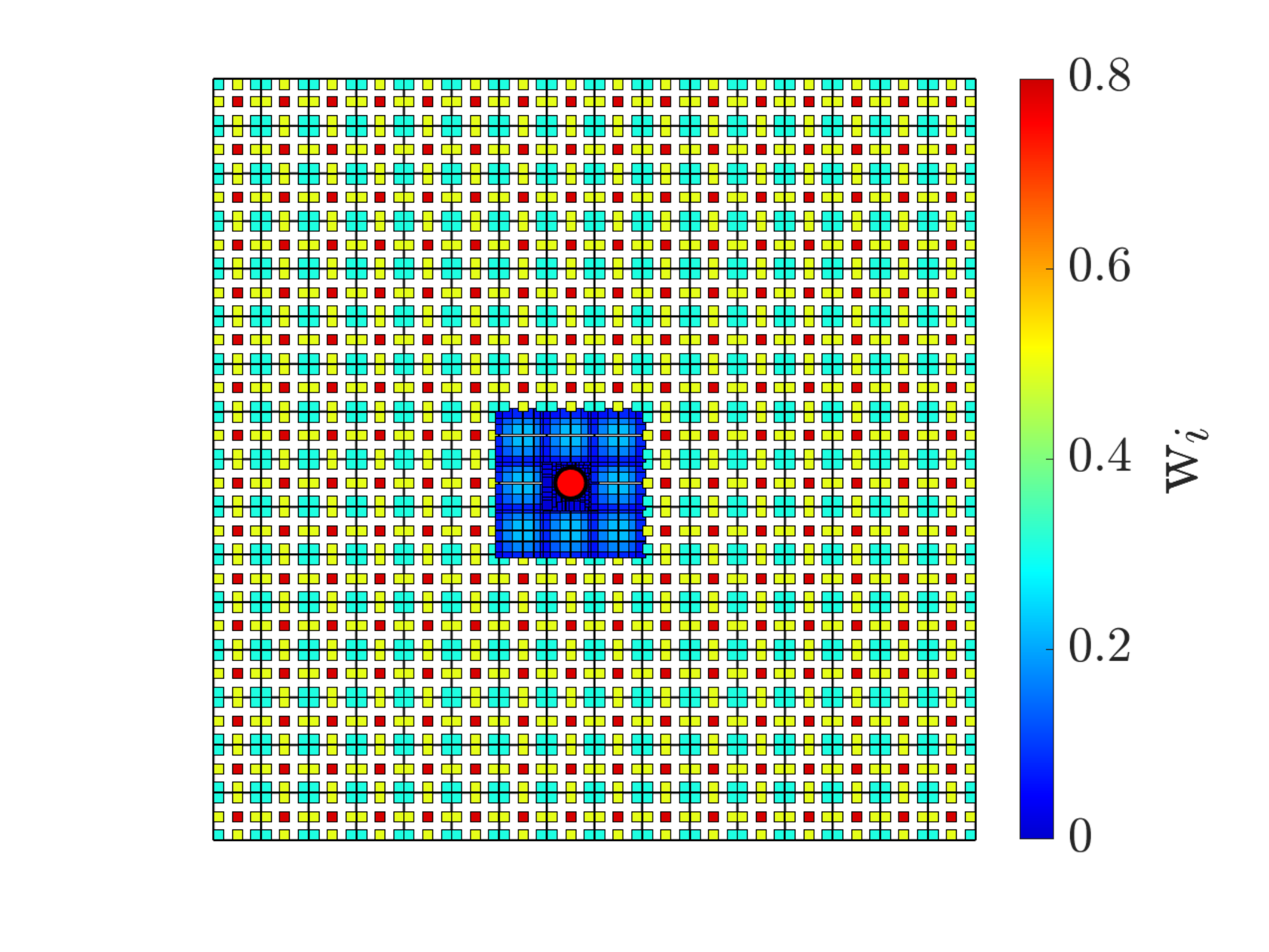}}
	\put(0.01,0){\includegraphics[trim = 90 50 130 20, clip, width=.3\linewidth]{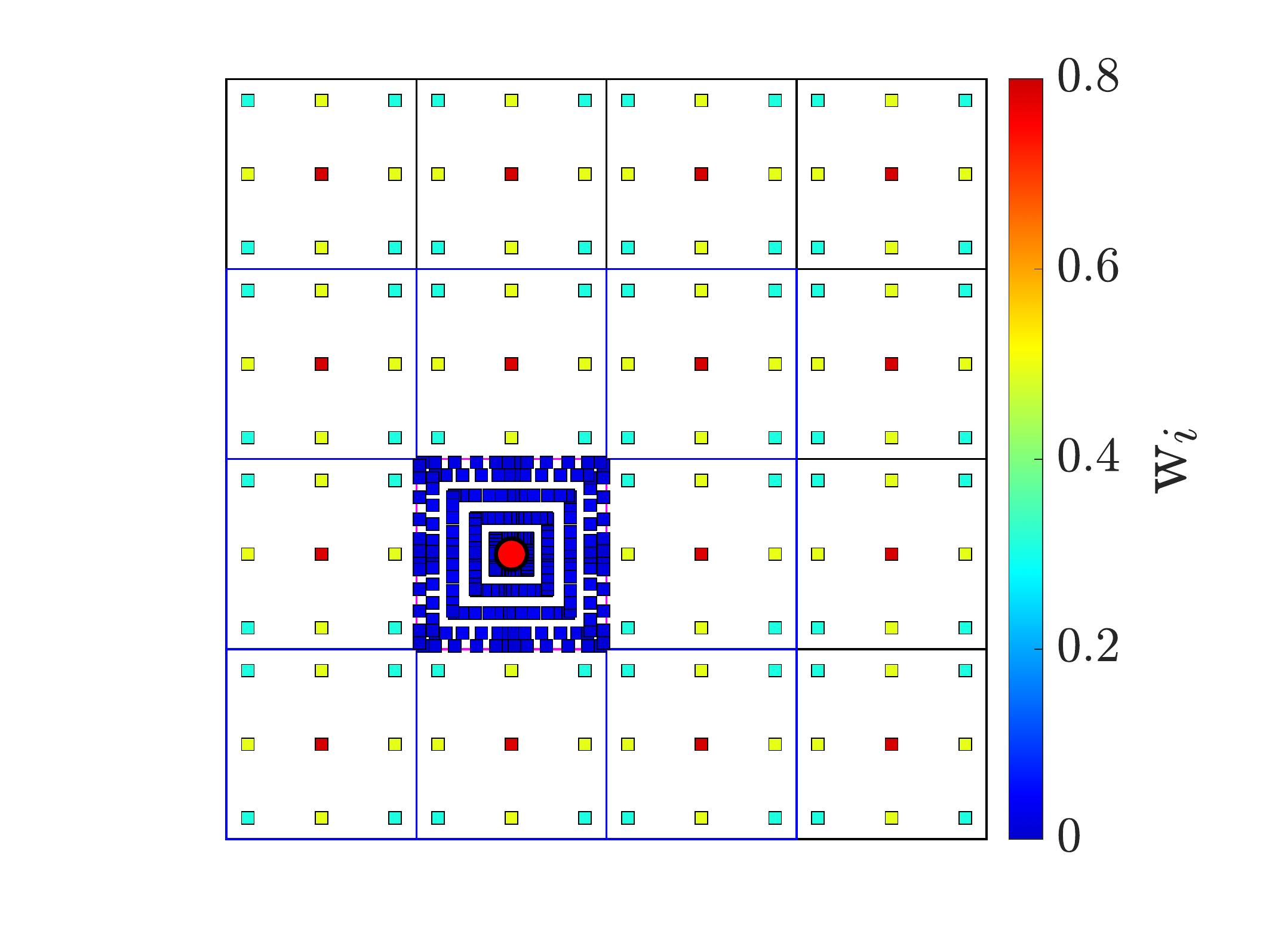}}
	\put(0.31,0){\includegraphics[trim = 90 50 130 20, clip, width=.3\linewidth ]{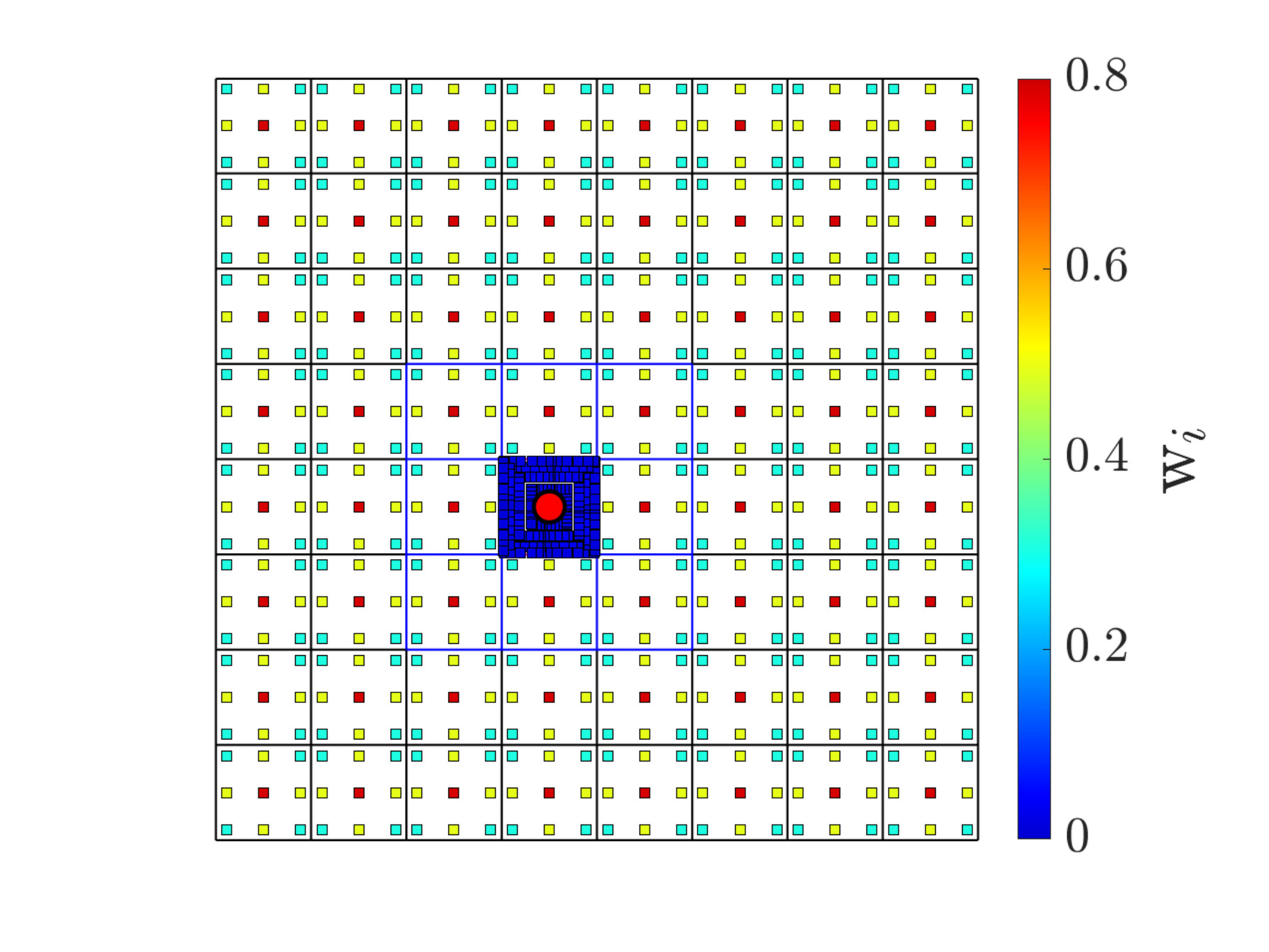}}
	\put(0.61,0){\includegraphics[trim = 90 50 130 20, clip, width=.3\linewidth ]{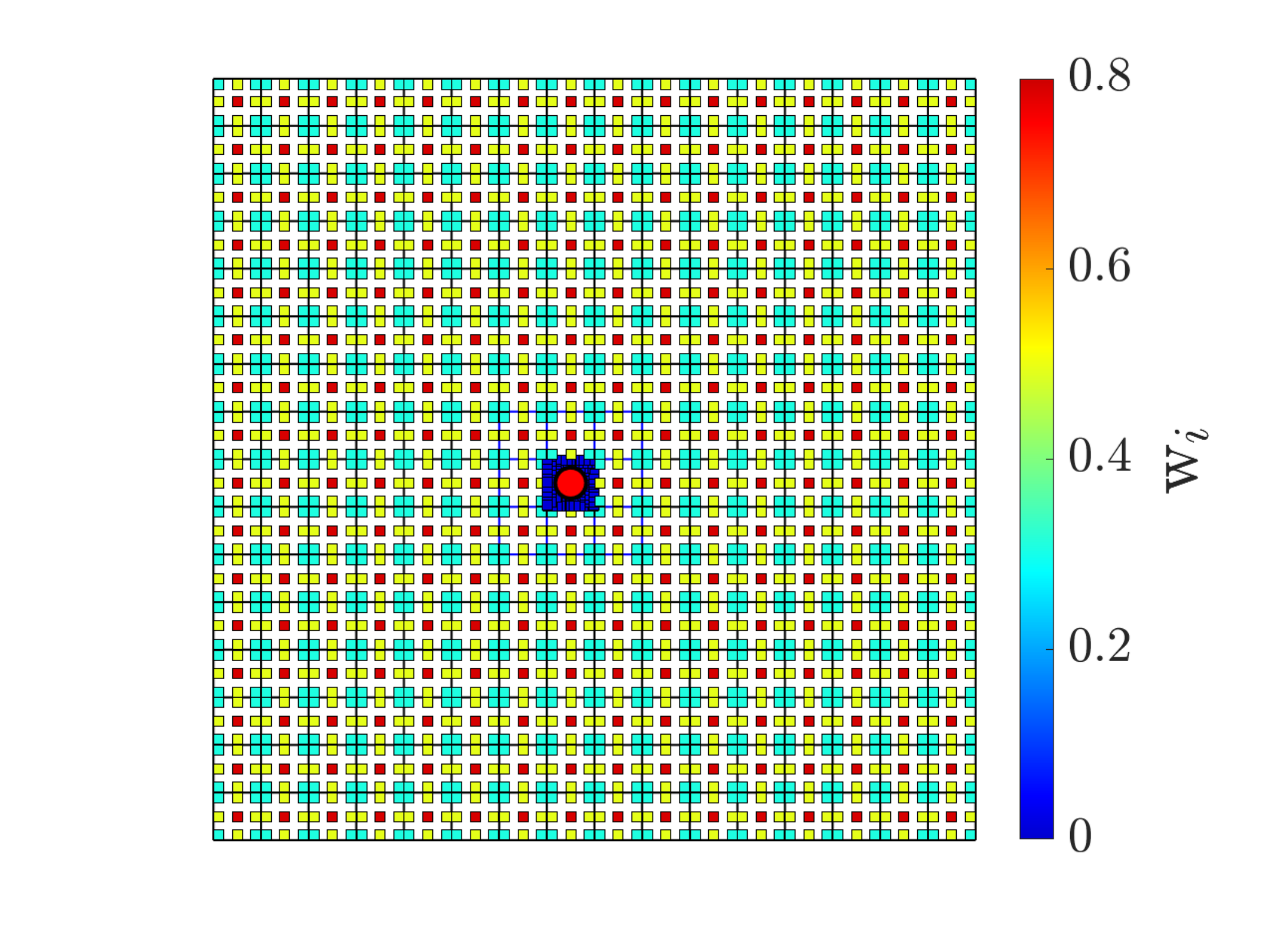}}
	\put(0.91,-0.01){\includegraphics[trim = 480 40 20 20, clip, width=.085\linewidth ]{figures/quad/scheme3/m3_nReg_3_nDuff_6/y153_w2.pdf}}
	\put(0.91,0.29){\includegraphics[trim = 480 40 20 20, clip, width=.085\linewidth ]{figures/quad/scheme3/m3_nReg_3_nDuff_6/y153_w2.pdf}}
	\put(0.91,0.59){\includegraphics[trim = 480 40 20 20, clip, width=.085\linewidth ]{figures/quad/scheme3/m3_nReg_3_nDuff_6/y153_w2.pdf}}
	\put(0.91,0.89){\includegraphics[trim = 480 40 20 20, clip, width=.085\linewidth ]{figures/quad/scheme3/m3_nReg_3_nDuff_6/y153_w2.pdf}}
	\put(0,.9){a.}\put(0,.6){b.}\put(0,.3){c.}\put(0,0){d.}
	\put(0.14,1.2){$\ell=1$}\put(0.43,1.2){$\ell=2$}\put(0.735,1.2){$\ell=3$}
	\end{picture}
\caption{\textit{Quadrature point locations and weights}: Hybrid quadrature on a biquadratic B-spline sheet of refinement level $\ell=1,2,3$. Location and weights of the quadrature points for G~(a.), DG~(b.), DGr~(c.) and for DGw~(d.) considering collocation point $\by_0$ and $n_0=3$.}\label{fig:app_w}
\end{figure}

This Appendix presents the quadrature point locations and weights for the new quadrature rule from Sec.~\ref{sec:quad_nearly_adjusted} (see Tables \ref{tab:Gauss} and \ref{tab:weights}) and for the new hybrid quadrature schemes from Sec.~\ref{sec:hybrid} (see Figs.~\ref{fig:app_rules} and \ref{fig:app_w}).

\begin{table}[h]
\small
\centering
\begin{tabular}{c|c|c|r|r|r} 
  $i$ & $\bxi^\mathrm{GL}_i$ & $\mrw^\mathrm{GL}_i$ & $\Delta \mrw_i^1$ & $\Delta \mrw_i^2$ & $\Delta \mrw_i^3$ \\ 
 \hline
1 & ($-\varphi$, $-\varphi$)	& 25/81	&	2.31938557634487e-03	& -3.06785320861036e-04	 &  2.91091742066762e-04	 \\
2 & (0, $-\varphi$)				& 40/81 &	-3.96194342728695e-04	& -7.65922184471912e-05	 &  -5.18107361679232e-04	 \\
3 & ($\varphi$, $-\varphi$) 	& 25/81 &	6.73693894590655e-04	& -8.35189243303947e-05	 &  3.01814399334244e-04	 \\
4 & ($-\varphi$, 0) 			& 40/81 &	-1.79631123847407e-03	& -7.65922184474688e-05	 &  2.78864715425431e-04	 \\
5 & (0, 0) 						& 64/81 &	-6.39888952746481e-03	& 9.40174607438005e-04	 &  -4.34093711446204e-04	 \\
6 & ($\varphi$, 0) 				& 40/81 &	2.35720567649134e-03	& -1.18292775950657e-04	 &  1.80583631911146e-04	 \\
7 & ($-\varphi$, $\varphi$) 	& 25/81 &	2.31938557634498e-03	& -8.35189243302836e-05	 &  -1.15997317956140e-04	 \\
8 & (0, $\varphi$ ) 			& 40/81 &	-3.96194342729583e-04	& -1.18292775950990e-04	 &  1.52292730246639e-05	 \\
9 & ($\varphi$, $\varphi$) 		& 25/81 &	6.73693894590766e-04	& -1.14956900993113e-04	 &  6.76167776109127e-06	 \\

\end{tabular}
\caption{\textit{Quadrature point locations and weights}: Locations $\bxi^\mathrm{GL}_i$ and weights $\mrw^\mathrm{GL}_i$ for classical Gauss-Legendre quadrature with $n_\mathrm{gp}^e= 3\times 3$, where $\varphi:= \sqrt{3/5}$. The first three sets of weight differences $\Delta \mrw_i^1$, $\Delta \mrw_i^2$ and $\Delta \mrw_i^3$ for Gauss-Legendre quadrature with adjusted weights (see Fig.~\ref{fig:DGw}).} \label{tab:Gauss}
\end{table}

\begin{table}[h]
\small
\centering
\begin{tabular}{c|r|r|r|r}
  $i$ & $\Delta \mrw_i^4$ & $\Delta \mrw_i^5$ &  $\Delta \mrw_i^6$ & $\Delta \mrw_i^7$ \\ 
 \hline
1 &  -6.54575192148843e-05 & 6.09298124450708e-04	 & 4.32446315561663e-04	 & -5.27414657617853e-05   \\
2 &  9.95820994587060e-06	 & -3.06041383078942e-03	 & -8.30579832303280e-05 &  1.80107109839511e-04   \\
3 &  -3.87801753661265e-05 & 1.56570492962638e-03	 & 4.32446315561719e-04	 & -1.42299297342940e-04   \\
4 &  9.95820994509344e-06	 & 2.47740158828003e-03	 & -2.18624664881650e-04 &  -1.98281169518255e-04  \\
5 &  1.74968939684494e-04	 & -1.36002345481279e-04 & -1.21107604659743e-03	 &  2.01683928024110e-04   \\
6 &  -1.01082187085044e-05 & 3.23573666747901e-04	 & -2.18624664881706e-04 &  -1.62236018943496e-05  \\
7 &  -3.87801753655159e-05 & -1.52026938222949e-03 & 2.20203796807650e-04	 &  6.72447376661167e-05   \\
8 &  -1.01082187095591e-05 & 1.03304298742501e-04	 & 3.91571815489933e-04	 &  4.48306966693846e-06   \\
9 &  -3.27884950691026e-05 & -1.93807081374009e-04 & 2.20203796807483e-04	 &  -4.43308431077871e-05  \\
\end{tabular}
\caption{\textit{Quadrature point locations and weights}: Weight differences $\Delta \mrw_i^4$, $\Delta \mrw_i^5$, $\Delta \mrw_i^6$ and $\Delta \mrw_i^7$ for Gauss quadrature with adjusted weights with $n_\mathrm{gp}^e= 3\times 3$ (see Fig.~\ref{fig:DGw}). The $j$-th set of adjusted weights is determined by $\mrw_i^j =\mrw_i^\mathrm{GL}+\Delta \mrw_i^j$ according to \eqref{eq:quad_adjusted_dw}.}\label{tab:weights}
\end{table}

\FloatBarrier

\bibliographystyle{own_abbrvnat}
\bibliography{bibliography}

\end{document}